%% file: kronecker-entropy.tex
\documentclass[11pt,a4paper]{article}
\usepackage[T1]{fontenc}
\usepackage{lmodern}
\usepackage[margin=27mm]{geometry}
\usepackage{amsthm,amsmath,amssymb}
\usepackage{stmaryrd}
\usepackage{mathtools,booktabs,enumitem}
\usepackage{microtype}
\usepackage{needspace}
\usepackage{graphicx,float}
\usepackage[dvipsnames]{xcolor}
\usepackage{hyperref}
\hypersetup{colorlinks=true,linkcolor=MidnightBlue,citecolor=MidnightBlue,urlcolor=MidnightBlue,pdftitle={On Douglas' conjecture for m-Kronecker quivers},pdfauthor={So Okada}}
\numberwithin{equation}{section}
\newtheorem{theorem}{Theorem}[section]

\newtheorem{corollary}[theorem]{Corollary}
\newtheorem{lemma}[theorem]{Lemma}
\newtheorem{proposition}[theorem]{Proposition}
\newtheorem{conjecture}[theorem]{Conjecture}
\newtheorem{definition}[theorem]{Definition}
\theoremstyle{definition}
\newtheorem{example}[theorem]{Example}
\theoremstyle{remark}

\newcommand\mdoubleplus{\ensuremath{\mathbin{+\mkern-10mu+}}}
\newcommand{\ccn}{\mdoubleplus}
\newcommand{\prf}[1]{\par{\emph{Proof of #1.}}}
\newcommand{\Z}{\mathbb Z}
\newcommand{\R}{\mathbb R}
\newcommand{\C}{\mathbb C}
\newcommand{\Q}{\mathbb Q}
\newcommand{\cA}{\mathcal A}
\newcommand{\cM}{\mathcal M}

\newcommand{\floor}[1]{\left\lfloor#1\right\rfloor}
\newcommand{\ceil}[1]{\left\lceil#1\right\rceil}
\setlist{itemsep=3pt,topsep=5pt}
\newlist{hypotheses}{enumerate}{1}
\setlist[hypotheses]{label=(\alph*),ref=(\alph*),leftmargin=*,font=\itshape}
\title{On Douglas' conjecture for $m$-Kronecker quivers%
\thanks{Some \TeX{} and Lean formalization files are available at\\
\url{https://github.com/so-okada/kronecker-euler-sources}.}}
\author{So Okada}
\date{\today{}}
\begin{document}
\maketitle
\begin{abstract}
We prove Douglas' conjecture for $m$-Kronecker quivers for each $m\geq3$
on the normalized logarithms of Euler characteristics of moduli
spaces.
Also, we state a log-concavity conjecture on Euler characteristics near slope one for the 3-Kronecker quiver.
\end{abstract}

\input{sections/introduction-and-statement-of-the-theorem}

\input{sections/finite-model-for-euler-characteristics-at-coprime-dimensions}

\input{sections/shell-paths-and-gradings}

\input{sections/run-edits}

\input{sections/rational-slopes}

\input{sections/irrational-slopes-for-m-ge-4}

\input{sections/uniform-lower-bound-near-slope-one-for-the-three-arrow-quiver}

\input{sections/central-upper-estimate-by-deletion}

\input{sections/three-arrow-case}

\input{sections/reflection-endpoints-and-the-central-value}

\input{sections/conjectural-log-concavity}

\input{sections/acknowledgements}

\bigskip
\noindent
So Okada\\
National Institute of Technology, Oyama College\\
Oyama Tochigi, Japan 323-0806\\
e-mail: \href{mailto:so.okada@gmail.com}{\nolinkurl{so.okada@gmail.com}}
\end{document}

%% file: sections/introduction-and-statement-of-the-theorem.tex
\section{Introduction and statement of the theorem}\label{sec:intro}
Let $m\in\Z_{\ge3}$. Let $K(m)$ be the quiver with two vertices and $m$ arrows from the first
vertex to the second. All quiver representations are over $\C$.
Let $\chi$ denote the ordinary topological Euler characteristic.
For positive coprime integers $a$ and $b$, let $\cM^s_{(a,b)}(K(m))$ denote the moduli space of stable representations of $K(m)$ of dimension $(a,b)$, where $a$ is the dimension at the source and $b$ at the target.
Let
\[
 E_m(a,b)=\chi\bigl(\cM^s_{(a,b)}(K(m))\bigr).
\]
Our stability requires
\begin{equation}\label{eq:stability}
 av-bu>0
\end{equation}
for every nonzero proper subrepresentation of dimension $(u,v)$.
For
\[
 \lambda_\pm=\frac{m\pm\sqrt{m^2-4}}2,
 \qquad I_m=[\lambda_-,\lambda_+],
\]
let
\[
 \cA_m=\{(a,b)\in\Z_{>0}^2:\gcd(a,b)=1,\
       \cM^s_{(a,b)}(K(m))\ne\emptyset\}.
\]

We prove Douglas' conjecture in the following
\cite[Conjecture 6.1]{Weist}
(see also~\cite{Okada}).
The study of black hole entropies motivates Douglas' conjecture~\cite{CS,Denef,DM,DFR,MPS,Weist}.

\begin{theorem}\label{thm:main}
For each integer $m\geq3$, there is a continuous function
$F_m:I_m\to\R_{\ge0}$ such that
\[
 F_m(\lambda_-)=F_m(\lambda_+)=0,\qquad
 0<F_m(r)<\infty\quad(\lambda_-<r<\lambda_+)
\]
and
\begin{equation}\label{eq:main-limit}
 \lim_{\substack{a\to\infty,\ (a,b)\in\cA_m\\\frac{b}{a}\to r\in I_m}}
 \frac{\log E_m(a,b)}a=F_m(r).
\end{equation}
\end{theorem}

By~\cite[Theorem 4.1]{Sch01} and~\cite[Theorem 6.1]{Sch92}, we have
\begin{equation}\label{eq:nonempty-domain}
 \cA_m=\{(a,b)\in\Z_{>0}^2:
       \gcd(a,b)=1,\ a^2+b^2-mab\le1\}.
\end{equation}
For $a^2+b^2-mab\le0$, see also~\cite[Theorems 2.7(1) and 4.1(1)]{Rei}.
Furthermore,~\eqref{eq:main-limit} is equivalent to the following: for every $r\in I_m$ and $\varepsilon>0$, there are $\delta>0$ and $A>0$
such that
\[
 \left|\frac{\log E_m(a,b)}a-F_m(r)\right|<\varepsilon
 \quad\text{if }(a,b)\in\cA_m,\ a>A,
       \left|\frac ba-r\right|<\delta.
\]

\subsection{The counting model and the comparison argument}
We use a counting formula for $mE_m(a,b)$~\cite{B} with quiver representations.
The finite model counts colored gradings on paths.
We construct maps between these sets of gradings and obtain a product estimate for their cardinalities.
The product estimate gives an exponential growth rate under repeated concatenation.

We first prove the limit in~\eqref{eq:main-limit} for rational $r$.
We then obtain the limit for irrational $r$ by comparison with rational values.
Finally, we prove continuity of $F_m$ and the endpoint values
$F_m(\lambda_-)=F_m(\lambda_+)=0$.

%% file: sections/finite-model-for-euler-characteristics-at-coprime-dimensions.tex
\section{A finite model for Euler characteristics at coprime dimensions}\label{sec:dictionary}
\subsection{Euler characteristics}\label{subsec:euler-characteristics}
By~\cite[Theorem 4.1 and Proposition 4.3]{King}
(see also~\cite[Proposition 3.5 and Section 3.5]{Rei08}), we identify
the stable moduli $\cM^s_{(a,b)}(K(m))$ with a quasi-projective complex
algebraic variety in its analytic topology.

By~\cite[Corollary 6.10]{Dimca}, an algebraic variety has the
homotopy type of a finite CW complex. Consequently, ordinary rational singular homology is finite-dimensional in each degree and vanishes in all sufficiently high degrees. Thus,
\[
 \chi(X)=\sum_{i\ge0}(-1)^i\dim_{\Q}H_i(X;\Q)
\]
is a finite sum for $X=\cM^s_{(a,b)}(K(m))$.

\subsection{Vertex exchange and reflection}
\begin{lemma}\label{lem:moduli-symmetries}
Let $a$ and $b$ be positive coprime integers.
Then, we have the following.

\begin{enumerate}
\item\label{c:moduli-symmetries-identities} We have
\begin{equation}\label{eq:symmetries}
 E_m(a,b)=E_m(b,a)=E_m(a,ma-b)
\end{equation}
provided $ma-b>0$.
\item\label{c:moduli-symmetries-invariants} We have the following properties for the transformations
from $(a,b)$ to $(b,a)$ and from $(a,b)$ to $(a,ma-b)$.
\begin{enumerate}[label=(\alph*),ref=\theenumi(\alph*)]
\item\label{c:moduli-coprimality} They preserve coprimality.
\item\label{c:moduli-quadratic} They preserve the quadratic form $mab-a^2-b^2$.
\end{enumerate}
\end{enumerate}
\end{lemma}
\begin{proof}
\prf{Claim~\ref{c:moduli-symmetries-identities}} Let us identify the dimension transformations that preserve the Euler characteristic.
By vertex exchange and the reflection isomorphism of
\cite[Proposition 4.3]{Weist}, we obtain~\eqref{eq:symmetries}.
On dimension vectors, the three steps are
\[
 (a,b)\xrightarrow{\text{exchange}}(b,a)
 \xrightarrow{\text{reflection}}(ma-b,a)
 \xrightarrow{\text{exchange}}(a,ma-b).
\]
Transpose duality exchanges
subrepresentations with quotients. Then, vertex exchange respects
\eqref{eq:stability}.

\prf{Claim~\ref{c:moduli-coprimality}}
We have $\gcd(b,a)=\gcd(a,b)$ and $\gcd(a,ma-b)=\gcd(a,b)$.

\prf{Claim~\ref{c:moduli-quadratic}} Let us check preservation of the quadratic form.
The form $mab-a^2-b^2$ is symmetric in $a$ and $b$. We have
\[
 ma(ma-b)-a^2-(ma-b)^2=mab-a^2-b^2.
\]
\end{proof}

\subsection{Connected colored gradings}
For $a\in\Z$, $b\in\Z$, and $n\in\Z$, let
\[
 [a,b]_{\Z}=\{j\in\Z:a\le j\le b\},\qquad [n]_{\Z}=[1,n]_{\Z}.
\]

We recall the notion of NE paths from~\cite{B} and fix the necessary notation.

\begin{definition}\label{def:word}
An \textit{NE path} is a finite lattice path $P$ in
$(\Z\times\R)\cup(\R\times\Z)$ such that $P$ starts from $(0,0)$ and
$P$ consists of unit east steps $E$ and north steps $N$.
For an NE path $P$ and $\mathbf t\in\Z^2$, we call $P+\mathbf t$ a \textit{translated NE path (tNE path)}.
\begin{enumerate}[font=\itshape]
\item\label{d:path-translation} For a tNE path $U$ and $\zeta\in\Z^2$, define
\[
 \operatorname{Trans}(U,\zeta)=U+\zeta.
\]
We call every path with no edges the \textit{empty path} and denote such a path by $\emptyset$.
\item Let $P$ be a tNE path.
\begin{enumerate}[label=(\alph*),font=\itshape]
\item Let $|P|$, $|P|_E$, and $|P|_N$ denote its numbers of edges,
horizontal edges, and vertical edges, respectively.
\item We label its horizontal edges $e_1,\ldots,e_{|P|_E}$ and its vertical
edges $n_1,\ldots,n_{|P|_N}$ in traversal order.
\item For $k\in[|P|]_{\Z}$, let $P_k$ denote its $k$-th edge.
\item For $P\ne\emptyset$, let $\operatorname{Start}_{\mathrm{pt}}(P)$ and
$\operatorname{End}_{\mathrm{pt}}(P)$ denote its starting and ending points, respectively.
We define its starting and ending edges by
\[
 \operatorname{Start}_{\mathrm{NE}}(P)=P_1,\qquad
 \operatorname{End}_{\mathrm{NE}}(P)=P_{|P|}.
\]
\end{enumerate}
\item\label{d:word-concatenation} For nonempty tNE paths $U$ and $V$, their
\textit{concatenation} $U\ccn V$ follows $U$ and then
\[
 \operatorname{Trans}\bigl(V,\operatorname{End}_{\mathrm{pt}}(U)-\operatorname{Start}_{\mathrm{pt}}(V)\bigr).
\]
We define $\emptyset\ccn V=V$ and $U\ccn\emptyset=U$.
For $p\in\Z$, $q\in\Z$, and tNE paths $U_i$ indexed by $i\in[p,q]_{\Z}$, define
\[
 \mathop{\ccn}\limits_{i=p}^{q}U_i=
 \begin{cases}
 U_p\ccn\cdots\ccn U_q,&p\le q,\\
 \emptyset,&p>q.
 \end{cases}
\]
For an integer $\nu\ge1$, we define the $\nu$-fold concatenation by
\[
 U^\nu=\underbrace{U\ccn\cdots\ccn U}_{\nu\text{ copies}}.
\]
We define $U^0=\emptyset$.
\item For an integer $j\ge1$, let $R_j=E\ccn N^j$ be an NE path.
We call a translate of $R_j$ a \textit{simple run} of vertical length $j$.
A \textit{run} is a concatenation of one or more simple runs.
\item For a tNE path $P$, a \textit{subpath} of $P$ is a path consisting of the edges
$P_{i+1},\ldots,P_j$ for some $i\in[0,|P|]_{\Z}$ and $j\in[0,|P|]_{\Z}$ with $i\le j$.
We denote this subpath by $P[i,j]$.
We define $P[i,i]=\emptyset$.
\item For tNE paths $A$ and $B$, let $A\subseteq B$ mean that $A$ is a subpath of $B$.
We call them \textit{disjoint} if they have no common edge.
\end{enumerate}
\end{definition}

\Needspace{5\baselineskip}
\begin{definition}\label{def:word-operations}
Let $U$ and $V$ be tNE paths.
\begin{enumerate}[font=\itshape]
\item For $i\in[0,|U|]_{\Z}$, we define the \textit{insertion}
\[
 \operatorname{Ins}(U,i,V)=U[0,i]\ccn V\ccn U[i,|U|].
\]
\item For a vertical edge $v=U_k$ of $U$, we define the \textit{deletion}
\[
 \operatorname{Del}(U,v)=U[0,k-1]\ccn\operatorname{Trans}\bigl(U[k,|U|],(0,-1)\bigr).
\]
\item Suppose $U\ne\emptyset$. Let $L$ and $R$ be NE paths. Let $A$ and $B$ be simple runs such that
$\operatorname{Start}_{\mathrm{pt}}(A)=\operatorname{Start}_{\mathrm{pt}}(B)=(0,0)$.
When $U=\operatorname{Trans}(L\ccn A\ccn B\ccn R,\operatorname{Start}_{\mathrm{pt}}(U))$, we define the \textit{exchange}
\[
 \operatorname{Exc}(U;|L|,A,B)
 =\operatorname{Trans}(L\ccn B\ccn A\ccn R,\operatorname{Start}_{\mathrm{pt}}(U)).
\]
\item Let $n\in\Z_{>0}$, let $U$ be a run, and let $A_1,\ldots,A_n$ be simple runs. Assume
\[
 U=\operatorname{Trans}(A_1\ccn\cdots\ccn A_n,\operatorname{Start}_{\mathrm{pt}}(U)),\qquad
 \operatorname{Start}_{\mathrm{pt}}(A_i)=(0,0)\quad(i\in[n]_{\Z}).
\]
For $j\in[0,n]_{\Z}$, we define the
\textit{cyclic rotation}
\[
 \operatorname{CRot}(U,j)=\operatorname{Trans}(A_{j+1}\ccn \cdots\ccn A_n\ccn A_1\ccn \cdots\ccn A_j,\operatorname{Start}_{\mathrm{pt}}(U)).
\]
\end{enumerate}
\end{definition}

The grading convention comes from~\cite{B}.

\begin{definition}\label{def:colored-grading}\label{def:balance}
Let $U$ be a tNE path.
\begin{enumerate}[font=\itshape]
\item Let $E(U)$, $N(U)$, and $\operatorname{NE}(U)$ be the sets of horizontal,
vertical, and all edges of $U$, respectively.
\item We refer to a function $g\colon \operatorname{NE}(U)\to\Z_{\ge0}$ as a
\textit{grading} of $U$. For each $v\in \operatorname{NE}(U)$, we call $g(v)$ the
\textit{grade} of $v$.
Let $G(U)$ denote the set of gradings of $U$.
For each integer $h\ge0$, we define
\[
 G_h(U)=\{g\in G(U):g(v)\le h\text{ for every }v\in\operatorname{NE}(U)\}.
\]
For $g\in G(U)$, $h\in G(U)$, and $e\in\operatorname{NE}(U)$, we define
\[
 \operatorname{Diff}_U(g,h;e)=g(e)-h(e).
\]
\item Consider $g\in G(U)$.
\begin{enumerate}[label=(\alph*),font=\itshape]
\item We define the
\textit{boundary balance function} of $g$
\[
 \operatorname{bBal}_{U,g}\colon
 [0,|U|]_{\Z}\to\Z
\]
such that
\[
 \operatorname{bBal}_{U,g}(i)
 =\#E(U[0,i])-\sum_{v\in N(U[0,i])}g(v).
\]
\item For $k\in[|U|]_{\Z}$, we define the \textit{edge balance} at the $k$-th edge $U_k$ of $U$ by
\[
 \operatorname{eBal}_{U,g}(U_k)=\operatorname{bBal}_{U,g}(k).
\]
\end{enumerate}
\item Consider $g\in G(U)$.
We use the convention $\binom{m}{r}=0$ for integers $r>m$.
We call a function $c\colon \operatorname{NE}(U)\to\Z_{>0}$ a \textit{color function} of $g$
if, for each $v\in \operatorname{NE}(U)$, its \textit{color} $c(v)$ satisfies
\[
 c(v)\in\left[\binom{m}{g(v)}\right]_{\Z}.
\]
Let $\operatorname{Col}(U,g,m)$ denote the set of color functions of $g$.
We call $(g,c)$ a \textit{colored grading} of $U$.
\item Let $g\in G(U)$ and $c\in\operatorname{Col}(U,g,m)$. Let $\mathbf t\in\Z^2$.
We define the \textit{translated grading} $g+\mathbf t$ and \textit{translated color function}
$c+\mathbf t$ on $\operatorname{Trans}(U,\mathbf t)$ by
\[
 (g+\mathbf t)(v+\mathbf t)=g(v),\qquad (c+\mathbf t)(v+\mathbf t)=c(v)
 \qquad(v\in\operatorname{NE}(U)).
\]
\item Let $A$ be a subpath of $U$, $g\in G(U)$, and $c\in\operatorname{Col}(U,g,m)$.
The \textit{restricted colored grading} on $A$ is
\[
 (g|_A,c|_A)=\bigl(g|_{\operatorname{NE}(A)},c|_{\operatorname{NE}(A)}\bigr).
\]
Let $g$ and $c$ also denote these restrictions when the ambient path is specified.
\item Let $j\in\Z_{\ge0}$. For $s\in[j]_{\Z}$, let $U_s$ be a tNE path and let $(g_s,c_s)$ be a colored grading of $U_s$.
Let $W=U_1\ccn\cdots\ccn U_j$.
For $s\in[j]_{\Z}$, define $\iota_s:\operatorname{NE}(U_s)\to\operatorname{NE}(W)$ by
\[
 \iota_s((U_s)_h)=W_{|U_1|+\cdots+|U_{s-1}|+h}
 \qquad(h\in[|U_s|]_{\Z}).
\]
The \textit{concatenated colored grading} $(g_1\ccn\cdots\ccn g_j,c_1\ccn\cdots\ccn c_j)$ is the pair $(g,c)$ such that
\[
 g(\iota_s(e))=g_s(e),\qquad c(\iota_s(e))=c_s(e)
 \qquad(s\in[j]_{\Z},\ e\in\operatorname{NE}(U_s)).
\]
\end{enumerate}
\end{definition}

\begin{example}\label{ex:colored-balance}
Let $m=3$. Let $U=E\ccn E\ccn N\ccn E\ccn N=e_1\ccn e_2\ccn n_1\ccn e_3\ccn n_2$.
Consider $g\in G(U)$ and $c\in\operatorname{Col}(U,g,m)$ such that
\[
\begin{array}{c|ccccc}
 x&e_1&e_2&n_1&e_3&n_2\\\hline
 g(x)&2&0&1&1&2\\
 c(x)&3&1&2&1&3
\end{array}
\]
Then, the boundary balances are
\[
\begin{array}{c|rrrrrr}
 i&0&1&2&3&4&5\\\hline
 \operatorname{bBal}_{U,g}(i)&0&1&2&1&2&0
\end{array}
\]

The edge balances are
\[
\begin{array}{c|rrrrr}
 x&e_1&e_2&n_1&e_3&n_2\\\hline
 \operatorname{eBal}_{U,g}(x)&1&2&1&2&0
\end{array}
\]
\end{example}

\begin{definition}\label{def:stack}
Let $U$ be a tNE path.
Consider $g\in G(U)$.
\begin{enumerate}[font=\itshape]
\item\label{d:grading-assignment} We define an \textit{h-precursor}
$P_{U,g}\colon N(U)\to E(U)$ of $g$ by the following conditions
for each $U_k\in N(U)$.
\begin{enumerate}[label=(\alph*),font=\itshape]
\item $P_{U,g}(U_k)\in E(U[0,k-1])$.
\item We have
\[
 \#\{w\in N(U[0,k-1]):P_{U,g}(w)=P_{U,g}(U_k)\}<g(P_{U,g}(U_k)).
\]
\item $P_{U,g}(U_k)$ is the last horizontal edge satisfying (a) and (b).
\end{enumerate}
For each $v\in N(U)$, we say that $P_{U,g}(v)$ is an
\textit{h-precursor} of $v$.
\item Suppose that an h-precursor $P_{U,g}$ exists.
Consider $U_i\in E(U)$ with
\[
 \#P_{U,g}^{-1}(\{U_i\})=g(U_i)>0.
\]
Let $j\in[|U|]_{\Z}$ be maximal such that
$U_j\in P_{U,g}^{-1}(\{U_i\})$.
\begin{enumerate}[label=(\alph*),font=\itshape]
\item We define the \textit{node}
\[
 \operatorname{Node}_{U,g}(U_i)=U[i-1,j].
\]
\item For each $w\in P_{U,g}^{-1}(\{U_i\})$, we call $\operatorname{Node}_{U,g}(U_i)$
the \textit{owner} of $w$.
\end{enumerate}
\item Suppose that an h-precursor $P_{U,g}$ exists.
Let $i\in[0,|U|]_{\Z}$. Let $x\in E(U[0,i])$.
We define the \textit{debt} of $x$ at $i$ by
\[
 \operatorname{Debt}_{U,g}(x,i)=g(x)-\#\{v\in N(U[0,i]):P_{U,g}(v)=x\}.
\]
\end{enumerate}
\end{definition}

\begin{lemma}\label{lem:predecessor-uniqueness}\label{lem:predecessor-existence}
Let $U$ be a tNE path.
Consider $g\in G(U)$.
Then, we have the following.

\begin{enumerate}
\item\label{c:predecessor-uniqueness} There is at most one h-precursor of $g$.
\item\label{c:predecessor-existence} The h-precursor $P_{U,g}$ exists if and only if
\begin{equation}\label{eq:predecessor-existence}
 \sum_{x\in E(U[0,i])}g(x)\ge\#N(U[0,i])
 \qquad(i\in[0,|U|]_{\Z}).
\end{equation}
\end{enumerate}
\end{lemma}
\begin{proof}
\prf{Claim~\ref{c:predecessor-uniqueness}} Let us prove uniqueness by induction on the edge index. Suppose that
$P\colon N(U)\to E(U)$ and $Q\colon N(U)\to E(U)$ are h-precursors of $g$.
We prove by induction on $k\in[0,|U|]_{\Z}$ that $P$ and $Q$
agree on $N(U[0,k])$.
For $k=0$, we have $N(U[0,0])=\emptyset$.

Let us prove the induction step.
Let $k\in[|U|]_{\Z}$. Suppose that $P$ and $Q$ agree on
$N(U[0,k-1])$. If $U_k$ is horizontal, then $N(U[0,k])=N(U[0,k-1])$.
Suppose that $U_k$ is vertical. For each $x\in E(U[0,k-1])$,
the induction hypothesis gives the common value
\begin{align*}
 \nu(x)&=\#\{w\in N(U[0,k-1]):P(w)=x\}\\
     &=\#\{w\in N(U[0,k-1]):Q(w)=x\}.
\end{align*}
Let
\[
 A_k=\{x\in E(U[0,k-1]):\nu(x)<g(x)\}.
\]
Conditions (a) and (b) of Definition~\ref{def:stack}(1) imply
that both $P(U_k)$ and $Q(U_k)$ belong to $A_k$.
Applying Definition~\ref{def:stack}(1)(c) to $P$ with $x=Q(U_k)$ and to $Q$ with
$x=P(U_k)$ shows that neither edge follows the other.
Hence, $P(U_k)=Q(U_k)$. Taking $k=|U|$ gives $P=Q$.

\prf{Claim~\ref{c:predecessor-existence}} Let us obtain the prefix inequalities from an h-precursor. Suppose that $P_{U,g}$ exists. Definition~\ref{def:stack}(1)(a) gives
$P_{U,g}(N(U[0,i]))\subseteq E(U[0,i])$. Definition~\ref{def:stack}(1)(b) gives
$\#P_{U,g}^{-1}(\{x\})\le g(x)$ for each $x\in E(U)$.
For $i\in[0,|U|]_{\Z}$, we therefore have
\[
 \#N(U[0,i])
 =\sum_{x\in E(U[0,i])}\#\{w\in N(U[0,i]):P_{U,g}(w)=x\}
 \le\sum_{x\in E(U[0,i])}g(x).
\]

Let us prove the converse by constructing $P_{U,g}$ in increasing edge order under~\eqref{eq:predecessor-existence}.
Suppose that we have defined its values on $N(U[0,k-1])$ and that
$U_k$ is vertical. For $x\in E(U[0,k-1])$, let
\[
 \nu_k(x)=\#\{w\in N(U[0,k-1]):P_{U,g}(w)=x\}.
\]
The induction hypothesis gives $\nu_k(x)\le g(x)$. Inequality~\eqref{eq:predecessor-existence} for $i=k$ gives
\[
 \sum_{x\in E(U[0,k-1])}\bigl(g(x)-\nu_k(x)\bigr)
 =\sum_{x\in E(U[0,k-1])}g(x)-\#N(U[0,k-1])\ge1.
\]
Thus, there is an $x$ with $\nu_k(x)<g(x)$. Among these edges, let
$P_{U,g}(U_k)$ be the last in edge order.
Induction defines $P_{U,g}(v)$ for every $v\in N(U)$.
\end{proof}

\begin{example}\label{ex:horizontal-predecessor}
Let $U=E\ccn E\ccn N\ccn N\ccn N=e_1\ccn e_2\ccn n_1\ccn n_2\ccn n_3$.
Consider $g\in G(U)$ such that
\[
 (g(e_1),g(e_2),g(n_1),g(n_2),g(n_3))=(2,1,0,0,0).
\]
Then, the h-precursor $P_{U,g}$ satisfies
\[
 P_{U,g}(n_1)=e_2,\qquad P_{U,g}(n_2)=P_{U,g}(n_3)=e_1.
\]
\end{example}

\begin{definition}\label{def:connected-grading}
Let $U$ be a tNE path.
\begin{enumerate}[font=\itshape]
\item Consider $g\in G(U)$.
We call $g$ \textit{connected} if the following conditions hold.
\begin{enumerate}[font=\itshape]
\item The h-precursor $P_{U,g}$ exists.
\item\label{d:connected-baseline} For every node $\operatorname{Node}_{U,g}(x)$
and every vertical edge $v$ owned by $\operatorname{Node}_{U,g}(x)$, we have
\begin{equation}\label{eq:connected-edge-balance}
 \operatorname{eBal}_{U,g}(v)\ge\operatorname{eBal}_{U,g}(x).
\end{equation}
We call this condition the \textit{node balance condition}.
\item $U=\operatorname{Node}_{U,g}(x)$ for some $x\in E(U)$.
\end{enumerate}
Let $\operatorname{CtdG}(U)$ denote the set of connected gradings of $U$.
For each integer $h\ge0$, let
\[
 \operatorname{CtdG}_h(U)=\operatorname{CtdG}(U)\cap G_h(U).
\]
\item For $g\in\operatorname{CtdG}(U)$, we define the \textit{surplus} $\operatorname{Surp}(U,g)$ by
\begin{equation}\label{eq:residual}
 \operatorname{Surp}(U,g)=\#E(U)-1-\sum_{v\in N(U)}g(v).
\end{equation}
\item Consider $g\in G(U)$ and $c\in\operatorname{Col}(U,g,m)$.
We call $(g,c)$ \textit{connected} if $g$ is connected.
For each integer $d\ge0$, we define
\[
 \mathcal C_d(U)=\left\{(g,c)\;\middle|\;
 \begin{array}{@{}l@{}}
 g\in\operatorname{CtdG}(U),\quad c\in\operatorname{Col}(U,g,m),\\
 \operatorname{Surp}(U,g)=d
 \end{array}\right\}.
\]
We define
\[
 C_d(U)=\#\mathcal C_d(U).
\]
\end{enumerate}
\end{definition}

\begin{lemma}\label{lem:word-reindexing}
Let $U$ be a tNE path. Let $\mathbf t\in\Z^2$ and $W=\operatorname{Trans}(U,\mathbf t)$.
For every integer $d\ge0$, we have
\[
 C_d(U)=C_d(W).
\]
\end{lemma}
\begin{proof}
Let us define the map
\[
 \mathcal C_d(U)\longrightarrow\mathcal C_d(W),\qquad
 (g,c)\longmapsto(g+\mathbf t,c+\mathbf t).
\]
Translation preserves edge order and each grade. Thus, for $v\in N(U)$,
\[
 P_{W,g+\mathbf t}(v+\mathbf t)=P_{U,g}(v)+\mathbf t.
\]
Each node is translated by $\mathbf t$, and its edge balances are unchanged.
The translated grading is therefore connected and has surplus $d$.
The color bounds are unchanged because the grades are unchanged.
Translation by $-\mathbf t$ defines the inverse map. Hence, the displayed map is a bijection.
\end{proof}

\begin{lemma}\label{lem:residual-nonnegative}
Let $U$ be a tNE path. Consider $g\in\operatorname{CtdG}(U)$.
Then, \mbox{$\operatorname{Surp}(U,g)\ge0$}.
\end{lemma}
\begin{proof}
Let us apply~\eqref{eq:connected-edge-balance} at $\operatorname{End}_{\mathrm{NE}}(U)$.
Let $x=\operatorname{Start}_{\mathrm{NE}}(U)$ and $v=\operatorname{End}_{\mathrm{NE}}(U)$.
By Definition~\ref{def:connected-grading}(1)(c),
$U=\operatorname{Node}_{U,g}(x)$.
Since $x=\operatorname{Start}_{\mathrm{NE}}(U)\in E(U)$, we have $\operatorname{eBal}_{U,g}(x)=1$.
Inequality~\eqref{eq:connected-edge-balance} gives
\[
\begin{aligned}
 \operatorname{Surp}(U,g)&=\#E(U)-1-\sum_{w\in N(U)}g(w)\\
  &=\operatorname{eBal}_{U,g}(v)-\operatorname{eBal}_{U,g}(x)\ge0.
\end{aligned}
\]
\end{proof}

\begin{lemma}\label{lem:colored-weight}
Let $U$ be a tNE path and $d\in\Z_{\ge0}$. Then,
\begin{equation}\label{eq:colored-weight}
\begin{aligned}
 C_d(U)&=\sum_{\substack{g\in\operatorname{CtdG}_m(U)\\\operatorname{Surp}(U,g)=d}}
           \#\operatorname{Col}(U,g,m)\\
 &=\sum_{\substack{g\in\operatorname{CtdG}_m(U)\\\operatorname{Surp}(U,g)=d}}
           \prod_{v\in \operatorname{NE}(U)}\binom{m}{g(v)}.
\end{aligned}
\end{equation}
\end{lemma}
\begin{proof}
Let us count the color functions for a fixed grading $g\in G(U)$.
By Definition~\ref{def:colored-grading}, choosing a color function is equivalent to
choosing one element of $[\binom m{g(v)}]_{\Z}$ for each edge $v$. Hence,
\[
 \operatorname{Col}(U,g,m)
 \simeq\prod_{v\in\operatorname{NE}(U)}\left[\binom m{g(v)}\right]_{\Z},\qquad
 \#\operatorname{Col}(U,g,m)=\prod_{v\in\operatorname{NE}(U)}\binom m{g(v)}.
\]
If $g(v)>m$ for some edge $v$, the corresponding factor is empty.
Thus, only gradings in $G_m(U)$ contribute to $C_d(U)$.
Partitioning $\mathcal C_d(U)$ according to its grading gives~\eqref{eq:colored-weight}.
\end{proof}

\begin{lemma}[Prefix counts for h-precursors]\label{lem:precursor-prefix-counts}
Let $U=U_1\ccn\cdots\ccn U_n$ be a tNE path. Let $g\in G(U)$. Suppose that
$P=P_{U,g}$ exists. For $e\in E(U)$ and $k\in[0,n]_{\Z}$, define
\[
 F_e=\{j\in[n]_{\Z}:U_j\in N(U),\ P(U_j)=e\},\qquad
 \operatorname{Num}(e,k)=\#(F_e\cap[k]_{\Z}).
\]
Then, we have the following.
\begin{enumerate}
\item\label{c:prefix-capacity} For $e\in E(U)$ and integers $0\le k\le\ell\le n$, we have
\[
 0=\operatorname{Num}(e,0)\le \operatorname{Num}(e,k)\le \operatorname{Num}(e,\ell)\le g(e).
\]
If $\operatorname{Num}(e,k)=g(e)$, then $F_e\subseteq[k]_{\Z}$ and
$\operatorname{Num}(e,\ell)=g(e)$ for $k\le\ell\le n$.
\item\label{c:assigned-prefix} Suppose that $U_s\in N(U)$ and $P(U_s)=U_r$.
Then, $r<s$. For $k\in[r,s-1]_{\Z}$, we have
\[
 \operatorname{Num}(U_r,k)\le \operatorname{Num}(U_r,s-1)\le g(U_r)-1,\qquad
 \operatorname{Debt}_{U,g}(U_r,k)\ge1.
\]
For each $U_k\in N(U)$ with $r<k\le s$, let $p_k$ be the unique index
satisfying $P(U_k)=U_{p_k}$. Then, $r\le p_k<k$.
\item\label{c:interior-prefix} Under the hypotheses of Claim~\ref{c:assigned-prefix},
let $r<i<s$ and $U_i\in E(U)$. Then,
\[
 \operatorname{Num}(U_i,s-1)=\operatorname{Num}(U_i,n)=g(U_i),\qquad
 F_{U_i}\subseteq[i+1,s-1]_{\Z}.
\]
If $g(U_i)>0$, the index $q_i=\max F_{U_i}$ satisfies
\[
 i<q_i<s,\qquad \operatorname{Node}_{U,g}(U_i)=U[i-1,q_i].
\]
\end{enumerate}
\end{lemma}
\begin{proof}
\prf{Claim~\ref{c:prefix-capacity}} Let us prove $\#F_e\le g(e)$.
The inclusions $F_e\cap[k]_{\Z}\subseteq F_e\cap[\ell]_{\Z}\subseteq F_e$
give monotonicity. If $F_e=\emptyset$, every count is zero.
Otherwise, let $q=\max F_e$. Definition~\ref{def:stack}(1)(b) gives
\[
 \#F_e-1=\operatorname{Num}(e,q-1)<g(e).
\]
Integrality gives $\#F_e\le g(e)$. If $\operatorname{Num}(e,k)=g(e)$, then
\[
 g(e)=\#(F_e\cap[k]_{\Z})\le\#F_e\le g(e).
\]
The two finite sets $F_e\cap[k]_{\Z}\subseteq F_e$ have equal cardinalities. Thus, they are equal, and $\operatorname{Num}(e,\ell)=g(e)$ for $k\le\ell\le n$.

\prf{Claim~\ref{c:assigned-prefix}} Let us bound the debt before an assigned vertical edge.
Conditions (a) and (b) of Definition~\ref{def:stack}(1) give
$r<s$ and $\operatorname{Num}(U_r,s-1)<g(U_r)$.
For $k\in[r,s-1]_{\Z}$, Claim~\ref{c:prefix-capacity} gives
\[
 \operatorname{Num}(U_r,k)\le \operatorname{Num}(U_r,s-1)\le g(U_r)-1.
\]
Hence, $\operatorname{Debt}_{U,g}(U_r,k)=g(U_r)-\operatorname{Num}(U_r,k)\ge1$.
For $r<k\le s$, the edge $U_r$ precedes $U_k$ and
$\operatorname{Num}(U_r,k-1)<g(U_r)$. Definition~\ref{def:stack}(1)(c) gives $r\le p_k$;
condition (a) gives $p_k<k$.

\prf{Claim~\ref{c:interior-prefix}} Let us prove $\operatorname{Num}(U_i,s-1)=g(U_i)$.
Claim~\ref{c:prefix-capacity} gives $\operatorname{Num}(U_i,s-1)\le g(U_i)$.
If $\operatorname{Num}(U_i,s-1)<g(U_i)$, Definition~\ref{def:stack}(1)(c) at $U_s$
would give $i\le r$, contrary to $r<i$.
Thus, $\operatorname{Num}(U_i,s-1)=g(U_i)$.
Claim~\ref{c:prefix-capacity} gives $F_{U_i}\subseteq[s-1]_{\Z}$ and $\operatorname{Num}(U_i,n)=g(U_i)$.
For each $t\in F_{U_i}$, we have $P(U_t)=U_i$. Definition~\ref{def:stack}(1)(a) gives $i<t$.
Hence, $F_{U_i}\subseteq[i+1,s-1]_{\Z}$.
If $g(U_i)>0$, then $F_{U_i}\ne\emptyset$ and its maximum $q_i$ satisfies
$i<q_i<s$. The definition of a node in Definition~\ref{def:stack}
gives $\operatorname{Node}_{U,g}(U_i)=U[i-1,q_i]$.
\end{proof}

\begin{lemma}[Restriction to an interval closed under h-precursors]\label{lem:precursor-closed-restriction}
Let $U=U_1\ccn\cdots\ccn U_n$ be a tNE path. Let $g\in G(U)$. Suppose that
$P=P_{U,g}$ exists. Let $a,b\in[0,n]_{\Z}$ satisfy $a<b$, and let $\ell=b-a$. Put
\[
 X=U[a,b],\qquad X_j=U_{a+j},\qquad g_X(X_j)=g(U_{a+j})
 \quad(j\in[\ell]_{\Z}).
\]
Suppose that $P(N(X))\subseteq E(X)$.
For $X_j\in N(X)$, let $p_j$ be the unique index satisfying
$P(U_{a+j})=U_{p_j}$.
Use $\operatorname{Num}(e,k)=\#(F_e\cap[k]_{\Z})$ from Lemma~\ref{lem:precursor-prefix-counts}.
Then, we have the following.
\begin{enumerate}
\item\label{c:closed-prefix-count} For $X_j\in N(X)$, we have $a<p_j<a+j$.
Define
\[
 Q\colon N(X)\longrightarrow E(X),\qquad Q(X_j)=X_{p_j-a}.
\]
For $X_i\in E(X)$ and $j\in[0,\ell]_{\Z}$,
\[
 \#\{k\in[j]_{\Z}:X_k\in N(X),\ Q(X_k)=X_i\}=\operatorname{Num}(U_{a+i},a+j).
\]
The map $Q$ is the h-precursor $P_{X,g_X}$.
\item\label{c:closed-node-correspondence} Let $X_i\in E(X)$ and $i<j\le\ell$.
Then, \begin{equation}\label{eq:closed-node-correspondence}
 \operatorname{Node}_{X,g_X}(X_i)=X[i-1,j]
 \quad\Longleftrightarrow\quad
 \operatorname{Node}_{U,g}(U_{a+i})=U[a+i-1,a+j].
\end{equation}
Each side of~\eqref{eq:closed-node-correspondence} asserts that its node exists and equals the specified subpath.
\item\label{c:closed-balance} For $j\in[0,\ell]_{\Z}$, we have
\[
 \operatorname{bBal}_{X,g_X}(j)
 =\operatorname{bBal}_{U,g}(a+j)-\operatorname{bBal}_{U,g}(a).
\]
For $j\in[\ell]_{\Z}$, we have
\[
 \operatorname{eBal}_{X,g_X}(X_j)
 =\operatorname{eBal}_{U,g}(U_{a+j})-\operatorname{bBal}_{U,g}(a).
\]
If $g$ satisfies the node balance condition, then $g_X$ satisfies that condition.
\end{enumerate}
\end{lemma}
\begin{proof}
\prf{Claim~\ref{c:closed-prefix-count}} Let us verify $a<p_j<a+j$ and the formula for $Q$.
For $X_j\in N(X)$, the hypothesis $P(N(X))\subseteq E(X)$ gives $a<p_j\le b$.
Definition~\ref{def:stack}(1)(a) gives $p_j<a+j$.
Hence, $X_{p_j-a}\in E(X)$ and $1\le p_j-a<j$.

Let us prove the equality of cardinalities in Claim~\ref{c:closed-prefix-count}.
For $X_i\in E(X)$ and $j\in[0,\ell]_{\Z}$, the map
\[
 \{k\in[j]_{\Z}:X_k\in N(X),\ Q(X_k)=X_i\}
 \longrightarrow F_{U_{a+i}}\cap[a+j]_{\Z},\qquad k\longmapsto a+k
\]
is a bijection. Its inverse is $t\mapsto t-a$: every $t\in F_{U_{a+i}}$ satisfies $t>a+i>a$ by Definition~\ref{def:stack}(1)(a).

Let us verify Definition~\ref{def:stack}(1)(b) and (c) for $Q$.
For $X_j\in N(X)$, apply the cardinality identity proved using $k\mapsto a+k$ with $i=p_j-a$ at boundary $j-1$. This gives
\[
 \#\{k\in[j-1]_{\Z}:X_k\in N(X),\ Q(X_k)=X_{p_j-a}\}
 =\operatorname{Num}(U_{p_j},a+j-1)<g(U_{p_j})=g_X(X_{p_j-a}),
\]
where the strict inequality follows from Definition~\ref{def:stack}(1)(b).
For $X_i\in E(X)$ with $i<j$, suppose
\[
 \#\{k\in[j-1]_{\Z}:X_k\in N(X),\ Q(X_k)=X_i\}<g_X(X_i).
\]
Then, $\operatorname{Num}(U_{a+i},a+j-1)<g(U_{a+i})$.
Definition~\ref{def:stack}(1)(c) gives $a+i\le p_j$. Thus, $i\le p_j-a$.
Claim~\ref{c:predecessor-uniqueness} of Lemma~\ref{lem:predecessor-existence}
gives $Q=P_{X,g_X}$.

\prf{Claim~\ref{c:closed-node-correspondence}} First, suppose $\operatorname{Node}_{X,g_X}(X_i)=X[i-1,j]$.
Then, $\#Q^{-1}(\{X_i\})=g_X(X_i)=g(U_{a+i})>0$.
Claim~\ref{c:closed-prefix-count} at boundary $\ell$ gives
$\operatorname{Num}(U_{a+i},b)=g(U_{a+i})$.
Claim~\ref{c:prefix-capacity} of Lemma~\ref{lem:precursor-prefix-counts} gives $F_{U_{a+i}}\subseteq[b]_{\Z}$.
For every $t\in F_{U_{a+i}}$, Definition~\ref{def:stack}(1)(a) gives $t>a+i>a$.
Hence,
\[
 F_{U_{a+i}}
 =\{a+k:k\in[\ell]_{\Z},\ X_k\in N(X),\ Q(X_k)=X_i\}.
\]
The maximum of $\{k\in[\ell]_{\Z}:X_k\in N(X),\ Q(X_k)=X_i\}$ is $j$.
Thus, $\max F_{U_{a+i}}=a+j$ and $\operatorname{Node}_{U,g}(U_{a+i})=U[a+i-1,a+j]$.

Second, suppose $\operatorname{Node}_{U,g}(U_{a+i})=U[a+i-1,a+j]$.
Then, $F_{U_{a+i}}\subseteq[a+i+1,a+j]_{\Z}$, $\#F_{U_{a+i}}=g(U_{a+i})>0$, and $\max F_{U_{a+i}}=a+j$.
Claim~\ref{c:closed-prefix-count} gives
\[
 \{k\in[\ell]_{\Z}:X_k\in N(X),\ Q(X_k)=X_i\}
 =\{t-a:t\in F_{U_{a+i}}\}.
\]
This set has cardinality $g_X(X_i)$ and maximum $j$.
Definition~\ref{def:stack}(2) gives $\operatorname{Node}_{X,g_X}(X_i)=X[i-1,j]$.

\prf{Claim~\ref{c:closed-balance}} Let us prove the two balance identities in Claim~\ref{c:closed-balance}.
For $k\in[n]_{\Z}$, Definition~\ref{def:colored-grading}(3)(a) gives
\[
 c_k:=\operatorname{bBal}_{U,g}(k)-\operatorname{bBal}_{U,g}(k-1)
 =\begin{cases}1,&U_k\in E(U),\\-g(U_k),&U_k\in N(U).\end{cases}
\]
Since $X_k=U_{a+k}$ and $g_X(X_k)=g(U_{a+k})$ for $k\in[\ell]_{\Z}$, for $j\in[0,\ell]_{\Z}$ we obtain
\[
 \operatorname{bBal}_{X,g_X}(j)=\sum_{k=a+1}^{a+j}c_k
 =\operatorname{bBal}_{U,g}(a+j)-\operatorname{bBal}_{U,g}(a).
\]
For $j=0$, the sum is empty and both sides are zero.
For $j\in[\ell]_{\Z}$, Definition~\ref{def:colored-grading}(3)(b) gives
\[
 \operatorname{eBal}_{X,g_X}(X_j)=\operatorname{bBal}_{X,g_X}(j),\qquad
 \operatorname{eBal}_{U,g}(U_{a+j})=\operatorname{bBal}_{U,g}(a+j).
\]
Substitution proves the edge-balance identity in Claim~\ref{c:closed-balance}.
Suppose that $g$ satisfies the node balance condition.
Let $Y=\operatorname{Node}_{X,g_X}(X_i)=X[i-1,j]$ and let $X_k$ be owned by $Y$.
Claims~\ref{c:closed-prefix-count} and~\ref{c:closed-node-correspondence} give
\[
 \operatorname{Node}_{U,g}(U_{a+i})=U[a+i-1,a+j],\qquad P(U_{a+k})=U_{a+i}.
\]
Inequality~\eqref{eq:connected-edge-balance} for $g$ gives
\[
 \operatorname{eBal}_{X,g_X}(X_k)-\operatorname{eBal}_{X,g_X}(X_i)
 =\operatorname{eBal}_{U,g}(U_{a+k})-\operatorname{eBal}_{U,g}(U_{a+i})\ge0.
\]
This is the node balance condition for $g_X$.
\end{proof}

\begin{lemma}[Structure of nodes]\label{lem:node-structure}
Let $U$ be a tNE path. Consider $g\in G(U)$. Suppose that the h-precursor $P_{U,g}$ exists.
Let $x=U_r\in E(U)$. Consider a node $X=\operatorname{Node}_{U,g}(x)$.
Let $v=\operatorname{End}_{\mathrm{NE}}(X)=U_s$.
Then, we have the following.

\begin{enumerate}
\item\label{c:node-internal} The following statements hold.
\begin{enumerate}[label=(\alph*),ref=\theenumi(\alph*)]
\item\label{c:node-predecessors} For every $w\in N(X)$, we have
$P_{U,g}(w)\in E(X)$.
\item\label{c:node-contained} For each $y\in E(X)$ with $g(y)>0$,
the node $\operatorname{Node}_{U,g}(y)$ exists and
$\operatorname{Node}_{U,g}(y)\subseteq X$.
\item\label{c:node-ends} For every $y\in E(X)\setminus\{x\}$ with $g(y)>0$,
$\operatorname{End}_{\mathrm{NE}}(\operatorname{Node}_{U,g}(y))$ precedes $v$.
\end{enumerate}
\item\label{c:node-nesting} For each pair of distinct edges $y\in E(U)$ and $z\in E(U)$
such that $\operatorname{Node}_{U,g}(y)$ and $\operatorname{Node}_{U,g}(z)$ exist,
these nodes are disjoint or
\[
 \operatorname{Node}_{U,g}(y)\subsetneq\operatorname{Node}_{U,g}(z)
 \quad\text{or}\quad
 \operatorname{Node}_{U,g}(z)\subsetneq\operatorname{Node}_{U,g}(y).
\]
\item\label{c:node-balances} Suppose that $g$ satisfies the node balance condition. The following statements hold.

\begin{enumerate}[label=(\alph*),ref=\theenumi(\alph*)]
\item\label{c:node-boundary-balances} For every $i\in[r,s]_{\Z}$, we have
\[
 \operatorname{bBal}_{U,g}(i)\ge\operatorname{eBal}_{U,g}(x).
\]
\item\label{c:node-descendant-balances} Let $y\in E(U)$ be such that
$\operatorname{Node}_{U,g}(y)\subsetneq X$.
Then, $\operatorname{eBal}_{U,g}(y)>\operatorname{eBal}_{U,g}(x)$.
\end{enumerate}
\item\label{c:node-residual} Suppose that $g$ satisfies the node balance condition. The following statements hold.
\begin{enumerate}[label=(\alph*),ref=\theenumi(\alph*)]
\item\label{c:node-restricted-connected} The grading $g|_{\operatorname{NE}(X)}$ is connected.
\item\label{c:node-residual-nonnegative} We have
$\operatorname{Surp}(X,g|_{\operatorname{NE}(X)})\ge0$.
\item\label{c:node-residual-identity} We have
\begin{equation}\label{eq:node-surplus-balance}
 \operatorname{eBal}_{U,g}(v)
 -\operatorname{bBal}_{U,g}(r-1)
 =\operatorname{Surp}(X,g|_{\operatorname{NE}(X)})+1.
\end{equation}
\end{enumerate}
\end{enumerate}
\end{lemma}
\begin{proof}
\prf{Claim~\ref{c:node-predecessors}} Let us locate the h-precursors inside the node.
Definition~\ref{def:stack}(2) gives $P_{U,g}(U_s)=U_r$.
For $U_k\in N(X)$, Claim~\ref{c:assigned-prefix} of
Lemma~\ref{lem:precursor-prefix-counts} gives $r\le p_k<k\le s$,
where $P_{U,g}(U_k)=U_{p_k}$. Thus, $P_{U,g}(U_k)\in E(X)$.

\prf{Claim~\ref{c:node-contained}} Let us prove containment of the smaller nodes.
If $y=x$, then $\operatorname{Node}_{U,g}(y)=X$.
For $y=U_i\in E(X)\setminus\{x\}$ with $g(y)>0$, we have $r<i<s$.
Claim~\ref{c:interior-prefix} of Lemma~\ref{lem:precursor-prefix-counts}
gives $q_i=\max F_{U_i}<s$ and
$\operatorname{Node}_{U,g}(y)=U[i-1,q_i]\subseteq X$.

\prf{Claim~\ref{c:node-ends}} Let us compare the node endpoints.
For $y=U_i\in E(X)\setminus\{x\}$ with $g(y)>0$,
Claim~\ref{c:interior-prefix} of Lemma~\ref{lem:precursor-prefix-counts} gives $q_i<s$.
Thus, the endpoint $U_{q_i}$ precedes $U_s$.

\prf{Claim~\ref{c:node-nesting}} Let us compare $\operatorname{Node}_{U,g}(y)$ and $\operatorname{Node}_{U,g}(z)$.
Interchange $y$ and $z$ if necessary so that $y$ precedes $z$.
Let $t$ be the endpoint of $\operatorname{Node}_{U,g}(y)$.
If $t$ precedes $z$, then $\operatorname{Node}_{U,g}(y)$ and $\operatorname{Node}_{U,g}(z)$ are disjoint.
Otherwise, $z$ lies in that node and is not its start.
Claims~\ref{c:node-contained} and~\ref{c:node-ends} give
$\operatorname{Node}_{U,g}(z)\subsetneq\operatorname{Node}_{U,g}(y)$.

\prf{Claim~\ref{c:node-boundary-balances}} Let us bound every boundary balance by the initial node balance.
Let $B_j=\operatorname{bBal}_{U,g}(j)$.
We prove $B_i\ge B_r$ by strong induction on $i\in[r,s]_{\Z}$.
The case $i=r$ is equality.
For $i>r$, assume $B_j\ge B_r$ for every $r\le j<i$.
If $U_i$ is horizontal, then $B_i=B_{i-1}+1\ge B_r+1$.
If $U_i$ is vertical, Claim~\ref{c:node-predecessors} gives
$P_{U,g}(U_i)=U_{p_i}$ for $r\le p_i<i$.
Definition~\ref{def:stack}(1)(b) implies $g(U_{p_i})>0$. Thus, Claim~\ref{c:node-contained} supplies its node.
The node balance condition and the induction hypothesis at $p_i$ give
\[
 B_i=\operatorname{eBal}_{U,g}(U_i)
 \ge\operatorname{eBal}_{U,g}(U_{p_i})=B_{p_i}\ge B_r.
\]
Since $B_r=\operatorname{eBal}_{U,g}(x)$, this proves Claim~\ref{c:node-boundary-balances}.

\prf{Claim~\ref{c:node-descendant-balances}} Let us prove the strict balance inequality for descendants.
Write $y=U_i$.
Proper containment gives $r<i\le s$.
Claim~\ref{c:node-boundary-balances} at $i-1$ gives
\[
 \operatorname{eBal}_{U,g}(y)=\operatorname{bBal}_{U,g}(i-1)+1
 \ge\operatorname{eBal}_{U,g}(x)+1>\operatorname{eBal}_{U,g}(x).
\]

\prf{Claim~\ref{c:node-restricted-connected}} Let us prove connectedness of the restricted grading.
Let $g_X=g|_{\operatorname{NE}(X)}$.
Claim~\ref{c:node-predecessors} gives $P_{U,g}(N(X))\subseteq E(X)$.
Apply Lemma~\ref{lem:precursor-closed-restriction} with $a=r-1$ and $b=s$.
Claim~\ref{c:closed-prefix-count} of Lemma~\ref{lem:precursor-closed-restriction} gives the h-precursor $P_{X,g_X}$.
Claim~\ref{c:closed-node-correspondence} of Lemma~\ref{lem:precursor-closed-restriction} gives $\operatorname{Node}_{X,g_X}(x)=X$.
Claim~\ref{c:closed-balance} of Lemma~\ref{lem:precursor-closed-restriction} gives the node balance condition for $g_X$.
Thus, $g_X$ is connected by Definition~\ref{def:connected-grading}(1).

\prf{Claim~\ref{c:node-residual-nonnegative}}
Claim~\ref{c:node-restricted-connected}
and Lemma~\ref{lem:residual-nonnegative} give
$\operatorname{Surp}(X,g|_{\operatorname{NE}(X)})\ge0$.

\prf{Claim~\ref{c:node-residual-identity}} Let us compute the surplus from the endpoint balances.
Claim~\ref{c:closed-balance} of
Lemma~\ref{lem:precursor-closed-restriction} and~\eqref{eq:residual} give
\[
 \operatorname{eBal}_{U,g}(v)-\operatorname{bBal}_{U,g}(r-1)
 =\#E(X)-\sum_{w\in N(X)}g(w)
 =\operatorname{Surp}(X,g|_{\operatorname{NE}(X)})+1.
\]
This proves~\eqref{eq:node-surplus-balance}.
\end{proof}

\begin{definition}\label{def:node-relations}
Let $U$ be a tNE path and $g\in G(U)$. Let $X$ and $Y$ be nodes of $(U,g)$.
\begin{enumerate}[font=\itshape]
\item We call $X$ the \textit{parent} of $Y$ and $Y$ a \textit{child} of $X$
if $Y\subsetneq X$
and there is no node $Z$ such that $Y\subsetneq Z\subsetneq X$.
\item We call $X$ an \textit{ancestor} of $Y$ and $Y$ a \textit{descendant} of $X$
if $Y\subsetneq X$.
\end{enumerate}
\end{definition}
\begin{lemma}\label{lem:parent-unique}
Let $U$ be a tNE path and $g\in G(U)$.
Each node of $(U,g)$ has at most one parent.
\end{lemma}
\begin{proof}
Let us prove that two parents of the same node are equal.
Suppose $X$ and $Z$ are parents of a node $Y$.
Both contain $Y$, so they are not disjoint.
Claim~\ref{c:node-nesting} of Lemma~\ref{lem:node-structure} gives $X\subseteq Z$ or $Z\subseteq X$.
If $X\ne Z$, interchanging their names if necessary gives
\[
 Y\subsetneq X\subsetneq Z.
\]
This contradicts the condition that $Z$ is a parent of $Y$ in Definition~\ref{def:node-relations}(1).
Hence, $X=Z$.
\end{proof}

\subsection{Compatibility and counting colored gradings}

\begin{definition}\label{def:cyclic-subpaths}
Let $U$ be a nonempty tNE path. Let $u$ and $v$ be edges of $U$.
\begin{enumerate}[font=\itshape]
\item Let $U\langle u,v\rangle$ denote the \textit{cyclic subpath} of $U$
read from $u$ in traversal order, returning from $\operatorname{End}_{\mathrm{NE}}(U)$
to $\operatorname{Start}_{\mathrm{NE}}(U)$,
and ending at the first occurrence of $v$. We include both end edges.
\item For $P=U\langle u,v\rangle$, let
$\operatorname{NE}(P)=\operatorname{NE}(U,u,v)$, $E(P)=E(U,u,v)$, and $N(P)=N(U,u,v)$
denote the sets of all, horizontal, and vertical edges of $U$ traversed by $P$, respectively.
\item We define the \textit{NE-path realization}
$\operatorname{NEReal}(U,u,v)$ of $U\langle u,v\rangle$ to be the NE path from $(0,0)$
with the same ordered sequence of east and north steps as $U\langle u,v\rangle$.
Each edge of $\operatorname{NEReal}(U,u,v)$ keeps its name and color, if any.
Let $X=\operatorname{NEReal}(U,u,v)$.
For $e\in\operatorname{NE}(U,u,v)$, let $\widetilde e\in\operatorname{NE}(X)$ be its realized edge.
For $g\in G(U)$, we define the \textit{transported grading} $g_{U,u,v}$ of $X$ by
\[
 g_{U,u,v}(\widetilde e)=g(e)\qquad(e\in\operatorname{NE}(U,u,v)).
\]
\item Let $x$ and $y$ be edges of $U$. We call $U\langle u,v\rangle$ a \textit{subpath} of $U\langle x,y\rangle$ and let
\[
 U\langle u,v\rangle\subseteq U\langle x,y\rangle
\]
if the edge sequence of $U\langle u,v\rangle$ occurs consecutively in the edge sequence of $U\langle x,y\rangle$.
\item Let $x$ and $y$ be edges of $U$. We call $U\langle u,v\rangle$ and $U\langle x,y\rangle$
\textit{disjoint} if
\[
 \operatorname{NE}(U,u,v)\cap\operatorname{NE}(U,x,y)=\varnothing.
\]
\end{enumerate}
\end{definition}

We define the notions of compatibility and shadows~\cite{B}
for tNE paths and fix some notation for later use.

\begin{definition}\label{def:compatible-grading}
Let $U$ be a tNE path and $g\in G(U)$.
\begin{enumerate}[font=\itshape]
\item\label{d:grading-differences} For edges $u$ and $t$ of $U$, we define
\[
\begin{aligned}
 \operatorname{NGrDiff}_{U,g}(u,t)
 &=\#N(U,u,t)-\sum_{x\in E(U,u,t)}g(x),\\
 \operatorname{EGrDiff}_{U,g}(u,t)
 &=\#E(U,u,t)-\sum_{x\in N(U,u,t)}g(x).
\end{aligned}
\]
We call $\operatorname{NGrDiff}_{U,g}(u,t)=0$ the \textit{NGrDiff-zero condition}.
We call $\operatorname{EGrDiff}_{U,g}(u,t)=0$ the \textit{EGrDiff-zero condition}.
\item\label{d:compatible} We call $g$ \textit{compatible} if for
every horizontal edge $u$ and vertical edge $v$ of $U$ with $g(u)g(v)>0$,
there is an edge $t\in \operatorname{NE}(U,u,v)$ satisfying one of the following conditions.
\begin{enumerate}[label=(\alph*),font=\itshape]
\item\label{d:compatibility-horizontal} $t$ is vertical, $t\ne v$,
and $\operatorname{NGrDiff}_{U,g}(u,t)=0$.
\item\label{d:compatibility-vertical} $t$ is horizontal, $t\ne u$,
and $\operatorname{EGrDiff}_{U,g}(t,v)=0$.
\end{enumerate}
We call such an edge $t$ a \textit{witness} for $(u,v)$.
\item\label{d:horizontal-shadow} Let $u$ be a horizontal edge of $U$.
If there is a vertical edge $t$ of $U$ such that
\[
 \operatorname{NGrDiff}_{U,g}(u,t)=0,
\]
let $t_{U,g}(u)$ be the first such edge in cyclic order starting at $u$.
Define the \textit{horizontal shadow} by
\[
 \operatorname{Sh}_{U,g}(u)=N(U,u,t_{U,g}(u)).
\]
If no vertical edge satisfies the NGrDiff-zero condition,
let $\operatorname{Sh}_{U,g}(u)=N(U)$.
\item\label{d:vertical-shadow} Let $v$ be a vertical edge of $U$.
If there is a horizontal edge $t$ of $U$ such that
\[
 \operatorname{EGrDiff}_{U,g}(t,v)=0,
\]
let $s_{U,g}(v)$ be the first such edge in reverse cyclic order starting at $v$.
We define the \textit{vertical shadow} $\operatorname{Sh}_{U,g}(v)$ of $v$ by
\[
 \operatorname{Sh}_{U,g}(v)=E(U,s_{U,g}(v),v).
\]
If no horizontal edge satisfies the EGrDiff-zero condition,
let $\operatorname{Sh}_{U,g}(v)=E(U)$.
For $C\subseteq\operatorname{NE}(U)$, we define
\[
 \operatorname{Sh}_{U,g}(C)=\bigcup_{e\in C}\operatorname{Sh}_{U,g}(e).
\]
\item\label{d:total-grades} We define the \textit{total vertical grade} and
\textit{total horizontal grade} of $g$ by
\[
 \operatorname{TNG}_{U,g}=\sum_{v\in N(U)}g(v),\qquad
 \operatorname{TEG}_{U,g}=\sum_{u\in E(U)}g(u),
\]
respectively.
\item\label{d:positive-edge-sets} We define the sets of positively graded horizontal and vertical edges by
\[
 E_g(U)=\{u\in E(U):g(u)>0\},\qquad
 N_g(U)=\{v\in N(U):g(v)>0\}.
\]
\item Let $\mathcal P_{U,g}$ be the set of
cyclic subpaths $U\langle u,t_{U,g}(u)\rangle$ for all $u\in E_g(U)$
with $t_{U,g}(u)$ existing.
Let $P\in\mathcal P_{U,g}$.
We call $P$ a \textit{maximal h-shadow subpath} of $g$ if there is no
$Q\in\mathcal P_{U,g}$ such that $P\subsetneq Q$.
\end{enumerate}
\end{definition}

We use the definitions in~\cite{B} for maximal Dyck paths,
the following function $\pi_{d_1,d_2}$, compatible gradings
with fixed totals, and tightness.

\begin{definition}\label{def:dyck-grading}
Let $d_1$ and $d_2$ be positive integers.
\begin{enumerate}[font=\itshape]
\item The \textit{maximal Dyck path} $\mathcal{D}=\mathcal{D}_{(d_1,d_2)}$ is the highest NE path
from $(0,0)$ to $(d_1,d_2)$ that stays on or below the line segment
from $(0,0)$ to $(d_1,d_2)$.
\item\label{d:dyck-coordinate} For $e\in\operatorname{NE}(\mathcal D)$, let $(x,y)$ be its left endpoint if $e\in E(\mathcal D)$ and its lower endpoint if $e\in N(\mathcal D)$.
We define $\pi_{d_1,d_2}\colon\operatorname{NE}(\mathcal{D})\to[0,d_1+d_2-1]_{\Z}$ by
\[
 \pi_{d_1,d_2}(e)=xd_2-yd_1.
\]
\item\label{d:compatible-grading-set} Let $p$ and $q$ be nonnegative integers.
Let $\operatorname{CG}(d_1,d_2,p,q)$ denote the set of compatible gradings $g\in G(\mathcal D)$ satisfying
\[
 \operatorname{TNG}_{\mathcal{D},g}=p,\qquad
 \operatorname{TEG}_{\mathcal{D},g}=q.
\]
\item\label{d:tight-grading} Let $g\in\operatorname{CG}(d_1,d_2,p,q)$.
We call $g$ \textit{tight} if there exists $\varepsilon\in\{-1,1\}$ satisfying
\[
 pd_2-qd_1=\varepsilon\gcd(p,q)
\]
and the corresponding support condition below \cite[Definition 4.8 and Remark 4.9]{B}.
\begin{enumerate}[label=(\alph*),font=\itshape]
\item If $\varepsilon=1$, then
$E_g(\mathcal{D})\subseteq\operatorname{Sh}_{\mathcal{D},g}(N_g(\mathcal{D}))$.
\item If $\varepsilon=-1$, then
$N_g(\mathcal{D})\subseteq\operatorname{Sh}_{\mathcal{D},g}(E_g(\mathcal{D}))$.
\end{enumerate}
We define
\[
 \operatorname{TCG}(d_1,d_2,p,q)
 =\{g\in\operatorname{CG}(d_1,d_2,p,q):g\text{ is tight}\}.
\]
\end{enumerate}
\end{definition}

\begin{lemma}[Horizontal defect as a sum of debts]\label{lem:node-debt-sum}
Let $U=U_1\ccn\cdots\ccn U_n$ be a tNE path. Let $h\in G(U)$. Suppose that
$P=P_{U,h}$ exists. Let $X=\operatorname{Node}_{U,h}(U_r)=U[r-1,s]$.
Use the counts $\operatorname{Num}(y,k)$ from Lemma~\ref{lem:precursor-prefix-counts} for the grading $h$.
For $k\in[r,s]_{\Z}$, put
\[
 H_{r,k}=E(U[r-1,k]),\qquad
 \Delta_r(k)=\operatorname{NGrDiff}_{U,h}(U_r,U_k).
\]
Then, we have the following.
\begin{enumerate}
\item\label{c:node-debt-identity} For $k\in[r,s]_{\Z}$, we have
\[
 \Delta_r(k)=-\sum_{y\in H_{r,k}}\bigl(h(y)-\operatorname{Num}(y,k)\bigr)
 =-\sum_{y\in H_{r,k}}\operatorname{Debt}_{U,h}(y,k).
\]
\item\label{c:node-debt-first-zero} We have
\[
 \Delta_r(k)\le-1\quad(r\le k<s),\qquad \Delta_r(s)=0.
\]
In particular, $t_{U,h}(U_r)=U_s$.
\end{enumerate}
\end{lemma}
\begin{proof}
\prf{Claim~\ref{c:node-debt-identity}} Let us count the vertical edges by their h-precursors.
For $U_j\in N(U[r-1,k])$, Claim~\ref{c:node-predecessors} of
Lemma~\ref{lem:node-structure} places $P(U_j)$ in $E(X)$.
Since $P(U_j)$ precedes $U_j$, it belongs to $H_{r,k}$.
For $y=U_i\in H_{r,k}$ and $t\in F_y$, Definition~\ref{def:stack}(1)(a) gives $t>i\ge r$.
Consequently,
\[
 \#N(U[r-1,k])=\sum_{y\in H_{r,k}}\operatorname{Num}(y,k).
\]
The subpath $U[r-1,k]$ does not cross from $U_n$ to $U_1$. Subtracting $\sum_{y\in H_{r,k}}h(y)$
proves Claim~\ref{c:node-debt-identity}.

\prf{Claim~\ref{c:node-debt-first-zero}} First, let us prove $\Delta_r(k)<0$ for $k\in[r,s-1]_{\Z}$.
Fix such an index $k$.
Claim~\ref{c:prefix-capacity} of Lemma~\ref{lem:precursor-prefix-counts} gives
\[
 h(y)-\operatorname{Num}(y,k)\ge0\qquad(y\in H_{r,k}).
\]
Claim~\ref{c:assigned-prefix} of Lemma~\ref{lem:precursor-prefix-counts}, applied to $P(U_s)=U_r$, gives
$h(U_r)-\operatorname{Num}(U_r,k)\ge1$.
Since $U_r\in H_{r,k}$, Claim~\ref{c:node-debt-identity} gives
\[
 \Delta_r(k)=-\sum_{y\in H_{r,k}}(h(y)-\operatorname{Num}(y,k))\le-1.
\]

Second, let us prove $\Delta_r(s)=0$.
For $y\in H_{r,s}$ with $h(y)>0$, Claim~\ref{c:node-contained} of Lemma~\ref{lem:node-structure} gives
\[
 \operatorname{Node}_{U,h}(y)\subseteq X,\qquad
 P^{-1}(\{y\})\subseteq N(X).
\]
The node definition gives $\#P^{-1}(\{y\})=h(y)$, so $\operatorname{Num}(y,s)=h(y)$.
For $y\in H_{r,s}$ with $h(y)=0$, Claim~\ref{c:prefix-capacity} of Lemma~\ref{lem:precursor-prefix-counts} gives $\operatorname{Num}(y,s)=0$.
Thus, Claim~\ref{c:node-debt-identity} gives
\[
 \Delta_r(s)=-\sum_{y\in H_{r,s}}(h(y)-\operatorname{Num}(y,s))=0.
\]

Let us identify the first vertical zero of NGrDiff from $U_r$.
The indices $r,\ldots,s$ occur in this order before the cyclic path returns to $U_1$.
For every $j\in[r,s]_{\Z}$ with $U_j\in N(U)$, Claim~\ref{c:node-debt-identity} gives
\[
 \operatorname{NGrDiff}_{U,h}(U_r,U_j)=\Delta_r(j).
\]
The strict inequalities for $j<s$ and equality at $j=s$ prove the claim.
\end{proof}

\begin{lemma}[An explicit witness and a first failing pair]\label{lem:explicit-compatibility-witness}
Let $U=U_1\ccn\cdots\ccn U_n$ be a tNE path. Let $h\in G(U)$. Suppose that
$P=P_{U,h}$ exists. Suppose that every horizontal edge of positive grade has a node.
For each $U_k\in N(U)$, let $p_k$ be the unique index satisfying $P(U_k)=U_{p_k}$.
Put $B_i=\operatorname{bBal}_{U,h}(i)$ for $i\in[0,n]_{\Z}$.
The following statements hold.
\begin{enumerate}
\item\label{c:explicit-balance-increments} For $i\in[n]_{\Z}$, we have
\[
 B_i-B_{i-1}=\begin{cases}1,&U_i\in E(U),\\-h(U_i),&U_i\in N(U).
 \end{cases}
\]
\item\label{c:explicit-positive-witness} Let $U_r\in E(U)$, $U_s\in N(U)$, and $r<s$.
Suppose that $h(U_r)h(U_s)>0$ and $B_s\ge B_r$.
The set
\[
 J_{r,s}=\{i\in[r+1,s-1]_{\Z}:B_i>B_s\}
\]
is nonempty. Its minimum $j_* =\min J_{r,s}$ satisfies
\[
 U_{j_*}\in E(U),\qquad B_{j_*-1}=B_s,\qquad B_{j_*}=B_s+1,
 \qquad \operatorname{EGrDiff}_{U,h}(U_{j_*},U_s)=0.
\]
Thus, $U_{j_*}$ is a compatibility witness for $(U_r,U_s)$.
\item\label{c:explicit-failing-pair} Define
\[
 S=\{k\in[n]_{\Z}:U_k\in N(U),\ B_k<B_{p_k}\}.
\]
If $S\ne\emptyset$, put $s=\min S$, $r=p_s$, and $z=B_s$.
Then,
\[
 h(U_r)>0,\qquad h(U_s)>0,\qquad r<s,\qquad
 B_i\ge B_r>z\quad(r\le i<s).
\]
For $r<j<s$, we have
\[
 \begin{cases}
 \operatorname{EGrDiff}_{U,h}(U_j,U_s)=z-B_{j-1}<0,&U_j\in E(U),\\
 \operatorname{NGrDiff}_{U,h}(U_r,U_j)<0,&U_j\in N(U).
 \end{cases}
\]
The pair $(U_r,U_s)$ has positive grades and no compatibility witness.
Its cyclic interval is the ordinary interval from $U_r$ to $U_s$.
\end{enumerate}
\end{lemma}
\begin{proof}
\prf{Claim~\ref{c:explicit-balance-increments}} Let us compute the balance change at each edge.
By Definition~\ref{def:colored-grading}(3)(a), a horizontal edge increases $\operatorname{bBal}_{U,h}$ by $1$. A vertical edge $U_i$ decreases it by $h(U_i)$.

\prf{Claim~\ref{c:explicit-positive-witness}} Let us determine $j_*=\min J_{r,s}$.
Since $h(U_s)>0$, we have $B_{s-1}=B_s+h(U_s)>B_s\ge B_r$.
Thus, $s-1>r$ and $s-1\in J_{r,s}$.
Minimality and $B_r\le B_s$ give
\[
 B_{j_*-1}\le B_s<B_{j_*}.
\]
For $U_i\in N(U)$, Claim~\ref{c:explicit-balance-increments} gives $B_i-B_{i-1}=-h(U_i)\le0$.
Hence, $U_{j_*}$ is horizontal and $B_{j_*}=B_{j_*-1}+1$.
Integrality gives $B_{j_*-1}=B_s$.
By Definition~\ref{def:compatible-grading}(1) and Claim~\ref{c:explicit-balance-increments}, summing $B_i-B_{i-1}$ for $i\in[j_*,s]_{\Z}$ gives
\[
 \operatorname{EGrDiff}_{U,h}(U_{j_*},U_s)=B_s-B_{j_*-1}=0.
\]
The bounds $r<j_*<s$ verify the witness location and $U_{j_*}\ne U_r$.

\prf{Claim~\ref{c:explicit-failing-pair}} Let us identify the node of $U_r$.
Definition~\ref{def:stack}(1)(a) and (b) gives
\[
 r<s,\qquad 0\le \operatorname{Num}(U_r,s-1)<h(U_r).
\]
Here $\operatorname{Num}(U_r,s-1)$ is the count in Lemma~\ref{lem:precursor-prefix-counts}.
Thus, $h(U_r)>0$. Let $t\in[r,n]_{\Z}$ satisfy
\[
 \operatorname{Node}_{U,h}(U_r)=U[r-1,t].
\]
Since $P(U_s)=U_r$, we have $s\le t$.

Let us prove $B_i\ge B_r$ for $i\in[r,s-1]_{\Z}$ by strong induction on $i$.
For $i=r$, the inequality is equality.
Fix $i\in[r+1,s-1]_{\Z}$ and assume
\[
 B_j\ge B_r\qquad(j\in[r,i-1]_{\Z}).
\]

First, let us prove the induction step when $U_i\in E(U)$.
Claim~\ref{c:explicit-balance-increments} and the induction hypothesis give
\[
 B_i=B_{i-1}+1\ge B_r+1>B_r.
\]

Second, let us prove the induction step when $U_i\in N(U)$.
Recall that $P(U_i)=U_{p_i}$.
Claim~\ref{c:node-predecessors} of Lemma~\ref{lem:node-structure} gives
$r\le p_i<i$. Since $i<s=\min S$, we have $i\notin S$. Therefore,
\[
 B_i\ge B_{p_i}\ge B_r.
\]
The inequality $B_{p_i}\ge B_r$ follows from the induction hypothesis, since $r\le p_i<i$.

Let us prove $h(U_s)>0$.
The definition of $S$ gives $B_r>B_s=z$.
Claim~\ref{c:explicit-balance-increments} and the bound at $i=s-1$ give
\[
 h(U_s)=B_{s-1}-B_s\ge B_r-z>0.
\]

First, let us exclude a horizontal interior witness.
Let $j\in[r+1,s-1]_{\Z}$ satisfy $U_j\in E(U)$.
Definition~\ref{def:compatible-grading}(1) and Claim~\ref{c:explicit-balance-increments} give
\[
 \operatorname{EGrDiff}_{U,h}(U_j,U_s)
 =B_s-B_{j-1}=z-B_{j-1}\le z-B_r<0.
\]

Second, let us exclude a vertical interior witness.
Let $j\in[r+1,s-1]_{\Z}$ satisfy $U_j\in N(U)$.
Since $r<j<s\le t$, Claim~\ref{c:node-debt-first-zero} of
Lemma~\ref{lem:node-debt-sum}, applied to the node $U[r-1,t]$, gives
\[
 \operatorname{NGrDiff}_{U,h}(U_r,U_j)<0.
\]

Let us exclude the endpoint witnesses.
Definition~\ref{def:compatible-grading}(\ref{d:compatible}) requires an
EGrDiff witness different from $U_r$ and an NGrDiff witness different from $U_s$.
Thus, $(U_r,U_s)$ has no compatibility witness.
\end{proof}

\begin{lemma}[Witnesses under cyclic realization and zero extension]\label{lem:cyclic-witness-extension}
Let $D=D_1\ccn\cdots\ccn D_N$ be a nonempty tNE path and let $q\in[N]_{\Z}$ and $\ell\in[N]_{\Z}$.
For an integer $a$, let $a\bmod N$ denote its remainder in $[0,N-1]_{\Z}$.
Define
\[
 \eta(j)=1+((q+j-2)\bmod N),\qquad
 \theta(k)=1+((k-q)\bmod N)\quad(j,k\in[N]_{\Z}).
\]
Put
\[
 U=\operatorname{NEReal}(D,D_q,D_{\eta(\ell)}),\qquad
 \iota(U_j)=D_{\eta(j)}\quad(j\in[\ell]_{\Z}).
\]
Let $h\in G(U)$. Suppose that $P_{U,h}$ exists and that every horizontal
edge of positive grade has a node. Let $g\in G(D)$ satisfy
$g(D_{\eta(j)})=h(U_j)$ for $j\in[\ell]_{\Z}$.
Then, we have the following.
\begin{enumerate}
\item\label{c:cyclic-witness-order} The following statements hold.
\begin{enumerate}[label=(\alph*),ref=\theenumi(\alph*)]
\item\label{c:cyclic-inverse-indices} The maps $\eta$ and $\theta$ are inverse permutations of $[N]_{\Z}$. For $j\in[N]_{\Z}$ and $k\in[N]_{\Z}$, we have
\[
 \theta(\eta(j))=j,\qquad \eta(\theta(k))=k.
\]
\item\label{c:cyclic-interval-order} We have $|U|=\ell$.
For $1\le i\le j\le\ell$, the map $\iota(U_k)=D_{\eta(k)}$ for $k\in[i,j]_{\Z}$ identifies $U[i-1,j]$ with $D\langle D_{\eta(i)},D_{\eta(j)}\rangle$ in edge order.
For $k\in[\ell]_{\Z}$, we have
\[
 \begin{gathered}
 U_k\in E(U)\ \Longleftrightarrow\ D_{\eta(k)}\in E(D),\qquad
 U_k\in N(U)\ \Longleftrightarrow\ D_{\eta(k)}\in N(D),\\
 h(U_k)=g(D_{\eta(k)}).
 \end{gathered}
\]
\item\label{c:cyclic-defect-identities} For $1\le i\le j\le\ell$, we have
\[
 \begin{aligned}
 \operatorname{NGrDiff}_{D,g}(D_{\eta(i)},D_{\eta(j)})
 &=\operatorname{NGrDiff}_{U,h}(U_i,U_j),\\
 \operatorname{EGrDiff}_{D,g}(D_{\eta(i)},D_{\eta(j)})
 &=\operatorname{EGrDiff}_{U,h}(U_i,U_j).
 \end{aligned}
\]
\item\label{c:cyclic-witness-bijection} Suppose that $U_i\in E(U)$, $U_j\in N(U)$, $i<j$, and $h(U_i)h(U_j)>0$.
For every $k\in[i,j]_{\Z}$, the edge $U_k$ is a witness for $(U_i,U_j)$ if and only if $D_{\eta(k)}$ is a witness for $(D_{\eta(i)},D_{\eta(j)})$.
Thus, $\iota$ restricts to a bijection between these two witness sets.
\end{enumerate}
\item\label{c:cyclic-reversed-witness} Suppose that $U_i\in E(U)$, $U_j\in N(U)$,
$j<i$, and $h(U_i)h(U_j)>0$. Define
\[
 F_i=\{k\in[\ell]_{\Z}:U_k\in N(U),\ P_{U,h}(U_k)=U_i\},
 \qquad t_i=\max F_i.
\]
Then, $U_{t_i}$ is a witness for $(U_i,U_j)$, and $D_{\eta(t_i)}$ is a witness
for $(D_{\eta(i)},D_{\eta(j)})$.
\item\label{c:cyclic-zero-extension} If $g$ is compatible, then $h$ is compatible.
Define
\[
 g_0(D_k)=\begin{cases}
 h(U_{\theta(k)}),&\theta(k)\le\ell,\\
 0,&\theta(k)>\ell.
 \end{cases}
\]
Then, $g_0$ is compatible if and only if $h$ is compatible.
\end{enumerate}
\end{lemma}
\begin{proof}
\prf{Claim~\ref{c:cyclic-inverse-indices}} Let us check that the two index maps are inverses.
Let $j\in[N]_{\Z}$ and $k\in[N]_{\Z}$.
The definitions of $\eta$ and $\theta$ give
\[
 \begin{aligned}
 \theta(\eta(j))&=1+((j-1)\bmod N)=j,\\
 \eta(\theta(k))&=1+((k-1)\bmod N)=k,
 \end{aligned}
\]
because $0\le j-1<N$ and $0\le k-1<N$.
Hence, $\theta=\eta^{-1}$ on $[N]_{\Z}$.

\prf{Claim~\ref{c:cyclic-interval-order}} Let us identify the edges of $U$.
We have $\eta(1)=q$ and
\[
 \eta(k+1)=\begin{cases}\eta(k)+1,&\eta(k)<N,\\1,&\eta(k)=N\end{cases}
 \qquad(k\in[N-1]_{\Z}).
\]
Claim~\ref{c:cyclic-inverse-indices} makes $\eta(1),\ldots,\eta(\ell)$ pairwise distinct.
Thus, the cyclic traversal from $D_q$ first reaches $D_{\eta(\ell)}$ after exactly $\ell$ edges, in the order $D_{\eta(1)},\ldots,D_{\eta(\ell)}$.
Definition~\ref{def:cyclic-subpaths} gives $|U|=\ell$ and preserves the direction of every edge under $\iota$.
For $1\le i\le j\le\ell$, the sequence $D_{\eta(i)},\ldots,D_{\eta(j)}$ has no repeated edge and ends at its first occurrence of $D_{\eta(j)}$.
Hence, it is precisely $D\langle D_{\eta(i)},D_{\eta(j)}\rangle$.
The hypothesis $g(D_{\eta(k)})=h(U_k)$ supplies grade preservation for each $k\in[\ell]_{\Z}$.

\prf{Claim~\ref{c:cyclic-defect-identities}} Let us compare the grading differences under the cyclic index map.
Let $1\le i\le j\le\ell$.
By Claim~\ref{c:cyclic-interval-order} and Definition~\ref{def:compatible-grading}(1), we have
\[
 \begin{aligned}
 \operatorname{NGrDiff}_{D,g}(D_{\eta(i)},D_{\eta(j)})
 &=\#\{k\in[i,j]_{\Z}:U_k\in N(U)\}
   -\sum_{\substack{k\in[i,j]_{\Z}\\U_k\in E(U)}}h(U_k)\\
 &=\operatorname{NGrDiff}_{U,h}(U_i,U_j),\\
 \operatorname{EGrDiff}_{D,g}(D_{\eta(i)},D_{\eta(j)})
 &=\#\{k\in[i,j]_{\Z}:U_k\in E(U)\}
   -\sum_{\substack{k\in[i,j]_{\Z}\\U_k\in N(U)}}h(U_k)\\
 &=\operatorname{EGrDiff}_{U,h}(U_i,U_j).
 \end{aligned}
\]

\prf{Claim~\ref{c:cyclic-witness-bijection}} Let us transport the compatibility witnesses in both directions.
Let $k\in[i,j]_{\Z}$.
By Definition~\ref{def:compatible-grading}(\ref{d:compatible}), the edge $U_k$ is a witness for $(U_i,U_j)$ precisely when one of the following holds:
\[
 \begin{aligned}
 &U_k\in N(U),\quad k\ne j,\quad \operatorname{NGrDiff}_{U,h}(U_i,U_k)=0;\\
 &U_k\in E(U),\quad k\ne i,\quad \operatorname{EGrDiff}_{U,h}(U_k,U_j)=0.
 \end{aligned}
\]
Claim~\ref{c:cyclic-inverse-indices} gives $k\ne j\Longleftrightarrow\eta(k)\ne\eta(j)$ and $k\ne i\Longleftrightarrow\eta(k)\ne\eta(i)$.
Claims~\ref{c:cyclic-interval-order} and~\ref{c:cyclic-defect-identities} therefore identify these two conditions with the two witness conditions for $D_{\eta(k)}$ and $(D_{\eta(i)},D_{\eta(j)})$.
Claim~\ref{c:cyclic-interval-order} includes every edge of $D\langle D_{\eta(i)},D_{\eta(j)}\rangle$.
The inverse on witnesses is $D_{\eta(k)}\mapsto U_{\theta(\eta(k))}=U_k$.

\prf{Claim~\ref{c:cyclic-reversed-witness}} Let us find a witness when the edge order is reversed.
The node at $U_i$ exists. Thus, $F_i$ is
nonempty and its maximum $t_i$ satisfies $i<t_i\le\ell$.
Claim~\ref{c:node-debt-first-zero} of Lemma~\ref{lem:node-debt-sum} gives
$\operatorname{NGrDiff}_{U,h}(U_i,U_{t_i})=0$.
Claim~\ref{c:cyclic-defect-identities} gives $\operatorname{NGrDiff}_{D,g}(D_{\eta(i)},D_{\eta(t_i)})=0$.

Let us verify that the two vertical edges lie strictly between the two edges of each pair.
In cyclic order on $U$, the distance from index $i$ to $j<i$ is $\ell-i+j$. We have
\[
 0<t_i-i\le\ell-i<\ell-i+j.
\]
In cyclic order on $D$, the distance from $\eta(i)$ to $\eta(j)$ is $N-i+j$. We have
\[
 t_i-i\le\ell-i\le N-i<N-i+j.
\]
The edges $U_{t_i}$ and $D_{\eta(t_i)}$ are vertical.
The inequalities for $t_i-i$ place them strictly before $U_j$ and $D_{\eta(j)}$ in the respective cyclic intervals.
The two NGrDiff identities and Definition~\ref{def:compatible-grading}(\ref{d:compatible}) therefore make them witnesses for $(U_i,U_j)$ and $(D_{\eta(i)},D_{\eta(j)})$, respectively.

\prf{Claim~\ref{c:cyclic-zero-extension}} Let us prove compatibility in both directions.
Suppose first that $g$ is compatible.
Let $U_i\in E(U)$ and $U_j\in N(U)$ satisfy $h(U_i)h(U_j)>0$.
Their images have positive $g$-grades.
If $i<j$, Claim~\ref{c:cyclic-witness-bijection} transports a witness from $D$ to $U$.
If $j<i$, Claim~\ref{c:cyclic-reversed-witness} supplies a witness in $U$.
The directions exclude $i=j$. Thus, $h$ is compatible.
Claim~\ref{c:cyclic-inverse-indices} gives $g_0(D_{\eta(j)})=h(U_j)$ for $j\le\ell$.
The implication $g\text{ compatible}\Longrightarrow h\text{ compatible}$ therefore applies to $g=g_0$.

Let us prove compatibility of the zero extension.
Conversely, let $h$ be compatible and let $(D_a,D_b)$ be a positive
horizontal--vertical pair for $g_0$.
Positivity gives the explicit indices $i=\theta(a)\le\ell$ and
$j=\theta(b)\le\ell$, and $h(U_i)h(U_j)>0$.
If $i<j$, Claim~\ref{c:cyclic-witness-bijection} transports a witness from $U$ to $D$.
If $j<i$, Claim~\ref{c:cyclic-reversed-witness} supplies an ambient witness.
Thus, $g_0$ is compatible. The argument includes $\ell=N$.
\end{proof}

\begin{lemma}[Nodes and horizontal shadows]\label{lem:compatibility}
Let $\mathcal{D}$ be a maximal Dyck path and $g\in G(\mathcal D)$.
Consider a maximal h-shadow subpath $\mathcal{D}\langle u_0,v_0\rangle$.
Let $U=\operatorname{NEReal}(\mathcal{D},u_0,v_0)$ and
$h=g_{\mathcal{D},u_0,v_0}$.
Then, we have the following.

\begin{enumerate}
\item\label{c:shadow-nodes} The following statements hold.
\begin{enumerate}[label=(\alph*),ref=\theenumi(\alph*)]
\item\label{c:shadow-predecessor} The map $P_{U,h}$ exists.
\item\label{c:shadow-node-identity} For each $x\in E(U)$ with $h(x)>0$,
the node $\operatorname{Node}_{U,h}(x)$ exists.
The edge $\operatorname{End}_{\mathrm{NE}}(\operatorname{Node}_{U,h}(x))$ corresponds to $t_{\mathcal{D},g}(x)$.
\item\label{c:shadow-whole-node} We have $\operatorname{Node}_{U,h}(u_0)=U$.
\end{enumerate}
\item\label{c:shadow-compatibility} The grading $h$ of $U$ is compatible
if and only if $h$ satisfies the node balance condition.
\item\label{c:shadow-extension} The following statements hold.
\begin{enumerate}[label=(\alph*),ref=\theenumi(\alph*)]
\item\label{c:shadow-zero-extension} Define $g_0\in G(\mathcal D)$ by
\[
 g_0(e)=\begin{cases}
 g(e),&e\in\operatorname{NE}(\mathcal D,u_0,v_0),\\
 0,&e\in\operatorname{NE}(\mathcal D)\setminus\operatorname{NE}(\mathcal D,u_0,v_0).
 \end{cases}
\]
Then, $g_0$ is compatible if and only if $h$ is connected.
\item\label{c:shadow-restricted-connected} If $g$ is compatible, then $h$ is connected.
\end{enumerate}
\end{enumerate}
\end{lemma}

\begin{proof}
\prf{Claim~\ref{c:shadow-predecessor}} Let us prove that $P_{U,h}$ exists.
Write $n=|U|$, $U_1=u_0$, $U_n=v_0$, and
\[
 A_k=\#N(U[0,k])-\sum_{y\in E(U[0,k])}h(y)\quad(k\in[0,n]_{\Z}).
\]
The realization preserves the ordered steps and grades from $u_0$ through $v_0$.
Definition~\ref{def:compatible-grading}(\ref{d:horizontal-shadow}) gives $A_n=0$ and $A_1=-h(u_0)<0$.
A horizontal step decreases $A_k$ by a nonnegative integer;
a vertical step increases it by one.
If $A_k\ge0$ for some $k\in[n-1]_{\Z}$, let $k$ be the least such index.
Since $A_1<0$, we have $k\ge2$.
The edge $U_k$ is vertical and $A_k=A_{k-1}+1=0$, contrary to the definition of $v_0$.
Thus, \begin{equation}\label{eq:shadow-component-balances}
 A_k<0\quad(1\le k<n),\qquad A_n=0.
\end{equation}
Since $A_0=0$, Claim~\ref{c:predecessor-existence} of
Lemma~\ref{lem:predecessor-existence} gives $P_{U,h}$.

\prf{Claim~\ref{c:shadow-node-identity}} Let us prove saturation of every fiber.
Use the counts $\operatorname{Num}(y,k)$ from Lemma~\ref{lem:precursor-prefix-counts}.
Claim~\ref{c:prefix-capacity} of Lemma~\ref{lem:precursor-prefix-counts} gives $\operatorname{Num}(y,n)\le h(y)$.
Each vertical edge has one h-precursor, so~\eqref{eq:shadow-component-balances} gives
\[
 \sum_{y\in E(U)}\bigl(h(y)-\operatorname{Num}(y,n)\bigr)
 =\sum_{y\in E(U)}h(y)-\#N(U)=0.
\]
Every summand is nonnegative. Hence, $\operatorname{Num}(y,n)=h(y)$ for every $y\in E(U)$.
For $y=U_r\in E(U)$ with $h(y)>0$, the fiber $P_{U,h}^{-1}(\{y\})$ is nonempty and saturated. Thus, its node exists; write $\operatorname{Node}_{U,h}(y)=U[r-1,s]$.
Claim~\ref{c:node-debt-first-zero} of Lemma~\ref{lem:node-debt-sum} gives
\[
 \operatorname{NGrDiff}_{U,h}(U_r,U_j)<0\quad(r\le j<s),\qquad
 \operatorname{NGrDiff}_{U,h}(U_r,U_s)=0.
\]
Claims~\ref{c:cyclic-interval-order} and~\ref{c:cyclic-defect-identities} of Lemma~\ref{lem:cyclic-witness-extension} transport $\operatorname{NGrDiff}_{U,h}(U_r,U_j)<0$ for $r\le j<s$ and $\operatorname{NGrDiff}_{U,h}(U_r,U_s)=0$ to $\mathcal D$.
Thus, $U_s$ corresponds to $t_{\mathcal D,g}(y)$.

\prf{Claim~\ref{c:shadow-whole-node}} Let us identify the node at the first edge.
Apply Claim~\ref{c:shadow-node-identity} to $u_0$.
Its first zero is $v_0$. Thus, $\operatorname{Node}_{U,h}(u_0)=U$.

\prf{Claim~\ref{c:shadow-compatibility}} Let us prove that node balance implies compatibility.
Suppose that $h$ satisfies the node balance condition.
Let $U_r\in E(U)$ and $U_s\in N(U)$ satisfy $h(U_r)h(U_s)>0$.
Let $t$ be the ending index of the node
\[
 \operatorname{Node}_{U,h}(U_r)=U[r-1,t].
\]

First, let us find a witness when $U_s\notin\operatorname{NE}(U[r-1,t])$.
Then, $s<r$ or $s>t$.
In either case, $U_t$ occurs before $U_s$ in cyclic order from $U_r$, and $U_t\ne U_s$.
Claim~\ref{c:node-debt-first-zero} of Lemma~\ref{lem:node-debt-sum} gives
\[
 \operatorname{NGrDiff}_{U,h}(U_r,U_t)=0.
\]
Thus, $U_t$ is a compatibility witness.

Second, let us find a witness when $U_s\in\operatorname{NE}(U[r-1,t])$.
Since $U_r$ is horizontal and $U_s$ is vertical, we have $r<s\le t$.
Claim~\ref{c:node-boundary-balances} of Lemma~\ref{lem:node-structure} gives
\[
 \operatorname{bBal}_{U,h}(s)\ge\operatorname{bBal}_{U,h}(r).
\]
Claim~\ref{c:explicit-positive-witness} of Lemma~\ref{lem:explicit-compatibility-witness}
therefore gives an EGrDiff witness for $(U_r,U_s)$.

Let us prove that compatibility implies node balance.
Suppose that the node balance condition fails.
For $U_k\in N(U)$, let $p_k$ satisfy $P_{U,h}(U_k)=U_{p_k}$.
Then, the set
\[
 S=\{k\in[|U|]_{\Z}:U_k\in N(U),\
       \operatorname{bBal}_{U,h}(k)<\operatorname{bBal}_{U,h}(p_k)\}
\]
is nonempty.
Claim~\ref{c:shadow-node-identity} gives a node for every horizontal edge of positive grade.
Claim~\ref{c:explicit-failing-pair} of Lemma~\ref{lem:explicit-compatibility-witness}, applied to $s=\min S$ and $r=p_s$, gives
\[
 h(U_r)h(U_s)>0
\]
and no compatibility witness for $(U_r,U_s)$.
Thus, $h$ is not compatible.

\prf{Claim~\ref{c:shadow-zero-extension}} Let us check compatibility of the zero extension.
Claim~\ref{c:cyclic-zero-extension} of
Lemma~\ref{lem:cyclic-witness-extension} gives compatibility of $g_0$
if and only if compatibility of $h$.
Claim~\ref{c:shadow-compatibility} identifies compatibility of $h$ with the node balance condition for $h$.
Claims~\ref{c:shadow-predecessor} and~\ref{c:shadow-whole-node} supply
the h-precursor $P_{U,h}$ and the identity $\operatorname{Node}_{U,h}(U_1)=U$.

\prf{Claim~\ref{c:shadow-restricted-connected}} Let us prove connectedness of the restriction.
Claim~\ref{c:cyclic-zero-extension} of
Lemma~\ref{lem:cyclic-witness-extension} sends compatibility of $g$ to compatibility of $h$.
Claim~\ref{c:shadow-compatibility} gives the node balance condition.
Claims~\ref{c:shadow-predecessor} and~\ref{c:shadow-whole-node} give connectedness.
\end{proof}

\begin{definition}\label{def:frame}
Let $b>a\ge2$ be coprime integers.
Let $t\in[a-1]_{\Z}$ be the unique integer satisfying $bt\equiv1\pmod a$.
We define
\begin{equation}\label{eq:ambient}
\begin{gathered}
 \mu_1(a,b)=a+t,\qquad \mu_2(a,b)=\frac{b\,\mu_1(a,b)-1}{a},\\
 \Gamma(a,b)=\frac{\mu_2(a,b)}{\mu_1(a,b)}=\frac ba-\frac1{a\,\mu_1(a,b)}.
\end{gathered}
\end{equation}
For $i\in[a]_{\Z}$, we define $\Lambda_i(a,b)=\floor{i\,\Gamma(a,b)}-\floor{(i-1)\Gamma(a,b)}$.
We define the run
\begin{equation}\label{eq:frame-word}
 \mathcal{U}(a,b)=R_{\Lambda_1(a,b)}\ccn \cdots\ccn R_{\Lambda_a(a,b)}\ccn R_1.
\end{equation}
\end{definition}

\begin{lemma}[Decomposition into maximal horizontal shadows]\label{lem:frame-shadow-decomposition}
Let $b>a\ge2$ be coprime integers, and let $\mu_j=\mu_j(a,b)$ for $j\in[2]_{\Z}$.
Let $\mathcal D=\mathcal D_{(\mu_1,\mu_2)}$ and $g\in\operatorname{TCG}(\mu_1,\mu_2,a,b)$.
The maximal h-shadow subpaths of $g$ are pairwise disjoint paths
$\mathcal D\langle u_i,v_i\rangle$, indexed by $i\in[M]_{\Z}$ for some $M\ge1$.
Every edge outside these subpaths has grade zero.

For $i\in[M]_{\Z}$, let
\[
 U_i=\operatorname{NEReal}(\mathcal D,u_i,v_i),\qquad
 g_i=g_{\mathcal D,u_i,v_i},\qquad
 q_i=\#N(U_i),\qquad p_i=\#E(U_i)-1.
\]
Then, $g_i$ is connected, and
\[
 \operatorname{TEG}_{U_i,g_i}=q_i,\qquad
 \operatorname{TNG}_{U_i,g_i}\le p_i,\qquad
 \sum_{i=1}^{M}q_i=b,\qquad
 \sum_{i=1}^{M}\operatorname{TNG}_{U_i,g_i}=a.
\]
\end{lemma}
\begin{proof}
Let us fix the ambient dimensions and the total grades.
Let $n=\mu_1+\mu_2$ and write $\mathcal D=D_1\ccn\cdots\ccn D_n$.
Equation~\eqref{eq:ambient} and Definition~\ref{def:dyck-grading}(\ref{d:compatible-grading-set}) give
\[
 \mu_1\ge a+1,\qquad \mu_2\ge b+1,\qquad
 b\mu_1-a\mu_2=1,\qquad
 \sum_{u\in E(\mathcal D)}g(u)=b,\qquad
 \sum_{v\in N(\mathcal D)}g(v)=a.
\]

Let us define the sums used to choose a cyclic cut.
For $i\in\Z_{>0}$, define
\[
 \overline i=1+((i-1)\bmod n),\qquad
 \delta_i=
 \begin{cases}
 g(D_{\overline i}),&D_{\overline i}\in E(\mathcal D),\\
 -1,&D_{\overline i}\in N(\mathcal D).
 \end{cases}
\]
For $t\in\Z_{\ge0}$, define $A(t)=\sum_{i=1}^{t}\delta_i$.
Then,
\[
 A(0)=0,\qquad A(n)=b-\mu_2<0,\qquad
 A(t+n)=A(t)+A(n)\quad(t\in\Z_{\ge0}).
\]

Let us choose the cut.
Choose $c\in[0,n-1]_{\Z}$ satisfying
$A(c)=\min\{A(j):j\in[0,n-1]_{\Z}\}$.
For $i\in[n]_{\Z}$, define $\eta(i)=1+((c+i-1)\bmod n)$, and let
\[
 W=\operatorname{NEReal}(\mathcal D,D_{\eta(1)},D_{\eta(n)}),
 \qquad h=g_{\mathcal D,D_{\eta(1)},D_{\eta(n)}}.
\]
The edge $W_i$ corresponds to $D_{\eta(i)}$, and
$h(W_i)=g(D_{\eta(i)})$.

Let us compute the sums after the cut.
For $t\in[0,n]_{\Z}$, define
\[
 B(t)=A(c+t)-A(c)
 =\sum_{u\in E(W[0,t])}h(u)-\#N(W[0,t]).
\]
Thus, $B(0)=0$ and $B(n)=A(n)$.
If $c+t<n$, minimality of $A(c)$ gives $B(t)\ge0>B(n)$.
If $c+t\ge n$, it gives
\[
 B(t)=A(c+t-n)-A(c)+A(n)\ge A(n)=B(n).
\]
Therefore, $B(n)\le B(t)$ for every $t\in[0,n]_{\Z}$.

Let us check the effect of the cut on compatibility and shadows.
The map $W_i\mapsto D_{\eta(i)}$ preserves directions, grades, and cyclic edge order.
For $r,s\in[n]_{\Z}$, the sums in Definition~\ref{def:compatible-grading}(\ref{d:grading-differences}) therefore give
\[
 \begin{aligned}
 \operatorname{NGrDiff}_{W,h}(W_r,W_s)
 &=\operatorname{NGrDiff}_{\mathcal D,g}(D_{\eta(r)},D_{\eta(s)}),\\
 \operatorname{EGrDiff}_{W,h}(W_r,W_s)
 &=\operatorname{EGrDiff}_{\mathcal D,g}(D_{\eta(r)},D_{\eta(s)}).
 \end{aligned}
\]
Hence, $h$ is compatible.
Since $a\mu_2-b\mu_1=-1$, tightness gives
$N_h(W)\subseteq\operatorname{Sh}_{W,h}(E_h(W))$.

Let us show that every positive horizontal grade closes before the end of $W$.
Let $r\in[n]_{\Z}$ satisfy $W_r\in E(W)$ and $h(W_r)>0$.
We have
\[
 B(r)>B(r-1),\qquad B(n)\le B(r-1).
\]

Let us find the first return to $B(r-1)$.
A horizontal edge increases $B$ by its nonnegative grade, and a vertical edge decreases $B$ by $1$.
Thus, the first index $t\in[r,n]_{\Z}$ with $B(t)\le B(r-1)$ satisfies
\[
 r<t,\qquad W_t\in N(W),\qquad B(t)=B(r-1),\qquad
 B(j)>B(r-1)\quad(r\le j<t).
\]
Write $t(r)=t$.
Since $\operatorname{NGrDiff}_{W,h}(W_r,W_j)=B(r-1)-B(j)$ for $r\le j\le t(r)$,
we have $t_{W,h}(W_r)=W_{t(r)}$.

Let us determine the grade of a vertical edge after a minimum.
Let $s\in[0,n-1]_{\Z}$ satisfy
$B(s)\le B(j)$ for every $j\in[0,s]_{\Z}$.
Suppose that $W_{s+1}$ is vertical and $h(W_{s+1})>0$.
The shadow support condition gives an $r\in[s]_{\Z}$ such that
\[
 W_r\in E_h(W),\qquad s+1\le t(r).
\]
The strict inequality at $j=s$ then gives $B(s)>B(r-1)$, contrary to the choice of $s$.
Therefore, $h(W_{s+1})=0$.

Let us divide $W$ into zero-grade edges and complete horizontal shadows.
Start at the character boundary $s=0$.
At each boundary $s<n$, maintain
$B(s)\le B(j)$ for $j\in[0,s]_{\Z}$.
If $h(W_{s+1})=0$, pass to $s+1$; the value of $B$ either stays fixed or decreases by $1$.
Thus, $B(s+1)\le B(j)$ for every $j\in[0,s+1]_{\Z}$.

Let us give the step for a positive grade.
If $h(W_{s+1})>0$, then $W_{s+1}\in E(W)$: under $B(s)\le B(j)$ for $j\in[0,s]_{\Z}$, we proved $W_{s+1}\in N(W)\Longrightarrow h(W_{s+1})=0$.
Record the subpath $W[s,t(s+1)]$ and pass to $t(s+1)$.
On $W[s,t(s+1)]$, we have
\[
 B(t(s+1))=B(s),\qquad
 B(j)>B(s)\quad(s<j<t(s+1)).
\]
Thus, $B(t(s+1))\le B(j)$ for every $j\in[0,t(s+1)]_{\Z}$ in the positive-grade case.

Let us record the resulting intervals.
Each step increases $s$, so the procedure ends at $s=n$.
Let the recorded subpaths be $W[r_i-1,t_i]$ for $i\in[M]_{\Z}$, in their order in $W$.
They satisfy
\[
 1\le r_i<t_i\le n,\qquad t_i=t(r_i),\qquad
 t_i<r_{i+1}\quad(1\le i<M).
\]
Every edge outside these subpaths has grade zero.
Since $\sum_{u\in E(W)}h(u)=b>0$, we have $M\ge1$.

Let us prove that these are exactly the maximal h-shadow subpaths.
Let $W_r\in E_h(W)$ belong to $W[r_i-1,t_i]$.
The first-return inequalities for $r_i$ give
$B(t_i)=B(r_i-1)\le B(r-1)$.
Since $B(r)>B(r-1)$ and each decrease of $B$ is $1$, the first return from $r$ satisfies
\[
 r<t(r)\le t_i.
\]
Thus, every h-shadow subpath starting at a positive horizontal grade is contained in a recorded subpath.
The recorded subpaths are disjoint, and each is itself an h-shadow subpath.
They are therefore exactly the maximal h-shadow subpaths.

Let us transport these subpaths to $\mathcal D$.
For $i\in[M]_{\Z}$, let
\[
 u_i=D_{\eta(r_i)},\qquad v_i=D_{\eta(t_i)},\qquad
 U_i=\operatorname{NEReal}(\mathcal D,u_i,v_i),\qquad
 g_i=g_{\mathcal D,u_i,v_i}.
\]
The cyclic index map identifies $U_i$ with $W[r_i-1,t_i]$ in edge order and preserves all grades.
Consequently, $\mathcal D\langle u_i,v_i\rangle$ for $i\in[M]_{\Z}$ are the maximal h-shadow subpaths of $g$.
Every edge of $\mathcal D$ outside them has grade zero.

Let us compute the totals on these subpaths.
For $i\in[M]_{\Z}$, define
$q_i=\#N(U_i)$ and $p_i=\#E(U_i)-1$.
The equality $B(t_i)=B(r_i-1)$ gives
$\operatorname{TEG}_{U_i,g_i}=q_i$.
Therefore,
\[
 \sum_{i=1}^{M}q_i=b,\qquad
 \sum_{i=1}^{M}\operatorname{TNG}_{U_i,g_i}=a.
\]

Let us bound the total vertical grade of each subpath.
Claim~\ref{c:shadow-restricted-connected} of Lemma~\ref{lem:compatibility} makes each $g_i$ connected.
Lemma~\ref{lem:residual-nonnegative} gives
\[
 p_i-\operatorname{TNG}_{U_i,g_i}
 =\operatorname{Surp}(U_i,g_i)\ge0\quad(i\in[M]_{\Z}).
\]
Summing over $i\in[M]_{\Z}$ gives $\sum_{i=1}^{M}p_i\ge a$.
\end{proof}

\Needspace{14\baselineskip}
\begin{lemma}[The initial path and its horizontal shadow]\label{lem:frame-shadow-position}
Let $b>a\ge2$ be coprime integers, and let $\mu_j=\mu_j(a,b)$ for $j\in[2]_{\Z}$.
Let $\mathcal D=\mathcal D_{(\mu_1,\mu_2)}$.

\begin{enumerate}
\item\label{c:frame-initial-path} The path $\mathcal U(a,b)$ is an initial subpath of $\mathcal D$, and
\[
 \sum_{i=1}^{a}\Lambda_i(a,b)=b-1,\qquad
 \#E(\mathcal U(a,b))=a+1,\qquad \#N(\mathcal U(a,b))=b.
\]
\item\label{c:frame-unique-shadow} Let $g\in\operatorname{TCG}(\mu_1,\mu_2,a,b)$.
The unique maximal h-shadow subpath of $g$ is the initial subpath $\mathcal U(a,b)$ of $\mathcal D$.
Every edge outside $\mathcal U(a,b)$ has grade zero.
The restriction $g|_{\mathcal U(a,b)}$ is connected and has surplus zero.
\end{enumerate}
\end{lemma}
\begin{proof}
\prf{Claim~\ref{c:frame-initial-path}} Let us identify the initial subpath.
Equation~\eqref{eq:ambient} gives $\mu_1\ge a+1$ and $\Gamma(a,b)>1$.
The first $a$ simple runs of $\mathcal D$ have vertical lengths $\Lambda_1(a,b),\ldots,\Lambda_a(a,b)$, and
\[
 \sum_{i=1}^{a}\Lambda_i(a,b)
 =\floor{a\Gamma(a,b)}
 =\floor{b-\frac1{\mu_1}}=b-1.
\]
Since $\Gamma(a,b)>1$, the next simple run starts with a horizontal edge followed by a vertical edge.
Appending these two edges gives $\mathcal U(a,b)$ by~\eqref{eq:frame-word}.
Its horizontal and vertical edge counts are $a+1$ and $b$.

\prf{Claim~\ref{c:frame-unique-shadow}} This shadow localization is the case
$p=a$, $q=b$, and $(d_1,d_2)=(\mu_1,\mu_2)$ of
\cite[Lemmas~4.12 and~4.13]{B}. We give a direct proof.
Let us fix the maximal h-shadow subpaths.
Use $M$, $u_i$, $v_i$, $U_i$, $g_i$, $p_i$, and $q_i$ from Lemma~\ref{lem:frame-shadow-decomposition}.
Let $n=\mu_1+\mu_2$ and write $\mathcal D=D_1\ccn\cdots\ccn D_n$.
Equation~\eqref{eq:ambient} gives $b\mu_1-a\mu_2=1$ and $\gcd(\mu_1,\mu_2)=1$.

Let us prove that the edge coordinates are distinct.
For $j\in[n]_{\Z}$, let $(x_j,y_j)$ be the starting point of $D_j$.
Then, $x_j+y_j=j-1$, so
\[
 \pi_{\mu_1,\mu_2}(D_j)
 =\mu_2(j-1)-ny_j.
\]
If $\pi_{\mu_1,\mu_2}(D_i)=\pi_{\mu_1,\mu_2}(D_j)$ for $i,j\in[n]_{\Z}$,
then $n$ divides $\mu_2(i-j)$.
Since $\gcd(\mu_2,n)=\gcd(\mu_1,\mu_2)=1$ and $|i-j|<n$, we obtain $i=j$.
In particular, $D_1$ is the unique edge with coordinate $0$.

Let us bound the number of subpaths.
Lemma~\ref{lem:frame-shadow-decomposition} gives $\sum_{i=1}^{M}q_i=b$ and $\sum_{i=1}^{M}p_i\ge a$.
Moving past a horizontal edge changes $\pi_{\mu_1,\mu_2}$ by $\mu_2$, and moving past a vertical edge changes it by $-\mu_1$.
The coordinate increments $+\mu_2$ on horizontal edges and $-\mu_1$ on vertical edges also apply from $D_n$ to $D_1$, since their total over $\mathcal D$ is $\mu_1\mu_2-\mu_2\mu_1=0$.
Since $u_i$ is horizontal and $v_i$ is vertical, Definition~\ref{def:dyck-grading}(\ref{d:dyck-coordinate}) gives
\begin{equation}\label{eq:shadow-endpoint-bound}
\begin{aligned}
 q_i\mu_1-p_i\mu_2
 &=\pi_{\mu_1,\mu_2}(u_i)+\mu_1+\mu_2-\pi_{\mu_1,\mu_2}(v_i)\\
 &\ge\pi_{\mu_1,\mu_2}(u_i)+1\ge1
 \qquad(i\in[M]_{\Z}).
\end{aligned}
\end{equation}
Summing~\eqref{eq:shadow-endpoint-bound} gives
\begin{equation}\label{eq:shadow-component-count}
 1\le M\le b\mu_1-\mu_2\sum_{i=1}^{M}p_i
 \le b\mu_1-a\mu_2=1.
\end{equation}
Hence, $M=1$.

Let us determine the endpoints of this subpath.
We have $q_1=b$ and $p_1\ge a$.
Equations~\eqref{eq:shadow-endpoint-bound} and~\eqref{eq:shadow-component-count} give
\[
 1\le b\mu_1-p_1\mu_2
 =1-(p_1-a)\mu_2\le1.
\]
Thus, $p_1=a$, and
\[
 \begin{gathered}
 1=\pi_{\mu_1,\mu_2}(u_1)+\mu_1+\mu_2-\pi_{\mu_1,\mu_2}(v_1),\\
 \pi_{\mu_1,\mu_2}(u_1)=0,\qquad
 \pi_{\mu_1,\mu_2}(v_1)=\mu_1+\mu_2-1.
 \end{gathered}
\]
\Needspace{5\baselineskip}
Let us identify the subpath with $\mathcal U(a,b)$.
Coordinate injectivity gives $u_1=D_1$.
The subpath has $a+1$ horizontal edges and $b$ vertical edges.
Claim~\ref{c:frame-initial-path} identifies the first $a+b+1$ edges of $\mathcal D$ with $\mathcal U(a,b)$.
Hence,
\[
 \mathcal D\langle u_1,v_1\rangle=\mathcal U(a,b),\qquad
 \operatorname{Surp}(\mathcal U(a,b),g_1)=p_1-a=0.
\]
\end{proof}

\Needspace{14\baselineskip}
\begin{lemma}[Restriction and zero extension]\label{lem:frame-grading-bijection}
Let $b>a\ge2$ be coprime integers, and let $\mu_j=\mu_j(a,b)$ for $j\in[2]_{\Z}$.
Let $\mathcal D=\mathcal D_{(\mu_1,\mu_2)}$ and $V=\mathcal U(a,b)$.
Restriction defines a bijection
\[
 \operatorname{TCG}(\mu_1,\mu_2,a,b)
 \longrightarrow
 \{f\in\operatorname{CtdG}(V):\operatorname{Surp}(V,f)=0\},
 \qquad g\longmapsto g|_V.
\]
Its inverse sends $f$ to the grading $g\in G(\mathcal D)$ defined by
\[
 g(e)=
 \begin{cases}
 f(e),&e\in\operatorname{NE}(V),\\
 0,&e\in\operatorname{NE}(\mathcal D)\setminus\operatorname{NE}(V).
 \end{cases}
\]
\end{lemma}
\begin{proof}
Let us check that restriction has the stated values.
Claim~\ref{c:frame-unique-shadow} of Lemma~\ref{lem:frame-shadow-position} identifies $V$ with the unique maximal h-shadow subpath of each tight grading.
It makes the restricted grading connected with surplus zero.

Let us fix a grading to extend by zero.
Let $f\in\operatorname{CtdG}(V)$ satisfy $\operatorname{Surp}(V,f)=0$, and let $g$ be its zero extension.
Claim~\ref{c:frame-initial-path} of Lemma~\ref{lem:frame-shadow-position} gives $\#E(V)=a+1$ and $\#N(V)=b$.
Let $h_0=V_1$ and $v_*=V_{a+b+1}$.

Let us compute the two total grades of $f$.
The definition of connectedness gives $\operatorname{Node}_{V,f}(h_0)=V$.
Claim~\ref{c:node-contained} of Lemma~\ref{lem:node-structure} gives a node for each $x\in E(V)$ with $f(x)>0$.
Thus, $\#P_{V,f}^{-1}(\{x\})=f(x)$ for each such $x$.
For $x\in E(V)$ with $f(x)=0$, the equality $\#P_{V,f}^{-1}(\{x\})=0$ follows from Claim~\ref{c:prefix-capacity} of Lemma~\ref{lem:precursor-prefix-counts}.

Let us sum these fiber sizes.
Each vertical edge has exactly one h-precursor, so
\[
 \operatorname{TEG}_{V,f}
 =\sum_{x\in E(V)}\#P_{V,f}^{-1}(\{x\})
 =\#N(V)=b.
\]
Since $\operatorname{Surp}(V,f)=0$ and $\#E(V)=a+1$, we also have $\operatorname{TNG}_{V,f}=a$.

Let us identify the horizontal shadow of $h_0$.
Claim~\ref{c:node-debt-first-zero} of Lemma~\ref{lem:node-debt-sum} gives
\[
 \operatorname{NGrDiff}_{V,f}(h_0,V_j)<0
 \quad(1\le j<a+b+1),\qquad
 \operatorname{NGrDiff}_{V,f}(h_0,v_*)=0.
\]
For $j\in[a+b+1]_{\Z}$, we have $g(V_j)=f(V_j)$.
Thus, $t_{\mathcal D,g}(h_0)=v_*$.

Let us prove that this h-shadow subpath is maximal.
Let $x\in E(V)$ satisfy $f(x)>0$.
Claim~\ref{c:node-contained} of Lemma~\ref{lem:node-structure} gives a node at $x$ contained in $V$.
Claim~\ref{c:node-debt-first-zero} of Lemma~\ref{lem:node-debt-sum} identifies the endpoint of this node with $t_{\mathcal D,g}(x)$.
Thus, every h-shadow subpath of $g$ is contained in $V=\mathcal D\langle h_0,v_*\rangle$.

Let us verify compatibility and tightness.
Claim~\ref{c:shadow-zero-extension} of Lemma~\ref{lem:compatibility}, applied to this maximal h-shadow subpath, makes $g$ compatible.
Its total grades are $a$ and $b$, and
\[
 N_g(\mathcal D)\subseteq N(V)
 =\operatorname{Sh}_{\mathcal D,g}(h_0),\qquad
 h_0\in E_g(\mathcal D).
\]
Definition~\ref{def:dyck-grading}(\ref{d:tight-grading})(b) gives
$g\in\operatorname{TCG}(\mu_1,\mu_2,a,b)$.

Let us verify that the two maps are inverse.
A tight grading is zero outside $V$ by Claim~\ref{c:frame-unique-shadow} of Lemma~\ref{lem:frame-shadow-position}.
Thus, extending its restriction by zero recovers that grading.
Conversely, restricting the zero extension of $f$ to $V$ gives $f$.
\end{proof}

\Needspace{14\baselineskip}
\begin{proposition}[Counting colored gradings on a run]\label{prop:frame}
Let $b>a\ge2$ be coprime integers.
Then, we have the following.

\begin{enumerate}
\item\label{c:frame-runs} The total vertical length of the first $a$ simple runs of $\mathcal U(a,b)$ is
\[
 \sum_{i=1}^{a}\Lambda_i(a,b)=b-1.
\]
\item\label{c:frame-dictionary} We have
\begin{equation}\label{eq:dictionary}
 C_0(\mathcal{U}(a,b))=mE_m(a,b).
\end{equation}
\end{enumerate}
\end{proposition}
\begin{proof}
\prf{Claim~\ref{c:frame-runs}} Apply Claim~\ref{c:frame-initial-path} of Lemma~\ref{lem:frame-shadow-position}.

\prf{Claim~\ref{c:frame-dictionary}} Let us apply the finite-enumeration formula.
Let $\mu_j=\mu_j(a,b)$ for $j\in[2]_{\Z}$, $\mathcal D=\mathcal D_{(\mu_1,\mu_2)}$, and $V=\mathcal U(a,b)$.
For incoming polynomials $(1+x)^m$, the finite enumeration in
\cite[Theorems 3.2 and 4.10]{B} and \cite[Theorem 12.1]{RW} gives
\[
 \sum_{g\in\operatorname{TCG}(\mu_1,\mu_2,a,b)}
       \prod_{w\in\operatorname{NE}(\mathcal D)}\binom{m}{g(w)}
 =mE_m(a,b).
\]

Let us count the color functions after restriction.
Lemma~\ref{lem:frame-grading-bijection} identifies the tight gradings with the connected gradings of $V$ having surplus zero.
For $g\in\operatorname{TCG}(\mu_1,\mu_2,a,b)$, Lemma~\ref{lem:frame-grading-bijection} gives $g(w)=0$ for $w\in\operatorname{NE}(\mathcal D)\setminus\operatorname{NE}(V)$. Since $\binom m0=1$,
\[
 \prod_{w\in\operatorname{NE}(\mathcal D)}\binom m{g(w)}
 =\prod_{w\in\operatorname{NE}(V)}\binom m{(g|_V)(w)}.
\]
Equation~\eqref{eq:colored-weight} therefore gives~\eqref{eq:dictionary}.
\end{proof}

\Needspace{10\baselineskip}
\begin{example}\label{ex:small-frame}
We have $\mu_1(2,5)=3$, $\mu_2(2,5)=7$, and
\[
 \mathcal U(2,5)=R_2^2\ccn R_1=E\ccn N^2\ccn E\ccn N^2\ccn E\ccn N.
\]
The run $\mathcal U(2,5)$ passes through $(2,4)$ and ends at $(3,5)$.
Since $4>\frac53\cdot2$, $\mathcal U(2,5)$ does not stay below the line segment
joining its endpoints.

For a connected grading $g$ of surplus zero, we have $g(n_1)=g(n_2)=0$.
The other grades have the following possibilities:
\[
\begin{array}{c|l}
 (g(e_1),g(e_2),g(e_3))&(g(n_3),g(n_4),g(n_5))\\\hline
 (3,2,0)&(0,0,2)\\
 (4,1,0)&(0,0,2),\ (0,1,1)\\
 (5,0,0)&(0,0,2),\ (1,0,1),\ (0,1,1).
\end{array}
\]
Therefore,
\[
\begin{aligned}
 C_0(R_2^2\ccn R_1)
 &=\binom m3\binom m2^2
   +m\binom m4\left(\binom m2+m^2\right)\\
 &\quad+\binom m5\left(\binom m2+2m^2\right).
\end{aligned}
\]
For $m=3$, only the first row contributes. Thus,
$C_0(R_2^2\ccn R_1)=9$ and $E_3(2,5)=3$.

\begin{figure}[H]
\centering
\includegraphics[height=66mm]{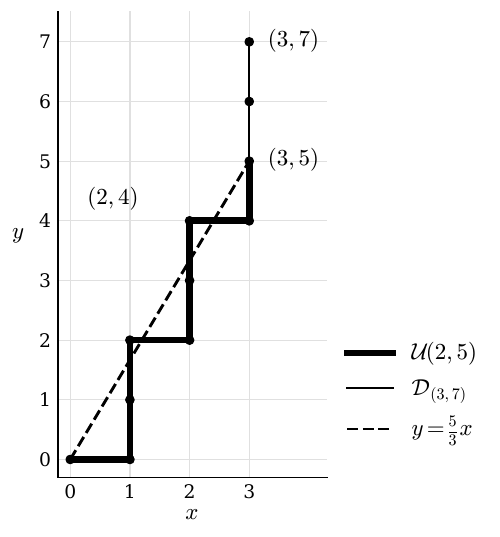}
\caption{The thick line shows $\mathcal U(2,5)$. The thin line continues $\mathcal U(2,5)$ to
$\mathcal D_{(3,7)}$. The dashed segment joins $(0,0)$ to $(3,5)$.}
\label{fig:example-run}
\end{figure}
\end{example}

%% file: sections/shell-paths-and-gradings.tex
\section{Shell paths and gradings}\label{sec:shells}
\subsection{Shell constructions and periodic growth}
\begin{definition}\label{def:shell}
Let $j\in\Z_{\ge0}$ and $\epsilon\in\Z_{\ge0}$ satisfy $j+1+\epsilon\le m$.
For $s\in[j]_{\Z}$, let $X_s$ be a tNE path.
Let $u\in[m-1]_{\Z}$ and $v\in[m-1]_{\Z}$.
Let $X=X_1\ccn \cdots\ccn X_j$.
For $j=0$, let $X=\emptyset$.
\begin{enumerate}[font=\itshape]
\item We define the \textit{shell path}
\[
 Y=\operatorname{Shell}(X_1,\ldots,X_j;u,v,\epsilon)=R_u\ccn X\ccn R_v.
\]
For $j=0$, use the notation $\operatorname{Shell}(;u,v,\epsilon)$.
\item For $s\in[j]_{\Z}$, let $g_s\in\operatorname{CtdG}_m(X_s)$.
Let $f$ be the concatenated grading on $X$.
For $j=0$, let $f=\emptyset$.
We define the \textit{shell grading} $g=\operatorname{GShell}_{f}(X_1,\ldots,X_j;u,v,\epsilon)$ on $Y$ by
\begin{equation}\label{eq:shell}
 g(Y_k)=
 \begin{cases}
 u+1,&k=1,\\
 f(X_{k-u-1}),&k\in[u+2,u+1+|X|]_{\Z},\\
 v-1,&k=u+|X|+2,\\
 j+1+\epsilon,&k=|Y|,\\
 0,&\text{otherwise}
 \end{cases}
 \qquad(k\in[|Y|]_{\Z}).
\end{equation}
For $j=0$, use the notation $\operatorname{GShell}_{\emptyset}(;u,v,\epsilon)$.
\item We call a shell path $\operatorname{Shell}(X_1,\ldots,X_j;u,v,\epsilon)$ \textit{unary} if $j=1$.
\item We call a shell path $\operatorname{Shell}(X_1,\ldots,X_j;u,v,\epsilon)$ \textit{binary} if $j=2$.
\item We call a shell grading $\operatorname{GShell}_{f}(X_1,\ldots,X_j;u,v,\epsilon)$ \textit{preserving} if $\epsilon=0$.
\item We call a shell grading $\operatorname{GShell}_{f}(X_1,\ldots,X_j;u,v,\epsilon)$ \textit{draining} if $j=\epsilon=1$.
\end{enumerate}
\end{definition}

\begin{definition}\label{def:recolor}
Let $U$ and $V$ be tNE paths. Let $S\subseteq\operatorname{NE}(U)$ and
$\iota\colon S\hookrightarrow\operatorname{NE}(V)$.
Let $g\in G_m(U)$, $\gamma\in\operatorname{Col}(U,g,m)$, and $h\in G_m(V)$.
We define the \textit{recoloring operator}
\[
 \operatorname{Recol}(U,V,\iota,h;g,\gamma)=\gamma'\in\operatorname{Col}(V,h,m)
\]
such that $\gamma'$ satisfies
\begin{equation}\label{eq:color-update}
 \gamma'(e)=
 \begin{cases}
 \gamma(x),&e=\iota(x),\ x\in S,\ h(e)=g(x),\\
 1,&\text{otherwise}
 \end{cases}
 \qquad(e\in\operatorname{NE}(V)).
\end{equation}
\end{definition}
\Needspace{7\baselineskip}
\begin{lemma}\label{lem:shell}
Let $j\in\Z_{\ge0}$ and $\epsilon\in\Z_{\ge0}$ satisfy $j+1+\epsilon\le m$. Let $u\in[m-1]_{\Z}$ and $v\in[m-1]_{\Z}$.
For $s\in[j]_{\Z}$, let $X_s$ be a tNE path and let $g_s\in\operatorname{CtdG}_m(X_s)$.
Let $X=X_1\ccn\cdots\ccn X_j$ and let $f$ be the concatenated grading on $X$.
For $j=0$, let $X=\emptyset$ and $f=\emptyset$.
Suppose that one of the following conditions holds.
\begin{hypotheses}
\item\label{hyp:shell-preserving} $j\in[0,2]_{\Z}$ and $\epsilon=0$.
\item\label{hyp:shell-draining} $j=\epsilon=1$ and $\operatorname{Surp}(X_1,g_1)\ge1$.
\end{hypotheses}
Let $Y=\operatorname{Shell}(X_1,\ldots,X_j;u,v,\epsilon)$ and $g=\operatorname{GShell}_{f}(X_1,\ldots,X_j;u,v,\epsilon)$.
Put
\[
 r_i=u+1+\sum_{s=1}^i|X_s|\qquad(i\in[0,j]_{\Z}).
\]
The empty sum is $0$.
Then, we have the following.

\begin{enumerate}
\item\label{c:shell-bounded} $g\in G_m(Y)$.
\item\label{c:shell-connected} $g\in\operatorname{CtdG}_m(Y)$.
\item\label{c:shell-residual} We have
\[
 \operatorname{Surp}(Y,g)=\sum_{i=1}^j\operatorname{Surp}(X_i,g_i)-\epsilon.
\]
\item\label{c:shell-recolor} Fix $X_1,\ldots,X_j$, $u$, $v$, and $\epsilon$.
Let $c_s\in\Z_{\ge0}$ for $s\in[j]_{\Z}$.
Assume $\sum_{s=1}^j c_s\ge\epsilon$.
For $(g_s,\gamma_s)\in\mathcal C_{c_s}(X_s)$, let $\eta$ be the concatenated coloring.
If $j=0$, let $\eta=\emptyset$.
We define
\[
 \iota\colon\operatorname{NE}(X)\hookrightarrow\operatorname{NE}(Y),
 \qquad\iota(X_k)=Y_{u+1+k}\quad(k\in[|X|]_{\Z}).
\]
Let
\begin{equation}\label{eq:shell-recolor}
 \gamma=\operatorname{Recol}(X,Y,\iota,g;f,\eta).
\end{equation}
Then, the map
\[
 \Phi_{\mathrm{sh}}\colon
 \prod_{s\in[j]_{\Z}}\mathcal C_{c_s}(X_s)
 \longrightarrow\mathcal C_{\sum_{s=1}^j c_s-\epsilon}(Y),
 \qquad \Phi_{\mathrm{sh}}\bigl((g_s,\gamma_s)_{s\in[j]_{\Z}}\bigr)=(g,\gamma)
\]
is injective. Also, we have
\[
 (g(\iota(X_k)),\gamma(\iota(X_k)))=(f(X_k),\eta(X_k))
 \qquad(k\in[|X|]_{\Z}).
\]
\item\label{c:shell-copy-precursor} For $i\in[j]_{\Z}$, $k\in[|X_i|]_{\Z}$, and $e=(X_i)_k$, define the copy of $e$ in $Y$ by
\begin{equation}\label{eq:shell-edge-copy}
 \widetilde e=Y_{r_{i-1}+k}.
\end{equation}
Then, we have
\begin{equation}\label{eq:shell-copy-precursor}
 P_{Y,g}(\widetilde w)=\widetilde{P_{X_i,g_i}(w)}
 \qquad(i\in[j]_{\Z},\ w\in N(X_i)).
\end{equation}
\item\label{c:shell-boundary-debts} For $i\in[0,j]_{\Z}$ and $e\in E(Y[0,r_i])$, we have
\[
 \operatorname{Debt}_{Y,g}(e,r_i)=
 \begin{cases}
 1,&e=Y_1,\\
 0,&e\ne Y_1.
 \end{cases}
\]
\item\label{c:shell-boundary-balances} We have
\[
 \operatorname{bBal}_{Y,g}(r_i)
 =1+i+\sum_{s=1}^i\operatorname{Surp}(X_s,g_s)
 \qquad(i\in[0,j]_{\Z}).
\]
\end{enumerate}
\end{lemma}
\begin{proof}
\prf{Claim~\ref{c:shell-bounded}} Let us verify that every shell grade lies in the prescribed range.
The bounds $u\in[m-1]_{\Z}$, $v\in[m-1]_{\Z}$, $j+1+\epsilon\le m$, and $g_i\in\operatorname{CtdG}_m(X_i)$ for $i\in[j]_{\Z}$ in
Definition~\ref{def:shell} give
\[
\begin{aligned}
 g_i(e)&\le m\qquad(i\in[j]_{\Z},\ e\in\operatorname{NE}(X_i)),\\
 \max\{u+1,v-1,j+1+\epsilon\}&\le m.
\end{aligned}
\]
Every edge of $Y$ inherits a grade from some $g_i$ or has grade
$0$, $u+1$, $v-1$, or $j+1+\epsilon$. Hence, $g\in G_m(Y)$.

\prf{Claim~\ref{c:shell-connected}} Let us verify the three conditions for connectedness in Definition~\ref{def:connected-grading}(1).
Let $y=\operatorname{Start}_{\mathrm{NE}}(Y)$ and $z=\operatorname{End}_{\mathrm{NE}}(Y)$.
We construct $P_{Y,g}$ in edge order, using Definition~\ref{def:stack}(1).
The grades $g(Y_1)=u+1$ and $g(Y_{k+1})=0$ for $k\in[u]_{\Z}$ give
\[
\begin{aligned}
 P_{Y,g}(Y_{k+1})&=y\qquad(k\in[u]_{\Z}),\\
 \operatorname{Debt}_{Y,g}(y,r_0)&=g(y)-u=(u+1)-u=1.
\end{aligned}
\]

Let us determine the h-precursors in each copy of $X_i$.
The subpath $Y[r_{i-1},r_i]$ is the translated copy of $X_i$ with the grading translated from $g_i$.
At each vertical edge $\widetilde w$ of this copy, the translated h-precursor from $X_i$ has positive debt and follows every horizontal edge before the copy.
Definition~\ref{def:stack}(1) therefore selects $\widetilde{P_{X_i,g_i}(w)}$ as the h-precursor of $\widetilde w$.

Let us compute the debt of $y$ after each copy of $X_i$.
Since $g_i$ is connected, Claim~\ref{c:node-contained} of Lemma~\ref{lem:node-structure}
shows that every horizontal edge of $X_i$ has debt zero at $|X_i|$.
Thus,
\[
 \operatorname{Debt}_{Y,g}(y,r_i)=1
 \qquad(i\in[j]_{\Z}).
\]

Let us compute the change in boundary balance across the copies.
The boundary balances satisfy
\[
 \operatorname{bBal}_{Y,g}(r_i)-\operatorname{bBal}_{Y,g}(r_{i-1})
 =|X_i|_E-\operatorname{TNG}_{X_i,g_i}
 =\operatorname{Surp}(X_i,g_i)+1
 \qquad(i\in[j]_{\Z}),
\]
and hence
\[
 \operatorname{bBal}_{Y,g}(r_j)-\operatorname{bBal}_{Y,g}(r_0)
 =\sum_{i=1}^j\operatorname{Surp}(X_i,g_i)+j.
\]

Let us prove $\operatorname{Node}_{Y,g}(y)=Y$.
The grades $g(Y_{r_j+1})=v-1$ and $g(Y_{r_j+1+k})=0$ for $k\in[v-1]_{\Z}$ give
\[
\begin{aligned}
 P_{Y,g}(Y_{r_j+1+k})&=Y_{r_j+1}\qquad(k\in[v-1]_{\Z}),\\
 \operatorname{bBal}_{Y,g}(r_j+v)-\operatorname{bBal}_{Y,g}(r_j)
 &=1-\sum_{k=1}^{v-1}g(Y_{r_j+1+k})=1,\\
 \operatorname{Debt}_{Y,g}(Y_{r_j+1},r_j+v)&=(v-1)-(v-1)=0.
\end{aligned}
\]
Immediately before $z$, the edge $y$ has debt one and every other horizontal edge has debt zero.
Hence,
\[
 P_{Y,g}(z)=y,\qquad
 \operatorname{Debt}_{Y,g}(y,|Y|)=0.
\]
This constructs $P_{Y,g}$ and gives $\operatorname{Node}_{Y,g}(y)=Y$.

Let us verify the node balance condition. For $i\in[j]_{\Z}$, $w\in N(X_i)$, and
$a=P_{X_i,g_i}(w)$, the identity $P_{Y,g}(\widetilde w)=\widetilde a$ and the connectedness of $g_i$ give
\[
\begin{aligned}
 &\operatorname{eBal}_{Y,g}(\widetilde w)-\operatorname{eBal}_{Y,g}(\widetilde a)\\
 &\qquad=\operatorname{eBal}_{X_i,g_i}(w)-\operatorname{eBal}_{X_i,g_i}(a)\ge0.
\end{aligned}
\]
For $Y_{k+1}$ with $k\in[u]_{\Z}$ and $Y_{r_j+1+k}$ with $k\in[v-1]_{\Z}$, we have
\[
\begin{aligned}
 \operatorname{eBal}_{Y,g}(Y_{k+1})&=1=\operatorname{eBal}_{Y,g}(y)
 &&(k\in[u]_{\Z}),\\
 \operatorname{eBal}_{Y,g}(Y_{r_j+1+k})&=\operatorname{eBal}_{Y,g}(Y_{r_j+1})
 &&(k\in[v-1]_{\Z}).
\end{aligned}
\]
By Lemma~\ref{lem:residual-nonnegative}, each $\operatorname{Surp}(X_i,g_i)$ is nonnegative.
Under condition~\ref{hyp:shell-preserving} or~\ref{hyp:shell-draining}, the edge balance at $z$ satisfies
\[
 \operatorname{eBal}_{Y,g}(z)=\operatorname{EGrDiff}_{Y,g}(y,z)
 =1+\sum_{i=1}^j\operatorname{Surp}(X_i,g_i)-\epsilon
 \ge1=\operatorname{eBal}_{Y,g}(y).
\]
Therefore, $g$ is connected. By Claim~\ref{c:shell-bounded},
$g\in\operatorname{CtdG}_m(Y)$.

\prf{Claim~\ref{c:shell-residual}}
By~\eqref{eq:residual},
\[
 \operatorname{Surp}(Y,g)=\operatorname{EGrDiff}_{Y,g}(y,z)-1
 =\sum_{i=1}^j\operatorname{Surp}(X_i,g_i)-\epsilon.
\]

\prf{Claim~\ref{c:shell-recolor}} Let us define the shell map and prove injectivity.
Claims~\ref{c:shell-connected} and~\ref{c:shell-residual} give
$(g,\gamma)\in\mathcal C_{\sum_{s=1}^j c_s-\epsilon}(Y)$.
Equation~\eqref{eq:shell} gives $g(\iota(X_k))=f(X_k)$.
Definition~\ref{def:recolor} then gives $\gamma(\iota(X_k))=\eta(X_k)$.
Restriction to each copy of $X_s$ recovers $(g_s,\gamma_s)$, proving injectivity.

\prf{Claim~\ref{c:shell-copy-precursor}} Let us compare the h-precursors in $X_i$ and its copy in $Y$.
Fix $i\in[j]_{\Z}$.
The copy map~\eqref{eq:shell-edge-copy} preserves edge order and grades.
For $w=(X_i)_k\in N(X_i)$, Definition~\ref{def:stack}(1) gives $\operatorname{Debt}_{X_i,g_i}(P_{X_i,g_i}(w),k-1)>0$.
Induction on $k$ preserves the debts of the horizontal edges in $X_i[0,k-1]$ under the copy map.
Every horizontal edge before the copy of $X_i$ precedes $\widetilde{P_{X_i,g_i}(w)}$.
Thus, Definition~\ref{def:stack}(1) selects $\widetilde{P_{X_i,g_i}(w)}$ as the h-precursor of $\widetilde w$.
Hence,~\eqref{eq:shell-copy-precursor} holds.

\prf{Claim~\ref{c:shell-boundary-debts}} Let us compute the debts at $r_i$.
The first $u$ vertical edges have h-precursor $Y_1$.
Hence,
\[
 \operatorname{Debt}_{Y,g}(Y_1,r_0)=g(Y_1)-u=1.
\]
By Claim~\ref{c:shell-copy-precursor}, every vertical edge in a copy of $X_s$ has its h-precursor in that copy.
Claim~\ref{c:node-contained} of Lemma~\ref{lem:node-structure} gives debt zero at $|X_s|$ for every horizontal edge of $X_s$.
Thus, at $r_i$, the edge $Y_1$ has debt $1$ and every horizontal edge in the copies of $X_1,\ldots,X_i$ has debt $0$.

\prf{Claim~\ref{c:shell-boundary-balances}} The definition of $\operatorname{bBal}$ and~\eqref{eq:residual} give
\[
 \operatorname{bBal}_{Y,g}(r_i)
 =1+\sum_{s=1}^i\bigl(|X_s|_E-\operatorname{TNG}_{X_s,g_s}\bigr)
 =1+i+\sum_{s=1}^i\operatorname{Surp}(X_s,g_s)
 \qquad(i\in[0,j]_{\Z}).
\]
\end{proof}

\begin{definition}[Fixed simple-run grading]\label{def:simple-run-grading}
Let $j\in[m]_{\Z}$.
Define $g_j^{\mathrm{run}}\in G_m(R_j)$ and $\gamma_j^{\mathrm{run}}\in\operatorname{Col}(R_j,g_j^{\mathrm{run}},m)$ by
\[
 g_j^{\mathrm{run}}(\operatorname{Start}_{\mathrm{NE}}(R_j))=j,\qquad
 g_j^{\mathrm{run}}(v)=0\quad(v\in N(R_j)),\qquad
 \gamma_j^{\mathrm{run}}(e)=1\quad(e\in\operatorname{NE}(R_j)).
\]
\end{definition}

\begin{lemma}\label{lem:word-bounds}
Let $U$ be a run with $\nu$ simple runs. Let $L$ be an integer.
Assume the following.
\begin{hypotheses}
\item\label{hyp:grading-range} $L\in[m-1]_{\Z}$.
\item\label{hyp:grading-lengths} Each simple run of $U$ has vertical length at most $L$.
\end{hypotheses}
Then, we have
\[
 m^{\ceil{\frac{\nu}{2}}}\le C_0(U)\le 2^{m(L+1)\nu}.
\]
\end{lemma}
\begin{proof}
Let us prove the lower bound first for one and two simple runs.
By Lemma~\ref{lem:word-reindexing}, we translate $U$ so that $\operatorname{Start}_{\mathrm{pt}}(U)=(0,0)$.

First, let us prove the bound for $\nu=1$.
Let $U=R_j$ and let $g=g_j^{\mathrm{run}}$.
This grading is connected and satisfies $\operatorname{Surp}(U,g)=0$.
By assumptions~\ref{hyp:grading-range} and~\ref{hyp:grading-lengths},
$1\le j\le L\le m-1$. Thus, the horizontal edge has $\binom mj\ge m$ colors.

Second, let us prove the bound for $\nu=2$.
Let $U=R_u\ccn R_v$ and use
$g=\operatorname{GShell}_{\emptyset}(;u,v,0)$. Claims~\ref{c:shell-connected} and~\ref{c:shell-residual} of Lemma~\ref{lem:shell} give a connected grading of surplus zero.
The edge $\operatorname{End}_{\mathrm{NE}}(U)$ has grade one and $m$ colors.

Let us extend the lower bound by induction on $\nu$.
For $\nu\ge3$, let $U=R_u\ccn X\ccn R_v$ with $|X|_E=\nu-2$.
For each $(g,\gamma)\in\mathcal C_0(X)$, Claims~\ref{c:shell-connected}--\ref{c:shell-recolor} of Lemma~\ref{lem:shell}, followed by Lemma~\ref{lem:word-reindexing}, give a colored grading $(h,\delta)\in\mathcal C_0(U)$ such that
\begin{equation}\label{eq:word-bound-restriction}
 h(U_{u+1+k})=g(X_k),\qquad
 \delta(U_{u+1+k})=\gamma(X_k)\quad(k\in[|X|]_{\Z}).
\end{equation}
The edge $\operatorname{End}_{\mathrm{NE}}(U)$ has grade two and $\binom m2\ge m$ colors.
Equation~\eqref{eq:word-bound-restriction} recovers $(g,\gamma)$ from $(h,\delta)$.
Induction gives
\[
 C_0(U)\ge\binom m2 C_0(X)
 \ge m\,m^{\ceil{\frac{\nu-2}{2}}}=m^{\ceil{\frac{\nu}{2}}}.
\]

Let us prove the upper bound by summing the choices of grade and color on each edge.
By assumption~\ref{hyp:grading-lengths}, there are at most $(L+1)\nu$ edges. Hence,
\[
 C_0(U)\le\prod_{e\in\operatorname{NE}(U)}\sum_{r=0}^m\binom mr
 =2^{m|U|}\le2^{m(L+1)\nu}.
\]
\end{proof}

\begin{lemma}\label{lem:outer-shell-extension}
Let $X$ be a tNE path. Let $A$ and $B$ be runs with $P\ge1$ and $S\ge1$
simple runs, respectively.
Assume the following.
\begin{hypotheses}
\item\label{hyp:outer-shell-lengths} Each simple run of $A$ and $B$ has vertical length at most $m-1$.
\item\label{hyp:outer-shell-counts} $\max(P,S)\le2\min(P,S)$.
\end{hypotheses}
Then, for every integer $d\ge0$, there is an injective map
\[
 i_d\colon\mathcal C_d(X)\longrightarrow\mathcal C_d(A\ccn X\ccn B)
\]
with the following property. For each $(g,\gamma)\in\mathcal C_d(X)$ and $(h,\delta)=i_d(g,\gamma)$, we have
\begin{equation}\label{eq:outer-extension-restriction}
 h\bigl((A\ccn X\ccn B)_{|A|+k}\bigr)=g(X_k),\qquad
 \delta\bigl((A\ccn X\ccn B)_{|A|+k}\bigr)=\gamma(X_k)
 \quad(k\in[|X|]_{\Z}).
\end{equation}
\end{lemma}
\begin{proof}
Let us construct shell gradings that preserve the surplus $d$.
Let $(g,\gamma)\in\mathcal C_d(X)$. Let $u\in[m-1]_{\Z}$, $v\in[m-1]_{\Z}$, and $t\in[m-1]_{\Z}$.
Give $R_t$ the colored grading $(g_t^{\mathrm{run}},\gamma_t^{\mathrm{run}})$.
The shell paths
\begin{equation}\label{eq:outer-shells}
\begin{gathered}
 \operatorname{Shell}(X;u,v,0),\\
 \operatorname{Shell}(R_t,X;u,v,0),\\
 \operatorname{Shell}(X,R_t;u,v,0)
\end{gathered}
\end{equation}
add $(1,1)$, $(2,1)$, and $(1,2)$ simple runs before and after $X$, respectively.
For each shell path $Y$, let $h$ be its shell grading from~\eqref{eq:shell}.
Claims~\ref{c:shell-connected} and~\ref{c:shell-residual} of Lemma~\ref{lem:shell} give $h\in\operatorname{CtdG}_m(Y)$ and
\[
 \operatorname{Surp}(Y,h)=\operatorname{Surp}(X,g)=d.
\]

Let us choose the numbers of shell paths to add exactly $P$ simple runs before $X$ and $S$ after $X$.
By assumption~\ref{hyp:outer-shell-counts}, the following decomposition has nonnegative integer coefficients:
\begin{equation}\label{eq:outer-shell-counts}
 (P,S)=
 \begin{cases}
 (P-S)(2,1)+(2S-P)(1,1),& S\le P\le2S,\\
 (S-P)(1,2)+(2P-S)(1,1),& P\le S\le2P.
 \end{cases}
\end{equation}

Let us match the simple runs of $A$ and $B$ in order.
Let $a_1,\ldots,a_P$ and $b_1,\ldots,b_S$ be the ordered vertical lengths
of the simple runs of $A$ and $B$.
By assumption~\ref{hyp:outer-shell-lengths},
\[
 a_i\in[m-1]_{\Z}\quad(i\in[P]_{\Z}),\qquad
 b_i\in[m-1]_{\Z}\quad(i\in[S]_{\Z}).
\]

Let us specify one step of the construction.
Let $p\in[0,P]_{\Z}$ and $q\in[0,S]_{\Z}$.
Suppose that $Y$ is a translated copy of
$R_{a_{p+1}}\ccn \cdots\ccn R_{a_P}\ccn X\ccn R_{b_1}\ccn \cdots\ccn R_{b_q}$.
For the three shell operators in~\eqref{eq:outer-shells} applied to $Y$, use
\[
 (u,v)=(a_p,b_{q+1}),\qquad
 (u,t,v)=(a_{p-1},a_p,b_{q+1}),\qquad
 (u,t,v)=(a_p,b_{q+1},b_{q+2}),
\]
respectively.
The first choice requires $p\ge1$ and $q<S$.
The second requires $p\ge2$ and $q<S$.
The third requires $p\ge1$ and $q+2\le S$.

Let us complete the construction.
Construct the binary shell paths before the unary shell paths.
Starting at $(p,q)=(P,0)$, the numbers of shell paths
in~\eqref{eq:outer-shell-counts} give $(p,q)=(0,S)$.
At each step, define the shell grading by~\eqref{eq:shell} and its coloring by~\eqref{eq:shell-recolor}.
At $(p,q)=(0,S)$, $Y$ is a translated copy of $A\ccn X\ccn B$.
Lemma~\ref{lem:word-reindexing} transports the colored grading on $Y$ to $A\ccn X\ccn B$.

Let us verify~\eqref{eq:outer-extension-restriction} and prove that $i_d$ is injective.
By Claim~\ref{c:shell-recolor} of Lemma~\ref{lem:shell}, the map $i_d$ copies the grade and color of $X_k$ to $(A\ccn X\ccn B)_{|A|+k}$ for every $k\in[|X|]_{\Z}$.
This proves~\eqref{eq:outer-extension-restriction}.
For $(g_1,\gamma_1)\in\mathcal C_d(X)$ and $(g_2,\gamma_2)\in\mathcal C_d(X)$,
equation~\eqref{eq:outer-extension-restriction} gives
\[
 i_d(g_1,\gamma_1)=i_d(g_2,\gamma_2)
 \ \Longrightarrow\ 
 \begin{cases}
 g_1(X_k)=g_2(X_k),\\
 \gamma_1(X_k)=\gamma_2(X_k)
 \end{cases}
 \quad(k\in[|X|]_{\Z}).
\]
Hence, $(g_1,\gamma_1)=(g_2,\gamma_2)$ and $i_d$ is injective.
\end{proof}

\begin{lemma}[A shifted supermultiplicative limit]\label{lem:shifted-growth}
Let $(A_n)_{n\in\Z_{>0}}$ be a sequence of real numbers and let $C\ge0$.
Assume
\[
 A_n\ge1,\qquad \log A_n\le Cn\quad(n\in\Z_{>0}),\qquad
 A_{i+j+2}\ge A_iA_j\quad(i,j\in\Z_{>0}).
\]
Then, the following statements hold.
\begin{enumerate}
\item\label{c:shifted-growth-limit} The limit exists and satisfies
\[
 L:=\lim_{n\to\infty}\frac{\log A_n}{n}
   =\sup_{n\ge1}\frac{\log A_n}{n+2}\in[0,C].
\]
\item\label{c:shifted-growth-bound} For $n\in\Z_{>0}$, we have $A_n\le e^{(n+2)L}$.
\end{enumerate}
\end{lemma}
\begin{proof}
\prf{Claim~\ref{c:shifted-growth-limit}} Let us write $\nu=n+2$.
Let $B_\nu=\log A_{\nu-2}$ for $\nu\in\Z_{\ge3}$ and let
\[
 L_0=\sup_{\nu\ge3}\frac{B_\nu}{\nu}.
\]
We have $0\le B_\nu\le C\nu$ for $\nu\in\Z_{\ge3}$ and
$B_{u+v}\ge B_u+B_v$ for $u\in\Z_{\ge3}$ and $v\in\Z_{\ge3}$.
Thus, $L_0\in[0,C]$.
Fix $t\in\Z_{\ge3}$.
For $\nu\ge3$, let $\nu=kt+r$ for $k\in\Z_{\ge0}$ and $r\in[3,t+2]_{\Z}$.
Induction on $k$ gives
\[
 B_\nu\ge kB_t+B_r\ge kB_t.
\]
Since $\frac{k}{\nu}\to\frac1t$, we obtain
\[
 \liminf_{\nu\to\infty}\frac{B_\nu}{\nu}\ge\frac{B_t}{t}.
\]
Taking the supremum over $t$ gives a lower bound of $L_0$.
The definition of $L_0$ gives $\frac{B_\nu}{\nu}\le L_0$.
Hence, $\frac{B_\nu}{\nu}\to L_0$ and
\[
 \frac{\log A_n}{n}=\frac{n+2}{n}\frac{B_{n+2}}{n+2}\longrightarrow L_0.
\]
Replacing $\nu$ by $n+2$ in the definition of $L_0$ proves Claim~\ref{c:shifted-growth-limit}.

\prf{Claim~\ref{c:shifted-growth-bound}} Let us bound the error in the finite growth estimate.
Claim~\ref{c:shifted-growth-limit} gives $\log A_n\le(n+2)L$.
Exponentiation gives the assertion.
\end{proof}

\begin{proposition}[Periodic growth]\label{prop:periodic}
Let $W$ be a run with $q$ simple runs. Assume the following.
\begin{hypotheses}
\item\label{hyp:periodic-lengths} Every simple run of $W$ has vertical length at most $m-1$.
\item Let $A_j=C_0(W^j)$ for $j\ge1$.
\end{hypotheses}
Then, we have the following.

\begin{enumerate}
\item\label{c:periodic-supermultiplicativity} For all integers $i\ge1$ and $j\ge1$, we have
\[
 A_{i+j+2}\ge A_iA_j.
\]
\item\label{c:periodic-limit} We have
\begin{equation}\label{eq:periodic-limit}
 \lim_{j\to\infty}\frac{\log A_j}{j}
       =\sup_{j\ge1}\frac{\log A_j}{j+2}
       \in(0,\infty).
\end{equation}
\item\label{c:periodic-bound} For every $j\ge1$, we have
\[
 A_j\le\exp\left((j+2)\lim_{k\to\infty}\frac{\log A_k}{k}\right).
\]
\end{enumerate}
\end{proposition}
\begin{proof}
\prf{Claim~\ref{c:periodic-supermultiplicativity}} Let us construct a colored grading on $W^{i+j+2}$ from colored gradings on $W^i$ and $W^j$.
Let $i\ge1$ and $j\ge1$.
Let $a_1,\ldots,a_q$ be the vertical lengths of the simple runs of $W$.
Let $X=W^i$ and $Y=W^j$.
For $(g_1,\gamma_1)\in\mathcal C_0(X)$ and $(g_2,\gamma_2)\in\mathcal C_0(Y)$, form the shell paths
\[
\begin{aligned}
 Z_1&=\operatorname{Shell}(X,Y;a_q,a_1,0),\\
 Z_{s+1}&=\operatorname{Shell}(Z_s;a_{q-s},a_{s+1},0)
 \quad(s\in[q-1]_{\Z}).
\end{aligned}
\]
Give each $Z_s$ its shell grading~\eqref{eq:shell} and coloring~\eqref{eq:shell-recolor}.
Assumption~\ref{hyp:periodic-lengths} and Claims~\ref{c:shell-connected} and~\ref{c:shell-residual} of Lemma~\ref{lem:shell} give a connected colored grading of surplus zero on $Z_q$.
By Lemma~\ref{lem:word-reindexing}, translate $Z_q$ and its colored grading to
$Z=W\ccn X\ccn Y\ccn W=W^{i+j+2}$. Denote the resulting colored grading by $(h,\delta)$.
This defines the map
\[
 \Phi_{i,j}\colon\mathcal C_0(W^i)\times\mathcal C_0(W^j)
 \longrightarrow\mathcal C_0(W^{i+j+2}),\qquad
 \Phi_{i,j}((g_1,\gamma_1),(g_2,\gamma_2))=(h,\delta).
\]

Let us recover $(g_1,\gamma_1)$ and $(g_2,\gamma_2)$ to prove that $\Phi_{i,j}$ is injective.
Claim~\ref{c:shell-recolor} of Lemma~\ref{lem:shell} gives
\begin{equation}\label{eq:periodic-restrictions}
\begin{aligned}
 h(Z_{|W|+k})&=g_1(X_k),&
 \delta(Z_{|W|+k})&=\gamma_1(X_k)
 &&(k\in[|X|]_{\Z}),\\
 h(Z_{|W|+|X|+k})&=g_2(Y_k),&
 \delta(Z_{|W|+|X|+k})&=\gamma_2(Y_k)
 &&(k\in[|Y|]_{\Z}).
\end{aligned}
\end{equation}
Equation~\eqref{eq:periodic-restrictions} determines both $(g_1,\gamma_1)$ and $(g_2,\gamma_2)$ from $(h,\delta)$.
Hence, $\Phi_{i,j}$ is injective and $A_{i+j+2}\ge A_iA_j$.

\prf{Claim~\ref{c:periodic-limit}} Let us prove~\eqref{eq:periodic-limit}.
Lemma~\ref{lem:word-bounds} gives $A_j\ge1$ and $\log A_j\le m^2qj\log2$ for $j\ge1$.
Claim~\ref{c:periodic-supermultiplicativity} gives $A_{i+j+2}\ge A_iA_j$ for $i\in\Z_{>0}$ and $j\in\Z_{>0}$.
Claim~\ref{c:shifted-growth-limit} of Lemma~\ref{lem:shifted-growth} therefore gives~\eqref{eq:periodic-limit}.
The lower bound in Lemma~\ref{lem:word-bounds} also gives
\[
 \lim_{j\to\infty}\frac{\log A_j}{j}\ge\frac q2\log m>0.
\]

\prf{Claim~\ref{c:periodic-bound}}
Claim~\ref{c:shifted-growth-bound} of Lemma~\ref{lem:shifted-growth} proves Claim~\ref{c:periodic-bound}.
\end{proof}

\begin{definition}\label{def:periodic-growth-rate}
Let $W$ be a run whose simple runs have vertical lengths at most $m-1$.
We define
\[
 \ell(W)=\lim_{j\to\infty}\frac{\log C_0(W^j)}{j}.
\]
\end{definition}

\subsection{A uniform upper bound for finite dimension pairs}
\begin{lemma}[Reduction of dimension pairs]\label{lem:dimension-reduction}
Let $m\in\Z_{\ge3}$ and let $a\in\Z_{>0}$ and $b\in\Z_{>0}$ be coprime.
Let $Q=mab-a^2-b^2>0$.
Then, there exist coprime $u\in\Z_{>0}$ and $v\in\Z_{>0}$ such that
\[
\begin{gathered}
 u\le v\le\frac{mu}{2},\qquad muv-u^2-v^2=Q,\qquad E_m(u,v)=E_m(a,b),\\
 u+v\le\sqrt{\frac{m+2}{m-2}}\sqrt Q.
\end{gathered}
\]
\end{lemma}
\begin{proof}
Let us obtain $(u,v)$ by decreasing $x+y$.
Let
\[
 (x,y)=(\min(a,b),\max(a,b)).
\]
If $2y>mx$, replace $(x,y)$ by
\[
 (\min(x,mx-y),\max(x,mx-y)).
\]
Claims~\ref{c:moduli-symmetries-identities}, \ref{c:moduli-coprimality}, and~\ref{c:moduli-quadratic}
of Lemma~\ref{lem:moduli-symmetries} show that each replacement preserves
the Euler characteristic, coprimality, and $Q$.
The inequality $mxy-x^2-y^2=Q>0$ implies $y<mx$.
Thus, both new entries are positive.
Their sum changes by $mx-2y<0$.
A strictly decreasing sequence of positive integers is finite. Thus, the replacements end at a pair
$u\le v\le\frac{mu}{2}$.

Let us bound $u+v$.
We have $mv-2u\ge0$ and $mu-2v\ge0$.
Therefore,
\[
 (m+2)Q-(m-2)(u+v)^2=(mv-2u)(mu-2v)\ge0.
\]
Dividing by $m-2>0$ and taking square roots gives $u+v\le\sqrt{\frac{m+2}{m-2}}\sqrt Q$.
\end{proof}

\begin{proposition}\label{prop:finite-endpoint}
For every $(a,b)\in\cA_m$, let $\rho(a,b)=\frac{b}{a}$.
Then, we have
\begin{equation}\label{eq:finite-endpoint}
\begin{split}
 0\le\frac{\log E_m(a,b)}a
 &\le m\log2\sqrt{\frac{m+2}{m-2}}
             \sqrt{\max(m\rho(a,b)-1-\rho(a,b)^2,0)}\\
 &\hspace{12mm}+\frac{\log\left(\frac{2^m}{m}\right)}a.
\end{split}
\end{equation}
\end{proposition}
\begin{proof}
Let us separate the cases according to the sign of $Q=mab-a^2-b^2$.
First, suppose $Q>0$.
Lemma~\ref{lem:dimension-reduction} gives a coprime pair $(u,v)$ satisfying
\[
 u\le v\le\frac{mu}{2},\qquad E_m(u,v)=E_m(a,b),\qquad
 u+v\le\sqrt{\frac{m+2}{m-2}}\sqrt Q.
\]

Let us bound the Euler characteristic at the reduced pair $(u,v)$.
If $u\ge2$, we have $v>u$ because $u$ and $v$ are coprime. Then, Claim~\ref{c:frame-dictionary} of Proposition~\ref{prop:frame} applies.
Equation~\eqref{eq:ambient} gives $\Gamma(u,v)<\frac vu\le\frac m2$.
The floor-difference formula in Definition~\ref{def:frame} gives
\[
 \Lambda_i(u,v)\le\ceil{\Gamma(u,v)}\le\ceil{\frac m2}\le m-1
 \qquad(i\in[u]_{\Z}).
\]
The final simple run is $R_1$ by~\eqref{eq:frame-word}.
The lower bound in Lemma~\ref{lem:word-bounds}
with $L=m-1$ and $\nu=u+1$ gives $C_0(\mathcal U(u,v))\ge m$.

Let us obtain the upper bound when $u\ge2$.
The run $\mathcal U(u,v)$ has exactly $u+v+1$ edges.
At each edge, the number of choices of a grade and its color is
\[
 \sum_{k=0}^m\binom{m}{k}=2^m.
\]
Therefore,
\[
 C_0(\mathcal U(u,v))\le(2^m)^{u+v+1}=2^{m(u+v+1)}.
\]
Dividing by $m$ in~\eqref{eq:dictionary} gives
\begin{equation}\label{eq:finite-reduced}
 1\le E_m(u,v)\le\frac{2^{m(u+v+1)}}m.
\end{equation}
If $u=1$, the stability condition~\eqref{eq:stability} means that the images of the $m$ arrow maps span $\C^v$.
The moduli space $\cM^s_{(1,v)}(K(m))$ is the Grassmannian of $v$-dimensional quotients of
$\C^m$. Then, $E_m(1,v)=\binom mv$. The bounds in~\eqref{eq:finite-reduced} also hold for $u=1$.

Let us apply $u+v\le\sqrt{\frac{m+2}{m-2}}\sqrt Q$ from Lemma~\ref{lem:dimension-reduction}.
Taking logarithms in~\eqref{eq:finite-reduced} therefore yields
\begin{equation}\label{eq:finite-Q}
 0\le\log E_m(a,b)
 \le m\log2\sqrt{\frac{m+2}{m-2}Q}+\log\left(\frac{2^m}{m}\right).
\end{equation}

Second, suppose $Q=0$. No positive integral pair has this value, because $m\ge3$ gives
\[
 (m-1)^2<m^2-4<m^2.
\]

Third, suppose $Q<0$. Let us prove $Q=-1$ and $E_m(a,b)=1$.
The nonempty stable moduli space $\cM^s_{(a,b)}(K(m))$ has dimension
$1+Q$~\cite[part (2) of Proposition 4.3]{Weist}. Since $Q$ is a negative integer, we obtain $Q=-1$.

Let us prove $E_m(a,b)=1$ when $Q=-1$.
The stable locus is a nonempty open subset of $\operatorname{Hom}_{\C}(\C^a,\C^b)^m$ and is therefore irreducible.
The smooth geometric quotient $\cM^s_{(a,b)}(K(m))$ is irreducible and zero-dimensional. Hence, $\cM^s_{(a,b)}(K(m))$ consists of one point.
Replacing $Q$ by $\max(Q,0)$ in~\eqref{eq:finite-Q} gives a bound on all of $\cA_m$.
Dividing by $a$ gives~\eqref{eq:finite-endpoint}.
\end{proof}

%% file: sections/run-edits.tex
\section{Run-edits}\label{sec:surgery}
\subsection{Insertion and deletion}
\begin{definition}\label{sync:surgery:2}\leavevmode
Let $()$ denote the empty tuple.
\begin{enumerate}[font=\itshape]
\item Let $U$ be a run. Let $i\in[0,|U|]_{\Z}$.
Let $v\in N(U)$. Let $j\ge1$ be an integer.
We refer to any one of the operations
\[
 U\mapsto\operatorname{Ins}(U,i,N),\qquad
 U\mapsto\operatorname{Del}(U,v),\qquad
 U\mapsto\operatorname{Ins}(U,i,R_j),
\]
or an exchange of two adjacent simple runs as a \textit{run-edit} of $U$.
\item We refer to a run-edit that inserts or deletes a vertical edge as a
\textit{vertical run-edit}.
\item Let $U_0$ be a run and $\ell\in\Z_{\ge0}$.
For each $j\in[\ell]_{\Z}$, let $\omega_j$ be a run-edit of $U_{j-1}$ and define
\[
 U_j=\omega_j(U_{j-1}).
\]
We call $\Omega=(\omega_j)_{j\in[\ell]_{\Z}}$ a \textit{run-edit list} from $U_0$ to $U_\ell$.
We define
\[
 \operatorname{RunEResult}(\Omega,U_0)=U_\ell,
 \qquad \operatorname{RunEResult}((),U_0)=U_0.
\]
\end{enumerate}
\end{definition}

\begin{lemma}[H-precursors after insertion]\label{lem:insert-precursor-chain}
Let $Y=Y_1\ccn\cdots\ccn Y_n$ be a run, let $g\in\operatorname{CtdG}(Y)$. Let $q\in[1,n-1]_{\Z}$. Put
\[
 h_0=Y_1,\qquad W=\operatorname{Ins}(Y,q,N),\qquad w_*=W_{q+1}.
\]
Define the map $\tau:\operatorname{NE}(Y)\to\operatorname{NE}(W)$ by
\[
 \tau(Y_k)=
 \begin{cases}
 W_k,&k\in[q]_{\Z},\\
 W_{k+1},&k\in[q+1,n]_{\Z}.
 \end{cases}
\]
Define $g'\in G(W)$ by
\begin{equation}\label{eq:insertion-grading}
 g'(w_*)=0,\qquad
 g'(\tau(h_0))=g(h_0)+1,\qquad
 g'(\tau(e))=g(e)\quad(e\in\operatorname{NE}(Y)\setminus\{h_0\}).
\end{equation}
Let
\[
 H_q=\{r\in[q]_{\Z}:Y_r\in E(Y),\
          \operatorname{Debt}_{Y,g}(Y_r,q)>0\}.
\]
The following statements hold.

\begin{enumerate}
\item\label{c:insert-active-chain} The set $H_q$ contains $1$. Its decreasing enumeration is uniquely determined:
\[
H_q=\{r_0,\ldots,r_d\},\qquad r_0>\cdots>r_d=1.
\]
For $j\in[0,d]_{\Z}$, put
\[
a_j=Y_{r_j},\qquad A_j=\operatorname{Node}_{Y,g}(a_j),\qquad
t_j=\operatorname{End}_{\mathrm{NE}}(A_j)=Y_{s_j}.
\]
Then, \begin{equation}\label{eq:insertion-node-chain}
A_0\subsetneq A_1\subsetneq\cdots\subsetneq A_d=Y,\qquad
q<s_0<\cdots<s_d=n,
\end{equation}
and $A_{j+1}$ is the parent of $A_j$ for $j\in[0,d-1]_{\Z}$.

\item\label{c:insert-chain-precursor} Use the indices from Claim~\ref{c:insert-active-chain}. The h-precursor $P_{W,g'}$ exists and satisfies
\begin{equation}\label{eq:insertion-chain-map}
\begin{aligned}
P_{W,g'}(w_*)&=\tau(a_0),\\
P_{W,g'}(\tau(t_j))&=\tau(a_{j+1})
   &&(j\in[0,d-1]_{\Z}),\\
P_{W,g'}(\tau(v))&=\tau(P_{Y,g}(v))
   &&\bigl(v\in N(Y)\setminus\{t_j:j\in[0,d-1]_{\Z}\}\bigr).
\end{aligned}
\end{equation}
For $d=0$, the middle line has no instances.

\item\label{c:insert-chain-fibers} For $e\in E(Y)$, put $F_e=P_{Y,g}^{-1}(\{e\})$. For a set of edges $F$, write $\tau(F)=\{\tau(v):v\in F\}$. Define
\[
b_0=w_*,\qquad b_j=\tau(t_{j-1})\quad(j\in[d]_{\Z}).
\]
Then, \begin{equation}\label{eq:insertion-chain-fibers}
\begin{aligned}
P_{W,g'}^{-1}(\{\tau(a_j)\})
   &=\tau(F_{a_j}\setminus\{t_j\})\sqcup\{b_j\}
     &&(j\in[0,d-1]_{\Z}),\\
P_{W,g'}^{-1}(\{\tau(a_d)\})
   &=\tau(F_{a_d})\sqcup\{b_d\},\\
P_{W,g'}^{-1}(\{\tau(e)\})
   &=\tau(F_e)
     &&\bigl(e\in E(Y)\setminus\{a_0,\ldots,a_d\}\bigr).
\end{aligned}
\end{equation}
Consequently, every $e\in E(Y)$ satisfies
\[
\#P_{W,g'}^{-1}(\{\tau(e)\})=g'(\tau(e)),
\]
and
\[
\operatorname{Node}_{W,g'}(\tau(h_0))=W.
\]
\end{enumerate}
\end{lemma}
\begin{proof}
\prf{Claim~\ref{c:insert-active-chain}} Let us enumerate $H_q$ and the nodes $\operatorname{Node}_{Y,g}(Y_r)$ for $r\in H_q$. Connectedness gives $Y=\operatorname{Node}_{Y,g}(h_0)$. Every positive-grade horizontal edge has a node by Claim~\ref{c:node-contained} of Lemma~\ref{lem:node-structure}. If $e=Y_r\in E(Y[0,q])$ and its node ends at $Y_s$, then
\begin{equation}\label{eq:insertion-active-end}
 \operatorname{Debt}_{Y,g}(e,q)>0\quad\Longleftrightarrow\quad q<s.
\end{equation}

Let us order these nodes by inclusion.
We have $P_{Y,g}(Y_n)=h_0$ and $q<n$.
Equation~\eqref{eq:insertion-active-end} gives $1\in H_q$.
Every node indexed by $H_q$ contains the boundary at $q$.
Claim~\ref{c:node-nesting} of Lemma~\ref{lem:node-structure} therefore orders them by inclusion.
The order $r_0>\cdots>r_d=1$ gives
\[
 A_0\subsetneq A_1\subsetneq\cdots\subsetneq A_d=Y.
\]
Claim~\ref{c:node-ends} of Lemma~\ref{lem:node-structure} gives
\[
 q<s_0<\cdots<s_d=n.
\]

Let us prove that $A_{j+1}$ is the parent of $A_j$ for $j\in[0,d-1]_{\Z}$.
If a node $\operatorname{Node}_{Y,g}(Y_r)$ satisfied
\[
 A_j\subsetneq\operatorname{Node}_{Y,g}(Y_r)\subsetneq A_{j+1},
\]
then it would contain the boundary at $q$. Hence,
\[
 r\in H_q,\qquad r_{j+1}<r<r_j.
\]
This contradicts the consecutive enumeration of $H_q$.
Definition~\ref{def:node-relations}(1) proves the parent assertion.

\prf{Claim~\ref{c:insert-chain-precursor}} Let us verify Definition~\ref{def:stack}(1)(a)--(c) for $Q$. Define $Q:N(W)\to E(W)$ by~\eqref{eq:insertion-chain-map}. Define its candidate debt, before knowing that $Q$ is an h-precursor, by
\begin{equation}\label{eq:insertion-candidate-debt}
 D_Q(x,\ell)=g'(x)-\#\{v\in N(W[0,\ell]):Q(v)=x\},
 \qquad x\in E(W[0,\ell]).
\end{equation}
Write $\mathbf1_{\mathcal P}=1$ when $\mathcal P$ holds and $0$ otherwise. For $k\in[0,q]_{\Z}$ and $e\in E(Y[0,k])$, equations~\eqref{eq:insertion-grading}, \eqref{eq:insertion-chain-map}, and~\eqref{eq:insertion-candidate-debt} give
\begin{equation}\label{eq:insertion-debt-before}
 D_Q(\tau(e),k)=\operatorname{Debt}_{Y,g}(e,k)+\mathbf1_{e=h_0}.
\end{equation}
For $k\in[q,n]_{\Z}$ and $e\in E(Y[0,k])$, equations~\eqref{eq:insertion-grading}, \eqref{eq:insertion-chain-map}, and~\eqref{eq:insertion-candidate-debt} give
\begin{equation}\label{eq:insertion-debt-after}
 \begin{aligned}
 D_Q(\tau(e),k+1)
 &=\operatorname{Debt}_{Y,g}(e,k)+\mathbf1_{e=h_0}-\mathbf1_{e=a_0}\\
 &\quad+\sum_{\substack{j\in[0,d-1]_{\Z}\\s_j\le k}}
       (\mathbf1_{e=a_j}-\mathbf1_{e=a_{j+1}}).
 \end{aligned}
\end{equation}
$w_*$ consumes one unit at $a_0$. Changing the owner of $t_j$ returns one unit at $a_j$ and consumes one unit at $a_{j+1}$.

Let us verify $Q(\tau(Y_k))=\tau(P_{Y,g}(Y_k))$ for $Y_k\in N(Y)$ with $k\le q$, and $Q(w_*)=\tau(a_0)$. The debt of $h_0$ is positive at every boundary before $n$. By~\eqref{eq:insertion-debt-before}, adding one unit to this first horizontal edge does not change the last horizontal edge of positive debt. Thus, $Q$ agrees with the h-precursor rule before $w_*$. At boundary $q$, the last positive-debt edge is $a_0$. Thus, $Q(w_*)=\tau(a_0)$.

Let us verify $Q(\tau(Y_k))$ for $Y_k\in N(Y)$ with $k>q$. Put $s_{-1}=q$. For $j\in[0,d]_{\Z}$, $s_{j-1}<k\le s_j$, and $e\in E(Y[0,k-1])$, the sum in~\eqref{eq:insertion-debt-after} telescopes to
\begin{equation}\label{eq:insertion-debt-between}
 D_Q(\tau(e),k)
 =\operatorname{Debt}_{Y,g}(e,k-1)+\mathbf1_{e=h_0}-\mathbf1_{e=a_j}.
\end{equation}
Here boundary $k$ of $W$ is immediately before $\tau(Y_k)=W_{k+1}$.

First, suppose $j<d$ and $Y_k$ is vertical with $k<s_j$. If $P_{Y,g}(Y_k)=a_j$, then both $Y_k$ and $t_j$ have owner $A_j$. Thus,
\begin{equation}\label{eq:insertion-two-units}
 \operatorname{Debt}_{Y,g}(a_j,k-1)\ge2,
 \qquad D_Q(\tau(a_j),k)\ge1.
\end{equation}
If $P_{Y,g}(Y_k)\ne a_j$, its h-precursor follows $a_j$, because $a_j$ has positive debt until $s_j$. Equation~\eqref{eq:insertion-debt-between} leaves its debt unchanged. By~\eqref{eq:insertion-debt-between} and~\eqref{eq:insertion-two-units}, in both cases $D_Q(\tau(P_{Y,g}(Y_k)),k)>0$ and $D_Q(\tau(Y_r),k)=0$ for every horizontal edge $Y_r$ strictly between $P_{Y,g}(Y_k)$ and $Y_k$. Hence, $Q(\tau(Y_k))=\tau(P_{Y,g}(Y_k))$.

Second, suppose $j<d$ and $k=s_j$. We have $t_j=Y_{s_j}$, $P_{Y,g}(t_j)=a_j$, and $s_j=\max\{k:Y_k\in N(Y),\ P_{Y,g}(Y_k)=a_j\}$. Its parent node ends at $s_{j+1}>s_j$ by~\eqref{eq:insertion-node-chain}. Hence, \eqref{eq:insertion-debt-between} gives
\begin{equation}\label{eq:insertion-ending-debts}
\begin{aligned}
 \operatorname{Debt}_{Y,g}(a_j,s_j-1)&=1,
 &D_Q(\tau(a_j),s_j)&=0,\\
 \operatorname{Debt}_{Y,g}(a_{j+1},s_j-1)&\ge1,
 &D_Q(\tau(a_{j+1}),s_j)&\ge1.
\end{aligned}
\end{equation}

Let us exclude a later horizontal edge of positive candidate debt when $j<d$ and $k=s_j$.
Let $r\in[r_{j+1}+1,s_j-1]_{\Z}\setminus\{r_j\}$ satisfy $Y_r\in E(Y)$.
If $g(Y_r)>0$, its node is contained in $A_j$, or is disjoint from $A_j$ and ends before $A_j$ starts.
Claim~\ref{c:insert-active-chain} excludes a node strictly between $A_j$ and $A_{j+1}$.
Claims~\ref{c:node-ends} and~\ref{c:node-nesting} of Lemma~\ref{lem:node-structure} give
\[
 \operatorname{Debt}_{Y,g}(Y_r,s_j-1)=0.
\]
If $g(Y_r)=0$, then $\operatorname{Debt}_{Y,g}(Y_r,s_j-1)=0$ by Claim~\ref{c:prefix-capacity} of Lemma~\ref{lem:precursor-prefix-counts} and Definition~\ref{def:stack}(3).
Equation~\eqref{eq:insertion-debt-between} therefore gives
\[
 D_Q(\tau(Y_r),s_j)=0.
\]
The equalities $D_Q(\tau(Y_r),s_j)=0$ for $r\in[r_{j+1}+1,s_j-1]_{\Z}\setminus\{r_j\}$ satisfying $Y_r\in E(Y)$, together with~\eqref{eq:insertion-ending-debts}, show that
$\tau(a_{j+1})$ is the last horizontal edge of positive candidate debt.
Thus, $Q(\tau(t_j))=\tau(a_{j+1})$.

Third, suppose $j=d$. The two indicator terms in~\eqref{eq:insertion-debt-between} cancel because $a_d=h_0$. Thus, $Q(\tau(Y_k))=\tau(P_{Y,g}(Y_k))$ for every $Y_k\in N(Y)$ with $s_{d-1}<k\le s_d$. Induction through the edges verifies all three conditions of Definition~\ref{def:stack}(1). Uniqueness in Claim~\ref{c:predecessor-uniqueness} of Lemma~\ref{lem:predecessor-uniqueness} gives $Q=P_{W,g'}$.

\prf{Claim~\ref{c:insert-chain-fibers}} Let us count the sets $P_{W,g'}^{-1}(\{\tau(e)\})$ for $e\in E(Y)$. For $j<d$, the fiber at $a_j$ loses $t_j$ and gains exactly $b_j$. The fiber at $a_d$ gains $b_d$ and loses no edge. Every other fiber is copied by $\tau$. These unions are disjoint: $w_*$ is new. For $j\in[d]_{\Z}$, we have $P_{Y,g}(t_{j-1})=a_{j-1}\ne a_j$. This proves~\eqref{eq:insertion-chain-fibers}. Connectedness of $g$ gives $\#F_e=g(e)$, including $F_e=\varnothing$ when $g(e)=0$. Therefore, $\#P_{W,g'}^{-1}(\{\tau(e)\})=g'(\tau(e))$ for every $e\in E(Y)$. The last edge of $P_{W,g'}^{-1}(\{\tau(h_0)\})$ is $\tau(Y_n)=W_{n+1}$. Thus, $\operatorname{Node}_{W,g'}(\tau(h_0))=W$. For $d=0$, this argument reads $P_{W,g'}^{-1}(\{\tau(h_0)\})=\tau(F_{h_0})\sqcup\{w_*\}$.
\end{proof}

\begin{lemma}[Balances after insertion]\label{lem:insert-balance-connected}
Use the hypotheses and construction of Lemma~\ref{lem:insert-precursor-chain}. Then, the following statements hold.

\begin{enumerate}
\item\label{c:insert-balances} For $\ell\in[0,n+1]_{\Z}$,
\[
\operatorname{bBal}_{W,g'}(\ell)=
\begin{cases}
\operatorname{bBal}_{Y,g}(\ell),&0\le\ell\le q,\\
\operatorname{bBal}_{Y,g}(\ell-1),&q+1\le\ell\le n+1.
\end{cases}
\]
Consequently,
\[
\operatorname{eBal}_{W,g'}(\tau(e))=\operatorname{eBal}_{Y,g}(e)
\quad(e\in\operatorname{NE}(Y)),\qquad
\operatorname{eBal}_{W,g'}(w_*)=\operatorname{bBal}_{Y,g}(q).
\]

\item\label{c:insert-connected-surplus} We have
\[
g'\in\operatorname{CtdG}(W),\qquad
\operatorname{Surp}(W,g')=\operatorname{Surp}(Y,g).
\]
If $m\in\Z_{\ge0}$, $g\in G_m(Y)$, and $g(h_0)+1\le m$, then $g'\in\operatorname{CtdG}_m(W)$.
\end{enumerate}
\end{lemma}
\begin{proof}
\prf{Claim~\ref{c:insert-balances}} Let us compute the two terms in the boundary balance.
Definition~\ref{def:colored-grading}(3)(a) gives
\[
 \operatorname{bBal}_{W,g'}(\ell)
 =\#E(W[0,\ell])-\sum_{v\in N(W[0,\ell])}g'(v)
 \qquad(\ell\in[0,n+1]_{\Z}).
\]

First, let us compute this expression for $\ell\in[0,q]_{\Z}$.
The insertion map $\tau$ gives
\[
 E(W[0,\ell])=\tau(E(Y[0,\ell])),\qquad
 N(W[0,\ell])=\tau(N(Y[0,\ell])).
\]
Equation~\eqref{eq:insertion-grading} gives $g'(\tau(v))=g(v)$ for every $v\in N(Y)$.
Hence,
\[
 \operatorname{bBal}_{W,g'}(\ell)
 =\#E(Y[0,\ell])-\sum_{v\in N(Y[0,\ell])}g(v)
 =\operatorname{bBal}_{Y,g}(\ell).
\]

Second, let us compute the expression for $\ell\in[q+1,n+1]_{\Z}$.
The prefix contains the new vertical edge $w_*$, and
\[
 \begin{aligned}
 E(W[0,\ell])&=\tau(E(Y[0,\ell-1])),\\
 N(W[0,\ell])&=\tau(N(Y[0,\ell-1]))\sqcup\{w_*\}.
 \end{aligned}
\]
Since $g'(w_*)=0$, we obtain
\[
 \operatorname{bBal}_{W,g'}(\ell)
 =\#E(Y[0,\ell-1])-\sum_{v\in N(Y[0,\ell-1])}g(v)-0
 =\operatorname{bBal}_{Y,g}(\ell-1).
\]

Let us evaluate the boundary balance after each edge.
For $k\in[q]_{\Z}$, we have $\tau(Y_k)=W_k$.
For $k\in[q+1,n]_{\Z}$, we have $\tau(Y_k)=W_{k+1}$.
Substituting $\operatorname{bBal}_{W,g'}(k)=\operatorname{bBal}_{Y,g}(k)$ for $k\le q$ and $\operatorname{bBal}_{W,g'}(k+1)=\operatorname{bBal}_{Y,g}(k)$ for $k>q$ gives
\[
 \operatorname{eBal}_{W,g'}(\tau(Y_k))
 =\operatorname{bBal}_{Y,g}(k)
 =\operatorname{eBal}_{Y,g}(Y_k)\qquad(k\in[n]_{\Z}).
\]
Since $w_*=W_{q+1}$, they also give
\[
 \operatorname{eBal}_{W,g'}(w_*)
 =\operatorname{bBal}_{W,g'}(q+1)=\operatorname{bBal}_{Y,g}(q).
\]

\prf{Claim~\ref{c:insert-connected-surplus}} Let us verify the node balance condition. Lemma~\ref{lem:insert-precursor-chain} gives the h-precursor, every positive-grade horizontal node, and the root node $W$. Let $v\in N(Y)$, write $a=P_{Y,g}(v)$. Let $b$ satisfy $P_{W,g'}(\tau(v))=\tau(b)$. Claim~\ref{c:insert-chain-precursor} of Lemma~\ref{lem:insert-precursor-chain} gives either $b=a$ or
\[
 \operatorname{Node}_{Y,g}(a)\subsetneq\operatorname{Node}_{Y,g}(b).
\]
The node balance condition for $g$, Claim~\ref{c:node-descendant-balances} of Lemma~\ref{lem:node-structure}, and Claim~\ref{c:insert-balances} give
\[
 \operatorname{eBal}_{W,g'}(\tau(v))
 =\operatorname{eBal}_{Y,g}(v)
 \ge\operatorname{eBal}_{Y,g}(a)
 \ge\operatorname{eBal}_{Y,g}(b)
 =\operatorname{eBal}_{W,g'}(\tau(b)).
\]
The inequality $\operatorname{eBal}_{Y,g}(a)\ge\operatorname{eBal}_{Y,g}(b)$ is equality when $a=b$.

Let us check the balance at the inserted edge.
For the inserted edge, Claim~\ref{c:insert-active-chain} of Lemma~\ref{lem:insert-precursor-chain} gives $r_0\le q<s_0$. Claim~\ref{c:node-boundary-balances} of Lemma~\ref{lem:node-structure} applied to $A_0=Y[r_0-1,s_0]$ yields
\[
 \operatorname{eBal}_{W,g'}(w_*)
 =\operatorname{bBal}_{Y,g}(q)
 \ge\operatorname{eBal}_{Y,g}(a_0)
 =\operatorname{eBal}_{W,g'}(P_{W,g'}(w_*)).
\]
Definition~\ref{def:connected-grading}(1) now gives connectedness.

Let us compute the surplus. We have
\[
 \#E(W)=\#E(Y),\qquad
 \sum_{v\in N(W)}g'(v)=\sum_{v\in N(Y)}g(v)+0.
\]
Substitution in Definition~\ref{def:connected-grading}(2) gives $\operatorname{Surp}(W,g')=\operatorname{Surp}(Y,g)$. If $g\in G_m(Y)$ and $g(h_0)+1\le m$, equation~\eqref{eq:insertion-grading} gives $g'(\tau(h_0))\le m$ and $g'(e)\le m$ for every $e\in\operatorname{NE}(W)\setminus\{\tau(h_0)\}$. Thus, $g'\in G_m(W)$.
\end{proof}

\begin{lemma}[Insertion]\label{lem:vertical-insert}
Let $M\in[m-2]_{\Z}$, $u\in[M]_{\Z}$, and $v\in[M]_{\Z}$.
Let $X$ be a run such that $\#N(A)\le M$ for every simple run $A\subseteq X$.
Let $c\in\Z_{\ge0}$.
Consider $(g,\gamma)\in\mathcal C_c(X)$.
Use the following construction.
\begin{hypotheses}
\item\label{hyp:insertion-shell} Let $Y=\operatorname{Shell}(X;u,v,0)$ with the shell grading $\widetilde g=\operatorname{GShell}_{g}(X;u,v,0)$.
\item\label{hyp:insertion-operation} Let $i\in[|X|]_{\Z}$.
Let
\[
 W=\operatorname{Ins}(Y,u+1+i,N)=R_u\ccn\operatorname{Ins}(X,i,N)\ccn R_v.
\]
Let $q=u+1+i$. Define $g'\in G(W)$ by
\[
 g'(W_k)=\begin{cases}
 u+2,&k=1,\\
 \widetilde g(Y_k),&k\in[2,q]_{\Z},\\
 0,&k=q+1,\\
 \widetilde g(Y_{k-1}),&k\in[q+2,|Y|+1]_{\Z}.
 \end{cases}
\]
\item\label{hyp:insertion-colors} Define
$\kappa\colon\operatorname{NE}(X)\hookrightarrow\operatorname{NE}(W)$ by
\begin{equation}\label{eq:insertion-edge-map}
 \kappa(X_k)=
 \begin{cases}
 W_{u+1+k},&k\in[i]_{\Z},\\
 W_{u+2+k},&k\in[i+1,|X|]_{\Z}.
 \end{cases}
\end{equation}
Let
\[
 \gamma'=\operatorname{Recol}(X,W,\kappa,g';g,\gamma).
\]
\end{hypotheses}
Then, we have the following.

\begin{enumerate}
\item\label{c:vertical-insert-connected} The grading $g'$ has the following properties.
\begin{enumerate}[label=(\alph*),ref=\theenumi(\alph*)]
\item\label{c:vertical-insert-connectivity} $g'\in\operatorname{CtdG}(W)$.
\item\label{c:vertical-insert-residual} We have
$\operatorname{Surp}(W,g')=\operatorname{Surp}(X,g)$.
\end{enumerate}
\item\label{c:vertical-insert-preserved} We have
\begin{equation}\label{eq:insertion-grade-color}
 g'(\kappa(X_k))=g(X_k),\qquad
 \gamma'(\kappa(X_k))=\gamma(X_k)
 \quad(k\in[|X|]_{\Z}).
\end{equation}
\item\label{c:vertical-insert-injective} Fix $X$, $i$, $u$, and $v$.
The construction~\ref{hyp:insertion-shell}--\ref{hyp:insertion-colors} gives an injective map
\[
 \Phi_c\colon\mathcal C_c(X)\longrightarrow\mathcal C_c(W)
\]
such that $\Phi_c(g,\gamma)=(g',\gamma')$ satisfies~\eqref{eq:insertion-grade-color}
for every $(g,\gamma)\in\mathcal C_c(X)$.
\end{enumerate}
\end{lemma}
\begin{proof}
\prf{Claim~\ref{c:vertical-insert-connectivity}} Let us apply the insertion lemmas to the shell grading.
Put $q=u+1+i$ and $h_0=\operatorname{Start}_{\mathrm{NE}}(Y)$.
The hypotheses give
\[
 1\le q<|Y|,\qquad \widetilde g(h_0)+1=u+2\le M+2\le m.
\]
Claim~\ref{c:shell-connected} of Lemma~\ref{lem:shell} gives
$\widetilde g\in\operatorname{CtdG}_m(Y)$.
Apply Lemma~\ref{lem:insert-precursor-chain} to $(Y,\widetilde g,q)$.
Claim~\ref{c:insert-connected-surplus} of Lemma~\ref{lem:insert-balance-connected}
gives $g'\in\operatorname{CtdG}_m(W)$.

\prf{Claim~\ref{c:vertical-insert-residual}}
Claim~\ref{c:insert-connected-surplus} of Lemma~\ref{lem:insert-balance-connected}
and Claim~\ref{c:shell-residual} of Lemma~\ref{lem:shell} give
\[
 \operatorname{Surp}(W,g')=\operatorname{Surp}(Y,\widetilde g)=\operatorname{Surp}(X,g).
\]

\prf{Claim~\ref{c:vertical-insert-preserved}} Let us verify~\eqref{eq:insertion-grade-color}.
By~\eqref{eq:insertion-edge-map}, the copy of $X_k$ in $W$ has index $u+1+k$ for $k\in[i]_{\Z}$ and $u+2+k$ for $k\in[i+1,|X|]_{\Z}$.
Assumptions~\ref{hyp:insertion-operation} and~\ref{hyp:insertion-colors} change no grade or color on these edges.
Thus, the map~\eqref{eq:insertion-edge-map} satisfies~\eqref{eq:insertion-grade-color}.

\prf{Claim~\ref{c:vertical-insert-injective}} Let us recover the input from the inserted grading.
Fix $X$, $i$, $u$, and $v$.
For $(g_1,\gamma_1)\in\mathcal C_c(X)$ and $(g_2,\gamma_2)\in\mathcal C_c(X)$, equation~\eqref{eq:insertion-grade-color} gives
\[
 \Phi_c(g_1,\gamma_1)=\Phi_c(g_2,\gamma_2)
 \ \Longrightarrow\ 
 (g_1(X_k),\gamma_1(X_k))=(g_2(X_k),\gamma_2(X_k))
 \quad(k\in[|X|]_{\Z}).
\]
Hence, $(g_1,\gamma_1)=(g_2,\gamma_2)$.
\end{proof}

\begin{lemma}[H-precursors after deletion]\label{lem:delete-precursors}
Let $Z$ be a tNE path. Let $g\in G(Z)$.
Suppose that the h-precursor $P_{Z,g}$ exists.
Let $v=Z_r\in N(Z)$ and put $x=P_{Z,g}(v)$.
Let $D=\operatorname{Del}(Z,v)$.
Define
\[
 \tau_v\colon\operatorname{NE}(Z)\setminus\{v\}\longrightarrow\operatorname{NE}(D),
 \qquad
 \tau_v(Z_k)=\begin{cases}D_k,&k<r,\\D_{k-1},&k>r.\end{cases}
\]
Define $g_v^-\in G(D)$ by
\[
 g_v^-(\tau_v(e))=
 \begin{cases}g(x)-1,&e=x,\\g(e),&e\ne x\end{cases}
 \qquad(e\in\operatorname{NE}(Z)\setminus\{v\}).
\]
Then, the h-precursor $P_{D,g_v^-}$ exists and satisfies
\begin{equation}\label{eq:delete-precursors}
 P_{D,g_v^-}(\tau_v(w))=\tau_v(P_{Z,g}(w))
 \qquad(w\in N(Z)\setminus\{v\}).
\end{equation}
\end{lemma}
\begin{proof}
Let us define $Q\colon N(D)\to E(D)$ by
\[
 Q(\tau_v(w))=\tau_v(P_{Z,g}(w))
 \qquad(w\in N(Z)\setminus\{v\}).
\]
The map $\tau_v$ preserves edge order and direction.
Condition (b) of Definition~\ref{def:stack}(1) gives $g(x)>0$.
Thus, $g_v^-$ has nonnegative grades.
For $Z_s\in N(Z)$ and $f\in E(Z[0,s-1])$, define
\[
 n_f(s)=\#\{w\in N(Z[0,s-1]):P_{Z,g}(w)=f\}.
\]
For $s\ne r$, define
\[
 n'_f(s)=\#\{w\in N(Z[0,s-1])\setminus\{v\}:Q(\tau_v(w))=\tau_v(f)\}.
\]
The definition of $Q$ gives
\[
 n'_f(s)=
 \begin{cases}
 n_f(s)-1,&s>r\text{ and }f=x,\\
 n_f(s),&s<r\text{ or }f\ne x.
 \end{cases}
\]

Let us verify that $Q$ satisfies the h-precursor conditions of Definition~\ref{def:stack}(1) at $\tau_v(Z_s)$ when $s<r$.
If $f=P_{Z,g}(Z_s)\ne x$, then $g_v^-(\tau_v(f))=g(f)$ and $n'_f(s)=n_f(s)$.
If $P_{Z,g}(Z_s)=x$, then
\[
 n_x(s)+1\le n_x(r)<g(x),\qquad n'_x(s)<g(x)-1.
\]
Thus, $Q(\tau_v(Z_s))$ satisfies conditions (a) and (b) of Definition~\ref{def:stack}(1).
For every $f\in E(Z[0,s-1])$, we have $n'_f(s)=n_f(s)$ and
$g_v^-(\tau_v(f))\le g(f)$.
Hence, every horizontal edge satisfying these two conditions after deletion
also satisfies them before deletion. Condition (c) follows.

Let us verify Definition~\ref{def:stack}(1)(a)--(c) for $Q$ at $\tau_v(Z_s)$ when $s>r$.
For $f\ne x$, we have $g_v^-(\tau_v(f))=g(f)$ and $n'_f(s)=n_f(s)$.
For $f=x$ and $s>r$, we have $v=Z_r\in N(Z[0,s-1])$ and $P_{Z,g}(v)=x$. Thus,
\[
 n_x(s)\ge1,\qquad
 n'_x(s)<g_v^-(\tau_v(x))
 \ \Longleftrightarrow\ n_x(s)-1<g(x)-1
 \ \Longleftrightarrow\ n_x(s)<g(x).
\]
Thus, the horizontal edges satisfying conditions (a) and (b) correspond under $\tau_v$.
The map $\tau_v$ preserves their order. Thus, condition (c) also holds.
Therefore, $Q$ is an h-precursor of $g_v^-$.
Claim~\ref{c:predecessor-uniqueness} of Lemma~\ref{lem:predecessor-uniqueness}
gives~\eqref{eq:delete-precursors}.
\end{proof}

\begin{lemma}[Deletion]\label{lem:vertical-delete}
Let $M\in[m-1]_{\Z}$.
Let $X$ be a run such that $\#N(A)\le M$ for every simple run $A\subseteq X$.
Let $v=X_r\in N(X)$. Let $c\in\Z_{\ge0}$.
Consider $(g,\gamma)\in\mathcal C_c(X)$.
Let $h=g(v)$. Let $u_0\in[M]_{\Z}$ and $v_0\in[M]_{\Z}$.
\begin{hypotheses}
\item\label{hyp:deletion-shell} Let $Y=\operatorname{Shell}(X;u_0,v_0,0)$ with the shell grading $\widetilde g=\operatorname{GShell}_{g}(X;u_0,v_0,0)$.
\item\label{hyp:deletion-edge} Let
\[
 W=\operatorname{Del}(Y,Y_{u_0+1+r})
 =R_{u_0}\ccn\operatorname{Del}(X,v)\ccn R_{v_0}.
\]
Let $q=u_0+1+r$. For each $x=X_j\in E(X[0,r-1])$ satisfying $g(x)>0$, define $g'_x\in G(W)$ by
\[
 g'_x(W_k)=\begin{cases}
 g(x)-1,&k=u_0+1+j,\\
 \widetilde g(Y_k),&k\in[1,q-1]_{\Z}\setminus\{u_0+1+j\},\\
 \widetilde g(Y_{k+1}),&k\in[q,|Y|-1]_{\Z}.
 \end{cases}
\]
\item\label{hyp:deletion-color} Define
$\kappa\colon\operatorname{NE}(X)\setminus\{v\}\hookrightarrow\operatorname{NE}(W)$ by
\[
 \kappa(X_k)=
 \begin{cases}
 W_{u_0+1+k},&k\in[r-1]_{\Z},\\
 W_{u_0+k},&k\in[r+1,|X|]_{\Z}.
 \end{cases}
\]
For each $x\in E(X[0,r-1])$ satisfying $g(x)>0$, define
\[
 \gamma'_x=\operatorname{Recol}(X,W,\kappa,g'_x;g,\gamma).
\]
\end{hypotheses}
Then, we have the following.

\begin{enumerate}
\item\label{c:vertical-delete-connected}\label{c:vertical-delete-connectivity}
There exists $y\in E(X[0,r-1])$ with $g(y)>0$ such that
$g'_y\in\operatorname{CtdG}(W)$.
\item\label{c:vertical-delete-residual} For every $x\in E(X[0,r-1])$ with $g(x)>0$, we have
\[
 \#E(W)-1-\sum_{w\in N(W)}g'_x(w)=\operatorname{Surp}(X,g)+h.
\]
\item\label{c:vertical-delete-preserved} For every $x\in E(X[0,r-1])$ with $g(x)>0$, we have
\begin{equation}\label{eq:deletion-grade-color}
\begin{aligned}
 g'_x(\kappa(x))&=g(x)-1,\\
 g'_x(\kappa(e))&=g(e),\qquad \gamma'_x(\kappa(e))=\gamma(e)
 \quad(e\in\operatorname{NE}(X)\setminus\{x,v\}).
\end{aligned}
\end{equation}
\item\label{c:vertical-delete-map} Fix $X$, $v$, $u_0$, and $v_0$.
Also, for each $(g,\gamma)\in\mathcal C_c(X)$, fix an edge $y(g,\gamma)$ satisfying Claim~\ref{c:vertical-delete-connected}.
The construction~\ref{hyp:deletion-shell}--\ref{hyp:deletion-color} gives a map
\[
 \Phi_c\colon\mathcal C_c(X)\longrightarrow
 \bigcup_{s\in[c,c+m]_{\Z}}\mathcal C_s(W)
\]
such that $\Phi_c(g,\gamma)=(g'_{y(g,\gamma)},\gamma'_{y(g,\gamma)})$.
The identities~\eqref{eq:deletion-grade-color} hold at $x=y(g,\gamma)$.
Also, we have
\[
 \operatorname{Surp}(W,g'_{y(g,\gamma)})=c+g(v).
\]
\end{enumerate}
\end{lemma}
\begin{proof}
\prf{Claim~\ref{c:vertical-delete-connectivity}} Let us choose $y\in E(X[0,r-1])$.
Let $y=P_{X,g}(v)$.
Definition~\ref{def:stack}(1) gives $y\in E(X[0,r-1])$ and $g(y)>0$.
Let $q=u_0+1+r$ and let $\widetilde y$ be the copy of $y$ in $Y$.
Claim~\ref{c:shell-connected} of Lemma~\ref{lem:shell} gives
$\widetilde g\in\operatorname{CtdG}_m(Y)$.
Claim~\ref{c:shell-copy-precursor} of Lemma~\ref{lem:shell} gives
$P_{Y,\widetilde g}(Y_q)=\widetilde y$.

Let us apply the deletion formula to the copied edge $Y_q$.
Let $\tau\colon\operatorname{NE}(Y)\setminus\{Y_q\}\to\operatorname{NE}(W)$
be the map $\tau_{Y_q}$ of Lemma~\ref{lem:delete-precursors}.
That lemma applies to $Y$, $\widetilde g$, and $Y_q$ and gives
\begin{equation}\label{eq:vertical-delete-precursor}
 P_{W,g'_y}(\tau(w))=\tau(P_{Y,\widetilde g}(w))
 \qquad(w\in N(Y)\setminus\{Y_q\}).
\end{equation}

Let us verify the node balance condition.
For $k\in[|Y|]_{\Z}\setminus\{q\}$, deletion of the vertical grade $h$ gives
\begin{equation}\label{eq:vertical-delete-balance}
 \operatorname{eBal}_{W,g'_y}(\tau(Y_k))
 =\operatorname{eBal}_{Y,\widetilde g}(Y_k)
 +\begin{cases}0,&k<q,\\h,&k>q.\end{cases}
\end{equation}
Let $w=Y_k\in N(Y)\setminus\{Y_q\}$ and let
$P_{Y,\widetilde g}(w)=Y_j$.
We have $j<k$. Thus,~\eqref{eq:vertical-delete-balance} gives
\[
\begin{aligned}
 &\operatorname{eBal}_{W,g'_y}(\tau(w))
 -\operatorname{eBal}_{W,g'_y}(\tau(Y_j))\\
 &\quad\ge\operatorname{eBal}_{Y,\widetilde g}(w)
 -\operatorname{eBal}_{Y,\widetilde g}(Y_j)\ge0.
\end{aligned}
\]
The inequality $\operatorname{eBal}_{Y,\widetilde g}(w)-\operatorname{eBal}_{Y,\widetilde g}(Y_j)\ge0$ follows from connectedness of $\widetilde g$ and
Claim~\ref{c:node-contained} of Lemma~\ref{lem:node-structure}.
Together with~\eqref{eq:vertical-delete-precursor}, this proves the node balance condition for $g'_y$.

Let us verify $\operatorname{Node}_{W,g'_y}(W_1)=W$.
The edges $Y_1$ and $Y_{|Y|}$ belong to the copies of $R_{u_0}$ and $R_{v_0}$, respectively.
Deleting $Y_q$ keeps $Y_1$ and $Y_{|Y|}$ and gives $g'_y(W_1)=\widetilde g(Y_1)$.
We have $P_{Y,\widetilde g}(Y_q)=\widetilde y\ne Y_1$.
Equation~\eqref{eq:vertical-delete-precursor} therefore gives $P_{W,g'_y}(\tau(w))=W_1$ for every $w\in P_{Y,\widetilde g}^{-1}(\{Y_1\})$.
Since $\operatorname{Node}_{Y,\widetilde g}(Y_1)=Y$, we obtain
\[
 \operatorname{Node}_{W,g'_y}(W_1)=W.
\]
Definition~\ref{def:connected-grading}(1) gives $g'_y\in\operatorname{CtdG}(W)$.
The grading $g'_y$ decreases only the positive grade $g(y)$.
Since $\widetilde g\in G_m(Y)$, we also have $g'_y\in G_m(W)$.

\prf{Claim~\ref{c:vertical-delete-residual}} Let us sum the vertical grades.
Let $x\in E(X[0,r-1])$ satisfy $g(x)>0$.
Deletion preserves the number of horizontal edges and removes vertical grade $h$.
By~\eqref{eq:residual} and Claim~\ref{c:shell-residual} of Lemma~\ref{lem:shell},
\[
 \#E(W)-1-\sum_{w\in N(W)}g'_x(w)
 =\#E(Y)-1-\left(\sum_{w\in N(Y)}\widetilde g(w)-h\right)
 =\operatorname{Surp}(X,g)+h.
\]

\prf{Claim~\ref{c:vertical-delete-preserved}} Let us identify the unchanged grades and colors.
Let $x\in E(X[0,r-1])$ with $g(x)>0$.
Assumption~\ref{hyp:deletion-edge} gives the grade identities
in~\eqref{eq:deletion-grade-color}.
For $e\in\operatorname{NE}(X)\setminus\{x,v\}$, we have
$g'_x(\kappa(e))=g(e)$.
Definition~\ref{def:recolor} and assumption~\ref{hyp:deletion-color}
therefore give $\gamma'_x(\kappa(e))=\gamma(e)$.

\prf{Claim~\ref{c:vertical-delete-map}} Let us prove
$\Phi_c(g,\gamma)\in\bigcup_{s\in[c,c+m]_{\Z}}\mathcal C_s(W)$.
Let $y=y(g,\gamma)$.
For the chosen $y$, Claim~\ref{c:vertical-delete-connectivity} gives
$g'_y\in\operatorname{CtdG}(W)$.
Definition~\ref{def:recolor} gives $\gamma'_y\in\operatorname{Col}(W,g'_y,m)$.
Definition~\ref{def:connected-grading}(2), Claim~\ref{c:vertical-delete-residual}, and $g(v)\in[0,m]_{\Z}$ give
\[
 c\le\operatorname{Surp}(W,g'_y)=c+g(v)\le c+m.
\]
Claim~\ref{c:vertical-delete-preserved} gives~\eqref{eq:deletion-grade-color} at $x=y$.
\end{proof}

\begin{lemma}[Insertion of a simple run]\label{lem:simple-run-insertion}
Let $Q$ be a tNE path, $c\in\Z_{\ge0}$, $q\in[1,|Q|-1]_{\Z}$, and $j\in[m]_{\Z}$.
Let $Q'=\operatorname{Ins}(Q,q,R_j)$.
Define $\tau\colon\operatorname{NE}(Q)\hookrightarrow\operatorname{NE}(Q')$ by
\begin{equation}\label{eq:run-insertion-map}
 \tau(Q_k)=
 \begin{cases}
 Q'_k,&k\in[q]_{\Z},\\
 Q'_{k+j+1},&k\in[q+1,|Q|]_{\Z}.
 \end{cases}
\end{equation}
Let $a_*=Q'_{q+1}$ and $v_t=Q'_{q+1+t}$ for $t\in[j]_{\Z}$.
For $(h,\eta)\in\mathcal C_c(Q)$, define $h'\in G_m(Q')$ by
\begin{equation}\label{eq:run-insertion-grades}
 h'=h|_{Q[0,q]}\ccn g_j^{\mathrm{run}}\ccn h|_{Q[q,|Q|]}.
\end{equation}
Let
\begin{equation}\label{eq:run-insertion-recolor}
 \eta'=\operatorname{Recol}(Q,Q',\tau,h';h,\eta).
\end{equation}
Then, the map
\[
 \Phi_{q,j,c}\colon\mathcal C_c(Q)\longrightarrow\mathcal C_{c+1}(Q'),
 \qquad \Phi_{q,j,c}(h,\eta)=(h',\eta')
\]
is well-defined and injective.
\end{lemma}
\begin{proof}
Let us determine the h-precursor after insertion of the simple run.
The new horizontal edge $a_*$ is the h-precursor of the $j$ new vertical edges.
The debt at $a_*$ decreases from $j$ to zero before the next original edge.
Definition~\ref{def:stack}(1) gives
\begin{equation}\label{eq:run-insertion-precursors}
\begin{aligned}
 P_{Q',h'}(v_t)&=a_*&&(t\in[j]_{\Z}),\\
 P_{Q',h'}(\tau(w))&=\tau(P_{Q,h}(w))&&(w\in N(Q)),\\
 \operatorname{Node}_{Q',h'}(a_*)&=Q'[q,q+j+1].
\end{aligned}
\end{equation}
For each $e\in E(Q)$, equation~\eqref{eq:run-insertion-precursors} gives
\[
 \#P_{Q',h'}^{-1}(\{\tau(e)\})=\#P_{Q,h}^{-1}(\{e\})=h(e)=h'(\tau(e)).
\]
For $a_0=\operatorname{Start}_{\mathrm{NE}}(Q)$, we therefore have
\[
 \operatorname{Node}_{Q',h'}(\tau(a_0))=Q'.
\]

Let us verify the node balance condition.
For $w=Q_k$, let $\varepsilon(w)=0$ when $k\le q$ and $\varepsilon(w)=1$ when $k>q$.
Equation~\eqref{eq:run-insertion-grades} gives
\begin{equation}\label{eq:run-insertion-balances}
\begin{aligned}
 \operatorname{eBal}_{Q',h'}(\tau(w))&=\operatorname{eBal}_{Q,h}(w)+\varepsilon(w),\\
 \operatorname{eBal}_{Q',h'}(a_*)&=\operatorname{eBal}_{Q',h'}(v_t)
 =\operatorname{bBal}_{Q,h}(q)+1\quad(t\in[j]_{\Z}).
\end{aligned}
\end{equation}
For $w\in N(Q)$ and $a=P_{Q,h}(w)$, the edge $a$ precedes $w$.
Hence, $\varepsilon(a)\le\varepsilon(w)$ and
\[
\begin{aligned}
 &\operatorname{eBal}_{Q',h'}(\tau(w))-\operatorname{eBal}_{Q',h'}(\tau(a))\\
 &\quad=\operatorname{eBal}_{Q,h}(w)-\operatorname{eBal}_{Q,h}(a)
       +\varepsilon(w)-\varepsilon(a)\ge0.
\end{aligned}
\]
Together with~\eqref{eq:run-insertion-precursors} and~\eqref{eq:run-insertion-balances}, this proves that $h'$ is connected.
Finally,
\begin{equation}\label{eq:run-insertion-surplus}
 |Q'|_E=|Q|_E+1,\qquad
 \operatorname{TNG}_{Q',h'}=\operatorname{TNG}_{Q,h},\qquad
 \operatorname{Surp}(Q',h')=\operatorname{Surp}(Q,h)+1.
\end{equation}

Let us recover $(h,\eta)$ from $\Phi_{q,j,c}(h,\eta)$.
Definition~\ref{def:recolor} and~\eqref{eq:run-insertion-grades} give
\[
 (h'(\tau(e)),\eta'(\tau(e)))=(h(e),\eta(e))
 \qquad(e\in\operatorname{NE}(Q)).
\]
Restriction to the image of $\tau$ therefore recovers $(h,\eta)$.
Hence, $\Phi_{q,j,c}$ is injective.
\end{proof}

\subsection{Composition of run-edits}
\begin{lemma}[Surplus reduction]\label{lem:surplus-reduction}
Let $m\in\Z_{\ge3}$ and $n\in\Z_{\ge0}$.
Let $Q$ be a tNE path.
Fix $\alpha=(a_j)_{j\in[n]_{\Z}}$ and
$\beta=(b_j)_{j\in[n]_{\Z}}$ with $a_j\in[m-1]_{\Z}$ and $b_j\in[m-1]_{\Z}$.
Define
\[
 Q_0=Q,\qquad Q_j=R_{a_j}\ccn Q_{j-1}\ccn R_{b_j}\quad(j\in[n]_{\Z}),
 \qquad \rho=\sum_{j=1}^n(a_j+1).
\]
Then, there is an injective map
\begin{equation}\label{eq:surplus-reduction-map}
 \Delta_{Q,\alpha,\beta}\colon
 \bigcup_{c\in[0,n]_{\Z}}\mathcal C_c(Q)\longrightarrow\mathcal C_0(Q_n)
\end{equation}
such that for each $(g,\gamma)\in\bigcup_{c\in[0,n]_{\Z}}\mathcal C_c(Q)$,
$(g',\gamma')=\Delta_{Q,\alpha,\beta}(g,\gamma)$ satisfies
\begin{equation}\label{eq:surplus-reduction-restriction}
 (g'((Q_n)_{\rho+k}),\gamma'((Q_n)_{\rho+k}))
 =(g(Q_k),\gamma(Q_k))\qquad(k\in[|Q|]_{\Z}).
\end{equation}
\end{lemma}
\begin{proof}
Let us construct the shell gradings and colorings.
Let $(g,\gamma)\in\bigcup_{c\in[0,n]_{\Z}}\mathcal C_c(Q)$.
Let $c=\operatorname{Surp}(Q,g)$ and $(g_0,\gamma_0)=(g,\gamma)$.
For $j\in[n]_{\Z}$, define
\[
 \epsilon_j(c)=\begin{cases}1,&j\le c,\\0,&j>c,\end{cases}
 \qquad
 g_j=\operatorname{GShell}_{g_{j-1}}(Q_{j-1};a_j,b_j,\epsilon_j(c)).
\]
Define $\iota_j\colon\operatorname{NE}(Q_{j-1})\hookrightarrow\operatorname{NE}(Q_j)$ by
\[
 \iota_j((Q_{j-1})_k)=(Q_j)_{a_j+1+k}
 \qquad(k\in[|Q_{j-1}|]_{\Z}).
\]
Let
\[
 \gamma_j=\operatorname{Recol}
 (Q_{j-1},Q_j,\iota_j,g_j;g_{j-1},\gamma_{j-1}).
\]

Let us prove that the final surplus is zero.
Claims~\ref{c:shell-connected}--\ref{c:shell-recolor} of Lemma~\ref{lem:shell} give, by induction,
\[
 (g_j,\gamma_j)\in\mathcal C_{\max(c-j,0)}(Q_j)
 \qquad(j\in[0,n]_{\Z}).
\]
For $j\le c$, we have $\epsilon_j(c)=1$ and
$\operatorname{Surp}(Q_{j-1},g_{j-1})=c-j+1\ge1$.
Thus, condition~\ref{hyp:shell-draining} of Lemma~\ref{lem:shell} applies.
For $j>c$, we have $\epsilon_j(c)=0$ and use condition~\ref{hyp:shell-preserving} of Lemma~\ref{lem:shell}.
Claim~\ref{c:shell-residual} of Lemma~\ref{lem:shell} gives
$\operatorname{Surp}(Q_j,g_j)=\operatorname{Surp}(Q_{j-1},g_{j-1})-\epsilon_j(c)$.
Since $c\le n$, define
\[
 \Delta_{Q,\alpha,\beta}(g,\gamma)=(g_n,\gamma_n)
 \in\mathcal C_0(Q_n).
\]

Let us recover $(g,\gamma)$ from $(g_n,\gamma_n)$.
Claim~\ref{c:shell-recolor} of Lemma~\ref{lem:shell} gives
\begin{equation}\label{eq:surplus-reduction-step-restriction}
 (g_j(\iota_j(e)),\gamma_j(\iota_j(e)))
 =(g_{j-1}(e),\gamma_{j-1}(e))
 \qquad(e\in\operatorname{NE}(Q_{j-1})).
\end{equation}
Composing~\eqref{eq:surplus-reduction-step-restriction} for $j\in[n]_{\Z}$ gives~\eqref{eq:surplus-reduction-restriction}.
Restriction to the copy of $Q$ therefore recovers $(g,\gamma)$. Thus,
$\Delta_{Q,\alpha,\beta}$ is injective on its whole domain.
\end{proof}

\begin{definition}\label{def:fiber-bound}
For $\nu\in\Z_{\ge0}$ and $D\in\Z_{\ge0}$, define
\[
 \operatorname{Fib}_m(\nu,D)=\binom{\nu+D}{D}2^{2mD}.
\]
\end{definition}

\begin{proposition}\label{prop:composition}
Let $M\in[m-1]_{\Z}$.
Let $U$ be a run with $\nu$ simple runs.
Let $K\in\Z_{\ge0}$, $A\in\Z_{\ge0}$, $D\in\Z_{\ge0}$, and $\ell=K+A+D$.
Assume the following.
\begin{hypotheses}
\item\label{hyp:comparison-edit-list} Let $\Omega=(\omega_s)_{s\in[\ell]_{\Z}}$ be a run-edit list starting at $U$.
Let
\[
 U_s=\operatorname{RunEResult}((\omega_1,\ldots,\omega_s),U)
 \quad(s\in[0,\ell]_{\Z}),\qquad V=\operatorname{RunEResult}(\Omega,U).
\]
Let $I_R$, $I_+$, and $I_-$ satisfy
\[
 [\ell]_{\Z}=I_R\sqcup I_+\sqcup I_-,\qquad \#I_R=K,\qquad \#I_+=A,\qquad \#I_-=D.
\]
For $s\in I_R\cup I_+$, fix $i_s\in[0,|U_{s-1}|]_{\Z}$.
For $s\in I_R$, fix $j_s\in\Z_{>0}$.
For $s\in I_-$, fix $v_s\in N(U_{s-1})$.
The run-edits satisfy
\[
 U_s=
 \begin{cases}
 \operatorname{Ins}(U_{s-1},i_s,R_{j_s}),&s\in I_R,\\
 \operatorname{Ins}(U_{s-1},i_s,N),&s\in I_+,\\
 \operatorname{Del}(U_{s-1},v_s),&s\in I_-.
 \end{cases}
\]
\item\label{hyp:comparison-run-lengths} For $s\in[0,\ell]_{\Z}$, the path $U_s$ is a run
and every simple run $B\subseteq U_s$ satisfies $\#N(B)\in[M]_{\Z}$.
\item\label{hyp:comparison-insertion-range} If $K+A>0$, then
$m\ge4$ and $M\le m-2$.
\item\label{hyp:comparison-shell-budget} Let $u\in[M]_{\Z}$ and $v\in[M]_{\Z}$.
Let $P$ be an integer satisfying
\begin{equation}\label{eq:layer-budget}
 P\ge K+A+(m+1)D.
\end{equation}
\item\label{hyp:comparison-target} Fix the tNE paths $W_{\mathrm L}$ and $W_{\mathrm R}$. Each path consists of $P$ simple runs.
Every simple run of these paths has vertical length in $[M]_{\Z}$.
Let
\begin{equation}\label{eq:target-shape}
 V_*=W_{\mathrm L}\ccn R_u\ccn V\ccn R_v\ccn W_{\mathrm R}.
\end{equation}
\end{hypotheses}
Then, we have the following.

\begin{enumerate}
\item\label{c:composition-target} There is a map
\[
 \Phi\colon\mathcal C_0(U)\longrightarrow\mathcal C_0(V_*)
\]
such that
\[
 \#\Phi^{-1}(\{\xi\})\le \operatorname{Fib}_m(\nu,D)
 \qquad(\xi\in\mathcal C_0(V_*)).
\]
\item\label{c:composition-bound} We have
\begin{equation}\label{eq:composition}
 C_0(V_*)\ge\frac{C_0(U)}{\operatorname{Fib}_m(\nu,D)}.
\end{equation}
\end{enumerate}
\end{proposition}
\begin{proof}
\prf{Claim~\ref{c:composition-target}} Let us construct $\Phi$.
Let $(g,\gamma)\in\mathcal C_0(U)$.
For $u$ and $v$ in assumption~\ref{hyp:comparison-shell-budget}, let
\[
 Y=\operatorname{Shell}(U;u,v,0)
 =R_u\ccn U\ccn R_v.
\]
Let $h=\operatorname{GShell}_{g}(U;u,v,0)$.
Define $\iota\colon\operatorname{NE}(U)\hookrightarrow\operatorname{NE}(Y)$ by
$\iota(U_k)=Y_{u+1+k}$ for $k\in[|U|]_{\Z}$ and let
\[
 \eta=\operatorname{Recol}(U,Y,\iota,h;g,\gamma).
\]
Claims~\ref{c:shell-connected} and~\ref{c:shell-residual} of Lemma~\ref{lem:shell} give $(h,\eta)\in\mathcal C_0(Y)$.
Suppose $P=0$.
Inequality~\eqref{eq:layer-budget} gives $K=A=D=\ell=0$.
Thus, $V=U$ and $W_{\mathrm L}=W_{\mathrm R}=\emptyset$.
The paths $Y$ and $V_*$ are nonempty translated copies.
Let
\[
 \mathbf t=\operatorname{Start}_{\mathrm{pt}}(V_*)-
             \operatorname{Start}_{\mathrm{pt}}(Y),\qquad
 \Phi(g,\gamma)=(h+\mathbf t,\eta+\mathbf t).
\]
Claim~\ref{c:shell-recolor} of Lemma~\ref{lem:shell} and Lemma~\ref{lem:word-reindexing}
show that $\Phi:\mathcal C_0(U)\to\mathcal C_0(V_*)$ is injective.
Its fibers therefore have size at most $1=\operatorname{Fib}_m(\nu,0)$.

Second, suppose $P\ge1$. Let us define $(Q_s,h_s,\eta_s)$ for $s\in[0,\ell]_{\Z}$.
Let $(Q_0,h_0,\eta_0)=(Y,h,\eta)$.
For $s\in[\ell]_{\Z}$, let us construct $(Q_s,h_s,\eta_s)$ recursively.
The run-edit $\omega_s$ acts on the copy of $U_{s-1}$ in $Q_{s-1}$.
Here $h_s=h_s(g,\gamma)$ and $\eta_s=\eta_s(g,\gamma)$.
Define $a_t\in[M]_{\Z}$ and $b_t\in[M]_{\Z}$ for $t\in[P]_{\Z}$ by
\[
\begin{aligned}
 W_{\mathrm L}&=\operatorname{Trans}(R_{a_P}\ccn \cdots\ccn R_{a_1},\operatorname{Start}_{\mathrm{pt}}(W_{\mathrm L})),\\
 W_{\mathrm R}&=\operatorname{Trans}(R_{b_1}\ccn \cdots\ccn R_{b_P},\operatorname{Start}_{\mathrm{pt}}(W_{\mathrm R})).
\end{aligned}
\]

Let us insert a simple run at step $s\in I_R$.
Let $Q_{s-1}=B_s^{\mathrm L}\ccn U_{s-1}\ccn B_s^{\mathrm R}$ and put $q_s=|B_s^{\mathrm L}|+i_s$.
The paths $B_s^{\mathrm L}$ and $B_s^{\mathrm R}$ contain the outer simple runs. Thus,
$q_s\in[1,|Q_{s-1}|-1]_{\Z}$.
Assumption~\ref{hyp:comparison-run-lengths} gives $j_s\in[M]_{\Z}\subseteq[m]_{\Z}$.
Apply Lemma~\ref{lem:simple-run-insertion} with
$Q=Q_{s-1}$, $q=q_s$, $j=j_s$, and $c=\operatorname{Surp}(Q_{s-1},h_{s-1})$.
Let $(h_s,\eta_s)$ be the image of $(h_{s-1},\eta_{s-1})$ under its map $\Phi_{q,j,c}$.
The resulting path is
\[
 Q_s=\operatorname{Ins}(Q_{s-1},q_s,R_{j_s})=B_s^{\mathrm L}\ccn U_s\ccn B_s^{\mathrm R}.
\]

Let us verify the hypotheses for the vertical run-edits.
At step $s\in I_+\cup I_-$, let
\[
 Q_{s-1}=B_s^{\mathrm L}\ccn U_{s-1}\ccn B_s^{\mathrm R},\qquad
 t_s=\#\bigl((I_+\cup I_-)\cap[s]_{\Z}\bigr).
\]
The paths $B_s^{\mathrm L}$ and $B_s^{\mathrm R}$ consist of the simple runs added around $U_{s-1}$.
Every simple run of $Q_{s-1}$ has vertical length in $[M]_{\Z}$
by assumptions~\ref{hyp:comparison-run-lengths}, \ref{hyp:comparison-shell-budget}, and~\ref{hyp:comparison-target}.

First, let us insert a vertical edge when $s\in I_+$.
Then, $A>0$, and assumption~\ref{hyp:comparison-insertion-range} gives
\[
 m\ge4,\qquad M\le m-2.
\]
Claims~\ref{c:vertical-insert-connected}--\ref{c:vertical-insert-injective}
of Lemma~\ref{lem:vertical-insert} therefore apply after the first $|B_s^{\mathrm L}|+i_s$ edges.

Second, let us delete a vertical edge when $s\in I_-$.
Let $v_s=(U_{s-1})_{r_s}$ and let
\[
 \widehat v_s=(Q_{s-1})_{|B_s^{\mathrm L}|+r_s}.
\]
Claims~\ref{c:vertical-delete-connected}--\ref{c:vertical-delete-map}
of Lemma~\ref{lem:vertical-delete} apply at $\widehat v_s$ because $M\le m-1$.
For each connected colored grading $(h,\eta)$ on $Q_{s-1}$, fix an edge
$y_s(h,\eta)$ satisfying Claim~\ref{c:vertical-delete-connected} of Lemma~\ref{lem:vertical-delete}
for $X=Q_{s-1}$, $v=\widehat v_s$, $u_0=a_{t_s}$, and $v_0=b_{t_s}$.
At step $s$, use $y_s(h_{s-1},\eta_{s-1})$.

Let us specify the shell path used at step $s\in I_+\cup I_-$.
Use the parameters $a_{t_s}$ and $b_{t_s}$ in the construction
\ref{hyp:insertion-shell}--\ref{hyp:insertion-colors} of Lemma~\ref{lem:vertical-insert}
or the construction \ref{hyp:deletion-shell}--\ref{hyp:deletion-color} of Lemma~\ref{lem:vertical-delete}, respectively.
The resulting path is
\[
 Q_s=R_{a_{t_s}}\ccn B_s^{\mathrm L}\ccn U_s\ccn B_s^{\mathrm R}\ccn R_{b_{t_s}}.
\]
Claim~\ref{c:vertical-insert-injective} of Lemma~\ref{lem:vertical-insert}
and Claim~\ref{c:vertical-delete-map} of Lemma~\ref{lem:vertical-delete}
supply the colored grading $(h_s,\eta_s)$.

Let us compute $\operatorname{Surp}(Q_\ell,h_\ell)$.
For $s\in[0,\ell]_{\Z}$, let $c_s=\operatorname{Surp}(Q_s,h_s)$.
The initial shell gives $c_0=0$.
Equation~\eqref{eq:run-insertion-surplus}, Claim~\ref{c:vertical-insert-residual} of Lemma~\ref{lem:vertical-insert},
and Claim~\ref{c:vertical-delete-residual} of Lemma~\ref{lem:vertical-delete} give
\[
 c_s-c_{s-1}=
 \begin{cases}
 1,&s\in I_R,\\
 0,&s\in I_+,\\
 h_{s-1}(\widehat v_s),&s\in I_-.
 \end{cases}
\]
Since $h_{s-1}(\widehat v_s)\in[0,m]_{\Z}$, summing gives
\[
 c_\ell=K+\sum_{s\in I_-}h_{s-1}(\widehat v_s)\le K+mD.
\]

Let us reduce the surplus to zero.
Let $n=P-A-D$.
Assumption~\ref{hyp:comparison-shell-budget} gives
\[
 n\ge K+mD\ge c_\ell.
\]
Apply Lemma~\ref{lem:surplus-reduction} with $Q=Q_\ell$ and
\[
 \alpha'=(a_{A+D+j})_{j\in[n]_{\Z}},\qquad
 \beta'=(b_{A+D+j})_{j\in[n]_{\Z}}.
\]
Let $Q'$ denote the path $Q_n$ in that application and let
\[
 (f,\delta)=\Delta_{Q_\ell,\alpha',\beta'}(h_\ell,\eta_\ell)
 \in\mathcal C_0(Q').
\]

Let us identify $Q'$ and define $\Phi$.
The paths $Q_\ell$ and $Q'$ are translated copies of
\[
 R_{a_{A+D}}\ccn \cdots\ccn R_{a_1}\ccn R_u\ccn V\ccn R_v\ccn R_{b_1}\ccn \cdots\ccn R_{b_{A+D}}
\]
and
\[
 R_{a_P}\ccn \cdots\ccn R_{a_1}\ccn R_u\ccn V\ccn R_v\ccn R_{b_1}\ccn \cdots\ccn R_{b_P},
\]
respectively. When $A+D=0$, the blocks indexed by $[A+D]_{\Z}$ are empty.
Let
\[
 \mathbf t=\operatorname{Start}_{\mathrm{pt}}(V_*)-
             \operatorname{Start}_{\mathrm{pt}}(Q').
\]
Then, $\operatorname{Trans}(Q',\mathbf t)=V_*$.
Using the translated grading and color function from
Definition~\ref{def:colored-grading}(5), define
\[
 \Phi(g,\gamma)=(f+\mathbf t,\delta+\mathbf t).
\]
Lemma~\ref{lem:word-reindexing} gives $\Phi(g,\gamma)\in\mathcal C_0(V_*)$.

Let us bound the fibers of $\Phi$ by recovering the colored gradings on $U$.
Let $D_U\subseteq N(U)$ be the set of edges of $U$ deleted by $\Omega$ in assumption~\ref{hyp:comparison-edit-list}.
For $e\in\operatorname{NE}(U)\setminus D_U$, let $\widetilde e\in\operatorname{NE}(V_*)$ be its image under the inclusion $U\to Q_0$, the edge maps $Q_{s-1}\to Q_s$ at steps $s=1,\ldots,\ell$, the inclusion $Q_\ell\to Q'$, and the translation by $\mathbf t$, in that order.
Fix $(h,\eta)\in\mathcal C_0(V_*)$.
Define $\widehat h\in G(U)$ by
\[
 \widehat h(e)=
 \begin{cases}
 h(\widetilde e),& e\notin D_U,\\
 0,& e\in D_U.
 \end{cases}
\]

Let us recover an input colored grading from its differences and deleted edges.
Let $e_1,\ldots,e_\nu$ be the horizontal edges of $U$ in path order.
For $(g,\gamma)\in\Phi^{-1}(\{(h,\eta)\})$, let
\[
 \delta_i=\operatorname{Diff}_U(g,\widehat h;e_i)
 =g(e_i)-h(\widetilde e_i)\qquad(i\in[\nu]_{\Z}).
\]
Claim~\ref{c:shell-recolor} of Lemma~\ref{lem:shell}, Claim~\ref{c:vertical-insert-preserved} of Lemma~\ref{lem:vertical-insert}, and Lemma~\ref{lem:surplus-reduction} preserve the grades and colors of the traced input edges.
Equation~\eqref{eq:run-insertion-grades} and Definition~\ref{def:recolor} give grade and color preservation for simple-run insertion.
Claim~\ref{c:vertical-delete-preserved} of Lemma~\ref{lem:vertical-delete} decreases one horizontal grade by $1$ and preserves all other surviving grades and colors.
Assumption~\ref{hyp:comparison-edit-list} specifies $D$ deletions.
Hence, $\delta_i\ge0$ for $i\in[\nu]_{\Z}$ and $\sum_{i=1}^{\nu}\delta_i\le D$.
If $\delta_i=0$, no deletion decreases the grade of the traced edge, so its color is preserved.
We recover $(g,\gamma)$ from:
\begin{enumerate}[label=(\roman*)]
\item the differences $\delta_1,\ldots,\delta_\nu$;
\item the colors $\gamma(e_i)$ for $i\in[\nu]_{\Z}$ such that $\delta_i>0$;
\item the pairs $(g(v),\gamma(v))$ for $v\in D_U$.
\end{enumerate}
We have
\[
\begin{aligned}
 g(e_i)&=h(\widetilde e_i)+\delta_i &&(i\in[\nu]_{\Z}),\\
 \gamma(e_i)&=\eta(\widetilde e_i) &&(\delta_i=0),\\
 g(v)&=h(\widetilde v),\qquad \gamma(v)=\eta(\widetilde v)
 &&(v\in N(U)\setminus D_U).
\end{aligned}
\]

Let us count the difference vectors, horizontal colors, and grade--color pairs in (i)--(iii).
The vector $(\delta_1,\ldots,\delta_\nu)$ has
$\binom{\nu+D}{D}$ possibilities. At most $D$
original horizontal colors change. Each has at most $2^m$ choices.
We delete at most $D$ original vertical edges. Each has
$\sum_{r=0}^m\binom mr=2^m$ possible grade--color pairs.
For each $\xi\in\mathcal C_0(V_*)$, we obtain
\[
 \#\Phi^{-1}(\{\xi\})\le\binom{\nu+D}{D}(2^m)^D(2^m)^D=\operatorname{Fib}_m(\nu,D).
\]

\prf{Claim~\ref{c:composition-bound}}
By Claim~\ref{c:composition-target},
\[
 C_0(U)=\sum_{\xi\in\mathcal C_0(V_*)}\#\Phi^{-1}(\{\xi\})
 \le \operatorname{Fib}_m(\nu,D)C_0(V_*).
\]
\end{proof}

\begin{definition}\label{def:entropy-loss-function}
For $x>0$, define
\[
 \psi_m(x)=(1+x)\log(1+x)-x\log x+2mx\log2.
\]
Let
\[
 \psi_m(0)=0.
\]
\end{definition}

\begin{lemma}\label{lem:entropy-loss}
Let $\nu\in\Z_{\ge0}$, $D\in\Z_{\ge0}$, $N\in\R_{>0}$, and $x\in\R_{\ge0}$ satisfy $\nu\le N$ and $D\le xN$.
Then,
\begin{equation}\label{eq:entropy-loss}
 \log \operatorname{Fib}_m(\nu,D)\le N\psi_m(x).
\end{equation}
\end{lemma}
\begin{proof}
First, suppose $x=0$. Then, $D=0$ and $\log \operatorname{Fib}_m(\nu,0)=0=N\psi_m(0)$.

Second, suppose $x>0$. Let us prove~\eqref{eq:entropy-loss}.
Let $t=\frac{x}{1+x}\in(0,1)$.
The binomial expansion gives
\[
 \binom{\nu+D}{D}t^D(1-t)^\nu
 \le\bigl(t+(1-t)\bigr)^{\nu+D}=1.
\]
Since $-\log t>0$ and $-\log(1-t)>0$, the bounds $\nu\le N$ and $D\le xN$ give
\[
\begin{aligned}
 \log \operatorname{Fib}_m(\nu,D)
 &\le-D\log t-\nu\log(1-t)+2mD\log2\\
 &\le-xN\log t-N\log(1-t)+2mxN\log2\\
 &=N\psi_m(x).
\end{aligned}
\]
\end{proof}

%% file: sections/rational-slopes.tex
\section{Rational slopes}\label{sec:rational}
\subsection{Integer slopes for \texorpdfstring{$m\ge4$}{m >= 4}}
\begin{proposition}\label{prop:integer}
Let $m$ and $k$ be integers and let $M=m-2$. Assume the following.
\begin{hypotheses}
\item\label{hyp:integer-arrows} $m\ge4$.
\item\label{hyp:integer-range} $k\in[M]_{\Z}$.
\end{hypotheses}
Let
\[
 \mathcal S_m=
 \left\{(a,b)\in\Z^2:\
 \begin{gathered}
 b>a\ge2,\quad \gcd(a,b)=1,\\
 |R|_N\in[M]_{\Z}\text{ for every simple run }R\subseteq\mathcal U(a,b)
 \end{gathered}\right\}.
\]
Then, we have the following.

\begin{enumerate}
\item\label{c:integer-limit} We have
\[
 \lim_{\substack{a\to\infty,\ (a,b)\in\cA_m\cap\mathcal S_m\\\frac{b}{a}\to k}}
 \frac{\log E_m(a,b)}a=\ell(R_k).
\]
\item\label{c:integer-endpoints} The following one-sided limits hold.
\begin{enumerate}[label=(\alph*),ref=\theenumi(\alph*)]
\item\label{c:integer-left-endpoint} We have
\[
 \lim_{\substack{a\to\infty,\ (a,b)\in\cA_m\\\frac ba\downarrow1}}
 \frac{\log E_m(a,b)}a=\ell(R_1).
\]
\item\label{c:integer-right-endpoint} We have
\[
 \lim_{\substack{a\to\infty,\ (a,b)\in\cA_m\\\frac ba\uparrow M}}
 \frac{\log E_m(a,b)}a=\ell(R_M).
\]
\end{enumerate}
\end{enumerate}
\end{proposition}
\begin{proof}
\prf{Claim~\ref{c:integer-limit}} Let us estimate the number of vertical run-edits needed to obtain a constant run.
Consider $(a,b)\in\cA_m\cap\mathcal S_m$ as $a\to\infty$ and $\frac ba\to k$.
For all sufficiently large $a$, the first $a$ simple runs of $\mathcal{U}(a,b)$ have
vertical lengths $k-1$ and $k$ or $k$ and $k+1$. The $(a+1)$-st simple run of $\mathcal U(a,b)$ is $R_1$. The total vertical run-edit distance to $R_k^{a+1}$ is
\[
 \Delta=|b-1-ka|+(k-1)=o(a).
\]

Let us prove the upper bound. Let $\nu=a+1$ and $P=(m+2)\Delta+1$. We change the run $\mathcal{U}(a,b)$ to the constant run $R_k^\nu$ by vertical run-edits. Every added outer simple run has vertical length $k$. Hence, the target is $R_k^{\nu+2P+2}$.
Assumptions~\ref{hyp:integer-arrows}--\ref{hyp:integer-range} and $(a,b)\in\mathcal S_m$ give the length and insertion hypotheses for Claim~\ref{c:composition-bound} of Proposition~\ref{prop:composition}.
Equations~\eqref{eq:composition} and~\eqref{eq:dictionary} give
\[
 mE_m(a,b)\le \operatorname{Fib}_m(\nu,D)\,C_0(R_k^{\nu+2P+2}),\qquad D\le\Delta.
\]

Let us check the error term as $a\to\infty$ and $\frac ba\to k$ through $\cA_m\cap\mathcal S_m$.
Set $x=\frac{\Delta}{\nu}$. We have $x\to0$ and $\frac\nu a\to1$.
Since $x\log x\to0$ as $x\downarrow0$, Definition~\ref{def:entropy-loss-function} gives $\psi_m(x)\to0$.
Equation~\eqref{eq:entropy-loss}, applied with $N=\nu$, gives
\[
 0\le\frac{\log\operatorname{Fib}_m(\nu,D)}a
 \le\frac\nu a\psi_m\left(\frac\Delta\nu\right)\longrightarrow0.
\]
Also, $P=o(\nu)$ gives $\frac{\nu+2P+2}{a}\to1$ and $\nu+2P+2\to\infty$.
Claim~\ref{c:periodic-limit} of Proposition~\ref{prop:periodic} therefore gives
\[
\begin{aligned}
 \frac{\log E_m(a,b)}a
 &\le\frac{\log \operatorname{Fib}_m(\nu,D)}a
 +\frac{\nu+2P+2}{a}\frac{\log C_0(R_k^{\nu+2P+2})}{\nu+2P+2}
 -\frac{\log m}{a},\\
 \limsup_{\substack{a\to\infty,\ (a,b)\in\cA_m\cap\mathcal S_m\\\frac{b}{a}\to k}}
 \frac{\log E_m(a,b)}a&\le\ell(R_k).
\end{aligned}
\]

Let us prove the lower bound. Let $A_1,\ldots,A_\nu$ be simple runs such that $\mathcal U(a,b)=A_1\ccn \cdots\ccn A_\nu$.
Let $l=\nu-2P-2$. Since $\frac P\nu\to0$, we have $l\ge1$ for all sufficiently large $a$.
We start from $R_k^l$. For the current path $Z$,
we use the run-edits $Z\mapsto\operatorname{Ins}(Z,i,N)$ at $i\in[|Z|]_{\Z}$
or $Z\mapsto\operatorname{Del}(Z,v)$ at a vertical edge $v$.
At most $\Delta$ such run-edits give a translated copy of $A_{P+2}\ccn \cdots\ccn A_{\nu-P-1}$.
We use $A_1,\ldots,A_{P+1}$ and $A_{\nu-P},\ldots,A_\nu$ for the outer simple runs. Equations~\eqref{eq:composition} and~\eqref{eq:dictionary} give
\[
 mE_m(a,b)\ge\frac{C_0(R_k^l)}{\operatorname{Fib}_m(l,D)},\qquad D\le\Delta.
\]
Since $\frac{l}{a}\to1$ and $D=o(l)$, Claim~\ref{c:periodic-limit} of Proposition~\ref{prop:periodic} and~\eqref{eq:entropy-loss} give
\[
\begin{aligned}
 \frac{\log E_m(a,b)}a
 &\ge\frac{l}{a}\frac{\log C_0(R_k^l)}l
 -\frac{\log \operatorname{Fib}_m(l,D)}a-\frac{\log m}a,\\
 \liminf_{\substack{a\to\infty,\ (a,b)\in\cA_m\cap\mathcal S_m\\\frac{b}{a}\to k}}
 \frac{\log E_m(a,b)}a&\ge\ell(R_k).
\end{aligned}
\]

\prf{Claim~\ref{c:integer-left-endpoint}} Let us prove the one-sided limit as $\frac ba\downarrow1$. We have $1<\frac ba<M$ for all sufficiently large $a$.
Since $b-a\ge1$, equation~\eqref{eq:ambient} gives $1<\Gamma(a,b)<M$.
The floor differences in~\eqref{eq:frame-word} lie in $[M]_{\Z}$.
The $(a+1)$-st simple run of $\mathcal U(a,b)$ is $R_1$.
Hence, $(a,b)\in\mathcal S_m$.
Thus, Claim~\ref{c:integer-limit} applies with $k=1$.

\prf{Claim~\ref{c:integer-right-endpoint}} Let us prove the one-sided limit as $\frac ba\uparrow M$. We have $1<\Gamma(a,b)<\frac ba<M$ for all sufficiently large $a$ by~\eqref{eq:ambient}.
The floor differences in~\eqref{eq:frame-word} and the appended $R_1$ have vertical lengths in $[M]_{\Z}$.
Hence, $(a,b)\in\mathcal S_m$. Claim~\ref{c:integer-limit} applies with $k=M$.
\end{proof}

\Needspace{16\baselineskip}
\subsection{Mechanical run sequences near a noninteger rational slope}
\begin{definition}\label{def:mechanical-paths}
Let $m\in\Z_{\ge3}$, $p\in\Z_{>0}$, and $q\in\Z_{>0}$ satisfy $p\ge q$ and $\gcd(p,q)=1$.
\begin{enumerate}[font=\itshape]
\item\label{d:rational-period} Define
\[
 M_{\mathrm{run}}(m)=\begin{cases}2,&m=3,\\m-2,&m\ge4.\end{cases}
\]
For $i\in\Z_{>0}$, define
\[
 w_i(p,q)=\floor{\frac{ip}{q}}-\floor{\frac{(i-1)p}{q}}.
\]
We define the \textit{rational run} and the \textit{rational run sequence} by
\[
 W_{p,q}=R_{w_1(p,q)}\ccn\cdots\ccn R_{w_q(p,q)},\qquad
 \mathcal R_{p,q}=(R_{w_i(p,q)})_{i\in\Z_{>0}}.
\]
For $j\in\Z_{>0}$ and $\alpha=\frac pq\in(1,M_{\mathrm{run}}(m))$, define
\begin{equation}\label{eq:rational-period}
 A_j(p,q)=C_0(W_{p,q}^j),\qquad \Upsilon_m(\alpha)=\frac{\ell(W_{p,q})}q.
\end{equation}
\item\label{d:mechanical-sequence} Let $s\ge1$ and $u\in[0,1)$ be real numbers. Define
\[
 \operatorname{MechRun}(s,u)=(R_{j_i})_{i\in\Z_{>0}},\qquad
 j_i=\floor{u+is}-\floor{u+(i-1)s}.
\]
We call $\operatorname{MechRun}(s,u)$ the
\textit{mechanical run sequence} of slope $s$ and intercept $u$.
\item Let $\alpha=\frac pq$.
For $h\in\Z_{\ge0}$, define
\[
 \operatorname{RunPhase}_q(h)=h-q\floor{\frac hq}.
\]
We call $\operatorname{RunPhase}_q(h)$ the \textit{run phase} at the start of the $(h+1)$-st simple run of $\mathcal R_{p,q}$.
\item Let $n\in\Z_{>0}$ and let $A_1,\ldots,A_n$ be simple runs.
Define $W=A_1\ccn \cdots\ccn A_n$.
Let $r\in\Z_{>0}$ and $\ell\in\Z_{>0}$ satisfy $r+\ell-1\le n$.
Define the \textit{$(r,\ell)$-run-factor} of $W$ by
\[
 \operatorname{RunFactor}_{r,\ell}(W)
 =W\left[\sum_{i=1}^{r-1}|A_i|,\ \sum_{i=1}^{r+\ell-1}|A_i|\right].
\]

Let $\mathcal A=(A_i)_{i\in\Z_{>0}}$ be a sequence of simple runs.
For $r\in\Z_{>0}$ and $\ell\in\Z_{>0}$, define
\[
 \operatorname{RunFactor}_{r,\ell}(\mathcal A)
 =\operatorname{RunFactor}_{r,\ell}(A_1\ccn \cdots\ccn A_{r+\ell-1}).
\]
\end{enumerate}
\end{definition}

\begin{lemma}\label{lem:rational-mechanical}
Let $m\in\Z_{\ge3}$.
Let $p\in\Z_{>0}$ and $q\in\Z_{>0}$ satisfy $\gcd(p,q)=1$, $q\ge2$, and
$1<\frac pq<M_{\mathrm{run}}(m)$.
Let $\alpha=\frac pq$.
Let $s\ge1$ and $u\in[0,1)$. Assume the following.
\begin{hypotheses}
\item\label{hyp:rational-mechanical-slope} $|s-\alpha|<q^{-2}$.
\end{hypotheses}
Then, we have the following.

\begin{enumerate}
\item\label{c:rational-mechanical-short} For every $r\in\Z_{>0}$ and $l\in[q-1]_{\Z}$, there exists $t\in\Z_{>0}$ such that the path
\[
 \operatorname{RunFactor}_{r,l}(\operatorname{MechRun}(s,u))
\]
is a translated copy of $\operatorname{RunFactor}_{t,l}(\mathcal R_{p,q})$.
\item\label{c:rational-mechanical-period} For $r\in\Z_{>0}$, let
$Z_r=\operatorname{RunFactor}_{r,q}(\operatorname{MechRun}(s,u))$.
If $|Z_r|_N=p$, then $Z_r$ is a translated copy of a cyclic rotation
of $W_{p,q}$ by simple runs.
\item\label{c:rational-mechanical-lengths} Let $k=\floor\alpha$.
Then, $w_i(p,q)\in[k,k+1]_{\Z}$ for $i\in[q]_{\Z}$ and $w_1(p,q)=k$.
\item\label{c:rational-mechanical-zero-phase} For $i\in[0,q-1]_{\Z}$, we have
\[
 \floor{is}=\floor{i\alpha}.
\]
\end{enumerate}
\end{lemma}
\begin{proof}
\prf{Claim~\ref{c:rational-mechanical-short}} Let us locate the floor thresholds within one period.
For $i\in[q-1]_{\Z}$, define
$T_i(t)=\ceil{\frac{ip}{q}}-it$ for $t\in\R$.
Let us prove
\begin{equation}\label{eq:mechanical-threshold-order}
 T_i(s)<T_j(s)\quad\Longleftrightarrow\quad T_i(\alpha)<T_j(\alpha)
 \qquad(i,j\in[q-1]_{\Z}).
\end{equation}
Since $\alpha=\frac pq$, the numbers $T_i(\alpha)$
are distinct elements of $\{\frac1q,\ldots,\frac{q-1}{q}\}$.
By assumption~\ref{hyp:rational-mechanical-slope},
\[
 |T_i(s)-T_i(\alpha)|=i|s-\alpha|<\frac1q,
 \qquad 0<T_i(s)<1.
\]
If $i\in[q-1]_{\Z}$ and $j\in[q-1]_{\Z}$ are distinct, then we also have
\[
 \left|(T_i(s)-T_j(s))-(T_i(\alpha)-T_j(\alpha))\right|
 =|i-j|\,|s-\alpha|<\frac1q\le|T_i(\alpha)-T_j(\alpha)|.
\]
Hence, equation~\eqref{eq:mechanical-threshold-order} holds.
For $t\in\{s,\alpha\}$ and $v\in[0,1)$,
\[
 \floor{v+it}=
 \begin{cases}
 \ceil{i\alpha}-1,&v<T_i(t),\\
 \ceil{i\alpha},&v\ge T_i(t).
 \end{cases}
\]
Thus, equation~\eqref{eq:mechanical-threshold-order} gives, for every $l\in[q-1]_{\Z}$,
\begin{equation}\label{eq:mechanical-prefix-vectors}
 \left\{(\floor{v+is})_{i\in[l]_{\Z}}:v\in[0,1)\right\}
 =\left\{(\floor{v+i\alpha})_{i\in[l]_{\Z}}:v\in[0,1)\right\}.
\end{equation}

Let us apply~\eqref{eq:mechanical-prefix-vectors} to $\operatorname{RunFactor}_{r,l}(\operatorname{MechRun}(s,u))$.
Let $v_r=u+(r-1)s-\floor{u+(r-1)s}$.
For $i\in[l]_{\Z}$,
\begin{equation}\label{eq:mechanical-factor-prefix}
\begin{aligned}
 \left|\operatorname{RunFactor}_{r,i}(\operatorname{MechRun}(s,u))\right|_N
 &=\floor{u+(r+i-1)s}-\floor{u+(r-1)s}\\
 &=\floor{v_r+is}.
\end{aligned}
\end{equation}
Every vector in the right-hand side of~\eqref{eq:mechanical-prefix-vectors} occurs at an intercept in
$\{0,\frac1q,\ldots,\frac{q-1}{q}\}$.
Since $\gcd(p,q)=1$,
\[
 \left\{(t-1)\alpha-\floor{(t-1)\alpha}:t\in[q]_{\Z}\right\}
 =\left\{0,\frac1q,\ldots,\frac{q-1}{q}\right\}.
\]
Equations~\eqref{eq:mechanical-prefix-vectors} and~\eqref{eq:mechanical-factor-prefix} therefore give
$t\in[q]_{\Z}$ such that
\[
 \left|\operatorname{RunFactor}_{r,i}(\operatorname{MechRun}(s,u))\right|_N
 =\left|\operatorname{RunFactor}_{t,i}(\mathcal R_{p,q})\right|_N
 \qquad(i\in[l]_{\Z}).
\]
Taking consecutive differences shows that
$\operatorname{RunFactor}_{r,l}(\operatorname{MechRun}(s,u))$
is a translated copy of $\operatorname{RunFactor}_{t,l}(\mathcal R_{p,q})$.

\prf{Claim~\ref{c:rational-mechanical-period}} Let us determine $Z_r$ from its first $q-1$ simple runs and its total vertical length.
Claim~\ref{c:rational-mechanical-short} gives $t\in\Z_{>0}$ such that the first $q-1$ simple runs
of $Z_r$ are a translated copy of $\operatorname{RunFactor}_{t,q-1}(\mathcal R_{p,q})$.
Since $q\alpha=p\in\Z$, we have
\[
 \left|\operatorname{RunFactor}_{t,q}(\mathcal R_{p,q})\right|_N
 =\floor{(t+q-1)\alpha}-\floor{(t-1)\alpha}=p.
\]
The assumption $|Z_r|_N=p$ therefore determines the last simple-run length.
Thus, $Z_r$ is a translated copy of $\operatorname{RunFactor}_{t,q}(\mathcal R_{p,q})$.
The identity $\floor{z+p}=\floor z+p$ for $z\in\R$ shows that the simple-run lengths of $\mathcal R_{p,q}$
are $q$-periodic. Hence, $\operatorname{RunFactor}_{t,q}(\mathcal R_{p,q})$ is a cyclic rotation of $W_{p,q}$ up to translation.

\prf{Claim~\ref{c:rational-mechanical-lengths}} Let us compute the periodic run lengths.
Let $\beta=\alpha-k\in[0,1)$.
For $i\in[q]_{\Z}$, we have
\[
 w_i(p,q)=k+\floor{i\beta}-\floor{(i-1)\beta}\in[k,k+1]_{\Z}.
\]
At $i=1$, we have $w_1(p,q)=k+\floor\beta-\floor0=k$.

\prf{Claim~\ref{c:rational-mechanical-zero-phase}} Let us bound the distance of $i\alpha$ from the integers.
For $i=0$, both floors are zero.
For $i\in[q-1]_{\Z}$, coprimality of $p$ and $q$ gives
\[
 \floor{i\alpha}+\frac1q\le i\alpha
 \le\floor{i\alpha}+1-\frac1q.
\]
Assumption~\ref{hyp:rational-mechanical-slope} gives
\[
 |is-i\alpha|<\frac{i}{q^2}<\frac1q.
\]
Hence, $\floor{i\alpha}<is<\floor{i\alpha}+1$, so $\floor{is}=\floor{i\alpha}$.
\end{proof}

\begin{definition}\label{sync:rational:264}
Let $m\in\Z_{\ge3}$, $p\in\Z_{>0}$, and $q\in\Z_{>0}$ satisfy $\gcd(p,q)=1$ and $q\ge2$.
Let $\alpha=\frac pq\in(1,M_{\mathrm{run}}(m))$.
Let $n\in\Z_{>0}$ and let $A_1,\ldots,A_n$ be simple runs.
Define $W=A_1\ccn \cdots\ccn A_n$.
For $r\in[n-q+1]_{\Z}$, define
\[
 Z_r=\operatorname{RunFactor}_{r,q}(W).
\]
We call $Z_r$ \textit{$\alpha$-bad} if $|Z_r|_N\ne p$.
Define
\[
 D_{\mathrm{bad}}(W,p,q)=\#\{r\in[n-q+1]_{\Z}:|Z_r|_N\ne p\}.
\]
\end{definition}

\begin{lemma}[Counting exceptional floor differences]\label{lem:floor-discrepancies}
Let $q\in\Z_{>0}$, $N\in\Z_{>0}$, $p\in\Z$, $s\in\R$, and $u\in\R$.
Assume $|qs-p|<1$.
Then,
\[
 \#\{i\in[0,N-1]_{\Z}:\floor{u+(i+q)s}-\floor{u+is}\ne p\}
 \le N|qs-p|+q.
\]
\end{lemma}
\begin{proof}
Let us express the exceptional floors as jumps of an integer sequence.
Let $\epsilon=qs-p$ and define
\[
 \delta_i=\floor{u+(i+q)s}-\floor{u+is}-p
 \qquad(i\in[0,N-1]_{\Z}).
\]
Since $p\in\Z$ and $|\epsilon|<1$,
\[
 \delta_i=\floor{u+is+\epsilon}-\floor{u+is}\in
 \begin{cases}
 [0,1]_{\Z},&\epsilon\ge0,\\
 [-1,0]_{\Z},&\epsilon<0.
 \end{cases}
\]
Thus, $\delta_i\ne0$ implies $\delta_i=1$ when $\epsilon\ge0$ and $\delta_i=-1$ when $\epsilon<0$.

Let us sum the floor differences in each residue class.
For $r\in[0,q-1]_{\Z}$, let
$n_r=\#\{i\in[0,N-1]_{\Z}:i\equiv r\pmod q\}$.
Let us compute $\sum_{t=0}^{n_r-1}\delta_{r+tq}$.
Telescoping gives
\[
\begin{aligned}
 \sum_{t=0}^{n_r-1}\delta_{r+tq}
 &=\floor{u+(r+n_rq)s}-\floor{u+rs}-n_rp,\\
 \left|\sum_{t=0}^{n_r-1}\delta_{r+tq}-n_r\epsilon\right|&<1.
\end{aligned}
\]
Since $\delta_i\in\{0,1\}$ for $\epsilon\ge0$ and $\delta_i\in\{-1,0\}$ for $\epsilon<0$, we obtain
\[
\begin{aligned}
 \#\{i\in[0,N-1]_{\Z}:\delta_i\ne0\}
 &=\sum_{r=0}^{q-1}\left|\sum_{t=0}^{n_r-1}\delta_{r+tq}\right|\\
 &\le\sum_{r=0}^{q-1}(n_r|\epsilon|+1)=N|\epsilon|+q.
\end{aligned}
\]
\end{proof}

\begin{lemma}\label{lem:bad-factors}
Let $m\in\Z_{\ge3}$, $p\in\Z_{>0}$, and $q\in\Z_{>0}$ satisfy $\gcd(p,q)=1$ and $q\ge2$.
Let $\alpha=\frac pq\in(1,M_{\mathrm{run}}(m))$.
Let $b>a\ge2$ be coprime integers.
Then, we have the following.
\begin{enumerate}
\item\label{c:bad-factors-bound} Assume the following.
\begin{hypotheses}
\item\label{hyp:bad-factor-length} $a\ge q$.
\item\label{hyp:bad-factor-slope} $|\Gamma(a,b)-\alpha|<q^{-2}$.
\end{hypotheses}
Then,
\begin{equation}\label{eq:bad-windows}
 D_{\mathrm{bad}}(\mathcal U(a,b),p,q)\le |qb-pa|+q+2.
\end{equation}
\item\label{c:bad-factors-negligible} We have
\[
 D_{\mathrm{bad}}(\mathcal U(a,b),p,q)=o(a)
 \qquad\left(a\to\infty,\ \frac ba\to\alpha,\ b>a,\ \gcd(a,b)=1\right).
\]
\end{enumerate}
\end{lemma}
\begin{proof}
\prf{Claim~\ref{c:bad-factors-bound}} Let us bound the number of $\alpha$-bad run-factors contained in the first $a$ simple runs of $\mathcal U(a,b)$.
Assumptions~\ref{hyp:bad-factor-length} and~\ref{hyp:bad-factor-slope} give $a-q+1\ge1$ and $|q\Gamma(a,b)-p|<1$.
Lemma~\ref{lem:floor-discrepancies} with $N=a-q+1$, $s=\Gamma(a,b)$, and $u=0$ gives
\begin{equation}\label{eq:frame-floor-exceptions}
\begin{aligned}
 &\#\left\{i\in[0,a-q]_{\Z}:\floor{(i+q)\Gamma(a,b)}-\floor{i\Gamma(a,b)}\ne p\right\}\\
 &\qquad\le(a-q+1)|q\Gamma(a,b)-p|+q.
\end{aligned}
\end{equation}
For $i\in[0,a-q]_{\Z}$, Definition~\ref{def:frame} gives
\[
 \left|\operatorname{RunFactor}_{i+1,q}(\mathcal U(a,b))\right|_N
 =\floor{(i+q)\Gamma(a,b)}-\floor{i\Gamma(a,b)}.
\]
By~\eqref{eq:frame-word}, the $(a+1)$-st simple run of $\mathcal U(a,b)$ is $R_1$.
For $r\in[a-q+2]_{\Z}$, the path $\operatorname{RunFactor}_{r,q}(\mathcal U(a,b))$
contains the $(a+1)$-st simple run if and only if $r=a-q+2$.
Hence, equation~\eqref{eq:frame-floor-exceptions} gives
\[
\begin{aligned}
 D_{\mathrm{bad}}(\mathcal U(a,b),p,q)
 &\le 1+(a-q+1)|q\Gamma(a,b)-p|+q\\
 &\le aq|\Gamma(a,b)-\alpha|+q+1.
\end{aligned}
\]
Assumption~\ref{hyp:bad-factor-length} and~\eqref{eq:ambient} give $\mu_1(a,b)\ge a+1>q$ and
\[
 aq|\Gamma(a,b)-\alpha|
 =\left|qb-pa-\frac q{\mu_1(a,b)}\right|
 \le |qb-pa|+\frac q{\mu_1(a,b)}.
\]
Since $\frac q{\mu_1(a,b)}<1$, inequality~\eqref{eq:bad-windows} follows.

\prf{Claim~\ref{c:bad-factors-negligible}} Let us verify the range of the ambient slope.
Equation~\eqref{eq:ambient} gives
\[
 |\Gamma(a,b)-\alpha|\le\left|\frac ba-\alpha\right|+\frac1{a(a+1)}\longrightarrow0.
\]
For all sufficiently large $a$, Claim~\ref{c:bad-factors-bound} applies and gives
\[
 0\le\frac{D_{\mathrm{bad}}(\mathcal U(a,b),p,q)}a
 \le\left|q\frac ba-p\right|+\frac{q+2}a\longrightarrow0.
\]
\end{proof}

\begin{lemma}\label{lem:separated-subwords}
Let $m\in\Z_{\ge3}$, $p\in\Z_{>0}$, and $q\in\Z_{>0}$ satisfy $\gcd(p,q)=1$ and $q\ge2$.
Let $\alpha=\frac pq\in(1,M_{\mathrm{run}}(m))$.

\begin{enumerate}
\item\label{c:separated-subwords-rotations} For every integer $J\ge4$, each cyclic rotation $V$ of $W_{p,q}^J$
by simple runs satisfies $C_0(V)\ge A_{J-3}(p,q)$.
\item\label{c:separated-subwords-construction} Let $h\ge1$ be an integer.
Let $X_1,\ldots,X_h$ be runs.
Let $B_0,\ldots,B_h$ each consist of three simple runs. Assume the following.
\begin{hypotheses}
\item\label{hyp:separated-lengths} Every simple run of $B_0,\ldots,B_h$ has vertical length in $[M_{\mathrm{run}}(m)]_{\Z}$.
\end{hypotheses}
Then,
\begin{equation}\label{eq:separated-subwords}
 C_0(B_0\ccn X_1\ccn B_1\ccn \cdots\ccn X_h\ccn B_h)\ge\prod_{i=1}^h C_0(X_i).
\end{equation}
\end{enumerate}
\end{lemma}
\begin{proof}
\prf{Claim~\ref{c:separated-subwords-rotations}} Let us apply Lemma~\ref{lem:outer-shell-extension} to a rotated periodic path.
Let $t=J-3$ and let $j\in[0,q-1]_{\Z}$.
Let $U_j=\operatorname{CRot}(W_{p,q},j)$.
Let $P=2q-j$ and $S=q+j$.
Let $X=\operatorname{RunFactor}_{P+1,qt}(U_j^{t+3})$. Then,
\[
 X=\operatorname{Trans}(W_{p,q}^t,\operatorname{Start}_{\mathrm{pt}}(X)),\qquad
 P+qt+S=q(t+3).
\]
Since $P+j=2q$, we have $\operatorname{RunPhase}_q(P+j)=0$.
Also, $\max(P,S)\le2\min(P,S)$.
Claim~\ref{c:rational-mechanical-lengths} of Lemma~\ref{lem:rational-mechanical} gives
$|R|_N\in[M_{\mathrm{run}}(m)]_{\Z}$ for every simple run $R$ of $W_{p,q}$.
Lemma~\ref{lem:outer-shell-extension} gives
\[
 C_0(U_j^{t+3})\ge C_0(W_{p,q}^t)=A_{J-3}(p,q).
\]

\prf{Claim~\ref{c:separated-subwords-construction}} Let us specify the shell paths.
Let $B_i=R_{b_{i,1}}\ccn R_{b_{i,2}}\ccn R_{b_{i,3}}$ for $i\in[0,h]_{\Z}$.
For $r\in[h]_{\Z}$ and $s\in[h]_{\Z}$ satisfying $r\le s$ and
$(u,v)\in\{(2,2),(2,1),(1,3),(3,2)\}$, define
\[
 T_{uv}(r,s)=
 \left(\mathop{\ccn}\limits_{j=4-u}^{3}R_{b_{r-1,j}}\right)
 \ccn X_r\ccn\left(\mathop{\ccn}\limits_{i=r+1}^{s}\bigl(B_{i-1}\ccn X_i\bigr)\right)
 \ccn\left(\mathop{\ccn}\limits_{j=1}^{v}R_{b_{s,j}}\right).
\]
The iterated concatenations use increasing index order.
We compare paths up to translation using Lemma~\ref{lem:word-reindexing}.
For $r=s$, we have
\[
\begin{aligned}
 T_{22}(r,r)&=\operatorname{Shell}
 (\operatorname{Shell}(X_r;b_{r-1,3},b_{r,1},0);b_{r-1,2},b_{r,2},0),\\
 T_{21}(r,r)&=\operatorname{Shell}(R_{b_{r-1,3}},X_r;b_{r-1,2},b_{r,1},0),\\
 T_{13}(r,r)&=\operatorname{Shell}(X_r,R_{b_{r,1}}\ccn R_{b_{r,2}};b_{r-1,3},b_{r,3},0),\\
 T_{32}(r,r)&=\operatorname{Shell}(T_{21}(r,r);b_{r-1,1},b_{r,2},0).
\end{aligned}
\]
For $r<s$, we have
\[
\begin{aligned}
 T_{22}(r,s)&=\operatorname{Shell}
 (\operatorname{Shell}(X_r;b_{r-1,3},b_{r,1},0),T_{21}(r+1,s);b_{r-1,2},b_{s,2},0),\\
 T_{21}(r,s)&=\operatorname{Shell}(T_{13}(r,s-1),X_s;b_{r-1,2},b_{s,1},0),\\
 T_{13}(r,s)&=\operatorname{Shell}(X_r,T_{32}(r+1,s);b_{r-1,3},b_{s,3},0),\\
 T_{32}(r,s)&=\operatorname{Shell}
 (T_{22}(r,s-1),\operatorname{Shell}(X_s;b_{s-1,3},b_{s,1},0);b_{r-1,1},b_{s,2},0).
\end{aligned}
\]

Let us count the resulting colored gradings.
Assumption~\ref{hyp:separated-lengths} gives $b_{i,j}\in[M_{\mathrm{run}}(m)]_{\Z}\subseteq[m-1]_{\Z}$.
Lemma~\ref{lem:word-bounds} gives $C_0(R_{b_{i,j}})\ge1$ and
$C_0(R_{b_{i,1}}\ccn R_{b_{i,2}})\ge1$.
Claim~\ref{c:shell-recolor} of Lemma~\ref{lem:shell} gives an injective map
into colored gradings of surplus zero at each shell step.
Induction on $s-r$ therefore gives
\[
 C_0(T_{uv}(r,s))\ge\prod_{i=r}^{s}C_0(X_i).
\]
Finally,
\[
 \operatorname{Shell}(T_{22}(1,h);b_{0,1},b_{h,3},0)
 =B_0\ccn X_1\ccn B_1\ccn \cdots\ccn X_h\ccn B_h
\]
up to translation. Claim~\ref{c:shell-recolor} of Lemma~\ref{lem:shell}
and Lemma~\ref{lem:word-reindexing} give~\eqref{eq:separated-subwords}.
\end{proof}

\begin{theorem}[Noninteger rational all-sequence limit]\label{thm:rational}
Let $m\in\Z_{\ge3}$, $p\in\Z_{>0}$, and $q\in\Z_{>0}$ satisfy $\gcd(p,q)=1$ and $q\ge2$.
Let $\alpha=\frac pq\in(1,M_{\mathrm{run}}(m))$.
Then,
\begin{equation}\label{eq:rational-limit}
 \lim_{\substack{a\to\infty,\ (a,b)\in\cA_m\\\frac{b}{a}\to\alpha}}
 \frac{\log E_m(a,b)}a=\frac{\ell(W_{p,q})}q.
\end{equation}
\end{theorem}
\begin{proof}
Let us construct the initial factor and the separated periodic factors.
Let $k=\floor\alpha$.
In the limit in~\eqref{eq:rational-limit}, we have
$a\ge q$ and $|\Gamma(a,b)-\alpha|<q^{-2}$ for all sufficiently large $a$.
\prf{the upper bound}
First, we divide $\mathcal{U}(a,b)$ into periodic pieces.
Equation~\eqref{eq:frame-word} and Claims~\ref{c:rational-mechanical-lengths} and~\ref{c:rational-mechanical-zero-phase}
of Lemma~\ref{lem:rational-mechanical} give
\[
 \operatorname{RunFactor}_{1,1}(\mathcal U(a,b))
 =\operatorname{RunFactor}_{1,1}(W_{p,q})=R_k,\qquad
 \operatorname{RunPhase}_q(0)=0.
\]
Let $0=c_0<c_1<\cdots<c_d=a$ satisfy
\[
 \{c_1,\ldots,c_{d-1}\}
 =\{r\in[a-q+1]_{\Z}:|\operatorname{RunFactor}_{r,q}(\mathcal U(a,b))|_N\ne p\}.
\]
For $i\in[d]_{\Z}$, let
\[
 X_i=\operatorname{RunFactor}_{c_{i-1}+1,c_i-c_{i-1}}(\mathcal U(a,b)).
\]
Each $X_i$ contains no $\alpha$-bad run-factor. We have
$d\le D_{\mathrm{bad}}(\mathcal U(a,b),p,q)+1$.
By Claim~\ref{c:rational-mechanical-short} of Lemma~\ref{lem:rational-mechanical}, a piece with $l<q$ simple runs
is a translated copy of $\operatorname{RunFactor}_{t,l}(\mathcal R_{p,q})$ for some $t\in\Z_{>0}$.
In a piece with at least $q$ simple runs, every run-factor consisting of $q$ simple runs has
vertical length $p$. Subtracting the sums for consecutive run-factors gives
\[
 0=\sum_{i=j+1}^{j+q}\Lambda_i(a,b)-\sum_{i=j}^{j+q-1}\Lambda_i(a,b)
 =\Lambda_{j+q}(a,b)-\Lambda_j(a,b).
\]
Thus, $\Lambda_{j+q}(a,b)=\Lambda_j(a,b)$ for $i\in[d]_{\Z}$ and $c_{i-1}<j<j+q\le c_i$.
Claim~\ref{c:rational-mechanical-period} of Lemma~\ref{lem:rational-mechanical} identifies the first $q$ simple runs of each such piece with a cyclic rotation of $W_{p,q}$.

Second, let us align each piece with a run-factor of $\mathcal R_{p,q}$.
Let us prove that no insertion is needed before the first piece.
Claim~\ref{c:rational-mechanical-zero-phase} of Lemma~\ref{lem:rational-mechanical}, applied with $s=\Gamma(a,b)$, gives
\[
 \floor{i\Gamma(a,b)}=\floor{i\alpha}\qquad(i\in[q-1]_{\Z}).
\]
If the first piece has at least $q$ simple runs, its first $q$ lengths sum to $p$.
Thus, its $q$-th length also agrees with $w_q(p,q)$.
The identities $\Lambda_{j+q}(a,b)=\Lambda_j(a,b)$ for $1\le j<j+q\le c_1$ identify $X_1$ with $\operatorname{RunFactor}_{1,c_1}(\mathcal R_{p,q})$.
Consequently, no insertion is needed before $X_1$.

Let us proceed by induction through $X_2,\ldots,X_d$.
For the current piece $X$, choose the least $t\in[q]_{\Z}$ such that
$X$ is a translated copy of $\operatorname{RunFactor}_{t,|X|_E}(\mathcal R_{p,q})$.
Let $P\ccn X\ccn Z$ be the current path and let $h=|P|_E$.
The induction hypothesis is that the simple-run lengths of $P$ are $w_1(p,q),\ldots,w_h(p,q)$.
The identification of $X_1$ gives the base case. Define
\[
 j=\operatorname{RunPhase}_q\bigl(q+t-1-\operatorname{RunPhase}_q(h)\bigr)
 \in[0,q-1]_{\Z}.
\]
Then, $\operatorname{RunPhase}_q(h+j)=t-1$.
If $j>0$, insert the simple runs of $\operatorname{RunFactor}_{h+1,j}(\mathcal R_{p,q})$
after $P$, producing
\[
 \operatorname{Ins}\bigl(P\ccn X\ccn Z,|P|,\operatorname{RunFactor}_{h+1,j}(\mathcal R_{p,q})\bigr).
\]
If $j=0$, the path is $P\ccn X\ccn Z$.

By $q$-periodicity, the simple-run lengths of $X$ now agree with those at positions $h+j+1,\ldots,h+j+|X|_E$ in $\mathcal R_{p,q}$.
Thus, the prefix through $X$ has simple-run lengths $w_1(p,q),\ldots,w_{h+j+|X|_E}(p,q)$, which proves the induction step.

Let us complete the final period before $R_1$.
Let $h$ be the number of simple runs before the final $R_1$ of $\mathcal U(a,b)$ after these insertions.
Let $j=\operatorname{RunPhase}_q\bigl(q-\operatorname{RunPhase}_q(h)\bigr)$.
If $j>0$, insert $\operatorname{RunFactor}_{h+1,j}(\mathcal R_{p,q})$ before $R_1$.
Then, $\operatorname{RunPhase}_q(h+j)=0$.
Let $K$ be the number of inserted simple runs. Let $n_0=\frac{a+K}{q}$.
The result is
\[
 W_{p,q}^{n_0}\ccn R_1,\qquad
 K\le(q-1)(D_{\mathrm{bad}}(\mathcal U(a,b),p,q)+2).
\]
For each insertion, use~\eqref{eq:run-insertion-map},
\eqref{eq:run-insertion-grades}, and~\eqref{eq:run-insertion-recolor}.
The resulting grading $g$ satisfies $\operatorname{Surp}(W_{p,q}^{n_0}\ccn R_1,g)=K$.

Third, let us add the outer simple runs.
Apply Lemma~\ref{lem:outer-shell-extension} with $(P,S)=(q,q-1)$
to add $W_{p,q}$ before $W_{p,q}^{n_0}\ccn R_1$ and $R_{w_2(p,q)}\ccn \cdots\ccn R_{w_q(p,q)}$ after $W_{p,q}^{n_0}\ccn R_1$.
The target is
\begin{equation}\label{eq:rational-initial-extension}
 Q_0=W_{p,q}^{n_0+1}\ccn R_1\ccn R_{w_2(p,q)}\ccn \cdots\ccn R_{w_q(p,q)}.
\end{equation}
If $k=1$, we have $Q_0=W_{p,q}^{n_0+2}$. For $m=3$, the condition $1<\alpha<2$ gives $k=1$.

Let us obtain the initial run of length $k$.
For $k>1$, let us replace the factor $R_1$ in~\eqref{eq:rational-initial-extension} by $R_k$.
We have $m\ge4$. Apply Claims~\ref{c:vertical-insert-connected}--\ref{c:vertical-insert-injective} of Lemma~\ref{lem:vertical-insert} $k-1$ times.
For each application, use $q$ unary shell paths to add one copy of $W_{p,q}$ on each side.
After the $j$th vertical-edge insertion, the path is
\[
 Q_j=W_{p,q}^{n_0+j+1}\ccn R_{j+1}\ccn R_{w_2(p,q)}\ccn \cdots\ccn R_{w_q(p,q)}\ccn W_{p,q}^j
 \qquad(j\in[k-1]_{\Z}).
\]
The transported shell grading $h_j$ satisfies $\operatorname{Surp}(Q_j,h_j)=K$.

Let us reduce the surplus to zero.
We have $Q_{k-1}=W_{p,q}^{n_0+2k}$.
Let $H=\ceil{\frac Kq}$.
Apply Lemma~\ref{lem:surplus-reduction} with $Q=W_{p,q}^{n_0+2k}$, $n=qH$, and
\[
 a_j=w_{q-\operatorname{RunPhase}_q(j-1)}(p,q),\qquad
 b_j=w_{1+\operatorname{RunPhase}_q(j-1)}(p,q)
 \qquad(j\in[qH]_{\Z}).
\]
The final path is a translated copy of $W_{p,q}^{n_0+2k+2H}$.

Let us prove injectivity of the composition.
At each simple-run insertion, equations~\eqref{eq:run-insertion-grades} and~\eqref{eq:run-insertion-recolor} recover the input grading and coloring.
At each vertical-edge insertion, Claim~\ref{c:vertical-insert-injective} of Lemma~\ref{lem:vertical-insert} gives injectivity.
Lemma~\ref{lem:outer-shell-extension} gives injectivity of each outer extension.
Lemma~\ref{lem:word-reindexing} and the injection~\eqref{eq:surplus-reduction-map} give
\begin{equation}\label{eq:rational-upper}
 mE_m(a,b)\le A_N(p,q),\qquad
 N=\frac{a+K}q+2k+2\ceil{\frac{K}{q}}.
\end{equation}

Let us obtain the upper limit from~\eqref{eq:rational-upper}.
By Claim~\ref{c:bad-factors-negligible} of Lemma~\ref{lem:bad-factors}, we have $K=o(a)$ with $q$ fixed. Hence, $\frac{N}{a}\to\frac1q$. Claim~\ref{c:periodic-limit} of Proposition~\ref{prop:periodic} gives
\[
 \begin{aligned}
 \frac{\log E_m(a,b)}a&\le\frac{N}{a}\frac{\log A_N(p,q)}{N}-\frac{\log m}{a},\\
 \limsup_{\substack{a\to\infty,\ (a,b)\in\cA_m\\\frac{b}{a}\to\alpha}}
 \frac{\log E_m(a,b)}a&\le\frac{\ell(W_{p,q})}q.
 \end{aligned}
\]

\prf{the lower bound}
Let us select separated subpaths on which to count gradings.
Let $J\ge4$ be an integer and let
\[
 h=\floor{\frac{a-2}{Jq+3}}
\]
if $a$ is sufficiently large and $h\ge1$. Define
\[
\begin{aligned}
 X_i&=\operatorname{RunFactor}_{4+(i-1)(Jq+3),Jq}(\mathcal U(a,b))
       &&(i\in[h-1]_{\Z}),\\
 X_h&=\operatorname{RunFactor}_{4+(h-1)(Jq+3),\,a-5-(h-1)(Jq+3)}(\mathcal U(a,b)),\\
 B_i&=\operatorname{RunFactor}_{i(Jq+3)+1,3}(\mathcal U(a,b))
       &&(i\in[0,h-1]_{\Z}),\\
 B_h&=\operatorname{RunFactor}_{a-1,3}(\mathcal U(a,b)).
\end{aligned}
\]
Then, $\mathcal U(a,b)=B_0\ccn X_1\ccn B_1\ccn \cdots\ccn X_h\ccn B_h$ and $|X_h|_E\ge Jq$.

Let us count colored gradings on these subpaths.
Let
\[
 \mathcal I=\{h\}\cup\{i\in[h-1]_{\Z}:D_{\mathrm{bad}}(X_i,p,q)>0\}.
\]
Then, $\#\mathcal I\le D_{\mathrm{bad}}(\mathcal U(a,b),p,q)+1$.
We call $X_i$ exceptional for $i\in\mathcal I$.
For $i\in[h]_{\Z}\setminus\mathcal I$, the path $X_i$ is a translated copy of $\operatorname{RunFactor}_{t,Jq}(\mathcal R_{p,q})$ for some $t\in\Z_{>0}$. By Claim~\ref{c:separated-subwords-rotations} of Lemma~\ref{lem:separated-subwords}, such a subpath has at least $A_{J-3}(p,q)$ connected colored gradings of surplus zero.
Each exceptional subpath $X$ has a connected colored grading $(g,\gamma)$ with
$\operatorname{Surp}(X,g)=0$ by Lemma~\ref{lem:word-bounds}.
Claim~\ref{c:separated-subwords-construction} of Lemma~\ref{lem:separated-subwords} gives
\[
 mE_m(a,b)\ge A_{J-3}(p,q)^{\max(0,h-D_{\mathrm{bad}}(\mathcal U(a,b),p,q)-1)}.
\]

Let us take the lower limit and then let $J$ tend to infinity.
If $J$ and $q$ are fixed, then~\eqref{eq:bad-windows} gives
\[
 \frac{h}{a}\longrightarrow\frac1{Jq+3},\qquad
 \frac{D_{\mathrm{bad}}(\mathcal U(a,b),p,q)+1}{a}\longrightarrow0
 \quad\left(a\to\infty,\ \frac ba\to\alpha,\ (a,b)\in\cA_m\right).
\]
We then let $J\to\infty$ with $q$ fixed. By Claim~\ref{c:periodic-limit} of Proposition~\ref{prop:periodic}, the lower bound becomes
\[
 \liminf_{\substack{a\to\infty,\ (a,b)\in\cA_m\\\frac{b}{a}\to\alpha}}
 \frac{\log E_m(a,b)}a
 \ge\lim_{J\to\infty}\frac{\log A_{J-3}(p,q)}{Jq+3}
 =\frac{\ell(W_{p,q})}q.
\]
\end{proof}

%% file: sections/irrational-slopes-for-m-ge-4.tex
\section{Irrational slopes for \texorpdfstring{$m\ge4$}{m >= 4}}\label{sec:irrational}
\subsection{Existence of the limit}

\begin{lemma}[Vertical length comparison]\label{lem:vertical-length-comparison}
Let $M\in\Z_{>0}$ and $L\in\Z_{>0}$. Let $X_1,\ldots,X_L$ and $Y_1,\ldots,Y_L$ be simple runs.
Put $X=X_1\ccn \cdots\ccn X_L$, $Y=Y_1\ccn \cdots\ccn Y_L$, $\alpha_i=|X_i|_N$, and $\beta_i=|Y_i|_N$.
Assume $\alpha_i\in[M]_{\Z}$ and $\beta_i\in[M]_{\Z}$ for $i\in[L]_{\Z}$.
Then, there are $n\in\Z_{\ge0}$ and a run-edit list $\Omega=(\omega_1,\ldots,\omega_n)$ satisfying the following.
\begin{enumerate}
\item\label{c:vertical-comparison-target} The list $\Omega$ starts at $X$ and satisfies
\[
 \operatorname{RunEResult}(\Omega,X)
 =\operatorname{Trans}\bigl(Y,\operatorname{Start}_{\mathrm{pt}}(X)-\operatorname{Start}_{\mathrm{pt}}(Y)\bigr).
\]
\item\label{c:vertical-comparison-count} The list $\Omega$ consists of exactly
\[
 n=\sum_{i=1}^L|\alpha_i-\beta_i|
\]
vertical run-edits.
\item\label{c:vertical-comparison-lengths} After each run-edit of $\Omega$, every simple-run vertical length belongs to $[M]_{\Z}$.
\end{enumerate}
\end{lemma}
\begin{proof}
\prf{Claim~\ref{c:vertical-comparison-target}} Let us construct $\Omega$.
Process the simple runs in the order $i=1,\ldots,L$.
At the $i$th simple run, let $W$ be the current path and write this simple run as $W[s,s+t+1]$.
Here $s\in[0,|W|-1]_{\Z}$ and $t\in\Z_{>0}$.

First, let us increase the vertical length when $\alpha_i<\beta_i$.
For each $t\in[\alpha_i,\beta_i-1]_{\Z}$ in increasing order, replace $W$ by
\[
 \operatorname{Ins}(W,s+t+1,N).
\]
This gives $\beta_i-\alpha_i$ insertions and final vertical length $\beta_i$.

Second, let us decrease the vertical length when $\beta_i<\alpha_i$.
For each $t=\alpha_i,\alpha_i-1,\ldots,\beta_i+1$, replace $W$ by
\[
 \operatorname{Del}(W,W_{s+t+1}).
\]
Since $t>\beta_i\ge1$, every deletion has $t\ge2$.
This gives $\alpha_i-\beta_i$ deletions and final vertical length $\beta_i$.

Let us identify the resulting path.
When $\alpha_i=\beta_i$, make no change at the $i$th simple run.
The resulting vertical lengths are $\beta_1,\ldots,\beta_L$.
Thus, the resulting path is $\operatorname{Trans}(Y,\operatorname{Start}_{\mathrm{pt}}(X)-\operatorname{Start}_{\mathrm{pt}}(Y))$.

\prf{Claim~\ref{c:vertical-comparison-count}} Let us count the vertical run-edits.
For each $i\in[L]_{\Z}$, the construction in Claim~\ref{c:vertical-comparison-target} uses exactly $|\alpha_i-\beta_i|$ vertical run-edits. Summing over $i\in[L]_{\Z}$ gives $n=\sum_{i=1}^L|\alpha_i-\beta_i|$.

\prf{Claim~\ref{c:vertical-comparison-lengths}} Let us bound the intermediate vertical lengths.
During the run-edits of $\Omega$ on the $i$th simple run, its vertical length $t$ satisfies
\[
 1\le\min(\alpha_i,\beta_i)\le t\le\max(\alpha_i,\beta_i)\le M.
\]
For $j\in[L]_{\Z}\setminus\{i\}$, these run-edits do not change the vertical length of the $j$th simple run.
\end{proof}

\Needspace{15\baselineskip}
\begin{definition}\label{def:admissible-mechanical}
Let $M\in\Z_{>0}$, $r\in\R$, $\eta>0$, $\delta>0$, and $B>0$.
Let $l\in\Z_{>0}$, $T\in\Z_{>0}$, $\sigma_1\in[1,\infty)$, $\sigma_2\in[1,\infty)$,
$u_1\in[0,1)$, and $u_2\in[0,1)$. Let $A_1,\ldots,A_l$ be simple runs.
We call $(l,T,\sigma_1,\sigma_2,u_1,u_2,(A_1,\ldots,A_l))$
\textit{$(M,r,\eta,\delta,B)$-admissible} when the following conditions hold.
\begin{hypotheses}
\item\label{hyp:mechanical-source} We have
\[
 |A_i|_N=\floor{u_1+i\sigma_1}-\floor{u_1+(i-1)\sigma_1}
 \quad(i\in[l-1]_{\Z}),\qquad |A_l|_N\in[M]_{\Z}.
\]
\item\label{hyp:mechanical-proximity} We have
\[
 |\sigma_i-r|<\delta\quad(i\in[2]_{\Z}),\qquad T\ge(1+\eta)l+B.
\]
\end{hypotheses}
\end{definition}

\begin{lemma}[Fractional parts near a rational slope]\label{lem:fractional-part-estimates}
Let $p\in\Z_{>0}$ and $q\in\Z_{>0}$ satisfy $\gcd(p,q)=1$.
Let $\sigma\in\R$, $w\in\R$, and $0<\eta<\frac14$. Assume
\[
 |\sigma-\frac pq|<\frac2{q^2}.
\]
For $x\in\R$, let $\{x\}=x-\floor x$.
Let $0\le x_1\le\cdots\le x_q<1$ be the nondecreasing rearrangement of
$(\{w+j\sigma\})_{j\in[0,q-1]_{\Z}}$. Let $x_{q+1}=x_1+1$.
Then, the following hold.
\begin{enumerate}
\item\label{c:fractional-part-gaps} We have
\[
 x_{i+1}-x_i\le\frac5q\qquad(i\in[q]_{\Z}).
\]
\item\label{c:fractional-part-count} We have
\[
 \#\{j\in[0,q-1]_{\Z}:\{w+j\sigma\}\in[0,2\eta]\cup[1-2\eta,1)\}
 \le4\eta q+10.
\]
\end{enumerate}
\end{lemma}
\begin{proof}
\prf{Claim~\ref{c:fractional-part-gaps}} Let us compare the fractional parts with a rational grid.
For $j\in[0,q-1]_{\Z}$, put $t_j=jp-q\floor{w+j\sigma}$.
Then,
\[
 \left|\{w+j\sigma\}-\left(w+\frac{t_j}{q}\right)\right|
 =j\left|\sigma-\frac pq\right|<\frac2q.
\]
The integers $t_j$ represent distinct residue classes modulo $q$ because
$t_j\equiv jp\pmod q$ and $\gcd(p,q)=1$.
Thus, each $t\in\Z$ has the form $t=t_j+nq$ for some
$j\in[0,q-1]_{\Z}$ and $n\in\Z$.
Suppose $x_{i+1}-x_i>\frac5q$.
The interval $(x_i+\frac2q,x_{i+1}-\frac2q)$ has length greater than $\frac1q$.
The interval contains $w+\frac tq$ for some $t\in\Z$.
Let $t=t_j+nq$.
The number $\{w+j\sigma\}+n$ belongs to $(x_i,x_{i+1})$,
contradicting the choice of $x_i$ and $x_{i+1}$.

\prf{Claim~\ref{c:fractional-part-count}} Let us count the grid points near the two endpoints.
If
$\{w+j\sigma\}\in[0,2\eta]\cup[1-2\eta,1)$, then
\[
 w+\frac{t_j}{q}\in
 \left[-\frac2q,2\eta+\frac2q\right]
 \cup\left[1-2\eta-\frac2q,1+\frac2q\right].
\]
Each interval has length $2\eta+\frac4q$ and contains at most $2\eta q+5$
numbers of the form $w+\frac tq$, $t\in\Z$.
Hence, at most $4\eta q+10$ indices satisfy the condition in Claim~\ref{c:fractional-part-count}.
\end{proof}

\Needspace{20\baselineskip}
\begin{lemma}[Phase alignment and prefix discrepancies]\label{lem:phase-prefix-discrepancy}
Let $p_0\in\Z_{>0}$, $q_0\in\Z_{>0}$, $L\in\Z_{>0}$, $0<\eta<\frac14$, $\sigma_1\in[1,\infty)$, and $\sigma_2\in[1,\infty)$. Assume
\[
 \gcd(p_0,q_0)=1,\qquad
 q_0>\frac{20}{\eta},\qquad q_0\le\eta L,
\]
\[
 \left|\sigma_a-\frac{p_0}{q_0}\right|<\frac{2}{q_0^2}
 \quad(a\in[2]_{\Z}),\qquad
 |\sigma_1-\sigma_2|<\frac{\eta}{2L}.
\]
Let $u_1\in[0,1)$, $u_2\in[0,1)$, and $b\in\Z_{\ge0}$. For $j\in[0,q_0-1]_{\Z}$, define
\[
 v_j=\{u_2+(b+j)\sigma_2\},\qquad
 J=\{j\in[0,q_0-1]_{\Z}:|u_1-v_j|<\eta\}.
\]
The following statements hold.

\begin{enumerate}
\item\label{c:phase-selection} The set $J$ is nonempty. Define the uniquely specified index and phase
\[
j_* =\min J,\qquad v=v_{j_*}.
\]
Then, $0\le j_*<q_0$ and $|u_1-v|<\eta$.

\item\label{c:phase-prefix-mass} Use $v$ from Claim~\ref{c:phase-selection}. Define
\[
d_i=\lfloor u_1+i\sigma_1\rfloor-\lfloor v+i\sigma_2\rfloor
\quad(i\in[0,L]_{\Z}).
\]
We have $d_0=0$. For $i\in[L]_{\Z}$,
\[
d_i\in\{-1,0,1\},\qquad
d_i\ne0\Longrightarrow
\{u_1+i\sigma_1\}\in[0,2\eta]\cup[1-2\eta,1).
\]
Moreover,
\[
\sum_{i=1}^L|d_i|
\le\left(4\eta+\frac{10}{q_0}\right)L+q_0
\le\frac{11}{2}\eta L.
\]

\item\label{c:phase-run-variation} Define, for $i\in[L]_{\Z}$,
\[
\alpha_i=\lfloor u_1+i\sigma_1\rfloor-\lfloor u_1+(i-1)\sigma_1\rfloor,
\qquad
\beta_i=\lfloor v+i\sigma_2\rfloor-\lfloor v+(i-1)\sigma_2\rfloor.
\]
Then,
\[
\alpha_i-\beta_i=d_i-d_{i-1},\qquad
\sum_{i=1}^L|\alpha_i-\beta_i|
\le2\sum_{i=1}^L|d_i|-|d_L|
\le11\eta L.
\]
\end{enumerate}
\end{lemma}
\begin{proof}
\prf{Claim~\ref{c:phase-selection}} Let us prove that $J$ is nonempty. Let $x_1\le\cdots\le x_{q_0}$ be the ordered values of $v_0,\ldots,v_{q_0-1}$. Let $x_{q_0+1}=x_1+1$.
Claim~\ref{c:fractional-part-gaps} of Lemma~\ref{lem:fractional-part-estimates} with phase $u_2+b\sigma_2$ gives
\[
 x_{a+1}-x_a\le\frac5{q_0}<\eta
 \quad(a\in[q_0]_{\Z}).
\]
If $x_a\le u_1\le x_{a+1}$ for some $a<q_0$, then $|u_1-x_a|\le\frac5{q_0}<\eta$. If $u_1\le x_1$, use
\[
 x_1-u_1\le x_1\le x_1+1-x_{q_0}\le\frac5{q_0}.
\]
If $u_1\ge x_{q_0}$, use
\[
 u_1-x_{q_0}<1-x_{q_0}\le x_1+1-x_{q_0}\le\frac5{q_0}.
\]
Every $x_i$ equals $v_j$ for some $j\in[0,q_0-1]_{\Z}$. Hence, $J\ne\varnothing$. Since $j_*=\min J\in J$, we have $0\le j_*<q_0$ and $|u_1-v_{j_*}|<\eta$.

\prf{Claim~\ref{c:phase-prefix-mass}} Let us locate every nonzero discrepancy. Since $u_1\in[0,1)$ and $v\in[0,1)$, we have $d_0=0$. For $i\in[L]_{\Z}$, put
\[
 a_i=u_1+i\sigma_1,\qquad b_i=v+i\sigma_2.
\]
Claim~\ref{c:phase-selection} gives $|u_1-v|<\eta$.
Together with $|\sigma_1-\sigma_2|<\frac{\eta}{2L}$ and $i\le L$, this gives
\[
 |a_i-b_i|
 \le|u_1-v|+i|\sigma_1-\sigma_2|
 <\eta+\frac{\eta}{2}=\frac{3\eta}{2}<1.
\]
Hence, the integer $d_i=\lfloor a_i\rfloor-\lfloor b_i\rfloor$ belongs to $\{-1,0,1\}$. If $d_i=1$, then $b_i<\lfloor a_i\rfloor$. Thus,
\[
 \{a_i\}=a_i-\lfloor a_i\rfloor<a_i-b_i<\frac{3\eta}{2}<2\eta.
\]
If $d_i=-1$, then $b_i\ge\lfloor a_i\rfloor+1$. Thus,
\[
 \{a_i\}\ge1-(b_i-a_i)>1-\frac{3\eta}{2}>1-2\eta.
\]
Thus, $d_i\ne0$ implies $\{u_1+i\sigma_1\}\in[0,2\eta]\cup[1-2\eta,1)$.

Let us count these indices. Write
\[
 H=\left\lfloor\frac{L}{q_0}\right\rfloor,\qquad
 s=L-Hq_0\in[0,q_0-1]_{\Z}.
\]
Partition $[L]_{\Z}$ into
\[
 B_z=[zq_0+1,(z+1)q_0]_{\Z}
 \quad(z\in[0,H-1]_{\Z})
\]
and the final interval $[Hq_0+1,L]_{\Z}$, which has $s$ elements. For $i=zq_0+1+j$, $j\in[0,q_0-1]_{\Z}$, let
\[
 w_z=u_1+(zq_0+1)\sigma_1.
\]
Then, $\{u_1+i\sigma_1\}=\{w_z+j\sigma_1\}$. Claim~\ref{c:fractional-part-count} of Lemma~\ref{lem:fractional-part-estimates} bounds the number of nonzero $d_i$ in $B_z$ by $4\eta q_0+10$. Since $|d_i|$ equals $1$ precisely at the nonzero discrepancies,
\[
 \sum_{i=1}^L|d_i|
 \le H(4\eta q_0+10)+s
 \le\left(4\eta+\frac{10}{q_0}\right)L+q_0.
\]
The inequalities $\frac{10}{q_0}<\frac\eta2$ and $q_0\le\eta L$ give the bound $\frac{11}{2}\eta L$. The estimate $\sum_{i=1}^L|d_i|\le H(4\eta q_0+10)+s$ also holds when $s=0$ and $[Hq_0+1,L]_{\Z}=\varnothing$.

\prf{Claim~\ref{c:phase-run-variation}} Let us express $\alpha_i-\beta_i$ as $d_i-d_{i-1}$ for $i\in[L]_{\Z}$. Expanding the definitions gives $\alpha_i-\beta_i=d_i-d_{i-1}$. Since $d_0=0$, the triangle inequality gives
\[
 \begin{aligned}
 \sum_{i=1}^L|\alpha_i-\beta_i|
 &\le\sum_{i=1}^L|d_i|+\sum_{i=1}^L|d_{i-1}|\\
 &=2\sum_{i=1}^L|d_i|-|d_L|\\
 &\le11\eta L.
 \end{aligned}
\]
The inequality $2\sum_{i=1}^L|d_i|-|d_L|\le11\eta L$ follows from Claim~\ref{c:phase-prefix-mass}.
\end{proof}

\begin{lemma}[Mechanical run-edits]\label{lem:mechanical-edits}
Let $m\in\Z_{\ge4}$ and $M=m-2$.
Let $r\in(1,M)$ be irrational and let $0<\eta<\frac14$.
Then, there are constants $\delta_m(r,\eta)>0$ and $B_m(r,\eta)>0$ depending only on $m$, $r$, and $\eta$ with the following property.
For each $(M,r,\eta,\delta_m(r,\eta),B_m(r,\eta))$-admissible tuple
$(l,T,\sigma_1,\sigma_2,u_1,u_2,(A_1,\ldots,A_l))$ and $U_0=A_1\ccn\cdots\ccn A_l$,
there are $n\in\Z_{\ge0}$, $\zeta\in\Z^2$, and a run-edit list $\Omega=(\omega_1,\ldots,\omega_n)$ with the following properties.
\begin{enumerate}
\item\label{c:mechanical-target} The list $\Omega$ starts at $U_0$ and satisfies
\[
 \operatorname{RunEResult}(\Omega,U_0)
 =\operatorname{Trans}\bigl(\operatorname{RunFactor}_{1,T}(\operatorname{MechRun}(\sigma_2,u_2)),\zeta\bigr).
\]
\item\label{c:mechanical-counts} $\Omega$ consists of $K=T-l$ insertions of simple runs,
$I$ vertical-edge insertions, and $D$ vertical-edge deletions. Their numbers satisfy
\begin{equation}\label{eq:mechanical-bound}
 n=K+I+D,\qquad I,D\in\Z_{\ge0},\qquad I+D\le12\eta l+B_m(r,\eta).
\end{equation}
\item\label{c:mechanical-lengths} For $U_j=\operatorname{RunEResult}((\omega_1,\ldots,\omega_j),U_0)$ $(j\in[0,n]_{\Z})$,
every simple run $R\subseteq U_j$ satisfies $|R|_N\in[M]_{\Z}$.
\end{enumerate}
\end{lemma}
\begin{proof}
\prf{Claim~\ref{c:mechanical-target}} Let us choose $p_0$, $q_0$, $L$, $\delta_m(r,\eta)$, and $B_m(r,\eta)$.
Put $k=\floor r$.
Choose relatively prime $p_0\in\Z_{>0}$ and $q_0\in\Z_{>0}$ such that
\begin{equation}\label{eq:mechanical-approximation}
 q_0>\frac{20}{\eta},\qquad |r-\frac{p_0}{q_0}|<\frac1{q_0^2}.
\end{equation}
To obtain such a fraction, partition $[0,1)$ into $Q$ intervals of length $\frac1Q$.
Two of the $Q+1$ numbers $\{jr\}$, $j\in[0,Q]_{\Z}$, lie in the same interval.
Subtracting their indices gives integers $p$ and $q$ such that
\[
 1\le q\le Q,\qquad |qr-p|<\frac1Q,\qquad
 \left|r-\frac pq\right|<\frac1{qQ}\le\frac1{q^2}.
\]
Let $\frac{p'}{q'}$ be $\frac pq$ in lowest terms. Then, $1\le q'\le q$ and $\left|r-\frac{p'}{q'}\right|<\frac1{q^2}\le\frac1{q'^2}$.
These reduced denominators are unbounded as $Q\to\infty$: otherwise, only finitely many
reduced fractions in $[r-1,r+1]$ could occur, and their positive distances from the irrational number $r$
could not tend to zero.
Thus, choose a reduced fraction with $q_0>\frac{20}{\eta}$ satisfying~\eqref{eq:mechanical-approximation}.
Choose $L\in\Z_{>0}$ such that $q_0\le\eta L$ and let
\begin{equation}\label{eq:mechanical-constants}
 \delta_m(r,\eta)=\frac12\min\left\{\frac1{q_0^2},\frac{\eta}{4L},r-k,k+1-r\right\},
 \qquad B_m(r,\eta)=(M-1)L+1.
\end{equation}
These constants depend only on $r$, $\eta$, and $M$.
Condition~\ref{hyp:mechanical-proximity} of Definition~\ref{def:admissible-mechanical} gives $k<\sigma_i<k+1$ for $i\in[2]_{\Z}$.
Thus, every simple run of either mechanical sequence has vertical length $k$ or $k+1$. Both lengths belong to $[M]_{\Z}$.

Let us compare the two slopes to the same rational approximation.
For $\sigma\in\{\sigma_1,\sigma_2\}$, equations~\eqref{eq:mechanical-approximation}--\eqref{eq:mechanical-constants} give
\[
 |\sigma-\frac{p_0}{q_0}|<\frac2{q_0^2}.
\]
Thus, Lemma~\ref{lem:fractional-part-estimates} applies with $p=p_0$ and $q=q_0$.

Let us partition the run into factors of length $L$ and one final factor.
Let $H=\floor{\frac{l-1}{L}}$.
Let us construct the part of $\Omega$ acting on the first $HL$ simple runs of $U_0$.
Put $c_0=0$.
For $z\in[0,H-1]_{\Z}$, suppose the run-edits on the first $zL$ simple runs have
produced the first $c_z$ simple runs of $\operatorname{MechRun}(\sigma_2,u_2)$.

Let us choose the target intercept for the factor indexed by $z\in[0,H-1]_{\Z}$.
Let $u=\{u_1+zL\sigma_1\}$.
Apply Claim~\ref{c:phase-selection} of Lemma~\ref{lem:phase-prefix-discrepancy}
with source phase $u$ and offset $b=c_z$.
Let $j_z=\min J$ and $v=\{u_2+(c_z+j_z)\sigma_2\}$ as in Claim~\ref{c:phase-selection} of Lemma~\ref{lem:phase-prefix-discrepancy}.
Then, $0\le j_z<q_0$ and $|u-v|<\eta$.
If $j_z>0$, insert
$\operatorname{RunFactor}_{c_z+1,j_z}(\operatorname{MechRun}(\sigma_2,u_2))$
immediately before the copy of $\operatorname{RunFactor}_{zL+1,L}(U_0)$.

Let us bound the number of vertical run-edits on the copy of $\operatorname{RunFactor}_{zL+1,L}(U_0)$.
Define
\[
 X=\operatorname{RunFactor}_{1,L}(\operatorname{MechRun}(\sigma_1,u)),\qquad
 Y=\operatorname{RunFactor}_{1,L}(\operatorname{MechRun}(\sigma_2,v)).
\]
Claim~\ref{c:phase-run-variation} of Lemma~\ref{lem:phase-prefix-discrepancy} gives
\[
 \sum_{i=1}^L
 \left|\,|\operatorname{RunFactor}_{i,1}(X)|_N
       -|\operatorname{RunFactor}_{i,1}(Y)|_N\,\right|\le11\eta L.
\]
By Lemma~\ref{lem:vertical-length-comparison}, there is a run-edit list $\Omega_z=(\omega_{z,1},\ldots,\omega_{z,n_z})$ satisfying $\operatorname{RunEResult}(\Omega_z,X)=Y$ and consisting of at most $11\eta L$
vertical run-edits.
Every simple-run vertical length produced by $\Omega_z$ belongs to $[M]_{\Z}$.
Let $c_{z+1}=c_z+j_z+L$.
For some $\zeta_z\in\Z^2$, the resulting first $c_{z+1}$ simple runs form
\[
 \operatorname{Trans}\bigl(\operatorname{RunFactor}_{1,c_{z+1}}(\operatorname{MechRun}(\sigma_2,u_2)),\zeta_z\bigr).
\]

Let us construct the run-edits on $A_{HL+1},\ldots,A_l$.
We have $l-HL\in[L]_{\Z}$.
For $i\in[l-HL]_{\Z}$, change the length of the copy of $A_{HL+i}$ to
$|\operatorname{RunFactor}_{c_H+i,1}(\operatorname{MechRun}(\sigma_2,u_2))|_N$.
Both lengths belong to $[M]_{\Z}$, including the length of $A_l$ by Condition~\ref{hyp:mechanical-source} of Definition~\ref{def:admissible-mechanical}.
Lemma~\ref{lem:vertical-length-comparison} gives at most $(M-1)L$ vertical run-edits
and keeps every length in $[M]_{\Z}$.
Now,
\[
 c_H+l-HL=l+\sum_{z=0}^{H-1}j_z\le l+Hq_0\le(1+\eta)l<T.
\]
Append
\[
 \operatorname{RunFactor}_{c_H+l-HL+1,\,T-c_H-l+HL}
       (\operatorname{MechRun}(\sigma_2,u_2)).
\]
Let $Q=\operatorname{RunEResult}(\Omega,U_0)$ and $\zeta=\operatorname{Start}_{\mathrm{pt}}(Q)$.
Then,
\[
 Q=\operatorname{Trans}\bigl(\operatorname{RunFactor}_{1,T}(\operatorname{MechRun}(\sigma_2,u_2)),\zeta\bigr).
\]

\prf{Claim~\ref{c:mechanical-counts}} Let us count the simple-run insertions and the vertical run-edits.
The construction in Claim~\ref{c:mechanical-target} uses exactly $K=T-l$ simple-run insertions. If $I$ and $D$ count its vertical-edge insertions and deletions, respectively, then $n=K+I+D$. The bound in~\eqref{eq:mechanical-bound} follows from
\[
 I+D\le11\eta HL+(M-1)L\le12\eta l+B_m(r,\eta).
\]
\prf{Claim~\ref{c:mechanical-lengths}} Every inserted simple run has length $k$ or $k+1$. Both belong to $[M]_{\Z}$.
Lemma~\ref{lem:vertical-length-comparison} keeps every intermediate vertical length in $[M]_{\Z}$ at each application in the construction of $\Omega$.
\end{proof}

\begin{lemma}\label{lem:approaching-sequence}
Let $m\in\Z_{\ge3}$.
For every $r\in I_m$, there is a sequence $((a_\nu,b_\nu))_{\nu\ge1}$ in $\cA_m$ satisfying
\begin{equation}\label{eq:approaching-sequence}
 \lim_{\nu\to\infty}a_\nu=\infty,\qquad
 \lim_{\nu\to\infty}\frac{b_\nu}{a_\nu}=r.
\end{equation}
\end{lemma}
\begin{proof}
Let us choose pairs in $\cA_m$ whose slopes tend to $r$.
Choose distinct reduced fractions $\frac{b_\nu}{a_\nu}\in I_m$ tending to $r$.
For each $A>0$, the bounded interval $I_m$ contains only finitely many reduced
fractions with denominator at most $A$. Thus, $a_\nu\to\infty$.
Since $\frac{b_\nu}{a_\nu}\in I_m$, we have
$a_\nu^2+b_\nu^2-ma_\nu b_\nu\le0$.
Equation~\eqref{eq:nonempty-domain} gives $(a_\nu,b_\nu)\in\cA_m$.
\end{proof}

\begin{theorem}\label{thm:irrational-large-m}
Let $m\in\Z_{\ge4}$ and $M=m-2$.
For every irrational $r\in(1,M)$, there is $\Upsilon_m(r)\in\R$ such that
\[
 \lim_{\substack{a\to\infty,\ (a,b)\in\cA_m\\\frac{b}{a}\to r}}
 \frac{\log E_m(a,b)}a=\Upsilon_m(r).
\]
\end{theorem}
\begin{proof}
Let us choose the comparison parameters for Lemma~\ref{lem:mechanical-edits} and Claim~\ref{c:composition-bound} of Proposition~\ref{prop:composition}.
Let
\[
 C=32(m+1),\qquad H=2C+16,
 \qquad 0<\eta<\frac1{2H},\qquad d=1-H\eta.
\]
Use $\delta_m(r,\eta)$ and $B_m(r,\eta)$ from Lemma~\ref{lem:mechanical-edits}.
Choose relatively prime $p\in\Z_{>0}$ and $q\in\Z_{>0}$ such that $q\ge2$,
$\alpha=\frac{p}{q}\in(1,M)$, and $|\alpha-r|<\delta_m(r,\eta)$.
We fix $\alpha$, $p$, and $q$ while $a\to\infty$, $(a,b)\in\cA_m$, and $\frac ba\to r$.

Let us verify~\eqref{eq:layer-budget} for Claim~\ref{c:composition-bound} of Proposition~\ref{prop:composition}.
For a run $Z$ of $\nu$ simple runs and sufficiently large $\nu$, let
\[
 P=\floor{C\eta\nu},\qquad T=\nu-2P-2,\qquad
 V=\operatorname{RunFactor}_{P+2,T}(Z).
\]
We have $P>0$ and $T>0$. Define
\[
\begin{aligned}
 W_{\mathrm L}&=\operatorname{RunFactor}_{1,P}(Z),&
 W_{\mathrm R}&=\operatorname{RunFactor}_{\nu-P+1,P}(Z),\\
 u&=|\operatorname{RunFactor}_{P+1,1}(Z)|_N,&
 v&=|\operatorname{RunFactor}_{\nu-P,1}(Z)|_N.
\end{aligned}
\]
Equation~\eqref{eq:target-shape} gives $Z=W_{\mathrm L}\ccn R_u\ccn V\ccn R_v\ccn W_{\mathrm R}$.
Let $l=d\nu+O(1)$. Then,
\[
 K=T-l,\qquad \lim_{\nu\to\infty}\frac K\nu=16\eta,\qquad
 T\ge(1+\eta)l+B_m(r,\eta)
\]
when $\nu$ is sufficiently large. If $I\in\Z_{\ge0}$ and $D\in\Z_{\ge0}$ satisfy $I+D\le13\eta\nu$, then eventually
\[
 K+I+(m+1)D\le(17+13(m+1))\eta \nu\le P.
\]
The inequality $(17+13(m+1))\eta\nu\le P$ uses $P=\floor{C\eta\nu}$ and
$(C-17-13(m+1))\eta\nu\ge1$ for sufficiently large $\nu$, since
$C=32(m+1)>17+13(m+1)$. Equation~\eqref{eq:entropy-loss} gives
\[
 \log \operatorname{Fib}_m(l,D)\le\nu\psi_m(13\eta).
\]

Let us verify admissibility in the two applications below.
Equation~\eqref{eq:ambient} gives $\Gamma(a,b)\to r$.
Thus, $1<\Gamma(a,b)<M$ and $|\Gamma(a,b)-r|<\delta_m(r,\eta)$ for all sufficiently large $a$.
The choice of $\alpha$ and the bound $T\ge(1+\eta)l+B_m(r,\eta)$ verify condition~\ref{hyp:mechanical-proximity} of Definition~\ref{def:admissible-mechanical} in both applications.

First, let us obtain a lower bound for $\frac1a\log E_m(a,b)$.
Let $(a,b)\in\cA_m$, $\nu=a+1$, and $j=\floor{\frac{d\nu}{q}}$.
Apply Lemma~\ref{lem:mechanical-edits} with
\[
\begin{gathered}
 U=W_{p,q}^j,\quad l=jq,\quad Z=\mathcal U(a,b),\quad
 V=\operatorname{RunFactor}_{P+2,T}(\mathcal U(a,b)),\\
 (\sigma_1,\sigma_2,u_1,u_2)
 =\bigl(\alpha,\Gamma(a,b),0,\{(P+1)\Gamma(a,b)\}\bigr),\\
 A_i=R_{w_i(p,q)}\qquad(i\in[l]_{\Z}).
\end{gathered}
\]
Definition~\ref{def:mechanical-paths} verifies condition~\ref{hyp:mechanical-source} of Definition~\ref{def:admissible-mechanical} for these source runs.
Lemma~\ref{lem:mechanical-edits} gives $I+D\le12\eta l+B_m(r,\eta)\le13\eta\nu$
for all sufficiently large $a$.
The last run $\operatorname{RunFactor}_{a+1,1}(\mathcal U(a,b))$
belongs to $W_{\mathrm R}$ and does not occur in $V$.
Claim~\ref{c:composition-bound} of Proposition~\ref{prop:composition}
and Lemma~\ref{lem:word-reindexing} apply with $V_*=\mathcal U(a,b)$.
Equations~\eqref{eq:dictionary}, \eqref{eq:composition}, and~\eqref{eq:entropy-loss} give
\[
 \frac{\log E_m(a,b)}a
 \ge\frac{j}{a}\frac{\log C_0(W_{p,q}^j)}j
 -\frac{\nu}{a}\psi_m(13\eta)-\frac{\log m}{a}.
\]
Since $\frac ja\to\frac dq$ and $\frac\nu a\to1$, Claim~\ref{c:periodic-limit} of Proposition~\ref{prop:periodic} and~\eqref{eq:rational-period} give
\begin{equation}\label{eq:irrational-lower}
 \liminf_{\substack{a\to\infty,\ (a,b)\in\cA_m\\\frac{b}{a}\to r}}
 \frac{\log E_m(a,b)}a
 \ge d\Upsilon_m(\alpha)-\psi_m(13\eta).
\end{equation}

Second, let us obtain an upper bound for $\frac1a\log E_m(a,b)$.
Let $(a,b)\in\cA_m$, $l=a+1$, $N=\ceil{\frac{l}{dq}}$, and $\nu=Nq$.
Apply Lemma~\ref{lem:mechanical-edits} with
\[
\begin{gathered}
 U=\mathcal U(a,b),\quad Z=W_{p,q}^N,\quad
 V=\operatorname{RunFactor}_{P+2,T}(W_{p,q}^N),\\
 (\sigma_1,\sigma_2,u_1,u_2)
 =\bigl(\Gamma(a,b),\alpha,0,\{(P+1)\alpha\}\bigr),\\
 A_i=R_{\Lambda_i(a,b)}\quad(i\in[a]_{\Z}),\qquad A_l=R_1.
\end{gathered}
\]
Definition~\ref{def:frame} verifies condition~\ref{hyp:mechanical-source} of Definition~\ref{def:admissible-mechanical}, including $|A_l|_N=1\in[M]_{\Z}$.
Lemma~\ref{lem:mechanical-edits} gives $I+D\le12\eta l+B_m(r,\eta)\le13\eta\nu$
for all sufficiently large $a$.
Claim~\ref{c:composition-bound} of Proposition~\ref{prop:composition}
and Lemma~\ref{lem:word-reindexing} apply with $V_*=W_{p,q}^N$.
Equations~\eqref{eq:composition}, \eqref{eq:dictionary}, and~\eqref{eq:entropy-loss} give
\[
 \frac{\log E_m(a,b)}a
 \le\frac{N}{a}\frac{\log C_0(W_{p,q}^N)}N
 +\frac{\nu}{a}\psi_m(13\eta)-\frac{\log m}{a}.
\]
Here, $\frac Na\to\frac1{dq}$ and $\frac\nu a\to\frac1d$. By Claim~\ref{c:periodic-limit} of Proposition~\ref{prop:periodic} and~\eqref{eq:rational-period},
\begin{equation}\label{eq:irrational-upper}
 \limsup_{\substack{a\to\infty,\ (a,b)\in\cA_m\\\frac{b}{a}\to r}}
 \frac{\log E_m(a,b)}a
 \le\frac{\Upsilon_m(\alpha)+\psi_m(13\eta)}d.
\end{equation}

Third, let us compare~\eqref{eq:irrational-lower} and~\eqref{eq:irrational-upper}. Lemma~\ref{lem:word-bounds} gives
$\Upsilon_m(\alpha)\le m(M+1)\log2$.
The interval in~\eqref{eq:irrational-lower}--\eqref{eq:irrational-upper} has width at most
\[
 m(M+1)\log2\left(\frac1d-d\right)
      +\psi_m(13\eta)\left(1+\frac1d\right).
\]
As $\eta\downarrow0$, we have $d=1-H\eta\to1$ and
\[
 \frac1d-d\longrightarrow0,\qquad
 \psi_m(13\eta)\longrightarrow0.
\]
The limit $\psi_m(13\eta)\to0$ as $\eta\downarrow0$ follows from $x\log x\to0$ as $x\downarrow0$ and Definition~\ref{def:entropy-loss-function}.

Let us prove that every sequence has the same limit.
For each $i\in[2]_{\Z}$, let $((a_{i,n},b_{i,n}))_{n\ge1}$ be a sequence in $\cA_m$ satisfying
\[
 \lim_{n\to\infty}a_{i,n}=\infty,\qquad
 \lim_{n\to\infty}\frac{b_{i,n}}{a_{i,n}}=r.
\]
Define
\[
 L_i=\liminf_{n\to\infty}\frac{\log E_m(a_{i,n},b_{i,n})}{a_{i,n}},\qquad
 U_i=\limsup_{n\to\infty}\frac{\log E_m(a_{i,n},b_{i,n})}{a_{i,n}}.
\]
Lemma~\ref{lem:word-bounds} and~\eqref{eq:dictionary} make these four numbers finite.

Let us prove $U_i\le L_j$ for $i,j\in[2]_{\Z}$.
For each fixed $\eta$, choose $\alpha$, $\delta_m(r,\eta)$, and $B_m(r,\eta)$ before taking either sequence limit.
Equations~\eqref{eq:irrational-lower} and~\eqref{eq:irrational-upper} give, for $i,j\in[2]_{\Z}$,
\[
 U_i-L_j\le m(M+1)\log2\left(\frac1d-d\right)
      +\psi_m(13\eta)\left(1+\frac1d\right).
\]
Letting $\eta\downarrow0$ gives $U_i\le L_j$.
Taking $i=j$ and using $L_i\le U_i$ gives $L_i=U_i$.
Taking $(i,j)=(1,2)$ and $(2,1)$ gives $L_1=U_1=L_2=U_2$.
Lemma~\ref{lem:approaching-sequence} supplies at least one such sequence.
Define $\Upsilon_m(r)$ to be its limit. The equality $L_1=U_1=L_2=U_2$ makes this definition independent of the chosen sequence.
\end{proof}

\subsection{Continuity of the limit}
\begin{lemma}\label{lem:continuity}
Let $m\in\Z_{\ge3}$.
Let $h_1\in\R$ and $h_2\in\R$ satisfy $h_1<h_2$. Let $J=[h_1,h_2]$.
Let $F:J\to\R$. Assume the following.
\begin{hypotheses}
\item\label{hyp:all-sequence-limit} For each $r\in J$ and every sequence $((a_n,b_n))_{n\ge1}$ in $\cA_m$ satisfying
\[
 \lim_{n\to\infty}a_n=\infty,\qquad
 \lim_{n\to\infty}\frac{b_n}{a_n}=r,
\]
we have
\[
 \lim_{n\to\infty}\frac{\log E_m(a_n,b_n)}{a_n}=F(r).
\]
\item\label{hyp:approaching-sequences} For each $r\in J$, there exists a sequence $((a_n,b_n))_{n\ge1}$ in $\cA_m$ such that
\[
 \lim_{n\to\infty}a_n=\infty,\qquad
 \lim_{n\to\infty}\frac{b_n}{a_n}=r.
\]
\end{hypotheses}
Then, we have the following.

\begin{enumerate}
\item\label{c:continuity-neighborhood} For every $r\in J$ and $\varepsilon>0$, there exist $A>0$ and $\delta>0$ such that
\begin{equation}\label{eq:continuity-neighborhood}
 \left|\frac{\log E_m(a,b)}a-F(r)\right|<\varepsilon
 \quad\left((a,b)\in\cA_m,\ a>A,\ \left|\frac ba-r\right|<\delta\right).
\end{equation}
\item\label{c:continuity-function} The function $F$ is continuous on $J$.
\item\label{c:continuity-endpoint} At the endpoints of $J$, we have the following.
\begin{enumerate}[label=(\alph*),ref=\theenumi(\alph*)]
\item\label{c:continuity-endpoint-neighborhood} For $r\in\{h_1,h_2\}$, the estimate in Claim~\ref{c:continuity-neighborhood} holds when additionally $\frac ba\in J$.
\item\label{c:continuity-endpoint-function} We have
\[
 \lim_{s\downarrow h_1}F(s)=F(h_1),\qquad
 \lim_{s\uparrow h_2}F(s)=F(h_2)\qquad(s\in J).
\]
\end{enumerate}
\end{enumerate}
\end{lemma}
\begin{proof}
\prf{Claim~\ref{c:continuity-neighborhood}} Let us prove~\eqref{eq:continuity-neighborhood} by contradiction.
If~\eqref{eq:continuity-neighborhood} fails at $r$, there is $\varepsilon_0>0$
such that, for every positive integer $\nu$, we find $(a_\nu,b_\nu)\in\cA_m$ with
\[
 a_\nu>\nu,\qquad
 \left|\frac{b_\nu}{a_\nu}-r\right|<\frac1\nu,\qquad
 \left|\frac{\log E_m(a_\nu,b_\nu)}{a_\nu}-F(r)\right|\ge\varepsilon_0.
\]
We have $a_\nu\to\infty$ and $\frac{b_\nu}{a_\nu}\to r$.
Assumption~\ref{hyp:all-sequence-limit} contradicts
$|\frac{\log E_m(a_\nu,b_\nu)}{a_\nu}-F(r)|\ge\varepsilon_0$.

\prf{Claim~\ref{c:continuity-function}} Let us obtain continuity from~\eqref{eq:continuity-neighborhood}.
Let $r\in J$ and $\varepsilon>0$.
Claim~\ref{c:continuity-neighborhood} gives $\delta>0$ and $A>0$ such that
\[
 \left|\frac{\log E_m(a,b)}a-F(r)\right|<\frac\varepsilon2
 \quad\left((a,b)\in\cA_m,\ a>A,\ \left|\frac ba-r\right|<\delta\right).
\]

Let us prove $|F(s)-F(r)|<\varepsilon$ when $s\in J$ and $|s-r|<\frac\delta2$.
Fix such an $s$.
By assumption~\ref{hyp:approaching-sequences}, choose a sequence $((a_\nu,b_\nu))_{\nu\ge1}$ in $\cA_m$ such that
\[
 a_\nu\longrightarrow\infty,\qquad \frac{b_\nu}{a_\nu}\longrightarrow s
 \qquad(\nu\longrightarrow\infty).
\]
Choose $\nu_0\in\Z_{>0}$ such that
\[
 a_\nu>A,\qquad
 \left|\frac{b_\nu}{a_\nu}-s\right|<\frac\delta2
 \qquad(\nu\ge\nu_0).
\]
For $\nu\ge\nu_0$, the triangle inequality gives
\[
 \left|\frac{b_\nu}{a_\nu}-r\right|
 \le\left|\frac{b_\nu}{a_\nu}-s\right|+|s-r|<\delta.
\]
Hence,~\eqref{eq:continuity-neighborhood} applies with $\varepsilon$ replaced by $\frac\varepsilon2$ for every $\nu\ge\nu_0$.

Let us take the limit as $\nu\to\infty$.
Assumption~\ref{hyp:all-sequence-limit} and continuity of the absolute value give
\[
 |F(s)-F(r)|
 =\lim_{\nu\to\infty}\left|\frac{\log E_m(a_\nu,b_\nu)}{a_\nu}-F(r)\right|
 \le\frac\varepsilon2<\varepsilon.
\]

\prf{Claim~\ref{c:continuity-endpoint-neighborhood}}
Claim~\ref{c:continuity-neighborhood} applies to $r\in\{h_1,h_2\}$ and pairs satisfying $\frac ba\in J$.

\prf{Claim~\ref{c:continuity-endpoint-function}} Let us identify the one-sided endpoint limits.
The hypothesis $h_1<h_2$ makes both one-sided domains nonempty.
By Claim~\ref{c:continuity-function}, for every $\varepsilon>0$ there are $\delta_1>0$ and $\delta_2>0$ such that
\[
\begin{aligned}
 |F(s)-F(h_1)|&<\varepsilon
 &&(s\in J,\ 0<s-h_1<\delta_1),\\
 |F(s)-F(h_2)|&<\varepsilon
 &&(s\in J,\ 0<h_2-s<\delta_2).
\end{aligned}
\]
These are the right limit at $h_1$ and the left limit at $h_2$, respectively.
\end{proof}

%% file: sections/uniform-lower-bound-near-slope-one-for-the-three-arrow-quiver.tex
\section{A uniform lower bound near slope one for the three-arrow quiver}\label{sec:central-lower}
\begin{definition}\label{def:central-euler-characteristic}
For positive coprime integers $a$ and $r$, let
\[
 E_r(a)=E_3(a,a+r).
\]
\end{definition}
From here to the end of the manuscript, let
\[
 \beta=\frac{256}{27}.
\]

\subsection{Graphs and quivers}\label{subsec:graphs-quivers}
\begin{definition}[Graphs]\label{def:graphs}
For a set $V$, define
\[
 \operatorname{Sym}^{2}(V)=(V\times V)/{\sim},\qquad
 (u,v)\sim(u',v')\ \Longleftrightarrow\
 (u',v')=(u,v)\ \text{or}\ (u',v')=(v,u).
\]
Let $\llbracket u,v\rrbracket$ denote the class of $(u,v)$.
For $U\subseteq V$, we identify $\operatorname{Sym}^{2}(U)$ with its image in $\operatorname{Sym}^{2}(V)$.
A \textit{graph} $G$ consists of a vertex set $\operatorname{Vert}(G)$,
an edge set $\operatorname{Edge}(G)$, and an endpoint map
\[
 \partial_G:\operatorname{Edge}(G)\longrightarrow
 \operatorname{Sym}^{2}(\operatorname{Vert}(G)).
\]
We call $G$ \textit{finite} if $\operatorname{Vert}(G)$ and $\operatorname{Edge}(G)$ are finite.
\begin{enumerate}[label=\arabic*.,ref=\arabic*,font=\itshape]
\item An edge $e$ is a \textit{loop at $v$} if $\partial_G(e)=\llbracket v,v\rrbracket$.
\item Edges $e$ and $e'$ are \textit{parallel} if $\partial_G(e)=\partial_G(e')$.
\item\label{d:graph-path-tuple} For $n\in\Z_{\ge0}$, consider a tuple
\[
 \gamma=(v_0,e_1,v_1,\ldots,e_n,v_n)
\]
such that $v_i\in\operatorname{Vert}(G)$ for $i\in[0,n]_{\Z}$ and
\[
 e_i\in\operatorname{Edge}(G),\qquad
 \partial_G(e_i)=\llbracket v_{i-1},v_i\rrbracket
 \qquad(i\in[n]_{\Z}).
\]
The tuple $\gamma$ starts at $v_0$ and ends at $v_n$.
We call $\gamma$ a \textit{simple path of length $n$} if $v_0,\ldots,v_n$ are pairwise distinct.
\item A \textit{cycle of length $n\in\Z_{\ge1}$} is a tuple
$(v_0,e_1,v_1,\ldots,e_n,v_n)$ such that $v_n=v_0$,
the vertices $v_0,\ldots,v_{n-1}\in\operatorname{Vert}(G)$ are pairwise distinct,
the edges $e_1,\ldots,e_n\in\operatorname{Edge}(G)$ are pairwise distinct, and
\[
 \partial_G(e_i)=\llbracket v_{i-1},v_i\rrbracket\qquad(i\in[n]_{\Z}).
\]
\end{enumerate}
\end{definition}

\begin{definition}[Graph notions]\label{def:graph-notions}
Let $G$ be a graph.
\begin{enumerate}[font=\itshape]
\item\label{d:graph-connected} We call $G$ \textit{connected} if, for every $\varnothing\ne U\subsetneq\operatorname{Vert}(G)$, there exist
$u\in U$, $v\in\operatorname{Vert}(G)\setminus U$, and $e\in\operatorname{Edge}(G)$ such that
\[
 \partial_G(e)=\llbracket u,v\rrbracket.
\]
\item We call $G$ \textit{simple} if it has no loops and $\partial_G$ is injective.
For a simple graph, we identify each edge with its image under $\partial_G$ and let
\[
 G=(\operatorname{Vert}(G),\operatorname{Edge}(G)),\qquad
 \operatorname{Edge}(G)\subseteq
 \{\llbracket u,v\rrbracket:u,v\in\operatorname{Vert}(G),\ u\ne v\}.
\]
We call a simple graph without cycles a \textit{forest} and a connected forest with at least one vertex a \textit{tree}.
\item Let $G'$ be a graph. A \textit{homomorphism} $\psi:G\to G'$ is a pair of maps
\[
 \psi_0:\operatorname{Vert}(G)\longrightarrow\operatorname{Vert}(G'),\qquad
 \psi_1:\operatorname{Edge}(G)\longrightarrow\operatorname{Edge}(G')
\]
such that, for every $e\in\operatorname{Edge}(G)$, $u\in\operatorname{Vert}(G)$, and $v\in\operatorname{Vert}(G)$,
\[
 \partial_G(e)=\llbracket u,v\rrbracket\quad\Longrightarrow\quad
 \partial_{G'}(\psi_1(e))=\llbracket \psi_0(u),\psi_0(v)\rrbracket.
\]
We call $\psi$ an \textit{isomorphism} if $\psi_0$ and $\psi_1$ are bijective.
We compose homomorphisms componentwise and let
$1_G=(1_{\operatorname{Vert}(G)},1_{\operatorname{Edge}(G)})$.
If $G$ and $G'$ are simple, we denote a homomorphism $\psi:G\to G'$ by its vertex map $\psi_0$.
\item Let $v\in\operatorname{Vert}(G)$. We define
\[
 \operatorname{AdjVert}_G(v)=\{u\in\operatorname{Vert}(G):
 \partial_G(e)=\llbracket v,u\rrbracket\text{ for some }e\in\operatorname{Edge}(G)\}.
\]
We call $u$ and $v$ \textit{adjacent} if $u\in\operatorname{AdjVert}_G(v)$.
We define the set of \textit{edges incident to $v$} by
\[
 \operatorname{IncEdge}_G(v)
 =\{e\in\operatorname{Edge}(G):\partial_G(e)=\llbracket v,u\rrbracket
 \text{ for some }u\in\operatorname{Vert}(G)\}.
\]
We define the \textit{degree of $v$} by
\[
\begin{aligned}
 \deg_G(v)
 &=\#\{e\in\operatorname{Edge}(G):\partial_G(e)=\llbracket v,u\rrbracket
       \text{ for some }u\ne v\}\\
 &\quad+2\#\{e\in\operatorname{Edge}(G):\partial_G(e)=\llbracket v,v\rrbracket\}.
\end{aligned}
\]
We define the set of \textit{leaves} by
\[
 \operatorname{Leaf}(G)=\{v\in\operatorname{Vert}(G):\deg_G(v)=1\}.
\]
\item\label{d:graph-subgraph} A graph $H$ is a \textit{subgraph} of $G$ if
\[
 \operatorname{Vert}(H)\subseteq\operatorname{Vert}(G),\qquad
 \operatorname{Edge}(H)\subseteq\operatorname{Edge}(G)
\]
and
\[
 \partial_H(e)=\partial_G(e)\qquad(e\in\operatorname{Edge}(H)).
\]
\item\label{d:graph-induced-subgraph} For $U\subseteq\operatorname{Vert}(G)$, the \textit{induced subgraph} $\operatorname{Subgraph}(G,U)$
has vertex set $U$ and edge set
\[
 \operatorname{Edge}(\operatorname{Subgraph}(G,U))
 =\{e\in\operatorname{Edge}(G):\partial_G(e)\in\operatorname{Sym}^{2}(U)\}.
\]
\item\label{d:graph-connected-components} Let $\operatorname{CtdComp}(G)$ denote the set of induced subgraphs
$\operatorname{Subgraph}(G,U)$ for nonempty $U\subseteq\operatorname{Vert}(G)$
that are maximal under inclusion among connected induced subgraphs of $G$.
We call an element of $\operatorname{CtdComp}(G)$ a \textit{connected component} of $G$.
\item\label{d:graph-edge-coloring} Let $C$ be a set.
An \textit{edge coloring} is a map $\kappa:\operatorname{Edge}(G)\to C$.
We call $(G,\kappa)$ a \textit{colored graph}.
\item A \textit{color-preserving homomorphism}
$(G,\kappa)\to(G',\kappa')$ is a graph homomorphism $\psi=(\psi_0,\psi_1):G\to G'$ such that
\[
 \kappa'\circ\psi_1=\kappa.
\]
We call it a \textit{color-preserving isomorphism} if $\psi_0$ and $\psi_1$ are bijective.
\item If $G$ is simple, we call $\kappa$ a \textit{proper edge coloring} if
\[
\begin{gathered}
 \llbracket u,v\rrbracket,\llbracket u,w\rrbracket\in\operatorname{Edge}(G),\qquad v\ne w\\
 \Longrightarrow\quad
 \kappa(\llbracket u,v\rrbracket)\ne\kappa(\llbracket u,w\rrbracket)\qquad(u,v,w\in\operatorname{Vert}(G)).
\end{gathered}
\]
\end{enumerate}
\end{definition}

\begin{definition}[Bipartite graphs]\label{def:bipartite-graphs}
Let $G$ be a graph. Let $I\subseteq\operatorname{Vert}(G)$ and $J\subseteq\operatorname{Vert}(G)$.
\begin{enumerate}[font=\itshape]
\item We call the pair $\mathbb G=(G,(I,J))$ a \textit{bipartite graph} if
\[
 \operatorname{Vert}(G)=I\sqcup J,\qquad
 \partial_G(\operatorname{Edge}(G))\subseteq
 \{\llbracket u,v\rrbracket:u\in I,\ v\in J\}.
\]
\item Let $\mathbb G=(G,(I,J))$ be a bipartite graph.
We call $I$ the \textit{source set} and $J$ the \textit{sink set} of $\mathbb G$.
Let $\operatorname{Src}_{\mathbb G}=I$ and $\operatorname{Snk}_{\mathbb G}=J$.
For $v\in\operatorname{Vert}(G)$, define
\[
\begin{aligned}
 \operatorname{SrcAdj}_{\mathbb G}(v)&=\operatorname{AdjVert}_G(v)\cap I,\\
 \operatorname{SnkAdj}_{\mathbb G}(v)&=\operatorname{AdjVert}_G(v)\cap J.
\end{aligned}
\]
\item Let $\mathbb G=(G,(I,J))$ be a bipartite graph. For a subgraph $H$ of $G$, define
\[
 \operatorname{SubSrc}_{\mathbb G}(H)=I\cap\operatorname{Vert}(H),\qquad
 \operatorname{SubSnk}_{\mathbb G}(H)=J\cap\operatorname{Vert}(H).
\]
\end{enumerate}
\end{definition}

\begin{lemma}[Rigidity of properly colored trees]\label{lem:colored-tree-rigidity}
Let $\mathbb T=(T,(I,J))$ be a bipartite graph such that $T$ is a finite tree.
Let $C$ be a set and let $\kappa:\operatorname{Edge}(T)\to C$ be a proper edge coloring.
Let $\phi=(\phi_0,\phi_1):T\to T$ be an automorphism satisfying
\[
 \kappa\circ\phi_1=\kappa,\qquad \phi_0(I)=I,\qquad \phi_0(J)=J.
\]
Then, $\phi=1_T$.
\end{lemma}
\begin{proof}
Let us use induction on $n=\#\operatorname{Vert}(T)$.
First, suppose $n=1$. The only vertex is fixed and there are no edges. Thus, $\phi=1_T$.

Second, suppose $n=2$. Let $I=\{u\}$ and $J=\{v\}$.
Since $\phi_0(I)=I$ and $\phi_0(J)=J$, we have $\phi_0(u)=u$ and $\phi_0(v)=v$.
The unique edge $\llbracket u,v\rrbracket$ is also fixed.

Third, suppose $n\ge3$. Let us remove the leaves and apply induction. Define
\[
 L=\operatorname{Leaf}(T),\qquad
 V'=\operatorname{Vert}(T)\setminus L,\qquad
 T'=\operatorname{Subgraph}(T,V').
\]
Let $(v_0,e_1,v_1,\ldots,e_\ell,v_\ell)$ be a longest simple path in $T$.
Then,
\[
 \ell\ge2,\qquad v_0,v_\ell\in L,\qquad
 \{v_1,\ldots,v_{\ell-1}\}\subseteq V'.
\]

Let us prove that $T'$ is connected.
For $u\in V'$ and $w\in V'$, let $P=(x_0,f_1,x_1,\ldots,f_r,x_r)$ be a simple path in $T$ from $x_0=u$ to $x_r=w$.
For $i\in[1,r-1]_{\Z}$, the distinct edges $f_i$ and $f_{i+1}$ are incident to $x_i$.
Hence, $\deg_T(x_i)\ge2$ and every vertex of $P$ belongs to $V'$.
Thus, $P$ lies in $T'$.
The graph $T'$ has no cycles. Therefore, $T'$ is a tree and
\[
 1\le\#\operatorname{Vert}(T')\le n-2.
\]

Let us apply the induction hypothesis to $T'$.
The automorphism $\phi$ preserves degrees. Hence, $\phi_0(L)=L$ and $\phi_0(V')=V'$.
Its restriction
\[
 \phi'=\bigl(\phi_0|_{V'},\phi_1|_{\operatorname{Edge}(T')}\bigr):T'\to T'
\]
preserves the coloring $\kappa|_{\operatorname{Edge}(T')}$ and the $(I\cap V',J\cap V')$-bipartition.
The induction hypothesis gives
\[
 \phi_0(v)=v\qquad(v\in V').
\]

Let us verify that $\phi_0(u)=u$ for each $u\in L$.
Its unique neighbor $v$ lies in $V'$ because $n\ge3$.
Let $e=\llbracket u,v\rrbracket$. Since $\phi_0(v)=v$,
\[
 e,\phi_1(e)\in\operatorname{IncEdge}_T(v),\qquad
 \kappa(\phi_1(e))=\kappa(e).
\]
Properness of $\kappa$ gives $\phi_1(e)=e$. Hence, $\phi_0(u)=u$.
Thus, $\phi_0=1_{\operatorname{Vert}(T)}$.
The injectivity of $\partial_T$ then gives $\phi_1=1_{\operatorname{Edge}(T)}$.
\end{proof}

\begin{definition}[Edge deletion and cut components]\label{def:cut-components}
Let $G$ be a graph.
\begin{enumerate}[font=\itshape]
\item\label{d:edge-deletion} For $e\in\operatorname{Edge}(G)$, we define a graph
\[
 \operatorname{DelEdge}(G,e)
 =\left(\operatorname{Vert}(G),\operatorname{Edge}(G)\setminus\{e\},
 \partial_G\big|_{\operatorname{Edge}(G)\setminus\{e\}}\right).
\]
\item\label{d:vertex-component} For $v\in\operatorname{Vert}(G)$, let $\operatorname{CtdComp}(G;v)$ be the connected component of $G$ containing $v$.
\item\label{d:sink-cut-component} Let $\mathbb S=(S,(I_S,J_S))$ be a bipartite graph such that $S$ is a tree.
Let $e=\llbracket x,b\rrbracket\in\operatorname{Edge}(S)$ for $x\in I_S$ and $b\in J_S$.
We define the \textit{sink-side cut component} by
\[
 \operatorname{SnkCtdComp}(\mathbb S;e)=\operatorname{CtdComp}(\operatorname{DelEdge}(S,e);b).
\]
\end{enumerate}
\end{definition}

\begin{definition}[Quivers]\label{def:quivers}
A \textit{quiver} is a triple
\[
 Q=(\operatorname{Vert}(Q),\operatorname{Edge}(Q),(s_Q,t_Q)).
\]
The sets $\operatorname{Vert}(Q)$ and $\operatorname{Edge}(Q)$ are its vertex and arrow sets.
The \textit{tail} and \textit{head} maps are
\[
 s_Q,t_Q:\operatorname{Edge}(Q)\longrightarrow\operatorname{Vert}(Q).
\]
We call $Q$ \textit{finite} if $\operatorname{Vert}(Q)$ and $\operatorname{Edge}(Q)$ are finite.
Let $\alpha:u\to v$ mean that $s_Q(\alpha)=u$ and $t_Q(\alpha)=v$.
If the map $\alpha\mapsto(s_Q(\alpha),t_Q(\alpha))$ is injective, we identify each arrow $\alpha$ with
$(s_Q(\alpha),t_Q(\alpha))$.
\end{definition}

\Needspace{6\baselineskip}
\begin{definition}[Directed paths and acyclicity]\label{def:directed-paths}
Let $Q$ be a quiver.
\begin{enumerate}[font=\itshape]
\item\label{d:quiver-directed-path} A \textit{directed path of length $n\in\Z_{\ge0}$} is a tuple
$(v_0,\alpha_1,v_1,\ldots,\alpha_n,v_n)$ such that
\[
 s_Q(\alpha_i)=v_{i-1},\qquad t_Q(\alpha_i)=v_i\qquad(i\in[n]_{\Z}).
\]
Here $v_i\in\operatorname{Vert}(Q)$ for $i\in[0,n]_{\Z}$ and
$\alpha_i\in\operatorname{Edge}(Q)$ for $i\in[n]_{\Z}$.
\item\label{d:quiver-directed-cycle} A directed path $(v_0,\alpha_1,v_1,\ldots,\alpha_n,v_n)$ is a \textit{directed cycle} if $n\ge1$, $v_n=v_0$, and $v_0,\ldots,v_{n-1}$ are pairwise distinct.
\item\label{d:quiver-acyclic} We call $Q$ \textit{acyclic} if it has no directed cycle.
\end{enumerate}
\end{definition}

\Needspace{6\baselineskip}
\begin{definition}[Quiver notions]\label{def:quiver-notions}
Let $Q$ be a quiver.
\begin{enumerate}[font=\itshape]
\item The \textit{underlying graph} $\operatorname{UnderG}(Q)$ of $Q$ has
\[
 \operatorname{Vert}(\operatorname{UnderG}(Q))=\operatorname{Vert}(Q),\qquad
 \operatorname{Edge}(\operatorname{UnderG}(Q))=\{\bar\alpha:\alpha\in\operatorname{Edge}(Q)\}.
\]
Here $\alpha\mapsto\bar\alpha$ is a bijection onto a copy of $\operatorname{Edge}(Q)$. The endpoint map of $\operatorname{UnderG}(Q)$ is
\[
 \partial_{\operatorname{UnderG}(Q)}(\bar\alpha)=\llbracket s_Q(\alpha),t_Q(\alpha)\rrbracket
 \qquad(\alpha\in\operatorname{Edge}(Q)).
\]
Let $\deg_Q(v)=\deg_{\operatorname{UnderG}(Q)}(v)$.
For $v\in\operatorname{Vert}(Q)$, we define the set of \textit{arrows incident to $v$} by
\[
\begin{aligned}
 \operatorname{IncEdgeQuiv}_Q(v)&:=\{\alpha\in\operatorname{Edge}(Q):\bar\alpha\in\operatorname{IncEdge}_{\operatorname{UnderG}(Q)}(v)\}\\
 &=\{\alpha\in\operatorname{Edge}(Q):s_Q(\alpha)=v\ \text{or}\ t_Q(\alpha)=v\}.
\end{aligned}
\]
We define $\operatorname{AdjVert}_Q(v)=\operatorname{AdjVert}_{\operatorname{UnderG}(Q)}(v)$.
We use $\operatorname{UnderG}(Q)$ to discuss trees and connected components of $Q$.
\item Let $Q$ and $Q'$ be quivers. A \textit{homomorphism} $\phi:Q\to Q'$ is a pair of maps
\[
 \phi_0:\operatorname{Vert}(Q)\longrightarrow\operatorname{Vert}(Q'),\qquad
 \phi_1:\operatorname{Edge}(Q)\longrightarrow\operatorname{Edge}(Q')
\]
satisfying
\[
 s_{Q'}\circ\phi_1=\phi_0\circ s_Q,\qquad
 t_{Q'}\circ\phi_1=\phi_0\circ t_Q.
\]
We call $\phi$ an \textit{isomorphism} if $\phi_0$ and $\phi_1$ are bijective.
We compose homomorphisms componentwise and let
$1_Q=(1_{\operatorname{Vert}(Q)},1_{\operatorname{Edge}(Q)})$.
The underlying graph homomorphism is
\[
 \operatorname{UnderG}(\phi)=(\phi_0,\bar\phi_1):
 \operatorname{UnderG}(Q)\longrightarrow\operatorname{UnderG}(Q'),\qquad
 \bar\phi_1(\bar\alpha)=\overline{\phi_1(\alpha)}.
\]
\item A quiver $Q'$ is a \textit{subquiver} of $Q$ if
\[
 \operatorname{Vert}(Q')\subseteq\operatorname{Vert}(Q),\qquad
 \operatorname{Edge}(Q')\subseteq\operatorname{Edge}(Q)
\]
and
\[
 s_{Q'}(\alpha)=s_Q(\alpha),\qquad t_{Q'}(\alpha)=t_Q(\alpha)
 \qquad(\alpha\in\operatorname{Edge}(Q')).
\]
For $U\subseteq\operatorname{Vert}(Q)$, the \textit{full subquiver} $\operatorname{SubQuiv}(Q,U)$
has vertex set $U$ and arrow set
\[
 \operatorname{Edge}(\operatorname{SubQuiv}(Q,U))
 =\{\alpha\in\operatorname{Edge}(Q):s_Q(\alpha),t_Q(\alpha)\in U\}.
\]
\end{enumerate}
\end{definition}

\begin{definition}[Bipartite quivers]\label{def:bipartite-quivers}
Let $Q$ be a quiver. Let $I\subseteq\operatorname{Vert}(Q)$ and $J\subseteq\operatorname{Vert}(Q)$.
\begin{enumerate}[font=\itshape]
\item We call the pair $\mathbb Q=(Q,(I,J))$ a \textit{bipartite quiver} if
\[
 \operatorname{Vert}(Q)=I\sqcup J,\qquad
 (s_Q,t_Q)(\operatorname{Edge}(Q))\subseteq I\times J.
\]
The underlying bipartite graph of $\mathbb Q$ is $(\operatorname{UnderG}(Q),(I,J))$.
\item Let $\mathbb Q=(Q,(I,J))$ and $\mathbb Q'=(Q',(I',J'))$ be bipartite quivers.
A quiver isomorphism $(\phi_0,\phi_1):Q\to Q'$ \textit{preserves the bipartitions} if
\[
 \phi_0(I)=I',\qquad \phi_0(J)=J'.
\]
\item Let $\mathbb Q=(Q,(I,J))$ be a bipartite quiver. We call $\mathbb Q$ \textit{simple} if the map
\[
 (s_Q,t_Q):\operatorname{Edge}(Q)\longrightarrow I\times J
\]
is injective.
\end{enumerate}
\end{definition}

See Remark~4.7(1) of~\cite{Weist} for arrow colors of bipartite quivers.

\begin{definition}[Colored quivers]\label{def:colored-quivers}
Let $Q$ be a quiver. Let $C$ be a set.
\begin{enumerate}[font=\itshape]
\item A \textit{colored quiver} with color set $C$ is a pair
\[
 \widetilde Q=(Q,c),\qquad c:\operatorname{Edge}(Q)\longrightarrow C.
\]
We define $\operatorname{QCol}(\widetilde Q)=c$ and use the notation
\[
\begin{gathered}
 \operatorname{Vert}(\widetilde Q)=\operatorname{Vert}(Q),\qquad
 \operatorname{Edge}(\widetilde Q)=\operatorname{Edge}(Q),\\
 s_{\widetilde Q}=s_Q,\qquad t_{\widetilde Q}=t_Q,\qquad
 \operatorname{UnderG}(\widetilde Q)=\operatorname{UnderG}(Q).
\end{gathered}
\]
Let
\[
 \deg_{\widetilde Q}(v)=\deg_Q(v),\qquad
 \operatorname{IncEdgeQuiv}_{\widetilde Q}(v)=\operatorname{IncEdgeQuiv}_Q(v)
 \qquad(v\in\operatorname{Vert}(Q)).
\]
For $v\in\operatorname{Vert}(Q)$, define
\[
 \operatorname{AdjVert}_{\widetilde Q}(v)=\operatorname{AdjVert}_Q(v).
\]
\item\label{d:colored-arrow} Let $\widetilde Q=(Q,c)$ be a colored quiver. Let $v\in\operatorname{Vert}(Q)$ and $d\in C$.
Define the set of \textit{$(v,d)$-colored arrows} by
\[
 \operatorname{ColArrow}_{\widetilde Q}(v,d)
 :=\{\alpha\in\operatorname{Edge}(Q):t_Q(\alpha)=v,\ c(\alpha)=d\}.
\]
\item\label{d:arrow-operators} Let $\widetilde Q=(Q,c)$ be a colored quiver.
\begin{enumerate}[label=(\alph*),ref=\alph*,font=\itshape]
\item\label{d:delete-arrow} For $\alpha\in\operatorname{Edge}(Q)$, we define
$\operatorname{DelArrow}(\widetilde Q,\alpha)$ by keeping $\operatorname{Vert}(Q)$,
replacing the arrow set by $\operatorname{Edge}(Q)\setminus\{\alpha\}$,
and restricting $s_Q$, $t_Q$, and $c$ to this set.
\item\label{d:add-arrow} For $x\in\operatorname{Vert}(Q)$, $v\in\operatorname{Vert}(Q)$, and $d\in C$, we define
$\operatorname{AddArrow}(\widetilde Q;x,v,d)$ by keeping $\operatorname{Vert}(Q)$
and adjoining a new arrow $\alpha_*$ to $\operatorname{Edge}(Q)$.
The source, target, and coloring maps extend those of $\widetilde Q$ by
\[
 s(\alpha_*)=x,\qquad t(\alpha_*)=v,\qquad c(\alpha_*)=d.
\]
\end{enumerate}
\item Let $\widetilde Q=(Q,c)$ and $\widetilde Q'=(Q',c')$ be colored quivers with color set $C$.
A \textit{color-preserving homomorphism} $\widetilde Q\to\widetilde Q'$
is a quiver homomorphism $(\phi_0,\phi_1):Q\to Q'$ such that
\[
 \operatorname{QCol}(\widetilde Q')\circ\phi_1
 =\operatorname{QCol}(\widetilde Q).
\]
We call this homomorphism a \textit{color-preserving isomorphism}
if $\phi_0$ and $\phi_1$ are bijective.
\item Let $\widetilde Q=(Q,c)$ and $\widetilde Q'=(Q',c')$ be colored quivers with color set $C$.
We call $\widetilde Q'$ a \textit{colored subquiver} of $\widetilde Q$ if $Q'$ is a subquiver of $Q$ and
\[
 c'=c\big|_{\operatorname{Edge}(Q')}.
\]
For $U\subseteq\operatorname{Vert}(Q)$, we define the \textit{full colored subquiver} on $U$ by
\[
 \operatorname{SubColQuiv}(\widetilde Q,U)=\left(\operatorname{SubQuiv}(Q,U),c\big|_{\operatorname{Edge}(\operatorname{SubQuiv}(Q,U))}\right).
\]
\item Let $\widetilde Q=(Q,c)$ be a colored quiver. Let $I\subseteq\operatorname{Vert}(Q)$ and $J\subseteq\operatorname{Vert}(Q)$.
We call the pair $\widetilde{\mathbb Q}=(\widetilde Q,(I,J))$ a \textit{colored bipartite quiver} if $\mathbb Q=(Q,(I,J))$ is a bipartite quiver.
We use
\[
\begin{gathered}
 \operatorname{Vert}(\widetilde{\mathbb Q})=\operatorname{Vert}(Q),\qquad
 \operatorname{Edge}(\widetilde{\mathbb Q})=\operatorname{Edge}(Q),\\
 s_{\widetilde{\mathbb Q}}=s_Q,\qquad t_{\widetilde{\mathbb Q}}=t_Q,\qquad
 \operatorname{QCol}(\widetilde{\mathbb Q})=c,\\
 \operatorname{UnderG}(\widetilde{\mathbb Q})=\operatorname{UnderG}(Q),\qquad
 \deg_{\widetilde{\mathbb Q}}(v)=\deg_Q(v)\quad(v\in\operatorname{Vert}(Q)).
\end{gathered}
\]
\item Let $\widetilde{\mathbb Q}=(\widetilde Q,(I,J))$ be a colored bipartite quiver.
We call $\widetilde{\mathbb Q}$ \textit{properly colored} if the restriction
\[
 c\big|_{\operatorname{IncEdgeQuiv}_Q(v)}:\operatorname{IncEdgeQuiv}_Q(v)\longrightarrow C
\]
is injective for every $v\in\operatorname{Vert}(Q)$.
\item\label{d:vacant-slot} Let $\widetilde{\mathbb Q}=(\widetilde Q,(I,J))$ be a colored bipartite quiver.
We define its set of \textit{colored-vacant slots} by
\[
 \operatorname{ColVacSlot}(\widetilde{\mathbb Q})
 :=\{(v,d)\in J\times C:
 \operatorname{ColArrow}_{\widetilde Q}(v,d)=\varnothing\}.
\]
We call $d$ the color of the colored-vacant slot $(v,d)$.
\end{enumerate}
\end{definition}

\subsection{Representations}\label{subsec:quiver-representations}
\begin{definition}[Representations of quivers]\label{def:quiver-representations}
Let $Q$ be a quiver.
\begin{enumerate}[font=\itshape]
\item A \textit{representation of $Q$ over $\C$} is a family
\[
 R=\bigl((R_v)_{v\in\operatorname{Vert}(Q)},
 (R_\alpha)_{\alpha\in\operatorname{Edge}(Q)}\bigr)
\]
such that each $R_v$ is a finite-dimensional complex vector space and each map
\[
 R_\alpha:R_{s_Q(\alpha)}\longrightarrow R_{t_Q(\alpha)}
 \qquad(\alpha\in\operatorname{Edge}(Q))
\]
is linear.
\item\label{d:representation-hom} Let $Q'$ be a quiver. Let $R$ and $R'$ be representations of $Q$ and $Q'$, respectively.
For a quiver homomorphism $\phi=(\phi_0,\phi_1):Q\to Q'$, a \textit{homomorphism}
$f:R\to R'$ \textit{over $\phi$} is a tuple $f=(f_v)_{v\in\operatorname{Vert}(Q)}$ of complex linear maps
$f_v:R_v\to R'_{\phi_0(v)}$ such that
\[
 f_{t_Q(\alpha)}\circ R_\alpha
 =R'_{\phi_1(\alpha)}\circ f_{s_Q(\alpha)}
 \qquad(\alpha\in\operatorname{Edge}(Q)).
\]
\item For a representation $R''$ of a quiver $Q''$ and a homomorphism $g:R'\to R''$
over $\psi:Q'\to Q''$, we define $g\circ f$ over $\psi\circ\phi$ by
\[
 (g\circ f)_v=g_{\phi_0(v)}\circ f_v\qquad(v\in\operatorname{Vert}(Q)).
\]
\item Let $\operatorname{Rep}_{\C}(Q)$ be the category of representations of $Q$ over $\C$ whose morphisms are homomorphisms over $1_Q$.
\item For $R\in\operatorname{Rep}_{\C}(Q)$, we define its \textit{dimension vector} and \textit{support} by
\[
\begin{aligned}
 \dim_{\C}(R)&=\bigl(\dim_{\C}(R_i)\bigr)_{i\in\operatorname{Vert}(Q)},\\
 \operatorname{RepSupp}(R)&=\{i\in\operatorname{Vert}(Q):R_i\ne0\}.
\end{aligned}
\]
The \textit{support quiver} is $\operatorname{SubQuiv}(Q,\operatorname{RepSupp}(R))$.
\item Let $R\in\operatorname{Rep}_{\C}(Q)$ and $R'\in\operatorname{Rep}_{\C}(Q)$.
An \textit{isomorphism} $\phi:R\to R'$ in $\operatorname{Rep}_{\C}(Q)$ is a tuple
\[
 \phi=(\phi_v)_{v\in\operatorname{Vert}(Q)}
\]
such that each $\phi_v:R_v\to R'_v$ is a complex linear isomorphism and
\[
 \phi_{t_Q(\alpha)}\circ R_\alpha
 =R'_\alpha\circ\phi_{s_Q(\alpha)}
 \qquad(\alpha\in\operatorname{Edge}(Q)).
\]
\item Let $R\in\operatorname{Rep}_{\C}(Q)$.
A \textit{subrepresentation} $W$ of $R$ consists of complex linear subspaces
$W_v\subseteq R_v$ for $v\in\operatorname{Vert}(Q)$ such that
\[
 R_\alpha(W_{s_Q(\alpha)})\subseteq W_{t_Q(\alpha)}
 \qquad(\alpha\in\operatorname{Edge}(Q)).
\]
Its arrow maps are
\[
 W_\alpha=R_\alpha\big|_{W_{s_Q(\alpha)}}:
 W_{s_Q(\alpha)}\longrightarrow W_{t_Q(\alpha)}.
\]
We denote the set of these subrepresentations by $\operatorname{SubRep}(R)$.
The zero subrepresentation has $W_v=0$ for every $v\in\operatorname{Vert}(Q)$.
We call $W$ \textit{proper} if $W\ne R$.
\end{enumerate}
\end{definition}

\Needspace{9\baselineskip}
\begin{definition}[Restriction and extension by zero]\label{def:rep-restriction-zero}
Let $Q$ be a quiver, let $U\subseteq\operatorname{Vert}(Q)$. Let $Q_U=\operatorname{SubQuiv}(Q,U)$.
\begin{enumerate}[font=\itshape]
\item Define
$\operatorname{Res}_{Q,U}:\operatorname{Rep}_{\C}(Q)\to\operatorname{Rep}_{\C}(Q_U)$
by keeping the spaces $R_v$ for $v\in U$, the maps $R_\alpha$ for $\alpha\in\operatorname{Edge}(Q_U)$, and the components $f_v$ for $v\in U$ of each homomorphism $f$.
\item Define
$\operatorname{ExtZero}_{Q,U}:\operatorname{Rep}_{\C}(Q_U)\to\operatorname{Rep}_{\C}(Q)$
by
\[
 (\operatorname{ExtZero}_{Q,U}R)_v=
 \begin{cases}R_v&v\in U,\\0&v\notin U,\end{cases}
 \qquad
 (\operatorname{ExtZero}_{Q,U}R)_\alpha=
 \begin{cases}R_\alpha&s_Q(\alpha),t_Q(\alpha)\in U,\\0&\text{otherwise}.\end{cases}
\]
For a homomorphism $f$, its component is $f_v$ on $U$ and the zero map outside $U$.
\end{enumerate}
\end{definition}

\begin{lemma}\label{lem:rep-restriction-zero}
Let $Q$ be a quiver, let $U\subseteq\operatorname{Vert}(Q)$. Let $Q_U=\operatorname{SubQuiv}(Q,U)$.
Then, we have the following.

\begin{enumerate}
\item\label{c:restriction-zero-left} Let $R\in\operatorname{Rep}_{\C}(Q_U)$.
The maps
\[
 (\varepsilon_R)_v=1_{R_v}\qquad(v\in U)
\]
define an isomorphism
\[
 \varepsilon_R:\operatorname{Res}_{Q,U}(\operatorname{ExtZero}_{Q,U}R)\longrightarrow R.
\]
\item\label{c:restriction-zero-right} Let $S\in\operatorname{Rep}_{\C}(Q)$ satisfy $\operatorname{RepSupp}(S)\subseteq U$.
The maps
\[
 (\eta_S)_v=
 \begin{cases}1_{S_v}&v\in U,\\0&v\notin U\end{cases}
\]
define an isomorphism
\[
 \eta_S:\operatorname{ExtZero}_{Q,U}(\operatorname{Res}_{Q,U}S)\longrightarrow S.
\]
\item\label{c:restriction-zero-subrep} Let $S$ satisfy the hypotheses of Claim~\ref{c:restriction-zero-right}.
The map
\[
 \rho_S:\operatorname{SubRep}(S)\longrightarrow
 \operatorname{SubRep}(\operatorname{Res}_{Q,U}S),\qquad
 W\longmapsto\operatorname{Res}_{Q,U}W
\]
is a bijection. Its inverse is
\[
 \sigma_S:\operatorname{SubRep}(\operatorname{Res}_{Q,U}S)\longrightarrow
 \operatorname{SubRep}(S),\qquad Z\longmapsto\widehat Z.
\]
For $Z\in\operatorname{SubRep}(\operatorname{Res}_{Q,U}S)$, define
\[
 \widehat Z_v=
 \begin{cases}Z_v&v\in U,\\0&v\notin U\end{cases}
 \qquad(v\in\operatorname{Vert}(Q)).
\]
For $W\in\operatorname{SubRep}(S)$ and $v\in U$,
\[
 \dim_{\C}(\rho_S(W)_v)=\dim_{\C}(W_v).
\]
\end{enumerate}
\end{lemma}
\begin{proof}
\prf{Claim~\ref{c:restriction-zero-left}} Let us compare the vertex spaces and arrow maps after restriction.
By Definition~\ref{def:rep-restriction-zero},
\[
\begin{aligned}
 (\operatorname{Res}_{Q,U}(\operatorname{ExtZero}_{Q,U}R))_v&=R_v
 &&(v\in U),\\
 (\operatorname{Res}_{Q,U}(\operatorname{ExtZero}_{Q,U}R))_\alpha&=R_\alpha
 &&(\alpha\in\operatorname{Edge}(Q_U)).
\end{aligned}
\]
Thus, the identity maps define $\varepsilon_R$ and its inverse.

\prf{Claim~\ref{c:restriction-zero-right}} Let us check the subrepresentation at the vertices outside $U$.
For $v\notin U$, we have $S_v=0$.
The map $(\eta_S)_v$ is therefore the unique isomorphism between zero vector spaces.
For $v\in U$, it is the identity on $S_v$.
Let $\widehat S=\operatorname{ExtZero}_{Q,U}(\operatorname{Res}_{Q,U}S)$.
For each $\alpha\in\operatorname{Edge}(Q)$, the equation
\[
 (\eta_S)_{t_Q(\alpha)}\circ\widehat S_\alpha
 =S_\alpha\circ(\eta_S)_{s_Q(\alpha)}
\]
holds as an equality of maps $\widehat S_{s_Q(\alpha)}\to S_{t_Q(\alpha)}$.
If $s_Q(\alpha)\in U$ and $t_Q(\alpha)\in U$, both sides equal $S_\alpha$. Otherwise, both sides are zero.
Hence, $\eta_S$ is an isomorphism of representations.

\prf{Claim~\ref{c:restriction-zero-subrep}} Let us verify $\rho_S(\widehat Z)=Z$ and $\widehat{\rho_S(W)}=W$.
Let $Z\in\operatorname{SubRep}(\operatorname{Res}_{Q,U}S)$ and $\alpha\in\operatorname{Edge}(Q)$.
If both endpoints of $\alpha$ lie in $U$, the inclusion
$S_\alpha(\widehat Z_{s_Q(\alpha)})\subseteq\widehat Z_{t_Q(\alpha)}$
is the subrepresentation condition for $Z$.
Otherwise, $S_{s_Q(\alpha)}=0$ or $S_{t_Q(\alpha)}=0$.
Then, $S_\alpha(\widehat Z_{s_Q(\alpha)})=0\subseteq\widehat Z_{t_Q(\alpha)}$.
For every $W\in\operatorname{SubRep}(S)$, we have $W_v=0$ for $v\notin U$.
Consequently,
\[
 \sigma_S\circ\rho_S=1_{\operatorname{SubRep}(S)},\qquad
 \rho_S\circ\sigma_S=1_{\operatorname{SubRep}(\operatorname{Res}_{Q,U}S)}.
\]
The dimension equality follows from $\rho_S(W)_v=W_v$ for $v\in U$.
\end{proof}

\begin{definition}[Euler and Tits forms]\label{def:euler-tits-forms}
Let $H$ be a finite acyclic quiver. For
\[
 d=(d_v)_{v\in\operatorname{Vert}(H)}\in\Z^{\operatorname{Vert}(H)},\qquad
 \delta=(\delta_v)_{v\in\operatorname{Vert}(H)}\in\Z^{\operatorname{Vert}(H)},
\]
we define the \textit{Euler form} and \textit{Tits quadratic form}
\[
\langle d,\delta\rangle_H=\sum_{v\in\operatorname{Vert}(H)}d_v\delta_v
 -\sum_{\alpha\in\operatorname{Edge}(H)}d_{s_H(\alpha)}\delta_{t_H(\alpha)},\qquad Q_H(d)=\langle d,d\rangle_H.
\]
\end{definition}

\begin{definition}[Representations of colored quivers]\label{def:colored-quiver-representations}
Let $C$ be a set, let $Q$ be a quiver, and let $\widetilde Q=(Q,c)$ be a colored quiver with color set $C$.
\begin{enumerate}[font=\itshape]
\item We define a \textit{representation $R$ of $\widetilde Q$ over $\C$} as a representation of $Q$ over $\C$.
\item Let $R$ be a representation of $\widetilde Q$ over $\C$. We define its \textit{colored support quiver} by
\[
 \operatorname{SuppColQuiv}(R):=\operatorname{SubColQuiv}(\widetilde Q,\operatorname{RepSupp}(R)).
\]
\item\begin{enumerate}[label=(\alph*),ref=\alph*,font=\itshape]
\item Let $\widetilde Q'=(Q',c')$ be a colored quiver with color set $C$.
Let $R$ and $R'$ be representations of $\widetilde Q$ and $\widetilde Q'$, respectively.
A \textit{color-preserving homomorphism} $\Phi:R\to R'$ is a tuple
\[
 \Phi=\bigl(\phi_0,\phi_1,(f_v)_{v\in\operatorname{Vert}(Q)}\bigr)
\]
such that $(\phi_0,\phi_1):\widetilde Q\to\widetilde Q'$ is a color-preserving quiver homomorphism
and the tuple $(f_v)_{v\in\operatorname{Vert}(Q)}$ is a homomorphism of representations over $(\phi_0,\phi_1)$.
\item For a composable color-preserving homomorphism
$\Psi=(\psi_0,\psi_1,(g_w)_{w\in\operatorname{Vert}(Q')})$, we define
\[
 \Psi\circ\Phi=
 \bigl(\psi_0\circ\phi_0,\psi_1\circ\phi_1,
 (g_{\phi_0(v)}\circ f_v)_{v\in\operatorname{Vert}(Q)}\bigr).
\]
The identity of $R$ is $(1_{\operatorname{Vert}(Q)},1_{\operatorname{Edge}(Q)},(1_{R_v})_{v\in\operatorname{Vert}(Q)})$.
\item Let $\operatorname{ColRep}_{\C}(\widetilde Q)$ denote the category of representations of $\widetilde Q$ and color-preserving homomorphisms.
For $R\in\operatorname{ColRep}_{\C}(\widetilde Q)$ and $R'\in\operatorname{ColRep}_{\C}(\widetilde Q)$, let
\[
 \operatorname{Hom}_{\operatorname{ColRep}_{\C}(\widetilde Q)}(R,R')
\]
denote the set of color-preserving homomorphisms $R\to R'$.
\item We call $\Phi$ a \textit{color-preserving isomorphism} if $\phi_0$ and $\phi_1$ are bijective
and every $f_v$ is a complex linear isomorphism.
\end{enumerate}
\item Let $R\in\operatorname{ColRep}_{\C}(\widetilde Q)$.
A \textit{colored subrepresentation} $W$ of $R$ is a representation of $\widetilde Q$
such that $W\in\operatorname{SubRep}(R)$.
We denote the set of these colored subrepresentations by $\operatorname{SubColRep}(R)$.
\item\label{d:total-dimensions} Let $(Q,(I,J))$ be a bipartite quiver. Let $R\in\operatorname{Rep}_{\C}(Q)$ have finite support $\operatorname{RepSupp}(R)$.
We define the \textit{total source dimension} and \textit{total sink dimension} of $R$ by
\[
 \operatorname{SrcDim}_{\C,(I,J)}(R)=\sum_{v\in I}\dim_{\C}R_v,\qquad
 \operatorname{SnkDim}_{\C,(I,J)}(R)=\sum_{v\in J}\dim_{\C}R_v.
\]
We define its \textit{total dimension pair} by
\[
 \operatorname{TotDim}_{\C,(I,J)}(R):=\bigl(\operatorname{SrcDim}_{\C,(I,J)}(R),\operatorname{SnkDim}_{\C,(I,J)}(R)\bigr).
\]
\end{enumerate}
\end{definition}

\begin{definition}[Bipartite stability]\label{def:bipartite-stability}
Let $\mathbb Q=(Q,(I,J))$ be a bipartite quiver.
Let $R\in\operatorname{Rep}_{\C}(Q)$ be nonzero and have finite support.
Let $(a,b)=\operatorname{TotDim}_{\C,(I,J)}(R)$.
We call $R$ \textit{stable} if every nonzero proper $W\in\operatorname{SubRep}(R)$ satisfies
\[
 a\operatorname{SnkDim}_{\C,(I,J)}(W)-b\operatorname{SrcDim}_{\C,(I,J)}(W)>0.
\]
\end{definition}

\begin{definition}[Arrow operators on representations]\label{def:representation-arrow-operators}
Let $\widetilde Q=(Q,c)$ be a colored quiver with color set $C$. Let $R\in\operatorname{Rep}_{\C}(Q)$.
\begin{enumerate}[font=\itshape]
\item\label{d:representation-delete-arrow} For $\alpha\in\operatorname{Edge}(Q)$, we define
$\operatorname{DelArrow}(R,\alpha)$ on $\operatorname{DelArrow}(\widetilde Q,\alpha)$
by keeping all vertex spaces and all arrow maps except $R_\alpha$.
\item\label{d:representation-add-arrow} For $x\in\operatorname{Vert}(Q)$, $v\in\operatorname{Vert}(Q)$, $d\in C$, and a linear map $f:R_x\to R_v$,
we define
\[
 \operatorname{AddArrow}(R;x,v,d,f)
\]
on $\operatorname{AddArrow}(\widetilde Q;x,v,d)$ by keeping the vertex spaces and
existing arrow maps and setting
\[
 (\operatorname{AddArrow}(R;x,v,d,f))_{\alpha_*}=f.
\]
\end{enumerate}
\end{definition}

\begin{definition}[Isomorphisms and markings of colored representations]\label{def:colored-representation-markings}
Let $C$ be a set. Let $\widetilde Q=(Q,c)$ and $\widetilde Q'=(Q',c')$ be colored quivers with color set $C$.
Let $\widetilde{\mathbb Q}=(\widetilde Q,(I,J))$ and $\widetilde{\mathbb Q}'=(\widetilde Q',(I',J'))$ be colored bipartite quivers such that $Q$ and $Q'$ are finite.
Let $R\in\operatorname{ColRep}_{\C}(\widetilde Q)$ and $R'\in\operatorname{ColRep}_{\C}(\widetilde Q')$.
\begin{enumerate}[font=\itshape]
\item\label{d:unmarked-colored-rep-iso} An isomorphism of $(\widetilde{\mathbb Q},R)$ and $(\widetilde{\mathbb Q}',R')$ is a color-preserving isomorphism $R\to R'$ that takes $(I,J)$ to $(I',J')$.
We denote the isomorphism class of $(\widetilde{\mathbb Q},R)$ by
$\operatorname{CIsomClass}_{\C}(\widetilde{\mathbb Q},R)$.
\item\label{d:ordered-colored-marking} Let $k\in\Z_{\ge0}$. An \textit{ordered marking of colored-vacant slots} of $R$ is a tuple
\[
 \mathbf s=((v_i,d_i))_{i\in[k]_{\Z}},\qquad
 (v_i,d_i)\in\operatorname{ColVacSlot}(\widetilde{\mathbb Q}).
\]
The marked representation is the pair $(R;\mathbf s)$.
Distinct marked colored-vacant slots means $(v_i,d_i)\ne(v_j,d_j)$ for $i\ne j$.
\item\label{d:marked-rep-category} Fix $\widetilde{\mathbb Q}$ and $\mathbf s$. We define the category
$\operatorname{MarkedColRep}_{\C}(\widetilde{\mathbb Q};\mathbf s)$ as follows.
Its objects are the pairs $(R;\mathbf s)$ with $R\in\operatorname{ColRep}_{\C}(\widetilde Q)$.
Its homomorphisms $(R;\mathbf s)\to(R';\mathbf s)$ are the homomorphisms
$\Phi=(\phi_0,\phi_1,(f_v)_v)\in\operatorname{Hom}_{\operatorname{ColRep}_{\C}(\widetilde Q)}(R,R')$
such that
\[
\begin{gathered}
 \phi_0(I)\subseteq I,\qquad
 \phi_0(J)\subseteq J,\\
 \phi_0(v_i)=v_i\qquad(i\in[k]_{\Z}).
\end{gathered}
\]
We use the identities and composition of $\operatorname{ColRep}_{\C}(\widetilde Q)$.
\item\label{d:marked-rep-iso} Let $\mathbf s'=((v'_i,d'_i))_{i\in[k]_{\Z}}$ be an ordered marking of $R'$.
An isomorphism $(\widetilde{\mathbb Q},R;\mathbf s)\to(\widetilde{\mathbb Q}',R';\mathbf s')$ is an isomorphism
$\Phi=(\phi_0,\phi_1,(f_v)_{v\in\operatorname{Vert}(Q)})$
of $(\widetilde{\mathbb Q},R)$ and $(\widetilde{\mathbb Q}',R')$ satisfying
\[
 (\phi_0(v_i),d_i)=(v'_i,d'_i)\qquad(i\in[k]_{\Z}).
\]
We denote its \textit{marked color-preserving isomorphism class} by
$\operatorname{MarkedCIsom}_{\C}(\widetilde{\mathbb Q},R;\mathbf s)$.
\end{enumerate}
\end{definition}

\subsection{Stable zero-sum representations and gluing}
From here to the end of this section, all colored quivers have color set $[0,2]_{\Z}$.

\begin{definition}[Zero-sum representations]\label{def:reflected-block}
Let $\mathbb U=(Q,(I_Q,J_Q))$ be a bipartite quiver such that $Q$ is finite.
We define the \textit{zero-sum representation} $\operatorname{ZS}(\mathbb U)$ of $Q$ as follows.
\begin{enumerate}[font=\itshape]
\item\label{d:reflected-kernel} For $u\in J_{Q}$, let $\operatorname{ZS}(\mathbb U)_u=\C$.
For $v\in I_Q$, put $A_v=\operatorname{IncEdgeQuiv}_Q(v)$.
On the space $\C^{A_v}$ of tuples $z=(z_\alpha)_{\alpha\in A_v}$, define
\begin{equation}\label{eq:zero-sum-map}
 \Sigma_v\colon\C^{A_v}\longrightarrow\C,\qquad
 \Sigma_v(z)=\sum_{\alpha\in A_v}z_\alpha,
\end{equation}
and let $\operatorname{ZS}(\mathbb U)_v=\ker\Sigma_v$.
\item\label{d:reflected-projections} For each arrow $\alpha:v\to u$ of $Q$, define
\begin{equation}\label{eq:zero-sum-arrow}
 \operatorname{ZS}(\mathbb U)_\alpha\colon
 \operatorname{ZS}(\mathbb U)_v\longrightarrow\operatorname{ZS}(\mathbb U)_u,
 \qquad z\longmapsto z_\alpha.
\end{equation}
\item Let $\widetilde Q=(Q,c)$ be a colored quiver.
Let $\widetilde{\mathbb Q}=(\widetilde Q,(I_Q,J_Q))$ be a colored bipartite quiver.
Let $\operatorname{ZS}(\widetilde{\mathbb Q})=\operatorname{ZS}(\mathbb U)$.
\end{enumerate}
\end{definition}

\begin{definition}[Colored source-trivalent trees and their quivers]\label{def:source-trivalent-tree}
Let $\mathbb T=(T,(I_T,J_T))$ be a bipartite graph such that $T$ is a finite tree.
\begin{enumerate}[font=\itshape]
\item We call $\mathbb T$ \textit{source-trivalent} if
\[
 \deg_T(u)=3\qquad(u\in I_T).
\]
\item A \textit{colored source-trivalent tree} is a pair $(\mathbb T,\kappa)$ such that $\mathbb T$ is source-trivalent and
\[
 \kappa\colon\operatorname{Edge}(T)\longrightarrow[0,2]_{\Z}
\]
is a proper edge coloring.
\item\label{d:reflected-quiver} For a colored source-trivalent tree $(\mathbb T,\kappa)$, define the colored bipartite quiver $Q(\mathbb T,\kappa)$ by
\[
\begin{aligned}
 I_{Q(\mathbb T,\kappa)}
 &=\{v\in J_{T}:\deg_T(v)\ge2\},\\
 J_{Q(\mathbb T,\kappa)}&=I_{T}.
\end{aligned}
\]
Its vertex set is $I_{Q(\mathbb T,\kappa)}\sqcup J_{Q(\mathbb T,\kappa)}$.
We define its arrows, tail and head maps, and colors by
\[
\begin{aligned}
 \operatorname{Edge}(Q(\mathbb T,\kappa))
 &=\{(v,u)\in I_{Q(\mathbb T,\kappa)}\times J_{Q(\mathbb T,\kappa)}:
 \llbracket u,v\rrbracket\in\operatorname{Edge}(T)\},\\
 s_{Q(\mathbb T,\kappa)}(v,u)&=v,\qquad t_{Q(\mathbb T,\kappa)}(v,u)=u,\\
 \operatorname{QCol}(Q(\mathbb T,\kappa))(v,u)&=\kappa(\llbracket u,v\rrbracket).
\end{aligned}
\]
\end{enumerate}
\end{definition}

\begin{lemma}[Source dimension of a subrepresentation]\label{lem:zero-sum-source-bound}
Let $(\mathbb T,\kappa)$ be a colored source-trivalent tree.
Let $\mathbb T=(T,(I_T,J_T))$.
Let $W\in\operatorname{SubColRep}(\operatorname{ZS}(Q(\mathbb T,\kappa)))$ satisfy $W\ne0$.
Let
\[
\begin{gathered}
 S=\{u\in I_{T}:\dim W_u=1\},\qquad n=\#S,\\
 u_W=\operatorname{SrcDim}_{\C,(I_{Q(\mathbb T,\kappa)},J_{Q(\mathbb T,\kappa)})}(W).
\end{gathered}
\]
Then, $n\ge1$ and $u_W\le n-1$.
\end{lemma}
\begin{proof}
Let us bound $\dim_{\C}W_v$ for $v\in I_{Q(\mathbb T,\kappa)}$.
Let
\[
 N=\{v\in J_{T}:\llbracket u,v\rrbracket\in\operatorname{Edge}(T)
       \text{ for some }u\in S\}.
\]
For $v\in J_{T}$, let
\[
 q_v=\#(S\cap\operatorname{SrcAdj}_{\mathbb T}(v)).
\]
Fix $v\in I_{Q(\mathbb T,\kappa)}$.
The arrows starting at $v$ are $(v,u)$ for $u\in\operatorname{SrcAdj}_{\mathbb T}(v)$.
For $z\in W_v$, equations~\eqref{eq:zero-sum-map} and~\eqref{eq:zero-sum-arrow} and the subrepresentation condition give
\begin{equation}\label{eq:zero-sum-subspace}
 \sum_{u\in\operatorname{SrcAdj}_{\mathbb T}(v)}z_{(v,u)}=0,
 \qquad z_{(v,u)}\in W_u
 \quad(u\in\operatorname{SrcAdj}_{\mathbb T}(v)).
\end{equation}
If $u\notin S$, then $W_u=0$ and $z_{(v,u)}=0$.
Thus, only the $q_v$ coordinates indexed by $u\in S$ can be nonzero. Their sum is zero by~\eqref{eq:zero-sum-subspace}.
If $q_v>0$, then $\dim_{\C}W_v\le q_v-1$. If $q_v=0$, then $W_v=0$.
If $n=0$, then $W_u=0$ for each $u\in J_{Q(\mathbb T,\kappa)}$ by the definition of $S$.
The bound on $\dim_{\C}W_v$ gives $W_v=0$ for each $v\in I_{Q(\mathbb T,\kappa)}$.
Then, $W=0$, contrary to the hypothesis. Hence, $n\ge1$.

Let us count the edges of $F=\operatorname{Subgraph}(T,S\sqcup N)$.
By Definition~\ref{def:graph-notions}(\ref{d:graph-induced-subgraph}),
\[
 \operatorname{Edge}(F)=
 \{\llbracket u,v\rrbracket\in\operatorname{Edge}(T):u\in S,\ v\in N\}.
\]
The graph $F$ is a forest because $T$ has no cycles.
Let $C=\#\operatorname{CtdComp}(F)\ge1$.
Each $u\in S$ has degree $3$ in $F$.
Each $v\in N\cap\operatorname{Leaf}(T)$ satisfies $q_v=1$.
For $v\in I_{Q(\mathbb T,\kappa)}\setminus N$, we have $W_v=0$.
Therefore,
\[
 3n=\#\operatorname{Edge}(F)=n+\#N-C,\qquad
 u_W\le\sum_{v\in N}(q_v-1)=3n-\#N=n-C\le n-1.
\]
\end{proof}

\begin{lemma}\label{lem:reflected-blocks}
Let $(\mathbb T,\kappa)$ be a colored source-trivalent tree.
Let $\mathbb T=(T,(I_T,J_T))$.
Assume $\#I_{T}=\nu\ge2$.
Let $H$ be the quiver underlying $Q(\mathbb T,\kappa)$ and let $d=\dim_{\C}(\operatorname{ZS}(Q(\mathbb T,\kappa)))$.
Then, we have the following.

\begin{enumerate}
\item\label{c:reflected-blocks-dimension} The representation $\operatorname{ZS}(Q(\mathbb T,\kappa))$ satisfies the following.
\begin{enumerate}[label=(\alph*),ref=\theenumi(\alph*)]
\item\label{c:reflected-blocks-dimension-vector} We have
\[
 \operatorname{TotDim}_{\C,(I_{Q(\mathbb T,\kappa)},J_{Q(\mathbb T,\kappa)})}(\operatorname{ZS}(Q(\mathbb T,\kappa)))=(\nu-1,\nu).
\]
\item\label{c:reflected-blocks-tits} We have $Q_H(d)=1$.
\end{enumerate}
\item\label{c:reflected-blocks-sinks} Each $u\in J_{Q(\mathbb T,\kappa)}$ satisfies
$\dim_{\C}(\operatorname{ZS}(Q(\mathbb T,\kappa))_u)=1$.
\item\label{c:reflected-blocks-sources} Each $v\in I_{Q(\mathbb T,\kappa)}$ satisfies
\[
 \bigl(\dim_{\C}(\operatorname{ZS}(Q(\mathbb T,\kappa))_v),\deg_{Q(\mathbb T,\kappa)}(v)\bigr)
 \in\{(1,2),(2,3)\}.
\]
\item\label{c:reflected-blocks-stable} The representation $\operatorname{ZS}(Q(\mathbb T,\kappa))$ is stable.
\item\label{c:reflected-blocks-generation} For each $u\in J_{Q(\mathbb T,\kappa)}$, we have
\[
 \sum_{\substack{\alpha\in\operatorname{Edge}(Q(\mathbb T,\kappa))\\t_{Q(\mathbb T,\kappa)}(\alpha)=u}}
 \operatorname{im}\bigl(\operatorname{ZS}(Q(\mathbb T,\kappa))_\alpha\bigr)
 =\operatorname{ZS}(Q(\mathbb T,\kappa))_u.
\]
\end{enumerate}
\end{lemma}
\begin{proof}
\prf{Claim~\ref{c:reflected-blocks-dimension-vector}} Let us compute the total dimension from the tree degrees.
Let $s=\#J_{T}$.
The identity $\#\operatorname{Edge}(T)=\#\operatorname{Vert}(T)-1$ gives
\[
 3\nu=\nu+s-1,\qquad s=2\nu+1.
\]
For $v\in J_{T}$, let $n_v=\deg_T(v)$.
By Definition~\ref{def:reflected-block}(\ref{d:reflected-kernel}), we obtain
\[
\begin{aligned}
 \operatorname{SrcDim}_{\C,(I_{Q(\mathbb T,\kappa)},J_{Q(\mathbb T,\kappa)})}(\operatorname{ZS}(Q(\mathbb T,\kappa)))
 &=\sum_{v\in J_{T}}(n_v-1)=3\nu-(2\nu+1)=\nu-1,\\
 \operatorname{SnkDim}_{\C,(I_{Q(\mathbb T,\kappa)},J_{Q(\mathbb T,\kappa)})}(\operatorname{ZS}(Q(\mathbb T,\kappa)))&=\#I_{T}=\nu.
\end{aligned}
\]
\prf{Claim~\ref{c:reflected-blocks-tits}} Claim~\ref{c:reflected-blocks-dimension-vector} gives
\[
Q_H(d)=\nu+\sum_{v\in J_{T}}\bigl((n_v-1)^2-n_v(n_v-1)\bigr)
 =\nu-\sum_{v\in J_{T}}(n_v-1)=1.
\]

\prf{Claim~\ref{c:reflected-blocks-sinks}} Let us count the sinks after deleting the leaves.
The sinks of $Q(\mathbb T,\kappa)$ belong to
$I_{T}$. By Definition~\ref{def:reflected-block}(\ref{d:reflected-kernel}),
$\operatorname{ZS}(Q(\mathbb T,\kappa))_u=\C$ for every such vertex $u$.

\prf{Claim~\ref{c:reflected-blocks-sources}} Let us count the sources after deleting the leaves.
By Definition~\ref{def:source-trivalent-tree}(2), the incident edges at each vertex have distinct colors in $[0,2]_{\Z}$. Hence, $n_v\le3$ for $v\in J_T$.
The source vertices of $Q(\mathbb T,\kappa)$ are
the vertices $v\in J_{T}$ satisfying $n_v\in[2,3]_{\Z}$.
For such a vertex, $\dim_{\C}(\operatorname{ZS}(Q(\mathbb T,\kappa))_v)=n_v-1$ and $\deg_{Q(\mathbb T,\kappa)}(v)=n_v$.

\prf{Claim~\ref{c:reflected-blocks-stable}} Let us prove stability by bounding each subrepresentation weight.
Let $W\in\operatorname{SubColRep}(\operatorname{ZS}(Q(\mathbb T,\kappa)))$ satisfy $0\ne W\ne\operatorname{ZS}(Q(\mathbb T,\kappa))$.
Claim~\ref{c:reflected-blocks-dimension-vector} gives the dimension pair $(\nu-1,\nu)$.
Define
\[
 S=\{u\in I_{T}:\dim W_u=1\},\qquad n=\#S,
\]
and let $u_W=\operatorname{SrcDim}_{\C,(I_{Q(\mathbb T,\kappa)},J_{Q(\mathbb T,\kappa)})}(W)$.
Lemma~\ref{lem:zero-sum-source-bound} gives $n\ge1$ and $u_W\le n-1$.
For $n<\nu$, we obtain
\[
 (\nu-1)n-\nu u_W\ge\nu-n>0.
\]
If $n=\nu$, properness gives $u_W\le\nu-2$. Hence,
\[
 (\nu-1)n-\nu u_W=\nu(\nu-1-u_W)\ge\nu>0.
\]

\prf{Claim~\ref{c:reflected-blocks-generation}} Let us verify surjectivity of the arrow maps in Definition~\ref{def:reflected-block}(\ref{d:reflected-projections}).
Since $T$ is connected and $\nu\ge2$, every $u\in I_{T}$
has a vertex $v\in\operatorname{SnkAdj}_{\mathbb T}(u)$ of degree $d\ge2$.
Let $e_1,\ldots,e_d$ be the standard basis of $\C^d$.
If $i$ and $j$ are distinct coordinate indices, then
\[
 e_i-e_j\in\ker\!\left(\C^d\xrightarrow{\;\sum\;}\C\right),
 \qquad (e_i-e_j)_i=1.
\]
Thus, each arrow map from this zero-sum space to $\C$ is surjective.
\end{proof}

\begin{lemma}\label{lem:vacant-slots}

Let $(\mathbb T,\kappa)$ be a colored source-trivalent tree.
Let $\mathbb T=(T,(I_T,J_T))$.
The map sending $(u,d)\in\operatorname{ColVacSlot}(Q(\mathbb T,\kappa))$ to the unique $v\in J_T$ satisfying
\[
 \deg_T(v)=1,\qquad \llbracket u,v\rrbracket\in\operatorname{Edge}(T),\qquad
 \kappa(\llbracket u,v\rrbracket)=d
\]
is a bijection onto $J_T\cap\operatorname{Leaf}(T)$.
\end{lemma}
\begin{proof}
Let us construct the inverse map from sink leaves to colored-vacant slots.
For $u\in I_{T}$, the three edges in $\operatorname{IncEdge}_T(u)$
have distinct colors in $[0,2]_{\Z}$.
The quiver $Q(\mathbb T,\kappa)$ keeps exactly the edges from $u$ to vertices of degree at least two.
Thus, $(u,d)\in\operatorname{ColVacSlot}(Q(\mathbb T,\kappa))$ determines the unique leaf
$v\in J_{T}$ such that $\llbracket u,v\rrbracket\in\operatorname{Edge}(T)$
and $\kappa(\llbracket u,v\rrbracket)=d$.
Conversely, let $v\in J_T$ be a leaf and let $u\in I_T$ satisfy $\operatorname{SrcAdj}_{\mathbb T}(v)=\{u\}$.
The inverse map is
\[
 J_T\cap\operatorname{Leaf}(T)\longrightarrow\operatorname{ColVacSlot}(Q(\mathbb T,\kappa)),
 \qquad v\longmapsto\bigl(u,\kappa(\llbracket u,v\rrbracket)\bigr).
\]
Uniqueness of the neighbor of $v$ and of the edge of color $d$ at $u$ shows that the two maps are inverse.
\end{proof}

\begin{definition}[Gluing zero-sum representations]\label{def:attachments}
Let $\mathbb S=(S,(I_S,J_S))$ be a source-trivalent tree.
For each $b\in J_S$, choose a colored source-trivalent tree $(\mathbb T_b,\kappa_b)$.
Let $\mathbb T_b=(T_b,(I_{T_b},J_{T_b}))$ for $b\in J_S$.
For each $e=\llbracket x,b\rrbracket\in\operatorname{Edge}(S)$ with
$x\in I_{S}$ and $b\in J_{S}$, choose
\[
 (v_e,d_e)\in\operatorname{ColVacSlot}(Q(\mathbb T_b,\kappa_b)),\qquad
 0\ne f_e:\C\longrightarrow\operatorname{ZS}(Q(\mathbb T_b,\kappa_b))_{v_e}.
\]
Let
\[
 \mathbf T=((\mathbb T_b,\kappa_b))_{b\in J_{S}},\qquad
 \mathbf p=((v_e,d_e))_{e\in\operatorname{Edge}(S)},\qquad
 \mathbf f=(f_e)_{e\in\operatorname{Edge}(S)}.
\]
We define the colored representation
\[
 R=\operatorname{Glue}_{\C}(\mathbb S;\mathbf T,\mathbf p,\mathbf f)
\]
by the following conditions.
\begin{enumerate}[font=\itshape]
\item Its quiver consists of disjoint copies of the $Q(\mathbb T_b,\kappa_b)$, a new vertex
$x$ for each $x\in I_{S}$, and one new arrow
\[
 \alpha_e:x\longrightarrow v_e,\qquad c(\alpha_e)=d_e
 \qquad(e=\llbracket x,b\rrbracket\in\operatorname{Edge}(S)).
\]
Each $v_e$ belongs to the copy of $Q(\mathbb T_b,\kappa_b)$ indexed by $b$.
The sources are the new vertices $x\in I_S$ and the copied vertices in $I_{Q(\mathbb T_b,\kappa_b)}$; the sinks are the copied vertices in $J_{Q(\mathbb T_b,\kappa_b)}$.
\item Its restriction to each $Q(\mathbb T_b,\kappa_b)$ is $\operatorname{ZS}(Q(\mathbb T_b,\kappa_b))$.
\item For each new vertex $x$, let $R_x=\C$.
For each new arrow $\alpha_e$, let $R_{\alpha_e}=f_e$.
\end{enumerate}
\end{definition}
\subsection{Flow sum equations}\label{subsec:flows-vertex-sums}
\begin{definition}[Flow sum equations]\label{def:flows-vertex-sums}
Let $\mathbb S=(S,(I_S,J_S))$ be a bipartite graph with a finite graph $S$. Let
\[
 \alpha=(\alpha_v\in\R)_{v\in I_{S}},\qquad
 \xi=(\xi_i\in\R)_{i\in J_{S}}.
\]
\begin{enumerate}[font=\itshape]
\item We define $\operatorname{Flow}(S)=\R^{\operatorname{Edge}(S)}$.
For $f\in\operatorname{Flow}(S)$, let $f_w=f(w)$.
We call
\[
\begin{aligned}
 \sum_{w\in\operatorname{IncEdge}_S(v)}f_w&=\alpha_v
 &&(v\in I_{S}),\\
 \sum_{w\in\operatorname{IncEdge}_S(i)}f_w&=\xi_i
 &&(i\in J_{S})
\end{aligned}
\]
the \textit{$(\alpha,\xi)$-flow sum equations on $S$}.
\item Define the set of \textit{$(\alpha,\xi)$-flows on $S$} by
\[
 \operatorname{Flow}(S,(\alpha,\xi))=
 \{f\in\operatorname{Flow}(S):f\text{ satisfies the }(\alpha,\xi)\text{-flow sum equations}\}.
\]
\item We call an $(\alpha,\xi)$-flow $f$ \textit{positive} if $f_w>0$ for every $w\in\operatorname{Edge}(S)$.
\end{enumerate}
\end{definition}

\begin{lemma}[Flows on a bipartite tree]\label{lem:tree-flow-cut}
Let $\mathbb S=(S,(I_S,J_S))$ be a bipartite graph such that $S$ is a finite tree.
Let $\alpha=(\alpha_v)_{v\in I_S}\in\R^{I_S}$ and
$\xi=(\xi_i)_{i\in J_S}\in\R^{J_S}$.
Assume
\begin{equation}\label{eq:tree-vertex-sum-balance}
 \sum_{v\in I_{S}}\alpha_v
 =\sum_{i\in J_{S}}\xi_i.
\end{equation}
Then, we have the following.

\begin{enumerate}
\item\label{c:tree-flow-existence} The set $\operatorname{Flow}(S,(\alpha,\xi))$ has exactly one element $f$.
\item\label{c:tree-flow-cut} Let $w=\llbracket v,i\rrbracket\in\operatorname{Edge}(S)$
with $i\in J_{S}$.
Let $S_w=\operatorname{SnkCtdComp}(\mathbb S;w)$.
Then,
\[
 f_w=\sum_{j\in \operatorname{SubSnk}_{\mathbb S}(S_w)}\xi_j
      -\sum_{u\in \operatorname{SubSrc}_{\mathbb S}(S_w)}\alpha_u.
\]
\item\label{c:tree-flow-positive-boundary} Assume $f_w>0$ for every $w\in\operatorname{Edge}(S)$.
Let $I\subseteq J_{S}$ and $J\subseteq I_{S}$.
Assume $\varnothing\ne I\sqcup J\subsetneq\operatorname{Vert}(S)$.
Assume $\operatorname{SnkAdj}_{\mathbb S}(v)\subseteq I$ for every $v\in J$.
Then,
\[
 \sum_{i\in I}\xi_i-\sum_{v\in J}\alpha_v>0.
\]
\end{enumerate}
\end{lemma}
\begin{proof}
\prf{Claim~\ref{c:tree-flow-existence}} Let us prove existence and uniqueness by induction on the number of vertices.

First, suppose $S$ has one vertex. If the vertex is a source $v$, equation~\eqref{eq:tree-vertex-sum-balance} gives $\alpha_v=0$.
If the vertex is a sink $i$, that equation gives $\xi_i=0$.
The empty map is the unique $(\alpha,\xi)$-flow.

Second, suppose $\#\operatorname{Vert}(S)\ge2$. Choose $z\in\operatorname{Leaf}(S)$ and let
\[
 S'=\operatorname{Subgraph}(S,\operatorname{Vert}(S)\setminus\{z\}).
\]
The graph $S'$ is a tree with one fewer vertex.

First, let us remove a source leaf $z=v\in I_S$. Let $\operatorname{SnkAdj}_{\mathbb S}(v)=\{i\}$ and $w=\llbracket v,i\rrbracket$.
The flow equation at $v$ forces $f_w=\alpha_v$.
On $S'$, define
\[
 \alpha'_u=\alpha_u\quad(u\in I_S\setminus\{v\}),\qquad
 \xi'_j=\begin{cases}\xi_i-\alpha_v,&j=i,\\\xi_j,&j\in J_S\setminus\{i\}.\end{cases}
\]
Both sums in~\eqref{eq:tree-vertex-sum-balance} decrease by $\alpha_v$.

Second, let us remove a sink leaf $z=i\in J_S$. Let $\operatorname{SrcAdj}_{\mathbb S}(i)=\{v\}$ and $w=\llbracket v,i\rrbracket$.
The flow equation at $i$ forces $f_w=\xi_i$.
On $S'$, define
\[
 \xi'_j=\xi_j\quad(j\in J_S\setminus\{i\}),\qquad
 \alpha'_u=\begin{cases}\alpha_v-\xi_i,&u=v,\\\alpha_u,&u\in I_S\setminus\{v\}.\end{cases}
\]
Both sums in~\eqref{eq:tree-vertex-sum-balance} decrease by $\xi_i$.

Let us extend the flow from $S'$ to $S$.
In either case, $(\alpha',\xi')$ satisfies~\eqref{eq:tree-vertex-sum-balance} on $S'$.
By induction, there is a unique $(\alpha',\xi')$-flow on $S'$.
The flow extends uniquely to $S$ by the forced value of $f_w$.

\prf{Claim~\ref{c:tree-flow-cut}} Let us compute the weight of the cut edge. Sum the sink equations over $S_w$ and subtract the source equations over $S_w$.
For every $e\in\operatorname{Edge}(S_w)$, the terms $f_e$ cancel. The edge $w$ contributes $f_w$.

\prf{Claim~\ref{c:tree-flow-positive-boundary}} Let us cancel the terms indexed by edges $\llbracket v,i\rrbracket$ with $v\in J$ and $i\in I$.
The assumption $\operatorname{SnkAdj}_{\mathbb S}(v)\subseteq I$ for every $v\in J$ gives
\[
 \sum_{i\in I}\xi_i-\sum_{v\in J}\alpha_v
 =\sum_{\substack{\llbracket v,i\rrbracket\in\operatorname{Edge}(S)\\
                  i\in I,\ v\notin J}}f_{\llbracket v,i\rrbracket}.
\]
Since $S$ is connected and $\varnothing\ne I\sqcup J\subsetneq\operatorname{Vert}(S)$,
there exist $i\in I$ and $v\notin J$ such that $\llbracket v,i\rrbracket\in\operatorname{Edge}(S)$.
Every summand is positive.
\end{proof}

\begin{lemma}[Weights of subrepresentations of a zero-sum representation]\label{lem:zero-sum-weight}
Let $(\mathbb T,\kappa)$ be a colored source-trivalent tree.
Let $\mathbb T=(T,(I_T,J_T))$, $Q=Q(\mathbb T,\kappa)$, and $Z=\operatorname{ZS}(Q)$.
Let $a_0\in\Z_{\ge0}$ satisfy
\[
 \operatorname{TotDim}_{\C,(I_Q,J_Q)}(Z)=(a_0,a_0+1).
\]
Let $0\ne W\in\operatorname{SubColRep}(Z)$. Let
\[
 (u,t)=\operatorname{TotDim}_{\C,(I_Q,J_Q)}(W).
\]
Let $y\in\R_{\ge0}$.
Then, the following statements hold.
\begin{enumerate}
\item\label{c:zero-sum-weight-bound} We have
\[
 y(t-u)-u\ge y-a_0.
\]
\item\label{c:zero-sum-weight-equality} We have
\[
 y(t-u)-u=y-a_0\quad\Longleftrightarrow\quad W=Z.
\]
\end{enumerate}
\end{lemma}
\begin{proof}
\prf{Claim~\ref{c:zero-sum-weight-bound}} Let us prove $t-u\ge1$ and $u\le a_0$.
Lemma~\ref{lem:zero-sum-source-bound} gives $t-u\ge1$.
For every $z\in I_Q$, the inclusion $W_z\subseteq Z_z$ gives
$\dim_\C W_z\le\dim_\C Z_z$.
Summing over $I_Q$ gives $u\le a_0$.
Therefore,
\[
 y(t-u)-u-(y-a_0)=y(t-u-1)+(a_0-u)\ge0.
\]

\prf{Claim~\ref{c:zero-sum-weight-equality}} Let us prove that a nonzero proper subrepresentation satisfies $u<a_0$.
Every $z\in J_Q=I_T$ has $Z_z=\C$.
Hence, $\#I_T=\#J_Q=a_0+1$.
If $a_0=0$, then $J_Q=\{z_0\}$ for one vertex $z_0$,
$Z_{z_0}=\C$. Every source space is zero.
The condition $W\ne0$ gives $W_{z_0}=\C$. Thus, $W=Z$.

Let us treat the case in which the subrepresentation has full source dimension.
Suppose that $a_0\ge1$ and $u=a_0$.
Then,
\[
 0=a_0-u=\sum_{z\in I_Q}
       \bigl(\dim_\C Z_z-\dim_\C W_z\bigr).
\]
Each summand is nonnegative.
Thus, $W_z=Z_z$ for every $z\in I_Q$.
Since $\#I_T=a_0+1\ge2$, Claim~\ref{c:reflected-blocks-generation}
of Lemma~\ref{lem:reflected-blocks} applies.
For every $z\in J_Q$, the subrepresentation condition gives
\[
 Z_z
 =\sum_{\substack{\beta\in\operatorname{Edge}(Q)\\t_Q(\beta)=z}}
       Z_\beta(Z_{s_Q(\beta)})
 =\sum_{\substack{\beta\in\operatorname{Edge}(Q)\\t_Q(\beta)=z}}
       Z_\beta(W_{s_Q(\beta)})
 \subseteq W_z\subseteq Z_z.
\]
Consequently, $W=Z$.
For $W\ne Z$, we therefore have $a_0-u>0$, and
\[
 y(t-u)-u-(y-a_0)=y(t-u-1)+(a_0-u)>0.
\]
For $W=Z$, substitution of $(u,t)=(a_0,a_0+1)$ gives equality.
\end{proof}

\begin{lemma}[Weights of subrepresentations of a gluing]\label{lem:gluing-subrepresentation-weights}
Let $\mathbb S=(S,(I_S,J_S))$ be a source-trivalent tree.
For $i\in J_S$, let $(\mathbb T_i,\kappa_i)$ be a colored source-trivalent tree, and let
\[
 Q_i=Q(\mathbb T_i,\kappa_i),\qquad Z_i=\operatorname{ZS}(Q_i).
\]
Let $a_i\in\Z_{\ge0}$ satisfy
\[
 \operatorname{TotDim}_{\C,(I_{Q_i},J_{Q_i})}(Z_i)=(a_i,a_i+1)
 \qquad(i\in J_S).
\]
For every $e=\llbracket v,i\rrbracket\in\operatorname{Edge}(S)$ with
$v\in I_S$ and $i\in J_S$, fix a colored-vacant slot and a nonzero complex linear map
\[
 p_e=(z_e,d_e)\in\operatorname{ColVacSlot}(Q_i),\qquad
 0\ne\phi_e:\C\longrightarrow(Z_i)_{z_e}.
\]
Let
\[
 \mathbf T=((\mathbb T_i,\kappa_i))_{i\in J_S},\qquad
 \mathbf p=(p_e)_{e\in\operatorname{Edge}(S)},\qquad
 \mathbf f=(\phi_e)_{e\in\operatorname{Edge}(S)},
\]
and let $R=\operatorname{Glue}_\C(\mathbb S;\mathbf T,\mathbf p,\mathbf f)$.
Use the vertex copies in Definition~\ref{def:attachments}. Let
\[
 I_R=I_S\sqcup\bigsqcup_{i\in J_S}I_{Q_i},\qquad
 J_R=\bigsqcup_{i\in J_S}J_{Q_i}.
\]
Let
\[
 c=\#I_S,\qquad k=c+1,\qquad
 a=\sum_{i\in J_S}a_i+c,\qquad x=\frac ak.
\]
Let $\alpha_v=x+1$ for $v\in I_S$ and $\xi_i=x-a_i$ for $i\in J_S$.
Fix $f\in\operatorname{Flow}(S,(\alpha,\xi))$.
Then, the following statements hold.
\begin{enumerate}
\item\label{c:gluing-weight-decomposition} Let $W\in\operatorname{SubColRep}(R)$.
For $i\in J_S$, let $W_i$ be the subrepresentation of $Z_i$ defined by
\[
 (W_i)_z=W_z\quad(z\in\operatorname{Vert}(Q_i)),\qquad
 (W_i)_\beta=R_\beta|_{W_{s_{Q_i}(\beta)}}
 \quad(\beta\in\operatorname{Edge}(Q_i)).
\]
Define
\[
 I=\{i\in J_S:W_i\ne0\},\qquad
 J=\{v\in I_S:W_v\ne0\}.
\]
\Needspace{5\baselineskip}
Let
\[
 (u_i,t_i)=\operatorname{TotDim}_{\C,(I_{Q_i},J_{Q_i})}(W_i)
 \quad(i\in I),\qquad
 (u,t)=\operatorname{TotDim}_{\C,(I_R,J_R)}(W).
\]
Then, $\operatorname{SnkAdj}_{\mathbb S}(v)\subseteq I$ for every $v\in J$. We have
\begin{equation}\label{eq:gluing-boundary-flow}
 \sum_{i\in I}(x-a_i)-\#J(x+1)
 =
 \sum_{\substack{e=\llbracket v,i\rrbracket\in\operatorname{Edge}(S)\\
                  v\in I_S\setminus J,\ i\in I}}f_e.
\end{equation}
Moreover,
\begin{equation}\label{eq:gluing-subobject-weight}
\begin{aligned}
 \frac{at-(a+k)u}{k}
 &=
 \sum_{\substack{e=\llbracket v,i\rrbracket\in\operatorname{Edge}(S)\\
                  v\in I_S\setminus J,\ i\in I}}f_e\\
 &\quad+\sum_{i\in I}
       \bigl(x(t_i-u_i-1)+(a_i-u_i)\bigr).
\end{aligned}
\end{equation}
\item\label{c:gluing-cut-subrepresentation} Let
$w=\llbracket v,i\rrbracket\in\operatorname{Edge}(S)$ with $v\in I_S$ and $i\in J_S$.
Let
\[
 S_w=\operatorname{SnkCtdComp}(\mathbb S;w),\qquad
 A_w=\operatorname{SubSnk}_{\mathbb S}(S_w),\qquad
 B_w=\operatorname{SubSrc}_{\mathbb S}(S_w).
\]
For each vertex $z$ of the quiver of $R$, define the subspace
\[
 W^w_z=
 \begin{cases}
 R_z,&z\in B_w\sqcup\bigsqcup_{j\in A_w}\operatorname{Vert}(Q_j),\\
 0,&\text{otherwise}.
 \end{cases}
\]
These subspaces and the restrictions of the arrow maps of $R$ define
a nonzero proper subrepresentation $W^w$.
For $(u_w,t_w)=\operatorname{TotDim}_{\C,(I_R,J_R)}(W^w)$, we have
\[
 u_w=\sum_{j\in A_w}a_j+\#B_w,\qquad
 t_w=\sum_{j\in A_w}(a_j+1),\qquad
 at_w-(a+k)u_w=kf_w.
\]
\end{enumerate}
\end{lemma}
\begin{proof}
\prf{Claim~\ref{c:gluing-weight-decomposition}} Let us determine which restrictions of $W$ are nonzero.
For $v\in J$, the inclusion $0\ne W_v\subseteq R_v=\C$ gives $W_v=\C$.
For $e=\llbracket v,i\rrbracket\in\operatorname{Edge}(S)$ with $v\in J$, we have
\[
 0\ne\phi_e(\C)=\phi_e(W_v)\subseteq W_{z_e}.
\]
Thus, $W_i\ne0$ and $i\in I$.
This proves $\operatorname{SnkAdj}_{\mathbb S}(v)\subseteq I$ for $v\in J$.

Let us sum the flow equations over $I$ and $J$.
They give
\[
 \sum_{i\in I}(x-a_i)-\#J(x+1)
 =
 \sum_{i\in I}\sum_{e\in\operatorname{IncEdge}_S(i)}f_e
 -
 \sum_{v\in J}\sum_{e\in\operatorname{IncEdge}_S(v)}f_e.
\]
An edge between $J$ and $I$ occurs in both sums.
There is no edge from $J$ to $J_S\setminus I$.
Subtracting the two finite sums proves~\eqref{eq:gluing-boundary-flow}.

Let us compute the stability weight from the vertex dimensions.
The disjoint vertex sets in the gluing give
\[
 u=\sum_{i\in I}u_i+\#J,\qquad t=\sum_{i\in I}t_i.
\]
Since $a=kx$ and $k>0$, we obtain
\[
\begin{aligned}
 \frac{at-(a+k)u}{k}
 &=x(t-u)-u\\
 &=\sum_{i\in I}\bigl(x(t_i-u_i)-u_i\bigr)-\#J(x+1)\\
 &=\sum_{i\in I}(x-a_i)-\#J(x+1)
   +\sum_{i\in I}\bigl(x(t_i-u_i-1)+(a_i-u_i)\bigr).
\end{aligned}
\]
Substitution of~\eqref{eq:gluing-boundary-flow} proves~\eqref{eq:gluing-subobject-weight}.

\prf{Claim~\ref{c:gluing-cut-subrepresentation}} Let us verify the subrepresentation condition at each arrow.
For $\beta\in\operatorname{Edge}(Q_j)$, both endpoint spaces of $W^w$ are full if $j\in A_w$ and zero if $j\notin A_w$.
For an added arrow $\alpha_e:y\to z_e$ indexed by $e=\llbracket y,j\rrbracket$, suppose $W^w_y\ne0$.
Then, $y\in B_w$.
The source endpoint $v$ of $w$ does not belong to $B_w$. Thus, $e\ne w$.
The edge $e$ therefore belongs to the component $S_w$ whenever $y\in B_w$.
Hence, $j\in A_w$ and
\[
 R_{\alpha_e}(W^w_y)=\phi_e(R_y)\subseteq R_{z_e}=W^w_{z_e}.
\]
If $W^w_y=0$, the subrepresentation condition holds because $\phi_e(0)=0$.

Let us check that the cut subrepresentation is nonzero and proper.
The index $i$ belongs to $A_w$. The representation $Z_i$ has sink dimension $a_i+1>0$.
Thus, $W^w\ne0$.
Also, $v\notin B_w$ gives $W^w_v=0\ne R_v=\C$. Thus, $W^w\ne R$.

Let us compute the dimension pair and its weight.
The full constituent spaces over $A_w$ and the full new source spaces over $B_w$ give
\[
 u_w=\sum_{j\in A_w}a_j+\#B_w,\qquad
 t_w=\sum_{j\in A_w}(a_j+1).
\]
Therefore,
\[
\begin{aligned}
 \frac{at_w-(a+k)u_w}{k}
 &=x(t_w-u_w)-u_w\\
 &=\sum_{j\in A_w}(x-a_j)-\sum_{y\in B_w}(x+1)=f_w.
\end{aligned}
\]
Summing the source and sink equations for $f$ gives~\eqref{eq:tree-vertex-sum-balance}.
The equality $\sum_{j\in A_w}(x-a_j)-\sum_{y\in B_w}(x+1)=f_w$ follows from Claim~\ref{c:tree-flow-cut} of Lemma~\ref{lem:tree-flow-cut}.
Multiplying by $k$ gives $at_w-(a+k)u_w=kf_w$.
\end{proof}

\begin{lemma}[Stability under gluing]\label{lem:flow-stability}
Let $\mathbb S=(S,(I_S,J_S))$ be a source-trivalent tree.
Let $c=\#I_{S}$.
For $i\in J_S$, let $(\mathbb T_i,\kappa_i)$ be a colored source-trivalent tree.
Let $\mathbf T=((\mathbb T_i,\kappa_i))_{i\in J_S}$.
Let $R=\operatorname{Glue}_{\C}(\mathbb S;\mathbf T,\mathbf p,\mathbf f)$ be a gluing of their zero-sum representations using colored-vacant slots $\mathbf p$ and nonzero complex linear maps $\mathbf f$.
For $i\in J_S$, let
\[
 \operatorname{TotDim}_{\C,(I_{Q(\mathbb T_i,\kappa_i)},J_{Q(\mathbb T_i,\kappa_i)})}
 (\operatorname{ZS}(Q(\mathbb T_i,\kappa_i)))=(a_i,a_i+1).
\]
Let
\[
 \widetilde Q_R=\operatorname{SuppColQuiv}(R),\qquad
 G_R=\operatorname{UnderG}(\widetilde Q_R).
\]
Let $\widetilde{\mathbb Q}_R=(\widetilde Q_R,(I_R,J_R))$ for
\[
\begin{aligned}
 I_R&=\operatorname{RepSupp}(R)\cap
 \left(I_S\sqcup\bigsqcup_{i\in J_S}I_{Q(\mathbb T_i,\kappa_i)}\right),\\
 J_R&=\operatorname{RepSupp}(R)\cap
 \bigsqcup_{i\in J_S}J_{Q(\mathbb T_i,\kappa_i)}.
\end{aligned}
\]
The disjoint unions use the vertex copies in $R$.
Assume the following conditions on the added arrows.
\begin{hypotheses}
\item\label{hyp:gluing-attachments} Each edge $\llbracket x,i\rrbracket$ of $S$ with $x\in I_{S}$ and $i\in J_{S}$ specifies one added $(v,d)$-colored arrow $x\to v$ for some $(v,d)\in\operatorname{ColVacSlot}(Q(\mathbb T_i,\kappa_i))$.
\item\label{hyp:gluing-slots} Each $(v,d)\in\operatorname{ColVacSlot}(Q(\mathbb T_i,\kappa_i))$ receives at most one added $(v,d)$-colored arrow. The three arrows starting at each $x\in I_{S}$ have distinct colors.
\end{hypotheses}
Let
\[
 k=c+1,\qquad a=\sum_{i\in J_{S}}a_i+c,\qquad x=\frac ak.
\]
For $\alpha_v=x+1$ and $\xi_i=x-a_i$, the $(\alpha,\xi)$-flow sum equations on $S$ are
\begin{equation}
\sum_{w\in\operatorname{IncEdge}_S(v)}f_w=x+1\quad(v\in I_{S}),\qquad
\sum_{w\in\operatorname{IncEdge}_S(i)}f_w=x-a_i\quad(i\in J_{S}).\label{flow}
\end{equation}
Then, we have the following.

\begin{enumerate}
\Needspace{6\baselineskip}
\item\label{c:flow-stability-tree} The following statements hold.
\begin{enumerate}[label=(\alph*),ref=\theenumi(\alph*)]
\item\label{c:flow-stability-support-tree} The graph $G_R$ is a tree.
\item\label{c:flow-stability-support-colors} The coloring $\operatorname{QCol}(\widetilde Q_R)$ induces a proper edge coloring of $G_R$.
\item\label{c:flow-stability-block-count} $\#J_{S}=2c+1$.
\item\label{c:flow-stability-total-dimension} We have $\operatorname{TotDim}_{\C,(I_R,J_R)}(R)=(a,a+k)$.
\end{enumerate}
\item\label{c:flow-stability-unique} The $(\alpha,\xi)$-flow sum equations~\eqref{flow} have a unique solution.
\item\label{c:flow-stability-positive} Let $f$ be the unique element of $\operatorname{Flow}(S,(\alpha,\xi))$ from Claim~\ref{c:flow-stability-unique}.
The representation $R$ is stable if and only if $f_w>0$ for every $w\in\operatorname{Edge}(S)$.
\end{enumerate}
\end{lemma}
\begin{proof}
\prf{Claim~\ref{c:flow-stability-support-tree}} Let us verify that $G_R$ is a tree.
Let $B=J_{S}$ and $b=\#B$.
Counting the edges of $S$ gives $3c=c+b-1$. Thus, $b=2c+1$.
For $i\in B$, let
\[
 G_i=\operatorname{UnderG}(\operatorname{SuppColQuiv}(\operatorname{ZS}(Q(\mathbb T_i,\kappa_i)))),\qquad
 v_i=\#\operatorname{Vert}(G_i),\qquad e_i=\#\operatorname{Edge}(G_i).
\]
Definition~\ref{def:reflected-block}(\ref{d:reflected-kernel}) gives $\#I_{T_i}=a_i+1\ge1$.
By Definitions~\ref{def:reflected-block} and~\ref{def:source-trivalent-tree}, $G_i$ is obtained from $T_i$ by deleting its sink leaves.
Every path between retained vertices avoids the deleted leaves. Thus, $G_i$ is a nonempty tree and $e_i=v_i-1$.
Replacing each $i\in B$ by $G_i$ and adding the arrows in assumption~\ref{hyp:gluing-attachments} gives the connected graph $G_R$.
Then,
\[
\begin{aligned}
 \#\operatorname{Vert}(G_R)&=c+\sum_{i\in B}v_i,\\
 \#\operatorname{Edge}(G_R)&=3c+\sum_{i\in B}e_i
 =3c+\sum_{i\in B}(v_i-1)
 =\#\operatorname{Vert}(G_R)-1.
\end{aligned}
\]
Hence, $G_R$ is a tree. Two colored-vacant slots sharing a vertex contribute two distinct added arrows to $3c$.
\prf{Claim~\ref{c:flow-stability-support-colors}} Assumption~\ref{hyp:gluing-slots} gives distinct colors at each new source.
At each old vertex, the old incident arrows have distinct colors because $\kappa_i$ is proper.
Each added arrow uses a colored-vacant slot by assumption~\ref{hyp:gluing-attachments}.
Its color is therefore unused by the old arrows at its target.
Assumption~\ref{hyp:gluing-slots} gives distinct colors among the added arrows at each old vertex.

\prf{Claim~\ref{c:flow-stability-block-count}} Let us count the sinks of the gluing tree.
Let $b=\#J_{S}$.
The identity $3c=c+b-1$ in the proof of Claim~\ref{c:flow-stability-support-tree} gives $b=2c+1$.

\prf{Claim~\ref{c:flow-stability-total-dimension}}
By Claim~\ref{c:flow-stability-block-count}, we obtain
\[
\begin{aligned}
 \operatorname{SrcDim}_{\C,(I_R,J_R)}(R)
 &=\sum_{i\in J_{S}}a_i+c=a,\\
 \operatorname{SnkDim}_{\C,(I_R,J_R)}(R)
 &=\sum_{i\in J_{S}}(a_i+1)=a-c+2c+1=a+k.
\end{aligned}
\]

\prf{Claim~\ref{c:flow-stability-unique}} Let us verify equality of the two total vertex weights.
Claim~\ref{c:flow-stability-block-count} gives
\[
\begin{aligned}
 \sum_{i\in J_{S}}\xi_i
 &=\sum_{i\in J_{S}}(x-a_i)\\
 &=(2c+1)x-(a-c)=c(x+1)
 =\sum_{v\in I_{S}}\alpha_v.
\end{aligned}
\]
Claim~\ref{c:tree-flow-existence} of Lemma~\ref{lem:tree-flow-cut} gives the unique $(\alpha,\xi)$-flow.

\prf{Claim~\ref{c:flow-stability-positive}} Let us obtain stability from positive flows.
Assume $f_e>0$ for every $e\in\operatorname{Edge}(S)$.
Let $0\ne W\ne R$ be a subrepresentation of $R$.
For $i\in J_S$, let
\[
 Q_i=Q(\mathbb T_i,\kappa_i),\qquad Z_i=\operatorname{ZS}(Q_i),
\]
and let $W_i$ be the restriction of $W$ to the copy of $Q_i$.
Define
\[
 I=\{i\in J_S:W_i\ne0\},\qquad
 J=\{v\in I_S:W_v\ne0\}.
\]
\Needspace{5\baselineskip}
Let
\[
 (u_i,t_i)=\operatorname{TotDim}_{\C,(I_{Q_i},J_{Q_i})}(W_i)\quad(i\in I),
 \qquad
 (u,t)=\operatorname{TotDim}_{\C,(I_R,J_R)}(W).
\]
We have $a_i\ge0$, $a\ge0$, $k>0$, and $x\ge0$.
Claim~\ref{c:zero-sum-weight-bound} of Lemma~\ref{lem:zero-sum-weight} gives
\[
 x(t_i-u_i-1)+(a_i-u_i)\ge0\qquad(i\in I).
\]
Claim~\ref{c:gluing-weight-decomposition}
of Lemma~\ref{lem:gluing-subrepresentation-weights} gives
$\operatorname{SnkAdj}_{\mathbb S}(v)\subseteq I$ for $v\in J$ and
the weight formula~\eqref{eq:gluing-subobject-weight}.
The dimension pair of $R$ is $(a,a+k)$ by Claim~\ref{c:flow-stability-total-dimension}.

First, suppose $I\sqcup J\subsetneq\operatorname{Vert}(S)$.
Since $W\ne0$, we have $I\sqcup J\ne\varnothing$.
Equation~\eqref{eq:gluing-boundary-flow} and
Claim~\ref{c:tree-flow-positive-boundary} of Lemma~\ref{lem:tree-flow-cut} give
\[
 \sum_{\substack{e=\llbracket v,i\rrbracket\in\operatorname{Edge}(S)\\
                  v\in I_S\setminus J,\ i\in I}}f_e
 =\sum_{i\in I}(x-a_i)-\#J(x+1)>0.
\]
Hence,~\eqref{eq:gluing-subobject-weight} gives $at-(a+k)u>0$.

Second, suppose $I\sqcup J=\operatorname{Vert}(S)$.
Then, $I=J_S$, $J=I_S$, and $W_v=R_v=\C$ for every $v\in I_S$.
Since $W\ne R$, there is $i_0\in J_S$ such that $0\ne W_{i_0}\ne Z_{i_0}$.
Claim~\ref{c:zero-sum-weight-equality} of Lemma~\ref{lem:zero-sum-weight}
and Claim~\ref{c:zero-sum-weight-bound} of Lemma~\ref{lem:zero-sum-weight} give
\[
 x(t_{i_0}-u_{i_0}-1)+(a_{i_0}-u_{i_0})>0.
\]
For $i\in I\setminus\{i_0\}$, Claim~\ref{c:zero-sum-weight-bound} of Lemma~\ref{lem:zero-sum-weight} gives $x(t_i-u_i-1)+(a_i-u_i)\ge0$.
The boundary-edge sum in~\eqref{eq:gluing-subobject-weight} is empty because $J=I_S$.
Thus, $at-(a+k)u>0$ in this case as well.
By Definition~\ref{def:bipartite-stability}, $R$ is stable.

Let us prove the converse. Suppose $R$ is stable.
Fix $w=\llbracket v,i\rrbracket\in\operatorname{Edge}(S)$ with $v\in I_S$ and $i\in J_S$.
Claim~\ref{c:gluing-cut-subrepresentation}
of Lemma~\ref{lem:gluing-subrepresentation-weights}
gives a nonzero proper subrepresentation $W^w$ of $R$.
For
\[
 (u_w,t_w)=\operatorname{TotDim}_{\C,(I_R,J_R)}(W^w),
\]
Claim~\ref{c:gluing-cut-subrepresentation} of Lemma~\ref{lem:gluing-subrepresentation-weights} gives
\[
 at_w-(a+k)u_w=kf_w.
\]
Stability gives $kf_w>0$.
Since $k>0$, we obtain $f_w>0$.
\end{proof}

\subsection{Exact counts for fixed colored-vacant-slot colors}
\begin{definition}\label{def:marked-block-counts}
Let $j\in\Z_{\ge1}$. Let $d\in[0,2]_{\Z}$ and $\epsilon\in[0,2]_{\Z}\setminus\{d\}$.
We consider all colored source-trivalent trees $(\mathbb T,\kappa)$ such that
\[
 \operatorname{TotDim}_{\C,(I_{Q(\mathbb T,\kappa)},J_{Q(\mathbb T,\kappa)})}(\operatorname{ZS}(Q(\mathbb T,\kappa)))=(j,j+1).
\]
\begin{enumerate}[font=\itshape]
\item We define $P_j$ to be the number of classes
\[
 \operatorname{MarkedCIsom}_{\C}(Q(\mathbb T,\kappa),\operatorname{ZS}(Q(\mathbb T,\kappa));((v,d))),
 \qquad (v,d)\in\operatorname{ColVacSlot}(Q(\mathbb T,\kappa)).
\]
\item We define $Q_j^{=}$ to be the number of classes
\[
 \operatorname{MarkedCIsom}_{\C}(Q(\mathbb T,\kappa),\operatorname{ZS}(Q(\mathbb T,\kappa));((v_1,d),(v_2,d)))
\]
for $(v_i,d)\in\operatorname{ColVacSlot}(Q(\mathbb T,\kappa))$ with $i\in[2]_{\Z}$ and $v_1\ne v_2$.
\item We define $Q_j^{\ne}$ to be the number of classes
\[
 \operatorname{MarkedCIsom}_{\C}(Q(\mathbb T,\kappa),\operatorname{ZS}(Q(\mathbb T,\kappa));((v_1,d),(v_2,\epsilon)))
\]
for $(v_1,d)\in\operatorname{ColVacSlot}(Q(\mathbb T,\kappa))$ and $(v_2,\epsilon)\in\operatorname{ColVacSlot}(Q(\mathbb T,\kappa))$.
\end{enumerate}
\end{definition}
\begin{definition}[Colored branches and their series]\label{def:colored-branches}
Let $d\in[0,2]_{\Z}$. Let $t$, $y_0$, $y_1$, and $y_2$ be commuting indeterminates.
\begin{enumerate}[font=\itshape]
\item\label{d:branch-set} Let $\mathbb B=(B,(I_B,J_B))$ be a bipartite graph such that $B$ is a finite tree.
We define $\operatorname{Branch}(\mathbb B;d)$ to be the set of pairs $(\kappa,r)$ satisfying the following conditions.
\begin{enumerate}[label=(\alph*),ref=\theenumi(\alph*),font=\itshape]
\item $\kappa:\operatorname{Edge}(B)\to[0,2]_{\Z}$ is a proper edge coloring.
\item $r\in\operatorname{Vert}(B)$ and $\kappa(e)\ne d$ for every $e\in\operatorname{IncEdge}_B(r)$.
\item For every $v\in I_{B}$,
\[
 \deg_B(v)=\begin{cases}2,&v=r,\\3,&v\ne r.\end{cases}
\]
\end{enumerate}
A \textit{branch with parent color $d$} is a tuple $\mathfrak b=(\mathbb B,\kappa,r)$ such that $(\kappa,r)\in\operatorname{Branch}(\mathbb B;d)$.
We call $r$ its \textit{root}.
\item\label{d:branch-types} We call $\mathfrak b$ a \textit{source branch} if $r\in I_{B}$ and a \textit{sink branch} if $r\in J_{B}$.
\item\label{d:branch-parent} For $\mathfrak b=(\mathbb B,\kappa,r)$, form $B^+$ by adjoining a new vertex $p$ in the part opposite to $r$ and the edge $\llbracket p,r\rrbracket$.
We call $p$ the \textit{parent vertex} and $\llbracket p,r\rrbracket$ the \textit{parent edge}.
Extend $\kappa$ to $\kappa^+:\operatorname{Edge}(B^+)\to[0,2]_{\Z}$ by setting
\[
 \kappa^+(e)=\kappa(e)\quad(e\in\operatorname{Edge}(B)),\qquad
 \kappa^+(\llbracket p,r\rrbracket)=d.
\]
\item\label{d:branch-isomorphism} Let $\mathfrak b=(\mathbb B,\kappa,r)$ and $\mathfrak b'=(\mathbb B',\kappa',r')$ be branches with parent color $d$.
Let $\mathbb B=(B,(I_B,J_B))$ and $\mathbb B'=(B',(I_{B'},J_{B'}))$.
An isomorphism $\mathfrak b\to\mathfrak b'$ is a graph isomorphism
$\psi=(\psi_0,\psi_1):B\to B'$ such that $\psi_0(I_B)=I_{B'}$, $\psi_0(J_B)=J_{B'}$, and
\[
 \psi_0(r)=r',\qquad \kappa'(\psi_1(e))=\kappa(e)\quad(e\in\operatorname{Edge}(B)).
\]
Let $\mathcal B_d^{\mathrm{src}}$ and $\mathcal B_d^{\mathrm{snk}}$ be the sets of isomorphism classes of source and sink branches with parent color $d$.
\item\label{d:branch-statistics} For $\mathfrak b=(\mathbb B,\kappa,r)$, define
\[
 n(\mathfrak b)=\#I_{B}.
\]
For $e\in[0,2]_{\Z}$, let $\ell_e(\mathfrak b)$ be the number of vertices
$v\in J_{B}$ such that $\deg_{B^+}(v)=1$ and the unique edge incident to $v$ in $B^+$ has color $e$ under $\kappa^+$.
\item\label{d:branch-weight} We define the weight of $\mathfrak b$ by
\[
 w(\mathfrak b)=t^{n(\mathfrak b)}\prod_{e=0}^2 y_e^{\ell_e(\mathfrak b)}.
\]
\item\label{d:tree-weight} Let $\mathcal T$ be the set of color-preserving, bipartition-preserving isomorphism classes of colored source-trivalent trees $(\mathbb T,\kappa)$.
For $\mathbb T=(T,(I_T,J_T))$, define $n(\mathbb T)=\#I_T$.
For $e\in[0,2]_{\Z}$, let $\ell_e(\mathbb T,\kappa)$ count the sink leaves of $T$ whose unique incident edge has color $e$.
Define
\[
 w(\mathbb T,\kappa)=t^{n(\mathbb T)}\prod_{e=0}^2 y_e^{\ell_e(\mathbb T,\kappa)}.
\]
\item\label{d:branch-series} We define
\[
 T_d=\sum_{[\mathfrak b]\in\mathcal B_d^{\mathrm{src}}}w(\mathfrak b),\qquad
 D_d=\sum_{[\mathfrak b]\in\mathcal B_d^{\mathrm{snk}}}w(\mathfrak b),\qquad
 U=\sum_{[(\mathbb T,\kappa)]\in\mathcal T}w(\mathbb T,\kappa).
\]
\end{enumerate}
\end{definition}

\Needspace{7\baselineskip}
\begin{lemma}[Branch-series identities]\label{lem:branch-series}
Let $\{d,\epsilon,f\}=[0,2]_{\Z}$.
Then, we have the following.
\begin{enumerate}
\item\label{c:slots-finite} For every $d\in[0,2]_{\Z}$ and $n\in\Z_{\ge0}$, there are finitely many classes in each of
$\mathcal B_d^{\mathrm{src}}$, $\mathcal B_d^{\mathrm{snk}}$, and $\mathcal T$ with source count $n$.
The weights $w(\mathfrak b)$ and $w(\mathbb T,\kappa)$ are invariant under branch isomorphisms and color-preserving, bipartition-preserving tree isomorphisms, respectively.
Thus, $T_d$, $D_d$, and $U$ belong to $\Z[y_0,y_1,y_2][[t]]$.
\item\label{c:branch-series-recursions} We have
\begin{equation}\label{eq:marked-branch-series}
D_d=y_d+T_\epsilon+T_f+T_\epsilon T_f,\qquad T_d=tD_\epsilon D_f.
\end{equation}
\item\label{c:branch-series-leaf} We have
\[
 \frac{\partial U}{\partial y_d}=T_d.
\]
\end{enumerate}
\end{lemma}
\begin{proof}
\prf{Claim~\ref{c:slots-finite}} Let us count the vertices for each source count $n$.
For $\mathfrak b=(\mathbb B,\kappa,r)$ with $\mathbb B=(B,(I_B,J_B))$ and $\#I_B=n$, we have
\[
 \#\operatorname{Edge}(B)=\sum_{v\in I_B}\deg_B(v)
 =\begin{cases}
 3n-1,&r\in I_B,\\
 3n,&r\in J_B.
 \end{cases}
\]
Since $B$ is a tree,
\[
 \#\operatorname{Vert}(B)=\#\operatorname{Edge}(B)+1
 =\begin{cases}
 3n,&r\in I_B,\\
 3n+1,&r\in J_B.
 \end{cases}
\]
For a source-trivalent tree $\mathbb T=(T,(I_T,J_T))$ with $\#I_T=n$, we obtain
\[
 \#\operatorname{Edge}(T)=\sum_{v\in I_T}3=3n,\qquad
 \#\operatorname{Vert}(T)=3n+1.
\]
For each of these vertex counts $N$, relabel the vertices by $[N]_{\Z}$.
There are finitely many choices of
\[
 I\subseteq[N]_{\Z},\qquad J=[N]_{\Z}\setminus I,\qquad
 E\subseteq\{\llbracket u,v\rrbracket:u,v\in[N]_{\Z},\ u\ne v\},\qquad
 \kappa:E\to[0,2]_{\Z}.
\]
For a branch, there are also finitely many choices of $r\in[N]_{\Z}$.
Every class at source count $n$ has such a representative.
Thus, each set of classes at source count $n$ is finite.

Let us verify that isomorphisms preserve the weights.
Let $\psi=(\psi_0,\psi_1):\mathfrak b\to\mathfrak b'$ be an isomorphism of branches with parent color $d$.
Let $\mathfrak b'=(\mathbb B',\kappa',r')$ and $\mathbb B'=(B',(I_{B'},J_{B'}))$.
Let $p$ and $p'$ be their parent vertices.
Extend $\psi$ to $\psi^+:B^+\to B'^+$ by
\[
\begin{aligned}
 \psi_0^+(v)&=\psi_0(v)\quad(v\in\operatorname{Vert}(B)),&
 \psi_0^+(p)&=p',\\
 \psi_1^+(g)&=\psi_1(g)\quad(g\in\operatorname{Edge}(B)),&
 \psi_1^+(\llbracket p,r\rrbracket)&=\llbracket p',r'\rrbracket.
\end{aligned}
\]
Since $\psi_0(r)=r'$ and both parent edges have color $d$, we have
\[
\begin{aligned}
 \deg_{B'^+}(\psi_0^+(v))&=\deg_{B^+}(v)
 &&(v\in\operatorname{Vert}(B^+)),\\
 \kappa'^+(\psi_1^+(g))&=\kappa^+(g)
 &&(g\in\operatorname{Edge}(B^+)).
\end{aligned}
\]
Also, $\psi_0(I_B)=I_{B'}$ and $\psi_0(J_B)=J_{B'}$.
Consequently,
\[
 n(\mathfrak b)=n(\mathfrak b'),\qquad
 \ell_e(\mathfrak b)=\ell_e(\mathfrak b')\quad(e\in[0,2]_{\Z}),
\]
and
\[
 w(\mathfrak b)
 =t^{n(\mathfrak b)}\prod_{e=0}^2y_e^{\ell_e(\mathfrak b)}
 =t^{n(\mathfrak b')}\prod_{e=0}^2y_e^{\ell_e(\mathfrak b')}
 =w(\mathfrak b').
\]

Let us verify that isomorphisms of colored bipartite trees preserve their weights.
For a color-preserving, bipartition-preserving tree isomorphism
$(\mathbb T,\kappa)\to(\mathbb T',\kappa')$, preservation of degrees and colors gives
\[
 n(\mathbb T)=n(\mathbb T'),\qquad
 \ell_e(\mathbb T,\kappa)=\ell_e(\mathbb T',\kappa')
 \quad(e\in[0,2]_{\Z}).
\]
Hence, $w(\mathbb T,\kappa)=w(\mathbb T',\kappa')$.
For each $n$, the coefficient of $t^n$ in $T_d$, $D_d$, and $U$ is therefore a finite sum of well-defined monomials in $\Z[y_0,y_1,y_2]$.

\prf{Claim~\ref{c:branch-series-recursions}} Let us specify the branches obtained by deleting the root.
Let $\mathfrak b=(\mathbb B,\kappa,r)$ have parent color $d$ and let $\mathbb B=(B,(I_B,J_B))$.
Let
\[
 A=\{\kappa(g):g\in\operatorname{IncEdge}_B(r)\}
 \subseteq\{\epsilon,f\}.
\]
For each $e\in A$, proper coloring gives a unique vertex $u_e$ satisfying
\[
 \llbracket r,u_e\rrbracket\in\operatorname{Edge}(B),\qquad
 \kappa(\llbracket r,u_e\rrbracket)=e.
\]
Let $B_e$ be the connected component of
$\operatorname{Subgraph}(B,\operatorname{Vert}(B)\setminus\{r\})$ containing $u_e$.
Define
\[
\begin{aligned}
 I_e&=I_B\cap\operatorname{Vert}(B_e),&
 J_e&=J_B\cap\operatorname{Vert}(B_e),\\
 \kappa_e&=\kappa|_{\operatorname{Edge}(B_e)},&
 \mathfrak b_e&=((B_e,(I_e,J_e)),\kappa_e,u_e).
\end{aligned}
\]
Since $B$ is a tree, its components after deleting $r$ are exactly the graphs $B_e$ for $e\in A$.
Each $\mathfrak b_e$ has parent color $e$.

First, let us treat a sink root. Let $r\in J_B$.
Then, $u_e\in I_e$ and
\[
 \deg_{B_e}(v)=
 \begin{cases}
 2,&v=u_e,\\
 3,&v\in I_e\setminus\{u_e\}.
 \end{cases}
\]
Thus, $\mathfrak b_e$ is a source branch.
The map
\[
 [\mathfrak b]\longmapsto(A,([\mathfrak b_e])_{e\in A})
\]
gives a bijection
\[
 \mathcal B_d^{\mathrm{snk}}
 \longrightarrow
 \bigsqcup_{A\subseteq\{\epsilon,f\}}
 \prod_{e\in A}\mathcal B_e^{\mathrm{src}}.
\]
To construct its inverse, choose disjoint representatives
$\mathfrak b_e=((B_e,(I_e,J_e)),\kappa_e,u_e)$ for $e\in A$ and a new sink $r$.
Let
\[
\begin{aligned}
 I_B&=\bigsqcup_{e\in A}I_e,&
 J_B&=\{r\}\sqcup\bigsqcup_{e\in A}J_e,\\
 \operatorname{Edge}(B)&=
 \bigsqcup_{e\in A}\operatorname{Edge}(B_e)
 \sqcup\{\llbracket r,u_e\rrbracket:e\in A\}.&&
\end{aligned}
\]
Let $\kappa|_{\operatorname{Edge}(B_e)}=\kappa_e$ and
$\kappa(\llbracket r,u_e\rrbracket)=e$ for $e\in A$.
These constructions are inverse on isomorphism classes.
For $A=\varnothing$, the branch consists of the sink $r$ and has weight $y_d$.
For $A\ne\varnothing$, the root $r$ is not a sink leaf in $B^+$. Thus,
\[
 n(\mathfrak b)=\sum_{e\in A}n(\mathfrak b_e),\qquad
 \ell_h(\mathfrak b)=\sum_{e\in A}\ell_h(\mathfrak b_e)
 \quad(h\in[0,2]_{\Z}).
\]
Hence,
\[
\begin{aligned}
 D_d
 &=y_d+\sum_{\varnothing\ne A\subseteq\{\epsilon,f\}}
   \prod_{e\in A}T_e\\
 &=y_d+T_\epsilon+T_f+T_\epsilon T_f.
\end{aligned}
\]

Second, let us treat a source root. Let $r\in I_B$.
Then, $A=\{\epsilon,f\}$, $u_e\in J_e$, and every vertex of $I_e$ has degree three in $B_e$.
Thus, $\mathfrak b_\epsilon$ and $\mathfrak b_f$ are sink branches.
The map
\[
 [\mathfrak b]\longmapsto([\mathfrak b_\epsilon],[\mathfrak b_f])
\]
gives a bijection
\[
 \mathcal B_d^{\mathrm{src}}
 \longrightarrow\mathcal B_\epsilon^{\mathrm{snk}}\times\mathcal B_f^{\mathrm{snk}}.
\]
Its inverse adjoins a new source $r$ to disjoint representatives of the two sink branches and adds the edges
$\llbracket r,u_\epsilon\rrbracket$ and $\llbracket r,u_f\rrbracket$ of colors $\epsilon$ and $f$.
The vertex parts are
\[
 I_B=\{r\}\sqcup I_\epsilon\sqcup I_f,\qquad
 J_B=J_\epsilon\sqcup J_f.
\]
We have
\[
\begin{aligned}
 n(\mathfrak b)&=1+n(\mathfrak b_\epsilon)+n(\mathfrak b_f),\\
 \ell_h(\mathfrak b)&=\ell_h(\mathfrak b_\epsilon)+\ell_h(\mathfrak b_f)
 \quad(h\in[0,2]_{\Z}).
\end{aligned}
\]
Therefore,
\[
 T_d=\sum_{[\mathfrak b_\epsilon]\in\mathcal B_\epsilon^{\mathrm{snk}}}
     \sum_{[\mathfrak b_f]\in\mathcal B_f^{\mathrm{snk}}}
     t\,w(\mathfrak b_\epsilon)w(\mathfrak b_f)
 =tD_\epsilon D_f.
\]

\prf{Claim~\ref{c:branch-series-leaf}} Let us specify the deletion of a marked sink leaf.
Let $\mathbb T=(T,(I_T,J_T))$ and let $v\in J_T$ be a sink leaf.
Let $u\in I_T$ be its unique neighbor and assume
\[
 \kappa(\llbracket u,v\rrbracket)=d.
\]
Let
\[
\begin{aligned}
 B&=\operatorname{Subgraph}(T,\operatorname{Vert}(T)\setminus\{v\}),\\
 \mathbb B&=(B,(I_T,J_T\setminus\{v\})),\qquad
 \kappa_B=\kappa|_{\operatorname{Edge}(B)}.
\end{aligned}
\]
Then, $\mathfrak b=(\mathbb B,\kappa_B,u)$ is a source branch with parent color $d$, since
\[
 \deg_B(u)=2,\qquad
 \deg_B(x)=3\quad(x\in I_T\setminus\{u\}),\qquad
 \kappa_B(g)\ne d\quad(g\in\operatorname{IncEdge}_B(u)).
\]
Thus, deletion defines the map on marked isomorphism classes
\[
 [((\mathbb T,\kappa);v)]\longmapsto[(\mathbb B,\kappa_B,u)].
\]

Let us construct the inverse map.
Let $(\mathbb B,\kappa_B,u)$ be a source branch with $\mathbb B=(B,(I_B,J_B))$.
Adjoin a new sink $v$ and let
\[
\begin{aligned}
 I_T&=I_B,&J_T&=J_B\sqcup\{v\},\\
 \operatorname{Edge}(T)&=\operatorname{Edge}(B)\sqcup\{\llbracket u,v\rrbracket\},&&\\
 \kappa|_{\operatorname{Edge}(B)}&=\kappa_B,&
 \kappa(\llbracket u,v\rrbracket)&=d.
\end{aligned}
\]
The marked leaf is $v$.
Both operations preserve isomorphisms and their compositions return the original marked tree or branch up to isomorphism.

Let us compute the weight after deleting the marked leaf.
Let $n=n(\mathbb T)$ and $\ell_e=\ell_e(\mathbb T,\kappa)$ for $e\in[0,2]_{\Z}$.
Deleting $v$ gives
\[
 n(\mathfrak b)=n,\qquad
 \ell_e(\mathfrak b)=
 \begin{cases}
 \ell_d-1,&e=d,\\
 \ell_e,&e\ne d.
 \end{cases}
\]
By Lemma~\ref{lem:colored-tree-rigidity}, the $\ell_d$ choices of $v$ give $\ell_d$ distinct marked isomorphism classes.
For $\ell_d\ge1$, their total branch weight is
\[
 \ell_d t^n y_d^{\ell_d-1}
 \prod_{e\in[0,2]_{\Z}\setminus\{d\}}y_e^{\ell_e}
 =\frac{\partial}{\partial y_d}
   \left(t^n\prod_{e=0}^2y_e^{\ell_e}\right).
\]
For $\ell_d=0$, both contributions are zero.
Claim~\ref{c:slots-finite} makes each coefficient a finite sum. Thus,
\[
 \frac{\partial U}{\partial y_d}
 =\sum_{[\mathfrak b]\in\mathcal B_d^{\mathrm{src}}}w(\mathfrak b)
 =T_d.\qedhere
\]
\end{proof}

\Needspace{7\baselineskip}
\begin{lemma}\label{slots}
Let $j\in\Z_{>0}$.
Then, the following hold.

\begin{subequations}\label{counts}
\begin{enumerate}
\item\label{c:slots-single} We have
\begin{equation}
P_j=\frac1{3j+4}\binom{4j+4}{j+1}.
\end{equation}
\item\label{c:slots-equal} We have
\begin{equation}
Q_j^{=}=2\binom{4j+2}{j-1}.
\end{equation}
\item\label{c:slots-distinct} We have
\begin{equation}
Q_j^{\ne}=\frac{2j+3}{2j}Q_j^{=}\ge Q_j^{=}.
\end{equation}
\end{enumerate}
\end{subequations}
\end{lemma}
\begin{proof}
\prf{Claim~\ref{c:slots-single}} Let us recover a colored tree from its quiver.
Let $(\mathbb T,\kappa)$ be a colored source-trivalent tree.
Let $\mathbb T=(T,(I_T,J_T))$ and assume $\#I_T=j+1$.
Put $Q=Q(\mathbb T,\kappa)$.
For each $(u,e)\in\operatorname{ColVacSlot}(Q)$, choose a new vertex $v_{u,e}$.
Choose these vertices distinct and outside $\operatorname{Vert}(Q)$.
Up to a color-preserving isomorphism, the tree is recovered by
\[
 \begin{aligned}
 I_T&=J_Q,\\
 J_T&=I_Q\sqcup\{v_{u,e}:(u,e)\in\operatorname{ColVacSlot}(Q)\},\\
 \operatorname{Edge}(T)
 &=\{\llbracket t_Q(\alpha),s_Q(\alpha)\rrbracket:
       \alpha\in\operatorname{Edge}(Q)\}\\
 &\quad\sqcup\{\llbracket u,v_{u,e}\rrbracket:
       (u,e)\in\operatorname{ColVacSlot}(Q)\},\\
 \kappa(\llbracket t_Q(\alpha),s_Q(\alpha)\rrbracket)
 &=\operatorname{QCol}(Q)(\alpha),\qquad
 \kappa(\llbracket u,v_{u,e}\rrbracket)=e.
 \end{aligned}
\]
Lemma~\ref{lem:vacant-slots} identifies the marked slot $(u,e)$ with the marked leaf $v_{u,e}$.

Let us construct the induced isomorphism of zero-sum representations.
Let $(\mathbb T',\kappa')$ be another colored source-trivalent tree.
Let $\mathbb T'=(T',(I_{T'},J_{T'}))$ and assume $\#I_{T'}=j+1$.
Put $Q'=Q(\mathbb T',\kappa')$ and let
$\phi=(\phi_0,\phi_1):Q\to Q'$ be a color-preserving isomorphism preserving the bipartitions.
For $v\in I_Q$, define
\[
 F_v:\operatorname{ZS}(Q)_v\longrightarrow
        \operatorname{ZS}(Q')_{\phi_0(v)},\qquad
 (F_vz)_{\phi_1(\alpha)}=z_\alpha
 \quad(\alpha\in\operatorname{IncEdgeQuiv}_Q(v)).
\]
The map $\phi_1$ is a bijection on the incident arrow sets. Hence,
\[
 \sum_{\alpha'\in\operatorname{IncEdgeQuiv}_{Q'}(\phi_0(v))}
 (F_vz)_{\alpha'}
 =\sum_{\alpha\in\operatorname{IncEdgeQuiv}_Q(v)}z_\alpha=0.
\]
For $u\in J_Q$, set $F_u=1_{\C}$.
For every arrow $\alpha:v\to u$, the identity $(F_vz)_{\phi_1(\alpha)}=z_\alpha$ gives
\[
 \operatorname{ZS}(Q')_{\phi_1(\alpha)}\circ F_v
 =F_u\circ\operatorname{ZS}(Q)_\alpha.
\]
The inverse maps are obtained from $\phi^{-1}$. Thus,
$(\phi_0,\phi_1,(F_v)_{v\in\operatorname{Vert}(Q)})$ is an isomorphism of the zero-sum representations.

Let us identify the marked isomorphism classes.
The extension of $\phi_0$ to the deleted leaves is
\[
 v_{u,e}\longmapsto v'_{\phi_0(u),e}.
\]
The map $v_{u,e}\mapsto v'_{\phi_0(u),e}$ preserves the marked leaves exactly when $\phi_0$ preserves the marked slots.
Conversely, an isomorphism of marked representations includes an isomorphism of the marked colored quivers by Definition~\ref{def:colored-representation-markings}(\ref{d:marked-rep-iso}).
The reconstruction therefore gives a bijection between the marked tree classes and the marked representation classes.

Let us express $P_j$ as a coefficient.
Let $d\in[0,2]_{\Z}$ and let $T=T_0|_{y_0=y_1=y_2=1}$.
Permuting the colors gives $T_d|_{y_0=y_1=y_2=1}=T$ for every $d\in[0,2]_{\Z}$.
Claims~\ref{c:branch-series-recursions} and~\ref{c:branch-series-leaf} of Lemma~\ref{lem:branch-series} give
\[
 T=t(1+T)^4,\qquad T(0)=0,\qquad P_j=[t^{j+1}]T.
\]

Let us extract the coefficient of $t^\nu$ for an integer $\nu\ge1$. Let $u$ be a formal indeterminate.
Lagrange inversion gives
\[
 [t^\nu]T
 =\frac1\nu[u^{\nu-1}](1+u)^{4\nu}
 =\frac1\nu\binom{4\nu}{\nu-1}
 =\frac1{3\nu+1}\binom{4\nu}{\nu}.
\]
Taking $\nu=j+1$ proves Claim~\ref{c:slots-single}.

\prf{Claim~\ref{c:slots-equal}} Let us compute $Q_j^{=}$.
Fix $d\in[0,2]_{\Z}$ and $\epsilon\in[0,2]_{\Z}\setminus\{d\}$.
By Claim~\ref{c:branch-series-leaf} of Lemma~\ref{lem:branch-series}, differentiating $T_d$ with respect to $y_d$ counts an ordered pair of distinct sink leaves of color $d$.
Lemma~\ref{lem:vacant-slots} identifies these pairs with the marked colored-vacant slots counted by $Q_j^{=}$.

Let us compute the derivatives in~\eqref{eq:marked-branch-series}.
Write
\[
 A=\left.\frac{\partial T_d}{\partial y_d}\right|_{y_0=y_1=y_2=1},\qquad
 B=\left.\frac{\partial T_\epsilon}{\partial y_d}\right|_{y_0=y_1=y_2=1}
 \quad(\epsilon\ne d).
\]
The two possible colors $\epsilon\ne d$ give the same $B$ by color symmetry.
Differentiating both recursions gives
\[
 A=\frac{2T}{1+T}(A+B),\qquad
 B=\frac{T}{(1+T)^2}+\frac{T}{1+T}(A+3B).
\]
The equation $A=\frac{2T}{1+T}(A+B)$ gives $(1-T)A=2TB$.
Substituting $(1-T)A=2TB$ into $B=\frac{T}{(1+T)^2}+\frac{T}{1+T}(A+3B)$ gives
\[
 B=\frac{T(1-T)}{(1+T)(1-3T)},\qquad
 A=\frac{2T^2}{(1+T)(1-3T)}.
\]
All denominators have constant term one and are invertible as formal power series.

Let us extract the coefficient that counts the two markings.
Differentiating $T=t(1+T)^4$ with respect to $t$ gives
\[
 tT'=\frac{T(1+T)}{1-3T},\qquad A=tT'\frac{2T}{(1+T)^2}.
\]
The derivative form of Lagrange inversion gives, for $\nu\ge2$,
\[
 [t^\nu]A=[u^{\nu-1}]\frac{2u}{(1+u)^2}(1+u)^{4\nu}
 =2\binom{4\nu-2}{\nu-2}.
\]
Taking $\nu=j+1$ gives $Q_j^{=}=2\binom{4j+2}{j-1}$.

\prf{Claim~\ref{c:slots-distinct}} Let us compute $Q_j^{\ne}$.
Fix $d\in[0,2]_{\Z}$ and $\epsilon\in[0,2]_{\Z}\setminus\{d\}$.
By Claim~\ref{c:branch-series-leaf} of Lemma~\ref{lem:branch-series}, differentiating $T_d$ with respect to $y_\epsilon$ counts an ordered pair of sink leaves of colors $d$ and $\epsilon$.
Lemma~\ref{lem:vacant-slots} identifies these pairs with the marked colored-vacant slots counted by $Q_j^{\ne}$.
The derivative calculation in Claim~\ref{c:slots-equal} and color symmetry give
\[
 \left.\frac{\partial T_d}{\partial y_\epsilon}\right|_{y_0=y_1=y_2=1}
 =\frac{T(1-T)}{(1+T)(1-3T)}
 =tT'\frac{1-T}{(1+T)^2}.
\]

Let us extract the coefficient for an integer $\nu\ge2$.
Lagrange inversion gives
\[
 [t^\nu]\left(tT'\frac{1-T}{(1+T)^2}\right)
 =[u^{\nu-1}](1-u)(1+u)^{4\nu-2}
 =\binom{4\nu-2}{\nu-1}-\binom{4\nu-2}{\nu-2}.
\]
Taking $\nu=j+1$ and using Claim~\ref{c:slots-equal}, we obtain
\[
 \frac{Q_j^{\ne}}{Q_j^{=}}
 =\frac12\left(\frac{\binom{4j+2}{j}}{\binom{4j+2}{j-1}}-1\right)
 =\frac12\left(\frac{3j+3}{j}-1\right)=\frac{2j+3}{2j}.
\]
\end{proof}

\Needspace{6\baselineskip}
\begin{definition}[Admissible pairs]\label{def:admissible-pairs}
Let $(p_j)_{j\in\Z_{>0}}$ and $(q_j)_{j\in\Z_{>0}}$ be real sequences and let $N\in\Z_{>0}$.
We call a pair $(A,B)\in\R_{>0}^2$ \textit{$(p_j,q_j,N)$-admissible} if
\[
\begin{aligned}
 p_j&\ge A\beta^j j^{-\frac32},\\
 q_j&\ge B\beta^j j^{-\frac12}
\end{aligned}
\qquad(j\in\Z_{\ge N}).
\]
\end{definition}

\begin{lemma}\label{lem:slot-bounds}
There exists a $(P_j,Q_j^{=},256)$-admissible pair $(A,B)$.
\end{lemma}
\begin{proof}
Let us obtain positive lower bounds for the normalized counts.
Claims~\ref{c:slots-single} and~\ref{c:slots-equal} of Lemma~\ref{slots} and Stirling's formula give
\[
 P_j\beta^{-j}j^{\frac{3}{2}}\longrightarrow
 \frac\beta3\sqrt{\frac{2}{3\pi}},\qquad
 Q_j^{=}\beta^{-j}j^{\frac{1}{2}}\longrightarrow
 \frac{32}{27}\sqrt{\frac{2}{3\pi}}.
\]

Let us include the finitely many smaller indices in these bounds.
Put
\[
 L_P=\frac\beta3\sqrt{\frac{2}{3\pi}},\qquad
 L_Q=\frac{32}{27}\sqrt{\frac{2}{3\pi}}.
\]
Choose $J>256$ such that, for every integer $j\ge J$,
\[
 P_j\beta^{-j}j^{\frac32}\ge\frac{L_P}{2},\qquad
 Q_j^{=}\beta^{-j}j^{\frac12}\ge\frac{L_Q}{2}.
\]
The explicit formulas in Lemma~\ref{slots} show that both counts are positive for $j\ge256$.
Hence, the finite minima
\[
\begin{aligned}
 A&=\min\left(\left\{\frac{L_P}{2}\right\}
       \cup\{P_j\beta^{-j}j^{\frac32}:256\le j<J\}\right)>0,\\
 B&=\min\left(\left\{\frac{L_Q}{2}\right\}
       \cup\{Q_j^{=}\beta^{-j}j^{\frac12}:256\le j<J\}\right)>0
\end{aligned}
\]
give a $(P_j,Q_j^{=},256)$-admissible pair.
\end{proof}

\subsection{A family of positive flows}
\begin{definition}[Labeled source-trivalent chains]\label{def:labeled-chain}
Let $c\in\Z_{>0}$. We define the \textit{labeled source-trivalent chain} $\operatorname{LSTrivCh}(c)$ as follows.
The vertices $x_t$ for $t\in[c]_{\Z}$, $z_t$ for $t\in[0,c]_{\Z}$, and $y_t$ for $t\in[c]_{\Z}$ are pairwise distinct.
\begin{enumerate}[font=\itshape]
\item Its bipartition is
\[
\begin{aligned}
 I_{\operatorname{LSTrivCh}(c)}&=\{x_t:t\in[c]_{\Z}\},\\
 J_{\operatorname{LSTrivCh}(c)}&=
 \{z_t:t\in[0,c]_{\Z}\}\sqcup\{y_t:t\in[c]_{\Z}\}.
\end{aligned}
\]
\item For $t\in[c]_{\Z}$, define
\[
\begin{aligned}
 \operatorname{LeftEdge}_c(t)&=\llbracket x_t,z_{t-1}\rrbracket,\\
 \operatorname{MidEdge}_c(t)&=\llbracket x_t,y_t\rrbracket,\\
 \operatorname{RightEdge}_c(t)&=\llbracket x_t,z_t\rrbracket.
\end{aligned}
\]
The edge set of $\operatorname{LSTrivCh}(c)$ consists of $\operatorname{LeftEdge}_c(t)$, $\operatorname{MidEdge}_c(t)$, and $\operatorname{RightEdge}_c(t)$ for $t\in[c]_{\Z}$.
We call them the left, middle, and right edges at $x_t$, respectively.
\end{enumerate}
\end{definition}

\begin{lemma}\label{lem:labeled-chain-tree}
For $c\in\Z_{>0}$, let $S=\operatorname{LSTrivCh}(c)$ and $\mathbb S=(S,(I_S,J_S))$.
The bipartite graph $\mathbb S$ is a source-trivalent tree.
Its sink sets satisfy
\[
 \begin{aligned}
 J_S\cap\operatorname{Leaf}(S)&=\{z_0,z_c\}\sqcup\{y_t:t\in[c]_{\Z}\},\\
 \{v\in J_S:\deg_S(v)=2\}&=\{z_t:t\in[c-1]_{\Z}\}.
 \end{aligned}
\]
\end{lemma}
\begin{proof}
First, suppose $c=1$. The graph is a tree with center $x_1$ and three leaves $z_0,y_1,z_1$.
Its only source has degree three, and it has no degree-two sinks.

Second, suppose $c\ge2$. Let us delete $z_0,z_c,y_1,\ldots,y_c$.
The remaining vertices and edges form the path
\[
 x_1,z_1,x_2,\ldots,z_{c-1},x_c.
\]
Each deleted vertex is joined to exactly one source by one edge.
Adding these vertices to the path preserves connectedness and creates no cycle.
Thus, $S$ is a tree. Its degrees are
\[
\begin{aligned}
 \deg_S(x_t)&=3 &&(1\le t\le c),\\
 \deg_S(z_t)&=2 &&(1\le t<c),\\
 \deg_S(z_0)=\deg_S(z_c)=\deg_S(y_t)&=1 &&(1\le t\le c).
\end{aligned}
\]
These formulas give the stated source degrees, sink leaves, and degree-two sinks.
\end{proof}

\begin{definition}[Chain flow family and source dimensions]\label{def:chain-flow-family}
Let $a\in\Z_{>0}$ and $k\in\Z_{>0}$ satisfy $k\ge2$.
Let $c=k-1$, $x=\frac ak$, and $S=\operatorname{LSTrivCh}(c)$.
For $f\in\operatorname{Flow}(S)$ and $t\in[c]_{\Z}$, define
\[
 l_t(f)=f_{\operatorname{LeftEdge}_c(t)},\qquad
 m_t(f)=f_{\operatorname{MidEdge}_c(t)},\qquad
 r_t(f)=f_{\operatorname{RightEdge}_c(t)}.
\]
Define $\mathcal F_{a,k}$ to be the set of $f\in\operatorname{Flow}(S)$ such that, for every $t\in[c]_{\Z}$,
\[
\begin{aligned}
 l_t(f)&\in[\frac{x}{256},\frac{2x}{256})\cap(tx+\Z),\\
 r_t(f)&\in(0,\frac{x}{256})\cap(-tx+\Z),\\
 m_t(f)&=x+1-l_t(f)-r_t(f).
\end{aligned}
\]
For $f\in\operatorname{Flow}(S)$ and $b\in J_S$, define the \textit{source dimension}
\[
 a_b(f)=x-\sum_{w\in\operatorname{IncEdge}_S(b)}f_w.
\]
\end{definition}

\begin{lemma}[Flow choices]\label{lem:flow-choices}
Let $a\in\Z_{>0}$ and $k\in\Z_{>0}$ satisfy $k\ge2$ and $a\ge100000k$.
Let $c=k-1$, $x=\frac ak$, and $S=\operatorname{LSTrivCh}(c)$.
Let
\[
 N=256,\qquad \gamma=1-\frac Nx,\qquad \lambda=\frac1{100}.
\]
Let $f\in\mathcal F_{a,k}$.
Let $l_t=l_t(f)$, $m_t=m_t(f)$, and $r_t=r_t(f)$ for $t\in[c]_{\Z}$.
Then, we have the following.

\begin{enumerate}
\item\label{c:flow-choices-dimensions} The flow $f$ and the source dimensions $(a_b(f))_{b\in J_{S}}$ satisfy the following.
\begin{enumerate}[label=(\alph*),ref=\theenumi(\alph*)]
\item\label{c:flow-choices-positive} We have $l_t>0$, $m_t>0$, and $r_t>0$ for $t\in[c]_{\Z}$.
\item\label{c:flow-choices-integral} For $b\in J_S$, we have $a_b(f)\in\Z$.
\item\label{c:flow-choices-sum} We have $\sum_{b\in J_S}a_b(f)=a-c$.
\item\label{c:flow-choices-range} For every $b\in J_S$, we have
\[
 256<a_b(f)<x.
\]
\end{enumerate}
\item\label{c:flow-choices-count} For $t\in[c]_{\Z}$, we have
\[
\begin{aligned}
 \#\left([\frac{x}{N},\frac{2x}{N})\cap(tx+\mathbb Z)\right)&\ge\frac{\gamma x}{N},\\
 \#\left((0,\frac{x}{N})\cap(-tx+\mathbb Z)\right)&\ge\frac{\gamma x}{N}.
\end{aligned}
\]
\item\label{c:flow-choices-distinct} For $f,f'\in\mathcal F_{a,k}$, we have
\[
 \bigl(a_b(f)=a_b(f')\text{ for every }b\in J_S\bigr)
 \quad\Longrightarrow\quad f=f'.
\]
\item\label{c:flow-choices-stable} Let $\mathbf T=((\mathbb T_b,\kappa_b))_{b\in J_S}$ be a family of colored source-trivalent trees satisfying
\[
 \operatorname{TotDim}_{\C,(I_{Q(\mathbb T_b,\kappa_b)},J_{Q(\mathbb T_b,\kappa_b)})}
 (\operatorname{ZS}(Q(\mathbb T_b,\kappa_b)))=(a_b(f),a_b(f)+1)
 \qquad(b\in J_S).
\]
Let $\mathbf p=((v_e,d_e))_{e\in\operatorname{Edge}(S)}$ be a family of colored-vacant slots and $\mathbf g=(g_e)_{e\in\operatorname{Edge}(S)}$ a family of nonzero complex linear maps such that
\[
 (v_e,d_e)\in\operatorname{ColVacSlot}(Q(\mathbb T_b,\kappa_b)),\qquad
 g_e:\C\longrightarrow\operatorname{ZS}(Q(\mathbb T_b,\kappa_b))_{v_e}
\]
for $e=\llbracket x,b\rrbracket\in\operatorname{Edge}(S)$, $x\in I_S$, and $b\in J_S$.
Let
\[
 R=\operatorname{Glue}_{\C}((S,(I_S,J_S));\mathbf T,\mathbf p,\mathbf g)
\]
satisfy hypotheses~\ref{hyp:gluing-attachments} and~\ref{hyp:gluing-slots} of Lemma~\ref{lem:flow-stability}.
Then, $R$ is stable.
\item\label{c:flow-choices-weight} For $t\in[c-1]_{\Z}$, we have
\begin{equation}\label{weight}
 (l_t+r_t-1)^{-\frac32}(x-r_t-l_{t+1})^{-\frac12}
 \ge\lambda N^2x^{-2}.
\end{equation}
The dimensions $a_{z_0}(f)$, $a_{z_c}(f)$, and $a_{y_c}(f)$ satisfy
\[
 (x-l_1)^{-\frac32}(x-r_c)^{-\frac32}(l_c+r_c-1)^{-\frac32}
 \ge x^{-\frac92}.
\]
\end{enumerate}
\end{lemma}
\begin{proof}
\prf{Claim~\ref{c:flow-choices-positive}} Let us prove positivity of all edge weights.
Since $f\in\mathcal F_{a,k}$, $l_t>0$, $r_t>0$, and $l_t+r_t<\frac{3x}{N}$. Hence,
\[
 m_t=x+1-l_t-r_t>x\left(1-\frac3N\right)+1>0.
\]

\prf{Claim~\ref{c:flow-choices-integral}} Let us prove $a_b(f)\in\Z$ for $b\in J_S$.
For $S=\operatorname{LSTrivCh}(c)$, the dimensions $a_b(f)$ satisfy
\[
\begin{aligned}
 a_{z_0}(f)&=x-l_1,& a_{z_c}(f)&=x-r_c,\\
 a_{y_t}(f)&=l_t+r_t-1&&(t\in[c]_{\Z}),\\
 a_{z_t}(f)&=x-r_t-l_{t+1}&&(t\in[c-1]_{\Z}).
\end{aligned}
\]
 The conditions $l_t-tx\in\Z$, $r_t+tx\in\Z$, and
$kx=a\in\Z$ give
\[
\begin{aligned}
 x-l_1&=-(l_1-x)\in\Z,\\
 x-r_c&=a-(r_c+cx)\in\Z,\\
 l_t+r_t-1&=(l_t-tx)+(r_t+tx)-1\in\Z,\\
 x-r_t-l_{t+1}&=-(r_t+tx)-(l_{t+1}-(t+1)x)\in\Z.
\end{aligned}
\]

\prf{Claim~\ref{c:flow-choices-sum}}
The dimensions $a_b(f)$ in Claim~\ref{c:flow-choices-integral} satisfy
\[
\begin{aligned}
 &(x-l_1)+(x-r_c)+\sum_{t=1}^c(l_t+r_t-1)
       +\sum_{t=1}^{c-1}(x-r_t-l_{t+1})\\
 &\qquad=(c+1)x-c=a-c.
\end{aligned}
\]

\prf{Claim~\ref{c:flow-choices-range}} Let us bound each source dimension between $256$ and $x$.
The bound $a\ge100000k$ gives $x\ge100000$.
First, for $t\in[c]_{\Z}$, the bounds $l_t\in[\frac{x}{N},\frac{2x}{N})$ and $r_t\in(0,\frac{x}{N})$ give
\[
 256<\frac{x}{N}-1<a_{y_t}(f)=l_t+r_t-1<\frac{3x}{N}-1<x.
\]
Second, for $t\in[0,c]_{\Z}$, the formulas in Claim~\ref{c:flow-choices-integral} give
\[
 256<x\left(1-\frac3N\right)<a_{z_t}(f)<x.
\]
For $t=0$, use $a_{z_0}(f)=x-l_1$ and $0<l_1<\frac{2x}{N}$.
For $t=c$, use $a_{z_c}(f)=x-r_c$ and $0<r_c<\frac{x}{N}$.
For $t\in[c-1]_{\Z}$, use $a_{z_t}(f)=x-r_t-l_{t+1}$ and $0<r_t+l_{t+1}<\frac{3x}{N}$.

\prf{Claim~\ref{c:flow-choices-count}} Let us count the integers after translating each interval.
For $t\in[c]_{\Z}$, we have
\[
 \begin{aligned}
 \#\left([\frac{x}{N},\frac{2x}{N})\cap(tx+\Z)\right)
 &=\left\lceil\frac{2x}{N}-tx\right\rceil
   -\left\lceil\frac{x}{N}-tx\right\rceil
 \ge\frac{x}{N}-1,\\
 \#\left((0,\frac{x}{N})\cap(-tx+\Z)\right)
 &=\left\lceil\frac{x}{N}+tx\right\rceil-\lfloor tx\rfloor-1
 \ge\frac{x}{N}-1.
 \end{aligned}
\]
Since $\frac{x}{N}-1=\frac{\gamma x}{N}$, both inequalities in Claim~\ref{c:flow-choices-count} follow.

\prf{Claim~\ref{c:flow-choices-distinct}} Let us recover the flow from its source dimensions.
Let $f\in\mathcal F_{a,k}$ and $\widetilde f\in\mathcal F_{a,k}$ satisfy $a_b(f)=a_b(\widetilde f)$ for every $b\in J_S$.
For $h\in\{f,\widetilde f\}$, we have
\[
 \sum_{w\in\operatorname{IncEdge}_S(v)}h_w=x+1\quad(v\in I_S),\qquad
 \sum_{w\in\operatorname{IncEdge}_S(b)}h_w=x-a_b(f)\quad(b\in J_S).
\]
Claim~\ref{c:tree-flow-existence} of Lemma~\ref{lem:tree-flow-cut} gives $f=\widetilde f$.

\prf{Claim~\ref{c:flow-choices-stable}} Let us prove stability.
Claim~\ref{c:flow-choices-sum} gives $\sum_{b\in J_S}a_b(f)+c=a$.
Definition~\ref{def:chain-flow-family} gives equations~\eqref{flow} for $f$ and the source dimensions $a_b(f)$.
Claim~\ref{c:flow-stability-unique} of Lemma~\ref{lem:flow-stability} therefore identifies $f$ as the unique flow used in Claim~\ref{c:flow-stability-positive} of Lemma~\ref{lem:flow-stability}.
Claim~\ref{c:flow-choices-positive} gives $f_w>0$ for every $w\in\operatorname{Edge}(S)$.
Claim~\ref{c:flow-stability-positive} of Lemma~\ref{lem:flow-stability}
gives stability of $R$.

\prf{Claim~\ref{c:flow-choices-weight}} Let us bound the products of inverse powers of the source dimensions.
For $t\in[c-1]_{\Z}$, we have $a_{y_t}(f)<\frac{3x}{N}$ and $a_{z_t}(f)<\frac{(N-1)x}{N}$.
Hence,
\[
(l_t+r_t-1)^{-\frac{3}{2}}
(x-r_t-l_{t+1})^{-\frac{1}{2}}
\ge\frac{N^2x^{-2}}{\sqrt{3^3(N-1)}}
\ge\lambda N^2x^{-2}.
\]
The inequality $\frac{N^2x^{-2}}{\sqrt{3^3(N-1)}}\ge\lambda N^2x^{-2}$ uses $\lambda=\frac1{100}$ and $3^3(N-1)=6885<100^2$.
We also have $a_{z_0}(f)^{-\frac32}a_{z_c}(f)^{-\frac32}a_{y_c}(f)^{-\frac32}\ge x^{-\frac92}$.
\end{proof}

\subsection{Covering quivers and dimension vectors}\label{subsec:cover-definitions}
\Needspace{10\baselineskip}
\begin{definition}[Coverings, embeddings, and deck transformations]\label{def:quiver-covering}
Let $Q$ and $\widehat Q$ be quivers. 
\begin{enumerate}[font=\itshape]
\item An \textit{embedding} $\iota=(\iota_0,\iota_1):Q\to\widehat Q$ consists of injective maps
\[
 \iota_0:\operatorname{Vert}(Q)\hookrightarrow\operatorname{Vert}(\widehat Q),\qquad
 \iota_1:\operatorname{Edge}(Q)\hookrightarrow\operatorname{Edge}(\widehat Q)
\]
such that, for every $\alpha\in\operatorname{Edge}(Q)$,
\[
 s_{\widehat Q}(\iota_1(\alpha))=\iota_0(s_Q(\alpha)),\qquad
 t_{\widehat Q}(\iota_1(\alpha))=\iota_0(t_Q(\alpha)).
\]
\item A quiver homomorphism $p=(p_0,p_1):\widehat Q\to Q$ is a \textit{covering}
if $p_0$ is surjective and, for every $v\in\operatorname{Vert}(\widehat Q)$, the restrictions
\[
\begin{aligned}
 p_1:s_{\widehat Q}^{-1}(\{v\})&\longrightarrow s_Q^{-1}(\{p_0(v)\}),\\
 p_1:t_{\widehat Q}^{-1}(\{v\})&\longrightarrow t_Q^{-1}(\{p_0(v)\})
\end{aligned}
\]
are bijective.
We call $p$ a \textit{universal covering} if $\operatorname{UnderG}(\widehat Q)$ is a tree.
\item A \textit{deck transformation} of $p$ is a quiver isomorphism
$\delta:\widehat Q\to\widehat Q$ such that
\[
 p_0\circ\delta_0=p_0,\qquad p_1\circ\delta_1=p_1.
\]
We denote the group of these transformations under composition by $\operatorname{Deck}(p)$.
\end{enumerate}
\end{definition}

\begin{definition}[Arrow traversals and lifts]\label{def:path-lift}
Let $Q$ be a quiver and $n\in\Z_{\ge0}$. A \textit{traversal of length $n$} is a tuple
$\gamma=(v_0,\alpha_1,\varepsilon_1,v_1,\ldots,\alpha_n,\varepsilon_n,v_n)$ such that
$v_i\in\operatorname{Vert}(Q)$ for $i\in[0,n]_{\Z}$ and
\[
 \alpha_i\in\operatorname{Edge}(Q),\qquad \varepsilon_i\in\{1,-1\},\qquad
 (s_Q(\alpha_i),t_Q(\alpha_i))=
 \begin{cases}(v_{i-1},v_i)&\varepsilon_i=1,\\(v_i,v_{i-1})&\varepsilon_i=-1\end{cases}
 \quad(i\in[n]_{\Z}).
\]
For a covering $p:\widehat Q\to Q$, a \textit{lift} of $\gamma$ is a traversal
$\widehat\gamma=(\widehat v_0,\widehat\alpha_1,\varepsilon_1,\widehat v_1,\ldots,\widehat\alpha_n,\varepsilon_n,\widehat v_n)$
satisfying $p_0(\widehat v_i)=v_i$ for $i\in[0,n]_{\Z}$ and $p_1(\widehat\alpha_i)=\alpha_i$ for $i\in[n]_{\Z}$.
\end{definition}

\begin{lemma}[Unique traversal lifting]\label{lem:path-lift}
Let $p:\widehat Q\to Q$ be a covering.
Let $n\in\Z_{\ge0}$ and let
$\gamma=(v_0,\alpha_1,\varepsilon_1,v_1,\ldots,\alpha_n,\varepsilon_n,v_n)$
be a traversal in $Q$.
Let $\widehat v_0\in\operatorname{Vert}(\widehat Q)$ satisfy $p_0(\widehat v_0)=v_0$.
Then, there is exactly one traversal
$\widehat\gamma=(\widehat v_0,\widehat\alpha_1,\varepsilon_1,\widehat v_1,\ldots,\widehat\alpha_n,\varepsilon_n,\widehat v_n)$ in $\widehat Q$ satisfying
\[
 p_0(\widehat v_i)=v_i\quad(i\in[0,n]_{\Z}),\qquad
 p_1(\widehat\alpha_i)=\alpha_i\quad(i\in[n]_{\Z}).
\]
\end{lemma}
\begin{proof}
Let us construct the lift one arrow at a time.
For a traversal of length zero, the lift consists of the specified vertex $\widehat v_0$.
For $i\ge1$, suppose that the lifted traversal through $\widehat v_{i-1}$ has been determined and
$p_0(\widehat v_{i-1})=v_{i-1}$.

First, suppose $\varepsilon_i=1$. The covering gives a bijection
\[
 p_1:s_{\widehat Q}^{-1}(\{\widehat v_{i-1}\})
       \longrightarrow s_Q^{-1}(\{v_{i-1}\}).
\]
Let $\widehat\alpha_i$ be the unique preimage of $\alpha_i$ and set
\[
 \widehat v_i=t_{\widehat Q}(\widehat\alpha_i).
\]
Then,
\[
 p_0(\widehat v_i)=t_Q(p_1(\widehat\alpha_i))=t_Q(\alpha_i)=v_i.
\]

Second, suppose $\varepsilon_i=-1$. Use the bijection
\[
 p_1:t_{\widehat Q}^{-1}(\{\widehat v_{i-1}\})
       \longrightarrow t_Q^{-1}(\{v_{i-1}\}).
\]
Let $\widehat\alpha_i$ be the unique preimage of $\alpha_i$ and set
\[
 \widehat v_i=s_{\widehat Q}(\widehat\alpha_i).
\]
Then,
\[
 p_0(\widehat v_i)=s_Q(p_1(\widehat\alpha_i))=s_Q(\alpha_i)=v_i.
\]

Let us prove uniqueness of the whole lift.
At each step, the appropriate covering bijection determines $\widehat\alpha_i$, and its other endpoint determines $\widehat v_i$.
Induction on $i$ proves existence and uniqueness through the last vertex.
\end{proof}

\begin{definition}[Lift of a traversal]\label{def:unique-path-lift}
Let $p:\widehat Q\to Q$ be a covering and let $\gamma$ be a traversal in $Q$ beginning at $v_0$.
For $\widehat v_0\in p_0^{-1}(\{v_0\})$, let
\[
 \operatorname{Lift}_p(\gamma;\widehat v_0)
\]
denote the unique lift of $\gamma$ beginning at $\widehat v_0$ from Lemma~\ref{lem:path-lift}.
\end{definition}

\begin{definition}[Supported dimension vectors]\label{def:cover-dimension-vectors}
Let $Q$ be a quiver.
\begin{enumerate}[font=\itshape]
\item A \textit{dimension vector} is a family $d=(d_v)_{v\in\operatorname{Vert}(Q)}$
of nonnegative integers. We define
\[
 \operatorname{VecSupp}(d)=\{v\in\operatorname{Vert}(Q):d_v>0\}.
\]
We call $d$ \textit{finitely supported} if $\operatorname{VecSupp}(d)$ is finite and let
\[
 \operatorname{DimVec}_{\mathrm{fin}}(Q)
 =\{d\in\Z_{\ge0}^{\operatorname{Vert}(Q)}:\#\operatorname{VecSupp}(d)<\infty\}.
\]
Its \textit{positive support quiver} is $\operatorname{SubQuiv}(Q,\operatorname{VecSupp}(d))$.
\item Let $p:\widehat Q\to Q$ be a covering and let $\delta\in\operatorname{Deck}(p)$.
For $d\in\operatorname{DimVec}_{\mathrm{fin}}(\widehat Q)$, define
\[
 (\delta\cdot d)_v=d_{\delta_0^{-1}(v)}
 \qquad(v\in\operatorname{Vert}(\widehat Q)).
\]
We call $d$ and $d'$ \textit{deck-equivalent} if $d'=\delta\cdot d$ for some $\delta\in\operatorname{Deck}(p)$.
We denote the resulting orbit by $[d]_p$.
\end{enumerate}
\end{definition}

\begin{definition}[Covering arrow coloring]\label{def:covering-arrow-coloring}
Fix a bijection $c_3:\operatorname{Edge}(K(3))\to[0,2]_{\Z}$ and a universal covering
$\varpi=(\varpi_0,\varpi_1):\widetilde K_3\to K(3)$ for the rest of this section.
\begin{enumerate}[font=\itshape]
\item Let $I_{K(3)}$ consist of the first vertex and $J_{K(3)}$ of the second vertex of $K(3)$.
We use the bipartition
\[
 I_{\widetilde K_3}=\varpi_0^{-1}(I_{K(3)}),\qquad
 J_{\widetilde K_3}=\varpi_0^{-1}(J_{K(3)}).
\]
\item We define the \textit{covering arrow coloring} by
\[
 \operatorname{CovCol}_{\varpi}:\operatorname{Edge}(\widetilde K_3)\longrightarrow[0,2]_{\Z},\qquad
 \operatorname{CovCol}_{\varpi}(\alpha)=c_3(\varpi_1(\alpha)).
\]
\end{enumerate}
\end{definition}

\begin{definition}[Cover dimension classes and stability weights]\label{def:cover-dimension-classes}
Let $a\in\Z_{>0}$ and $b\in\Z_{>0}$. Let
$I=I_{\widetilde K_3}$ and $J=J_{\widetilde K_3}$.
\begin{enumerate}[font=\itshape]
\item We define
\[
 \operatorname{CovDim}_{\varpi}(a,b)=
 \left\{d\in\operatorname{DimVec}_{\mathrm{fin}}(\widetilde K_3):
 \sum_{i\in I}d_i=a,\ \sum_{j\in J}d_j=b\right\}.
\]
\item For $e\in\operatorname{DimVec}_{\mathrm{fin}}(\widetilde K_3)$, define
\[
 \theta^{\varpi}_{a,b}(e)=a\sum_{j\in J}e_j-b\sum_{i\in I}e_i.
\]
\item For $k\in\Z_{>0}$, we define the sets of \textit{cover dimension classes}
and \textit{stable cover dimension classes} by
\[
\begin{aligned}
 \mathcal D_{a,k}&=\{[d]_{\varpi}:d\in\operatorname{CovDim}_{\varpi}(a,a+k)\},\\
 \mathcal D^{\mathrm{st}}_{a,k}&=\{[d]_{\varpi}\in\mathcal D_{a,k}:
 \cM^s_d(\widetilde K_3)\ne\varnothing\}.
\end{aligned}
\]
The moduli spaces here use the weight $\theta^{\varpi}_{a,a+k}$.
\end{enumerate}
\end{definition}

\subsection{Counting cover representations}
\begin{definition}\label{def:marked-chain-data}
Let $a\in\Z_{>0}$ and $k\in\Z_{>0}$ satisfy $k\ge2$ and $a\ge100000k$.
Let $c=k-1$ and $S=\operatorname{LSTrivCh}(c)$.
Let $\mathbb S=(S,(I_S,J_S))$.
Define
\[
 L=J_{S}\cap\operatorname{Leaf}(S),\qquad
 I=\{i\in J_{S}:\deg_S(i)=2\}.
\]
We call a tuple
\[
 M=(S,f,\tau,(B_i)_{i\in L\sqcup I})
\]
\textit{marked chain data} if the following conditions hold.
\begin{hypotheses}
\item\label{hyp:marked-flows} Let $f\in\mathcal F_{a,k}$ and let $a_i=a_i(f)$ for $i\in J_S$.
\item\label{hyp:marked-colors} For each $t\in[c]_{\Z}$, the map
$\tau_t\colon[3]_{\Z}\longrightarrow[0,2]_{\Z}$ is a bijection.
The indices $1$, $2$, and $3$ denote $\operatorname{LeftEdge}_c(t)$,
$\operatorname{MidEdge}_c(t)$, and $\operatorname{RightEdge}_c(t)$, respectively.
\item\label{hyp:marked-blocks} For each $i\in L\sqcup I$, let $(\mathbb T_i,\kappa_i)$ be a colored source-trivalent tree.
Let $\mathbb T_i=(T_i,(I_{T_i},J_{T_i}))$ for $i\in L\sqcup I$.
Require the following conditions.
\begin{enumerate}[label=(\roman*),font=\itshape]
\item For $i\in L$, let $\mathbf s_i=(s_{i,1})$ satisfy
\[
 s_{i,1}\in\operatorname{ColVacSlot}(Q(\mathbb T_i,\kappa_i)).
\]
For $i\in I$, let $\mathbf s_i=(s_{i,1},s_{i,2})$ satisfy
\[
 s_{i,j}\in\operatorname{ColVacSlot}(Q(\mathbb T_i,\kappa_i))\quad(j\in[2]_{\Z}),
 \qquad s_{i,1}\ne s_{i,2}.
\]
\item For $i\in L\sqcup I$, let
\[
 B_i=\operatorname{MarkedCIsom}_{\C}(Q(\mathbb T_i,\kappa_i),\operatorname{ZS}(Q(\mathbb T_i,\kappa_i));\mathbf s_i)
\]
and require
\[
 \operatorname{TotDim}_{\C,(I_{Q(\mathbb T_i,\kappa_i)},J_{Q(\mathbb T_i,\kappa_i)})}
 (\operatorname{ZS}(Q(\mathbb T_i,\kappa_i)))=(a_i,a_i+1).
\]
\item Define the colored-vacant slots $p_e$ for $e\in\operatorname{Edge}(S)$ by
\[
 p_e=s_{i,1}\quad(i\in L,\ e\in\operatorname{IncEdge}_S(i)),
\]
and, for $t\in[c-1]_{\Z}$, let
\[
 p_{\operatorname{RightEdge}_c(t)}=s_{z_t,1},\qquad
 p_{\operatorname{LeftEdge}_c(t+1)}=s_{z_t,2}.
\]
\item For $e\in\operatorname{Edge}(S)$, let $p_e=(v_e,d_e)$.
For every $t\in[c]_{\Z}$, require
\[
 d_{\operatorname{LeftEdge}_c(t)}=\tau_t(1),\qquad
 d_{\operatorname{MidEdge}_c(t)}=\tau_t(2),\qquad
 d_{\operatorname{RightEdge}_c(t)}=\tau_t(3).
\]
\end{enumerate}
\end{hypotheses}
\begin{enumerate}[font=\itshape]
\item\label{d:marked-chain-set} Let $\mathcal H_{a,k}$ be the set of tuples $M$ satisfying assumptions~\ref{hyp:marked-flows}--\ref{hyp:marked-blocks}.
For a family $\tau=(\tau_t)_{t\in[c]_{\Z}}$ of bijections $\tau_t:[3]_{\Z}\to[0,2]_{\Z}$, define
\[
 I_=(\tau)=\{z_t:t\in[c-1]_{\Z},\ \tau_t(3)=\tau_{t+1}(1)\},\qquad
 I_{\ne}(\tau)=I\setminus I_=(\tau).
\]
\item\label{d:marked-chain-support} For each $M\in\mathcal H_{a,k}$, fix a family
\[
 \rho_M=((\mathbb T_i,\kappa_i;\mathbf s_i))_{i\in J_S}
\]
satisfying assumption~\ref{hyp:marked-blocks} and defining the classes $B_i$ in $M$.
Let $\rho=(\rho_M)_{M\in\mathcal H_{a,k}}$.
For $M\in\mathcal H_{a,k}$, let
$\mathbf T_M^\rho=((\mathbb T_i,\kappa_i))_{i\in J_S}$,
$\mathbf p_M^\rho=(p_e)_{e\in\operatorname{Edge}(S)}$, and
$\mathbf 1=(1_{\C})_{e\in\operatorname{Edge}(S)}$.
Define
\[
 R_M^\rho=\operatorname{Glue}_{\C}(\mathbb S;\mathbf T_M^\rho,\mathbf p_M^\rho,\mathbf 1),\qquad
 \widetilde Q_M^\rho=\operatorname{SuppColQuiv}(R_M^\rho),\qquad
 d_M^\rho(v)=\dim_{\C}(R_M^\rho)_v.
\]
Use the tagged vertices $(i,v)$ for the copy of $Q(\mathbb T_i,\kappa_i)$ and distinct new vertices $x_t$ for $t\in[c]_{\Z}$.
Define
\[
\begin{aligned}
 I_M^\rho&=\operatorname{Vert}(\widetilde Q_M^\rho)\cap
 \left(\{x_t:t\in[c]_{\Z}\}\sqcup
 \bigsqcup_{i\in J_S}\{(i,v):v\in I_{Q(\mathbb T_i,\kappa_i)}\}\right),\\
 J_M^\rho&=\operatorname{Vert}(\widetilde Q_M^\rho)\cap
 \bigsqcup_{i\in J_S}\{(i,v):v\in J_{Q(\mathbb T_i,\kappa_i)}\},\\
 \widetilde{\mathbb Q}_M^\rho&=(\widetilde Q_M^\rho,(I_M^\rho,J_M^\rho)).
\end{aligned}
\]
Let $\iota_i^\rho:Q(\mathbb T_i,\kappa_i)\hookrightarrow\widetilde{\mathbb Q}_M^\rho$ denote the inclusion of the copy indexed by $i$.
\item\label{d:dimensioned-support-classes} An isomorphism between colored bipartite quivers $\widetilde{\mathbb G}=(G,(I_G,J_G))$ and $\widetilde{\mathbb G}'=(G',(I_{G'},J_{G'}))$ equipped with positive vertex dimensions $d_G$ and $d_{G'}$ is a color-preserving quiver isomorphism $\phi=(\phi_0,\phi_1)$ that preserves the bipartitions and satisfies
\[
 d_{G'}(\phi_0(v))=d_G(v)\qquad(v\in\operatorname{Vert}(G)).
\]
Define
\[
 \mathcal S_{a,k}^\rho=\{[(\widetilde{\mathbb Q}_M^\rho,d_M^\rho)]:M\in\mathcal H_{a,k}\},\qquad
 \pi_\rho:\mathcal H_{a,k}\longrightarrow\mathcal S_{a,k}^\rho,\quad
 \pi_\rho(M)=[(\widetilde{\mathbb Q}_M^\rho,d_M^\rho)].
\]

\item\label{d:connector-set} Let $X=(\widetilde{\mathbb G},d_G)$ represent a class in $\mathcal S_{a,k}^\rho$.
Let $\widetilde{\mathbb G}=(G,(I_G,J_G))$.
Define
\[
 C(X)=\{v\in I_G:d_G(v)=1,\ \deg_G(v)=3\}.
\]
\item\label{d:deleted-components} Define
\[
 \mathcal K(X)=\operatorname{CtdComp}\bigl(\operatorname{Subgraph}(\operatorname{UnderG}(G),\operatorname{Vert}(G)\setminus C(X))\bigr).
\]
\item\label{d:recovered-skeleton} Define the bipartite graph $\mathbb S_X=(S_X,(I_{S_X},J_{S_X}))$ with vertex set $C(X)\sqcup\mathcal K(X)$, identifying each part with its canonical image, by
\[
 I_{S_X}=C(X),\qquad
 J_{S_X}=\mathcal K(X).
\]
For $v\in C(X)$ and $K\in\mathcal K(X)$, define its edges by
\[
 \llbracket v,K\rrbracket\in\operatorname{Edge}(S_X)
 \ \Longleftrightarrow\
 \exists u\in\operatorname{Vert}(K),\quad
 \llbracket v,u\rrbracket\in\operatorname{Edge}(\operatorname{UnderG}(G)).
\]
\end{enumerate}
\end{definition}

\begin{lemma}[Recovery of a colored tree and its markings]\label{lem:retained-tree-recovery}
Let $(\mathbb T,\kappa)$ and $(\mathbb T',\kappa')$ be colored source-trivalent trees with $\#I_T\ge1$ and $\#I_{T'}\ge1$.
Let $Q=Q(\mathbb T,\kappa)$ and $Q'=Q(\mathbb T',\kappa')$.
Let $q=(q_0,q_1):Q\to Q'$ be a color-preserving isomorphism preserving the bipartitions.
Then, we have the following.
\begin{enumerate}
\item\label{c:retained-tree-extension} There is a unique color-preserving isomorphism
\[
 \Psi=(\Psi_0,\Psi_1):(\mathbb T,\kappa)\longrightarrow(\mathbb T',\kappa')
\]
preserving the bipartitions and satisfying
$\Psi_0(v)=q_0(v)$ for $v\in\operatorname{Vert}(Q)$.
\item\label{c:retained-marked-transport} Let $\ell\in\Z_{\ge0}$ and let
$\mathbf s=((u_j,d_j))_{j\in[\ell]_{\Z}}$ be an ordered marking of colored-vacant slots of $Q$.
Let $\mathbf s'=((q_0(u_j),d_j))_{j\in[\ell]_{\Z}}$.
Then, $\mathbf s'$ is an ordered marking of $Q'$ and
\[
 \operatorname{MarkedCIsom}_{\C}(Q,\operatorname{ZS}(Q);\mathbf s)
 =\operatorname{MarkedCIsom}_{\C}(Q',\operatorname{ZS}(Q');\mathbf s').
\]
\end{enumerate}
\end{lemma}
\begin{proof}
\prf{Claim~\ref{c:retained-tree-extension}} Let us define the images of the deleted leaves.
Let
\[
 L=\{v\in J_T:\deg_T(v)=1\},\qquad
 L'=\{v'\in J_{T'}:\deg_{T'}(v')=1\}.
\]
Connectedness, $\#I_T\ge1$, and $\#I_{T'}\ge1$ imply
\[
 \operatorname{Vert}(T)=\operatorname{Vert}(Q)\sqcup L,\qquad
 \operatorname{Vert}(T')=\operatorname{Vert}(Q')\sqcup L'.
\]
For $v\in L$, let $u\in I_T$ be its unique neighbor and let
$d=\kappa(\llbracket u,v\rrbracket)$.
Lemma~\ref{lem:vacant-slots} gives $(u,d)\in\operatorname{ColVacSlot}(Q)$.
Since $q$ preserves incidence and colors,
\[
 (u,d)\in\operatorname{ColVacSlot}(Q)
 \quad\Longleftrightarrow\quad
 (q_0(u),d)\in\operatorname{ColVacSlot}(Q').
\]
Lemma~\ref{lem:vacant-slots} gives a unique $v'\in L'$ satisfying
\[
 \llbracket q_0(u),v'\rrbracket\in\operatorname{Edge}(T'),\qquad
 \kappa'(\llbracket q_0(u),v'\rrbracket)=d.
\]
Define $\eta(v)=v'$.
Applying the construction of $\eta:L\to L'$ to $q^{-1}:Q'\to Q$ gives $\eta':L'\to L$ satisfying
\[
 \eta'(\eta(v))=v\quad(v\in L),\qquad
 \eta(\eta'(v'))=v'\quad(v'\in L').
\]
Define
\[
 \Psi_0(v)=
 \begin{cases}
 q_0(v)&v\in\operatorname{Vert}(Q),\\
 \eta(v)&v\in L,
 \end{cases}
 \qquad
 \Psi_1(\llbracket u,v\rrbracket)
 =\llbracket\Psi_0(u),\Psi_0(v)\rrbracket.
\]
Let $e=\llbracket u,v\rrbracket\in\operatorname{Edge}(T)$ with $u\in I_T$ and $v\in J_T\setminus L$.
The arrow $v\to u$ belongs to $Q$.
Since $q$ preserves incidence and colors, we have
\[
 \Psi_1(e)=\llbracket q_0(u),q_0(v)\rrbracket\in\operatorname{Edge}(T'),\qquad
 \kappa'(\Psi_1(e))=\kappa(e).
\]
For $v\in L$, the construction of $\eta(v)$ gives $\Psi_1(e)\in\operatorname{Edge}(T')$ and $\kappa'(\Psi_1(e))=\kappa(e)$.
The maps defined from $q^{-1}$ and $\eta'$ are inverse to $\Psi_0$ and $\Psi_1$.
Thus, $\Psi$ is a color-preserving, bipartition-preserving isomorphism extending $q$.
Any such isomorphism must send a leaf $v$ with neighbor $u$ and color $d$ to the unique leaf of color $d$ adjacent to $q_0(u)$.
Such an isomorphism therefore agrees with $\Psi$ on all vertices and edges.

\prf{Claim~\ref{c:retained-marked-transport}} Let us give the linear maps on the vertex spaces.
The colored-vacant-slot equivalence in Claim~\ref{c:retained-tree-extension} shows that $\mathbf s'$ is a marking of $Q'$.
For $v\in I_Q$, put
\[
 A_v=\operatorname{IncEdgeQuiv}_Q(v),\qquad
 A'_{q_0(v)}=\operatorname{IncEdgeQuiv}_{Q'}(q_0(v)).
\]
The quiver isomorphism $q$ restricts to a bijection
$q_1:A_v\to A'_{q_0(v)}$.
Define $h_u=\operatorname{id}_{\C}$ for $u\in J_Q$.
For $v\in I_Q$ and $z\in\operatorname{ZS}(Q)_v$, define
\begin{equation}\label{eq:zero-sum-transport}
 (h_v(z))_\beta=z_{q_1^{-1}(\beta)}
 \qquad(\beta\in A'_{q_0(v)}).
\end{equation}

Let us prove that $h_v$ is a linear isomorphism of the zero-sum spaces.
By~\eqref{eq:zero-sum-transport} and~\eqref{eq:zero-sum-map}, we have
\[
 \sum_{\beta\in A'_{q_0(v)}}(h_v(z))_\beta
 =\sum_{\alpha\in A_v}z_\alpha=0.
\]
Thus, $h_v(z)\in\operatorname{ZS}(Q')_{q_0(v)}$.
The map $h_v$ is linear. Its inverse satisfies
\[
 (h_v^{-1}(w))_\alpha=w_{q_1(\alpha)}
 \qquad(\alpha\in A_v,\ w\in\operatorname{ZS}(Q')_{q_0(v)}).
\]

Let us check the arrow maps and the markings.
For $\alpha:v\to u$ and $z\in\operatorname{ZS}(Q)_v$, equations~\eqref{eq:zero-sum-arrow} and~\eqref{eq:zero-sum-transport} give
\[
 \operatorname{ZS}(Q')_{q_1(\alpha)}(h_v(z))
 =(h_v(z))_{q_1(\alpha)}
 =z_\alpha
 =h_u(\operatorname{ZS}(Q)_\alpha(z)).
\]
Hence, $(q_0,q_1,(h_v)_{v\in\operatorname{Vert}(Q)})$ is a representation isomorphism.
For each $j\in[\ell]_{\Z}$, its image of the $j$th marking is $(q_0(u_j),d_j)$.
Hence, $(q_0,q_1,(h_v)_{v\in\operatorname{Vert}(Q)})$ is an isomorphism of the marked representations.
\end{proof}

\begin{lemma}[Recovery from a labeled support]\label{lem:marked-chain-recovery}
Let $a\in\Z_{>0}$ and $k\in\Z_{>0}$ satisfy $k\ge2$ and $a\ge100000k$.
Fix a family of representatives $\rho=(\rho_M)_{M\in\mathcal H_{a,k}}$.
Let $c=k-1$ and $S=\operatorname{LSTrivCh}(c)$.
Let $\mathbb S=(S,(I_S,J_S))$.
For $M=(S,f,\tau,(B_i)_i)\in\mathcal H_{a,k}$, let
$\rho_M=((\mathbb T_i,\kappa_i;\mathbf s_i))_{i\in J_S}$.
Let $\iota_i=\iota_i^\rho$ for $i\in J_S$.

Then, we have the following.
\begin{enumerate}
\item\label{c:marked-representative-independence} For any other family of representatives $\rho'$, we have
\[
 \mathcal S_{a,k}^\rho=\mathcal S_{a,k}^{\rho'},\qquad
 \pi_\rho=\pi_{\rho'}.
\]
\item\label{c:marked-recovery-unique} Assume the following.
\begin{hypotheses}
\item\label{hyp:recovery-support} $H\in\mathcal S_{a,k}^\rho$ and $X=(\widetilde{\mathbb G},d_G)$ represents $H$.
Let $\widetilde{\mathbb G}=(G,(I_G,J_G))$.
\item\label{hyp:recovery-labeling} $\lambda=(\lambda_0,\lambda_1):S\to S_X$ is an isomorphism preserving the bipartition.
\end{hypotheses}
Then, there is at most one $M\in\pi_\rho^{-1}(\{H\})$ admitting an isomorphism
$\phi:(\widetilde{\mathbb Q}_M^\rho,d_M^\rho)\to X$ such that
\begin{equation}\label{eq:recovery-matching}
\begin{aligned}
 \phi_0(x_t)&=\lambda_0(x_t)&& (t\in[c]_{\Z}),\\
 \phi_0\bigl(\iota_{i,0}(\operatorname{Vert}(Q(\mathbb T_i,\kappa_i)))\bigr)
 &=\operatorname{Vert}(\lambda_0(i))&& (i\in J_{S}).
\end{aligned}
\end{equation}
\item\label{c:marked-recovery-skeleton} Let $H\in\mathcal S_{a,k}^\rho$ and let $X=(\widetilde{\mathbb G},d_G)$ represent $H$.
For $M\in\pi_\rho^{-1}(\{H\})$, every isomorphism
$\phi:(\widetilde{\mathbb Q}_M^\rho,d_M^\rho)\to X$
induces an isomorphism $\lambda:S\to S_X$ preserving the bipartition and satisfying~\eqref{eq:recovery-matching}.
\end{enumerate}
\end{lemma}
\begin{proof}
\prf{Claim~\ref{c:marked-representative-independence}} Let us compare the two representatives of each $B_i$.
Fix $M=(S,f,\tau,(B_i)_i)\in\mathcal H_{a,k}$.
The equality of the marked classes gives isomorphisms
\[
 \Phi_i=(\phi_{i,0},\phi_{i,1},(h_{i,v})_v)
 \quad(i\in J_S)
\]
between the marked representations in $\rho_M$ and $\rho'_M$.
They preserve the bipartitions, the arrow colors, and the ordered marked colored-vacant slots.
Define maps on the glued vertices and arrows by
\[
\begin{aligned}
 \psi_0(x_t)&=x_t &&(t\in[c]_{\Z}),\\
 \psi_0(i,v)&=(i,\phi_{i,0}(v)) &&(i\in J_S,\ v\in\operatorname{Vert}(Q(\mathbb T_i,\kappa_i))),\\
 \psi_1(\alpha_e)&=\alpha'_e &&(e\in\operatorname{Edge}(S)),\\
 \psi_1(i,\alpha)&=(i,\phi_{i,1}(\alpha)) &&(i\in J_S,\ \alpha\in\operatorname{Edge}(Q(\mathbb T_i,\kappa_i))).
\end{aligned}
\]
For an added arrow $\alpha_e:x_t\to(i,v_e)$, preservation of its marked colored-vacant slot gives
$\alpha'_e:x_t\to(i,\phi_{i,0}(v_e))$ with the same color $d_e$.
Thus, $\psi=(\psi_0,\psi_1)$ is a color-preserving quiver isomorphism that preserves the bipartitions.
The maps $h_{i,v}$ give
$d_M^{\rho'}(\psi_0(i,v))=d_M^\rho(i,v)$.
At each $x_t$, both dimensions are one.
Restricting $\psi$ to the positive supports gives
\[
 (\widetilde{\mathbb Q}_M^\rho,d_M^\rho)
 \cong(\widetilde{\mathbb Q}_M^{\rho'},d_M^{\rho'}).
\]
Hence, $\pi_\rho(M)=\pi_{\rho'}(M)$ for every $M\in\mathcal H_{a,k}$.
Their image sets are therefore equal.

\prf{Claim~\ref{c:marked-recovery-unique}} Let us identify the components determined by $\lambda$.
Let $M\in\pi_\rho^{-1}(\{H\})$ and let $\phi:(\widetilde{\mathbb Q}_M^\rho,d_M^\rho)\to X$ be an isomorphism satisfying~\eqref{eq:recovery-matching}.
For $i\in J_{S}$, let
\[
\begin{gathered}
 K_i=\lambda_0(i),\qquad
 \widetilde Q_i=\operatorname{SubColQuiv}(G,\operatorname{Vert}(K_i)),\\
 I_i=I_G\cap\operatorname{Vert}(K_i),\qquad
 J_i=J_G\cap\operatorname{Vert}(K_i),\\
 \widetilde{\mathbb Q}_i=(\widetilde Q_i,(I_i,J_i)).
\end{gathered}
\]
Assumption~\ref{hyp:recovery-labeling} and~\eqref{eq:recovery-matching} give
$\widetilde{\mathbb Q}_i\cong Q(\mathbb T_i,\kappa_i)$ by a color-preserving isomorphism taking $(I_i,J_i)$ to $(I_{Q(\mathbb T_i,\kappa_i)},J_{Q(\mathbb T_i,\kappa_i)})$.

Let us recover the colored-vacant slots of the added arrows.
For $e=\llbracket x_t,i\rrbracket\in\operatorname{Edge}(S)$, there is a unique arrow
\[
 \alpha_e:\lambda_0(x_t)\longrightarrow u_e,\qquad
 u_e\in\operatorname{Vert}(K_i).
\]
Existence follows from the definition of $S_X$ and assumption~\ref{hyp:recovery-labeling}.
If there were two such arrows, their edges together with the path inside $K_i$ would form a cycle in $\operatorname{UnderG}(G)$.
This contradicts Claim~\ref{c:flow-stability-support-tree} of Lemma~\ref{lem:flow-stability}.
Define
\[
 d_e=\operatorname{QCol}(G)(\alpha_e),\qquad p_e=(u_e,d_e).
\]
Proper coloring of $G$ gives $p_e\in\operatorname{ColVacSlot}(\widetilde{\mathbb Q}_i)$.

Let us recover the color permutations and the ordered markings.
For $t\in[c]_{\Z}$, the color permutation is
\[
 \bigl(\tau_t(1),\tau_t(2),\tau_t(3)\bigr)
 =\bigl(d_{\operatorname{LeftEdge}_c(t)},
 d_{\operatorname{MidEdge}_c(t)},d_{\operatorname{RightEdge}_c(t)}\bigr).
\]
For $i\in L$, let $e_i$ be its unique incident edge. Its marking is $(p_{e_i})$.
For $t\in[c-1]_{\Z}$, its ordered marking is
\[
 \mathbf s_{z_t}
 =\bigl(p_{\operatorname{RightEdge}_c(t)},p_{\operatorname{LeftEdge}_c(t+1)}\bigr).
\]
These are exactly the orders required by assumption~\ref{hyp:marked-blocks} of Definition~\ref{def:marked-chain-data}.

Let us recover the marked zero-sum representations.
Suppose $M'=(S,f',\tau',(B'_i)_i)\in\pi_\rho^{-1}(\{H\})$ admits an isomorphism $\phi'$ satisfying~\eqref{eq:recovery-matching}.
Let $((\mathbb T'_i,\kappa'_i;\mathbf s'_i))_i=\rho_{M'}$ and let $\iota'_i$ be its inclusion maps.
For $i\in J_S$, let
\[
 Q_i=Q(\mathbb T_i,\kappa_i),\qquad Q'_i=Q(\mathbb T'_i,\kappa'_i),\qquad
 q_i=(\phi'\circ\iota'_i)^{-1}\circ(\phi\circ\iota_i):Q_i\longrightarrow Q'_i.
\]
Both maps in this composition are isomorphisms onto $\widetilde{\mathbb Q}_i$.
For each edge $e$ incident to $i$, let $(v_{i,e},d_e)$ and $(v'_{i,e},d_e)$ be the colored-vacant slots assigned to $e$ in the two representatives.
Uniqueness of $\alpha_e$ gives
\[
 (q_{i,0}(v_{i,e}),d_e)=(v'_{i,e},d_e).
\]
The ordered markings agree under $q_i$ because their orders are determined by $S$.

Let us recover each class $B_i$ for $i\in J_S$.
Assumption~\ref{hyp:marked-blocks} of Definition~\ref{def:marked-chain-data} gives
\[
 \#I_{T_i}=a_i(f)+1\ge1,\qquad
 \#I_{T'_i}=a_i(f')+1\ge1.
\]
Claim~\ref{c:retained-tree-extension} of Lemma~\ref{lem:retained-tree-recovery} recovers an isomorphism of these colored trees.
Claim~\ref{c:retained-marked-transport} of Lemma~\ref{lem:retained-tree-recovery} gives
\[
 B_i=B'_i\qquad(i\in J_S).
\]
Thus, $G$, $d_G$, and $\lambda$ determine every $B_i$.

Let us recover the flow coordinates.
For each $i\in J_{S}$, its source dimension is
\[
 a_i=\sum_{v\in I_i}d_G(v).
\]
For $e\in\operatorname{Edge}(S)$, put $S_e=\operatorname{SnkCtdComp}(\mathbb S;e)$.
Claim~\ref{c:tree-flow-cut} of Lemma~\ref{lem:tree-flow-cut} gives
\[
 f_e=\sum_{i\in \operatorname{SubSnk}_{\mathbb S}(S_e)}\left(\frac ak-a_i\right)
      -\left(\frac ak+1\right)\#\operatorname{SubSrc}_{\mathbb S}(S_e).
\]
In particular,
\[
 l_t=f_{\operatorname{LeftEdge}_c(t)},\qquad
 r_t=f_{\operatorname{RightEdge}_c(t)}\qquad(t\in[c]_{\Z}).
\]
Therefore, $G$, $d_G$, and $\lambda$ determine every entry of $M=(S,f,\tau,(B_i)_{i\in L\sqcup I})$.

\prf{Claim~\ref{c:marked-recovery-skeleton}} Let us identify the sources in $C(X)$.
Claim~\ref{c:flow-choices-range} of Lemma~\ref{lem:flow-choices} and assumption~\ref{hyp:marked-blocks} of Definition~\ref{def:marked-chain-data} give $\#I_{T_i}=a_i(f)+1\ge2$ for $i\in J_S$.
For each $t\in[c]_{\Z}$, the new source satisfies
\[
 d_M^\rho(x_t)=1,\qquad \deg_{\widetilde{\mathbb Q}_M^\rho}(x_t)=3.
\]
For $i\in J_S$ and $v\in I_{Q(\mathbb T_i,\kappa_i)}$, the added arrows do not meet the copied source $(i,v)$.
Thus,
\[
 \deg_{\widetilde{\mathbb Q}_M^\rho}(i,v)
 =\deg_{Q(\mathbb T_i,\kappa_i)}(v).
\]
Claim~\ref{c:reflected-blocks-sources} of Lemma~\ref{lem:reflected-blocks} gives
\[
 d_M^\rho(i,v)=1\quad\Longrightarrow\quad
 \deg_{\widetilde{\mathbb Q}_M^\rho}(i,v)=2.
\]
Consequently,
\[
 C(X)=\{\phi_0(x_t):t\in[c]_{\Z}\}.
\]

Let us identify the components after deleting $C(X)$.
Each $Q(\mathbb T_i,\kappa_i)$ is connected because it is obtained from $T_i$ by deleting its sink leaves.
By Definition~\ref{def:attachments}, deleting the new sources $x_t$ removes every added arrow and leaves exactly the disjoint copied quivers.
Thus, for each $i\in J_S$, there is a unique $K_i\in\mathcal K(X)$ satisfying
\[
 \operatorname{Vert}(K_i)=\phi_0\bigl(\iota_{i,0}(\operatorname{Vert}(Q(\mathbb T_i,\kappa_i)))\bigr).
\]
Define $\lambda_0(x_t)=\phi_0(x_t)$ for $t\in[c]_{\Z}$ and $\lambda_0(i)=K_i$ for $i\in J_S$.
These maps are bijections onto $C(X)$ and $\mathcal K(X)$, respectively.

Let us verify incidence.
The added arrows in Definition~\ref{def:attachments} and Definition~\ref{def:marked-chain-data}(\ref{d:recovered-skeleton}) give
\[
 \llbracket x_t,i\rrbracket\in\operatorname{Edge}(S)
 \quad\Longleftrightarrow\quad
 \llbracket\phi_0(x_t),K_i\rrbracket\in\operatorname{Edge}(S_X).
\]
Hence, $\lambda_0$ induces the required isomorphism $\lambda:S\to S_X$.
\end{proof}

\begin{lemma}[Marked chains and their supports]\label{lem:forget-markings}
Let $a\in\Z_{>0}$ and $k\in\Z_{>0}$ satisfy $k\ge2$ and $a\ge100000k$.
Fix a family of representatives $\rho=(\rho_M)_{M\in\mathcal H_{a,k}}$.
Let $c=k-1$ and $S=\operatorname{LSTrivCh}(c)$. Consider the map $\pi_\rho\colon\mathcal H_{a,k}\to\mathcal S_{a,k}^\rho$.
Let
\[
 L=J_S\cap\operatorname{Leaf}(S),\qquad I=\{i\in J_S:\deg_S(i)=2\}.
\]
Then, we have the following.
\begin{enumerate}
\item\label{c:forget-markings-eight} For $H\in\mathcal S_{a,k}^\rho$,
we have $\#\pi_\rho^{-1}(\{H\})\le8$.
\item\label{c:forget-markings-six} If $c=1$, then
$\#\pi_\rho^{-1}(\{H\})\le6$ for $H\in\mathcal S_{a,k}^\rho$.
\item\label{c:marked-color-count} There are $6^c$ families $\tau=(\tau_t)_{t\in[c]_{\Z}}$ such that each $\tau_t:[3]_{\Z}\to[0,2]_{\Z}$ is a bijection.
\item\label{c:marked-block-count} Let $f\in\mathcal F_{a,k}$.
Let $\tau=(\tau_t)_{t\in[c]_{\Z}}$ be a family of bijections $\tau_t:[3]_{\Z}\to[0,2]_{\Z}$.
Let $a_i=a_i(f)$ for $i\in J_S$.
The number of choices of $(B_i)_{i\in L\sqcup I}$ is
\[
 \prod_{i\in L}P_{a_i}
 \prod_{i\in I_=(\tau)}Q_{a_i}^{=}
 \prod_{i\in I_{\ne}(\tau)}Q_{a_i}^{\ne}.
\]
For every $(P_j,Q_j^{=},256)$-admissible pair $(A,B)$, this number is at least
\[
 \beta^{a-c}A^{c+2}B^{c-1}
 \prod_{i\in L}a_i^{-\frac32}
 \prod_{i\in I}a_i^{-\frac12}.
\]
\end{enumerate}
\end{lemma}
\begin{proof}
\prf{Claim~\ref{c:forget-markings-eight}} Let us bound the number of marked chains over one support.
Choose a representative $X=(\widetilde{\mathbb G},d_G)$ of $H$.
Let us count the bipartition-preserving isomorphisms
\[
 \mathcal I_X=\{\lambda:S\xrightarrow{\sim}S_X:
 \lambda_0(I_S)=I_{S_X},\ \lambda_0(J_S)=J_{S_X}\}.
\]

Let us associate an element of $\mathcal I_X$ to each $M\in\pi_\rho^{-1}(\{H\})$.
Choose an isomorphism
\[
 \phi_M:(\widetilde{\mathbb Q}_M^\rho,d_M^\rho)\longrightarrow X.
\]
Claim~\ref{c:marked-recovery-skeleton} of Lemma~\ref{lem:marked-chain-recovery} gives a bipartition-preserving isomorphism $\lambda_M:S\to S_X$ satisfying
\[
 \begin{aligned}
 (\lambda_M)_0(x_t)&=(\phi_M)_0(x_t)&& (t\in[c]_{\Z}),\\
 \operatorname{Vert}((\lambda_M)_0(i))
 &=(\phi_M)_0\bigl(\iota_{i,0}^\rho(\operatorname{Vert}(Q(\mathbb T_i,\kappa_i)))\bigr)
 &&(i\in J_S).
 \end{aligned}
\]
Thus, $\lambda_M\in\mathcal I_X$.
By Claim~\ref{c:marked-recovery-unique} of Lemma~\ref{lem:marked-chain-recovery}, each $\lambda\in\mathcal I_X$ determines at most one such $M$.
Hence,
\[
 \#\pi_\rho^{-1}(\{H\})\le\#\mathcal I_X.
\]

First, suppose $c\ge2$. Deleting the sink leaves of $S$ gives a path with endpoints $x_1,x_c$ and vertex set
\[
 \{x_t:t\in[c]_{\Z}\}\sqcup\{z_t:t\in[c-1]_{\Z}\}.
\]
An isomorphism sends $x_1$ to one of the two endpoints of the corresponding path in $S_X$.
For $\lambda\in\mathcal I_X$, choosing $\lambda_0(x_1)$ determines
$\lambda_0(x_t)$ for $t\in[c]_{\Z}$ and $\lambda_0(z_t)$ for $t\in[c-1]_{\Z}$.
The two sink leaves at $x_1$ are $z_0,y_1$, and those at $x_c$ are $y_c,z_c$.
Each pair has $2!$ possible permutations.
For $1<t<c$, the only sink leaf at $x_t$ is $y_t$, whose image is then determined.
Thus,
\[
 \#\pi_\rho^{-1}(\{H\})\le\#\mathcal I_X
 \le2\cdot2!\cdot2!=8.
\]

Second, suppose $c=1$. The only source is $x_1$, and its three sink neighbors are $z_0,y_1,z_1$.
An isomorphism is determined by a permutation of these sinks. Therefore,
\[
 \#\pi_\rho^{-1}(\{H\})\le\#\mathcal I_X\le3!=6\le8.
\]

\prf{Claim~\ref{c:forget-markings-six}}
The case $c=1$ in the proof of Claim~\ref{c:forget-markings-eight} gives
\[
 \#\pi_\rho^{-1}(\{H\})\le6.
\]

\prf{Claim~\ref{c:marked-color-count}} Let us count the color bijections independently at each source.
Each $\tau_t$ is a bijection between
sets of three elements. The choices are independent. Thus, their number is
$(3!)^c=6^c$.

\prf{Claim~\ref{c:marked-block-count}} Let us count the choices of zero-sum representations. Fix $f$ and $\tau$.
Definition~\ref{def:marked-block-counts}(1) gives $P_{a_i}$ choices for $i\in L$.
Parts~(2) and~(3) of that definition give $Q_{a_i}^{=}$ choices for $i\in I_=(\tau)$
and $Q_{a_i}^{\ne}$ choices for $i\in I_{\ne}(\tau)$.
The choices of zero-sum representations are independent. Thus, their counts multiply.

Let us check that the counts depend only on equality of the marked colors.
For colors $d,d'\in[0,2]_{\Z}$, choose a permutation $\sigma$ of $[0,2]_{\Z}$ such that $\sigma(d)=d'$.
For ordered pairs $(d,e)$ and $(d',e')$ satisfying
\[
 d=e\quad\Longleftrightarrow\quad d'=e',
\]
choose $\sigma$ such that $\sigma(d)=d'$ and $\sigma(e)=e'$.
Replacing $\kappa$ by $\sigma\circ\kappa$ gives a bijection of the corresponding marked classes.
Its inverse uses $\sigma^{-1}$. Thus, the counts apply to every family $\tau=(\tau_t)_{t\in[c]_{\Z}}$.

Let us obtain a uniform lower bound for the product of marked counts.
Let $(A,B)$ be a $(P_j,Q_j^{=},256)$-admissible pair.
Claim~\ref{c:slots-distinct} of Lemma~\ref{slots} gives $Q_j^{\ne}\ge Q_j^=$.
Claims~\ref{c:flow-choices-sum} and~\ref{c:flow-choices-range} of Lemma~\ref{lem:flow-choices}
give $\sum_i a_i=a-c$ and $a_i>256$.
The $(P_j,Q_j^{=},256)$-admissibility of $(A,B)$ gives the lower bound since $\#L=c+2$ and $\#I=c-1$.
\end{proof}

\begin{proposition}[Counting stable dimensioned supports]\label{prop:cover-count}
Let $k\ge2$ and $a>0$ be integers. Let $c=k-1$, $x=\frac ak$,
$N=256$, $\gamma=1-\frac Nx$, and $\lambda=\frac1{100}$.
Let $S=\operatorname{LSTrivCh}(c)$.
Assume the following.
\begin{hypotheses}
\item\label{hyp:cover-dimensions} $a\ge100000k$.
\item\label{hyp:cover-slot-bounds} Let $(A,B)$ be a $(P_j,Q_j^{=},256)$-admissible pair.
\end{hypotheses}
Fix a family of representatives $\rho=(\rho_M)_{M\in\mathcal H_{a,k}}$.
Then,
\begin{equation}
\#\mathcal S_{a,k}^\rho\ge\beta^a x^{-\frac{5}{2}}\frac{A^2}{8B\lambda N^2}
 \left(\frac{6AB\lambda\gamma^2}{\beta}\right)^c.\label{count}
\end{equation}
\end{proposition}
\begin{proof}
Let us count the marked chains for a fixed color family.
Assumption~\ref{hyp:cover-dimensions} gives $x\ge100000$.
Fix a family $\tau=(\tau_t)_{t\in[c]_{\Z}}$ of bijections
$\tau_t:[3]_{\Z}\to[0,2]_{\Z}$.
Claim~\ref{c:flow-choices-count} of Lemma~\ref{lem:flow-choices} gives
\[
 \#\mathcal F_{a,k}\ge\left(\frac{\gamma x}{N}\right)^{2c}.
\]
Claim~\ref{c:flow-choices-distinct} of Lemma~\ref{lem:flow-choices}
gives distinct labeled size lists for distinct $f\in\mathcal F_{a,k}$.
Apply Claim~\ref{c:marked-block-count} of Lemma~\ref{lem:forget-markings}
and Claim~\ref{c:flow-choices-weight} of Lemma~\ref{lem:flow-choices}.
The $c-1$ pairs in~\eqref{weight} and the leaves of sizes
$x-l_1$, $x-r_c$, and $l_c+r_c-1$ give
\[
\begin{aligned}
 &\#\{(S,f,\tau,(B_i)_i)\in\mathcal H_{a,k}\}\\
 &\quad\ge\left(\frac{\gamma x}{N}\right)^{2c}
   \beta^{a-c}A^{c+2}B^{c-1}
   (\lambda N^2x^{-2})^{c-1}x^{-\frac92}\\
 &\quad=\beta^{a-c}A^{c+2}B^{c-1}\gamma^{2c}x^{-\frac52}N^{-2}\lambda^{c-1}.
\end{aligned}
\]
For $c=1$, the product over pairs is empty and $B^{c-1}\lambda^{c-1}=1$.

Let us sum over color families and pass to isomorphism classes.
Claim~\ref{c:marked-color-count} of Lemma~\ref{lem:forget-markings}
gives $6^c$ choices of $\tau$.
Claim~\ref{c:forget-markings-eight} of Lemma~\ref{lem:forget-markings} gives
\[
\begin{aligned}
 \#\mathcal S_{a,k}^\rho
 &\ge\frac18\#\mathcal H_{a,k}\\
 &\ge\frac{6^c}{8}\beta^{a-c}A^{c+2}B^{c-1}
 \gamma^{2c}x^{-\frac52}N^{-2}\lambda^{c-1}.
\end{aligned}
\]
Substitute $\beta^{a-c}=\beta^a\beta^{-c}$ to obtain~\eqref{count}.
\end{proof}

\begin{lemma}[Traversals without immediate reversal in a tree]\label{lem:tree-traversal-distinct}
Let $G$ be a tree and let $n\in\Z_{\ge0}$.
Let $v_i\in\operatorname{Vert}(G)$ for $i\in[0,n]_{\Z}$ and suppose
\[
 e_i=\llbracket v_{i-1},v_i\rrbracket\in\operatorname{Edge}(G)
 \quad(i\in[n]_{\Z}),\qquad
 e_i\ne e_{i+1}\quad(i\in[n-1]_{\Z}).
\]
Then, $v_p\ne v_q$ for $0\le p<q\le n$.
\end{lemma}
\begin{proof}
Let us suppose $v_p=v_q$ for some $0\le p<q\le n$.
For $i\in[p,q]_{\Z}$, let $h_i$ be the length of the unique simple path from $v_p$ to $v_i$.
Then,
\[
 h_p=h_q=0,\qquad h_{p+1}=1.
\]
In particular, $q\ge p+2$.
Choose $j\in[p+1,q-1]_{\Z}$ satisfying
\[
 h_j=\max_{i\in[p,q]_{\Z}}h_i.
\]
The unique simple paths from $v_p$ to adjacent vertices differ by one edge.
Thus, $|h_{i+1}-h_i|=1$ for $p\le i<q$.
Maximality of $h_j$ gives
\[
 h_{j-1}=h_j-1=h_{j+1}.
\]
Both $v_{j-1}$ and $v_{j+1}$ are the predecessor of $v_j$ on the unique simple path from $v_p$ to $v_j$.
Consequently,
\[
 v_{j-1}=v_{j+1},\qquad
 e_j=\llbracket v_{j-1},v_j\rrbracket
 =\llbracket v_j,v_{j+1}\rrbracket=e_{j+1}.
\]
This contradicts $e_j\ne e_{j+1}$.
\end{proof}

\begin{lemma}[Embedding a colored tree]\label{lem:cover-tree-embedding}
Consider the universal covering
\[
 \varpi:\widetilde K_3\longrightarrow K(3)
\]
with arrow coloring $\operatorname{CovCol}_{\varpi}=c_3\circ\varpi_1$.
Assume the following.
\begin{hypotheses}
\item\label{hyp:embedding-tree} $(T,(I_T,J_T))$ is a bipartite quiver such that $T$ is finite and $\operatorname{UnderG}(T)$ is a tree.
The map $\kappa:\operatorname{Edge}(T)\to[0,2]_{\Z}$ is a proper coloring.
\item\label{hyp:embedding-root} $v_0\in\operatorname{Vert}(T)$ and $\widehat v_0\in\operatorname{Vert}(\widetilde K_3)$ satisfy
\[
 v_0\in I_T\quad\Longleftrightarrow\quad\widehat v_0\in I_{\widetilde K_3}.
\]
\end{hypotheses}
Then, there is a unique bipartition-preserving embedding $\iota:T\to\widetilde K_3$ such that
\[
 \iota_0(v_0)=\widehat v_0,\qquad
 \operatorname{CovCol}_{\varpi}(\iota_1(\alpha))=\kappa(\alpha)
 \quad(\alpha\in\operatorname{Edge}(T)).
\]
\end{lemma}
\begin{proof}
Let us define the vertex and arrow maps.
Define $q:T\to K(3)$ by $q_0(v)=1$ for $v\in I_T$ and $q_0(v)=2$ for $v\in J_T$, where $1$ and $2$ denote the source and sink vertices of $K(3)$, respectively. Let
\[
 q_1(\alpha)=c_3^{-1}(\kappa(\alpha)).
\]
For $v\in\operatorname{Vert}(T)$, let $\gamma_v$ be the unique simple path from $v_0$ to $v$, viewed as an arrow traversal.
Its projection $q\gamma_v$ starts at $\varpi_0(\widehat v_0)$ by assumption~\ref{hyp:embedding-root}.
Define $\iota_0(v)$ to be the final vertex of $\operatorname{Lift}_{\varpi}(q\gamma_v;\widehat v_0)$.
For $\alpha:u\to v$, define $\iota_1(\alpha)$ to be the unique lift of $q_1(\alpha)$ starting at $\iota_0(u)$.
The simple paths to adjacent vertices differ by one edge.
Uniqueness of their lifts gives
\[
 s_{\widetilde K_3}(\iota_1(\alpha))=\iota_0(u),\qquad
 t_{\widetilde K_3}(\iota_1(\alpha))=\iota_0(v).
\]
The equality $\varpi_1(\iota_1(\alpha))=q_1(\alpha)$ gives $\operatorname{CovCol}_{\varpi}(\iota_1(\alpha))=\kappa(\alpha)$.

Let us prove injectivity.
Let $u\in\operatorname{Vert}(T)$ and $v\in\operatorname{Vert}(T)$ satisfy $u\ne v$.
Let
\[
 (z_0,\bar\alpha_1,z_1,\ldots,\bar\alpha_n,z_n),\qquad
 z_0=u,\quad z_n=v,\quad n\ge1
\]
be the unique simple path from $u$ to $v$.
For $i\in[0,n]_{\Z}$, let $w_i=\iota_0(z_i)$.
For $i\in[n]_{\Z}$, let $\widehat\alpha_i=\iota_1(\alpha_i)$.
Assumption~\ref{hyp:embedding-tree} and the color identity give
\[
 \operatorname{CovCol}_{\varpi}(\widehat\alpha_i)
 \ne\operatorname{CovCol}_{\varpi}(\widehat\alpha_{i+1})
 \qquad(i\in[n-1]_{\Z}).
\]
Let $H=\operatorname{UnderG}(\widetilde K_3)$ and
\[
 e_i=\llbracket w_{i-1},w_i\rrbracket\in\operatorname{Edge}(H)
 \qquad(i\in[n]_{\Z}).
\]
The distinct colors give $\widehat\alpha_i\ne\widehat\alpha_{i+1}$ for $i\in[n-1]_{\Z}$.
Since $H$ has no parallel edges, $e_i\ne e_{i+1}$ for $i\in[n-1]_{\Z}$.
Lemma~\ref{lem:tree-traversal-distinct} gives
\[
 \iota_0(u)=w_0\ne w_n=\iota_0(v).
\]
If $\iota_1(\alpha)=\iota_1(\beta)$, their two endpoints agree by injectivity of $\iota_0$.
The tree $\operatorname{UnderG}(T)$ has no parallel edges by assumption~\ref{hyp:embedding-tree}.
Hence, $\alpha=\beta$.

Let us prove uniqueness.
For any embedding $\jmath:T\to\widetilde K_3$ satisfying $\jmath_0(v_0)=\widehat v_0$ and $\operatorname{CovCol}_{\varpi}\circ\jmath_1=\kappa$, the traversal $\jmath(\gamma_v)$ is the lift of $q\gamma_v$ from $\widehat v_0$.
Lemma~\ref{lem:path-lift} determines this lift and hence every vertex and arrow image.
\end{proof}

\Needspace{8\baselineskip}
\begin{lemma}[Deck transformations with a fixed vertex image]\label{lem:cover-deck-transitivity}
Consider the universal covering $\varpi:\widetilde K_3\to K(3)$.
Assume the following.
\begin{hypotheses}
\item\label{hyp:deck-vertices} $u\in\operatorname{Vert}(\widetilde K_3)$ and $v\in\operatorname{Vert}(\widetilde K_3)$.
\item\label{hyp:deck-fiber} $\varpi_0(u)=\varpi_0(v)$.
\end{hypotheses}
Then, there is a unique $\delta\in\operatorname{Deck}(\varpi)$ such that $\delta_0(u)=v$.
\end{lemma}
\begin{proof}
Let us construct $\delta$ by lifting paths.
For $w\in\operatorname{Vert}(\widetilde K_3)$, let $\gamma_w$ be the unique simple path from $u$ to $w$.
By assumption~\ref{hyp:deck-fiber}, the vertex $v$ lies over the initial vertex of $\varpi\gamma_w$. Lemma~\ref{lem:path-lift} gives its unique lift starting at $v$.
Define $\delta_0(w)$ to be its final vertex.
For an arrow $\alpha:w\to w'$, define $\delta_1(\alpha)$ to be the unique lift of $\varpi_1(\alpha)$ starting at $\delta_0(w)$.
Uniqueness in Lemma~\ref{lem:path-lift} gives its endpoint $\delta_0(w')$ and
\[
 \varpi_0\circ\delta_0=\varpi_0,\qquad
 \varpi_1\circ\delta_1=\varpi_1.
\]
Interchanging $u$ and $v$ gives a quiver homomorphism $\epsilon$ over $K(3)$.
The composites $\epsilon\circ\delta$ and $\delta\circ\epsilon$ fix $u$ and $v$, respectively.
Applied to the paths from these vertices, uniqueness of lifting gives
\[
 \epsilon\circ\delta=1_{\widetilde K_3},\qquad
 \delta\circ\epsilon=1_{\widetilde K_3}.
\]
Hence, $\delta\in\operatorname{Deck}(\varpi)$.
If $\delta'\in\operatorname{Deck}(\varpi)$ satisfies $\delta'_0(u)=v$, then $\delta'\gamma_w$ and $\delta\gamma_w$ both lift $\varpi\gamma_w$ from $v$ for every $w\in\operatorname{Vert}(\widetilde K_3)$.
Lemma~\ref{lem:path-lift} gives $\delta'_0(w)=\delta_0(w)$ for every $w$ and $\delta'_1(\alpha)=\delta_1(\alpha)$ for every arrow $\alpha$.
Thus, $\delta'=\delta$.
\end{proof}

\begin{lemma}[Dimension vectors from colored trees]\label{lem:cover-dimension-equivalence}
Consider the universal covering
\[
 \varpi:\widetilde K_3\longrightarrow K(3)
\]
with arrow coloring $\operatorname{CovCol}_{\varpi}=c_3\circ\varpi_1$.
Assume the following.
\begin{hypotheses}
\item\label{hyp:dimensioned-trees} $((T,\kappa),(I_T,J_T))$ and $((T',\kappa'),(I_{T'},J_{T'}))$ are properly colored bipartite quivers such that $T$ and $T'$ are finite and $\operatorname{UnderG}(T)$ and $\operatorname{UnderG}(T')$ are trees.
\item\label{hyp:dimensioned-embeddings} $\iota:T\hookrightarrow\widetilde K_3$ and $\iota':T'\hookrightarrow\widetilde K_3$ preserve the bipartitions and satisfy
\[
 \operatorname{CovCol}_{\varpi}\circ\iota_1=\kappa,\qquad
 \operatorname{CovCol}_{\varpi}\circ\iota'_1=\kappa'.
\]
\item\label{hyp:positive-tree-dimensions} $\rho:\operatorname{Vert}(T)\to\Z_{>0}$ and $\rho':\operatorname{Vert}(T')\to\Z_{>0}$.
\end{hypotheses}
Define finitely supported dimension vectors $d$ and $d'$ by
\[
\begin{aligned}
 d_{\iota_0(v)}&=\rho(v)&& (v\in\operatorname{Vert}(T)),&
 d_w&=0&& (w\notin\iota_0(\operatorname{Vert}(T))),\\
 d'_{\iota'_0(v')}&=\rho'(v')&& (v'\in\operatorname{Vert}(T')),&
 d'_w&=0&& (w\notin\iota'_0(\operatorname{Vert}(T'))).
\end{aligned}
\]
Then, $[d]_{\varpi}=[d']_{\varpi}$ if and only if there is a bipartition-preserving quiver isomorphism $\psi:T\to T'$ satisfying
\[
 \kappa'\circ\psi_1=\kappa,\qquad
 \rho'\circ\psi_0=\rho.
\]
\end{lemma}
\begin{proof}
Let us restrict a deck transformation to the positive supports.
Suppose $d'=\delta\cdot d$ for $\delta\in\operatorname{Deck}(\varpi)$.
Assumption~\ref{hyp:positive-tree-dimensions} gives
\[
 \operatorname{VecSupp}(d)=\iota_0(\operatorname{Vert}(T)),\qquad
 \operatorname{VecSupp}(d')=\iota'_0(\operatorname{Vert}(T')).
\]
The images of $\iota$ and $\iota'$ are full subquivers of $\widetilde K_3$: an additional edge between image vertices would form a cycle together with their path in the image tree.
Thus, $\delta$ restricts to an isomorphism from the image of $\iota$ to the image of $\iota'$.
The maps $\iota$ and $\iota'$ in assumption~\ref{hyp:dimensioned-embeddings} give an isomorphism $\psi:T\to T'$ such that
\[
 \iota'\circ\psi=\delta\circ\iota.
\]
Deck transformations preserve $\operatorname{CovCol}_{\varpi}$, $I_{\widetilde K_3}$, and $J_{\widetilde K_3}$.
The equality $d'=\delta\cdot d$ gives $\rho'(\psi_0(v))=\rho(v)$ for every $v$.

Let us prove the converse by extending the colored tree isomorphism to the covering.
Conversely, let $\psi:T\to T'$ be a bipartition-preserving quiver isomorphism satisfying $\kappa'\circ\psi_1=\kappa$ and $\rho'\circ\psi_0=\rho$.
By assumption~\ref{hyp:dimensioned-trees}, choose $v_0\in\operatorname{Vert}(T)$.
Bipartition preservation gives
\[
 \varpi_0(\iota_0(v_0))=\varpi_0(\iota'_0(\psi_0(v_0))).
\]
Lemma~\ref{lem:cover-deck-transitivity} gives $\delta\in\operatorname{Deck}(\varpi)$ satisfying $\delta_0(\iota_0(v_0))=\iota'_0(\psi_0(v_0))$.
Their arrow colors satisfy
\[
 \operatorname{CovCol}_{\varpi}\circ\delta_1\circ\iota_1
 =\kappa=\kappa'\circ\psi_1
 =\operatorname{CovCol}_{\varpi}\circ\iota'_1\circ\psi_1.
\]
Uniqueness in Lemma~\ref{lem:cover-tree-embedding} gives $\delta\circ\iota=\iota'\circ\psi$.
The identity $\rho'\circ\psi_0=\rho$ then gives $d'=\delta\cdot d$.
\end{proof}

\begin{lemma}[Finiteness of connected support classes]\label{lem:cover-trees-finite}
Consider the universal covering $\varpi:\widetilde K_3\to K(3)$.
Let $a\in\Z_{>0}$ and $b\in\Z_{>0}$.
There are finitely many classes in
\[
 \left\{[d]_{\varpi}:d\in\operatorname{CovDim}_{\varpi}(a,b),\
 \operatorname{SubQuiv}(\widetilde K_3,\operatorname{VecSupp}(d))
 \text{ is connected}\right\}.
\]
\end{lemma}
\begin{proof}
Let us bound each connected support after a deck transformation.
Let $d\in\operatorname{CovDim}_{\varpi}(a,b)$ and assume its positive support $D=\operatorname{VecSupp}(d)$ is connected.
Since every $d_v$ for $v\in D$ is a positive integer,
\[
 \#D\le\sum_{v\in D}d_v=a+b.
\]
The inequality $a>0$ gives a source in $D$.
Fix a source $o$ of $\widetilde K_3$.
Lemma~\ref{lem:cover-deck-transitivity} gives a representative of $[d]_{\varpi}$ with $o\in D$.
For a vertex $v$, define $\operatorname{dist}(o,v)$ as the number of edges in the unique path from $o$ to $v$.
Since $D$ is connected, its path from $o$ to each $v\in D$ has at most $\#D-1$ edges. Hence,
\[
 D\subseteq\{v:\operatorname{dist}(o,v)\le a+b-1\},\qquad
 d_v\in[a+b]_{\Z}\quad(v\in D).
\]
Every vertex of $\widetilde K_3$ has degree three.
Thus, $\{v\in\operatorname{Vert}(\widetilde K_3):\operatorname{dist}(o,v)\le a+b-1\}$ is finite.
There are finitely many sets $D\subseteq\{v:\operatorname{dist}(o,v)\le a+b-1\}$ and finitely many functions $D\to[a+b]_{\Z}$.
These finite choices contain a representative of every class in the statement.
\end{proof}

\begin{lemma}[Cover stability and localization]\label{lem:cover-localization}
Let $a$ and $k$ be integers. Assume the following.
\begin{hypotheses}
\item\label{hyp:localization-dimensions} $k\ge2$ and $a\ge100000k$.
\item\label{hyp:localization-coprimality} $\gcd(a,k)=1$.
\end{hypotheses}
Fix a family of representatives $\rho=(\rho_M)_{M\in\mathcal H_{a,k}}$.
Let $I=I_{\widetilde K_3}$ and $J=J_{\widetilde K_3}$.
Then, the following hold.
\begin{enumerate}
\item\label{c:cover-stability-weight} Let $R$ be a representation of dimension $\widetilde d\in\operatorname{CovDim}_{\varpi}(a,a+k)$.
For $W\in\operatorname{SubRep}(R)$, let
\[
 u=\sum_{i\in I}\dim W_i,\qquad v=\sum_{j\in J}\dim W_j.
\]
Then,
\[
 \theta^{\varpi}_{a,a+k}(\dim R)=0,\qquad
 \theta^{\varpi}_{a,a+k}(\dim W)=av-(a+k)u.
\]
Thus, the stability condition~\eqref{eq:stability} with $b=a+k$ is equivalent to
$\theta^{\varpi}_{a,a+k}(\dim W)>0$ for every proper nonzero subrepresentation $W$.
\item\label{c:cover-stability-connected} Let $R$ be a representation of dimension $\widetilde d\in\operatorname{CovDim}_{\varpi}(a,a+k)$.
If $R$ satisfies the stability condition in Claim~\ref{c:cover-stability-weight},
then the quiver $\operatorname{SubQuiv}(\widetilde K_3,\operatorname{VecSupp}(\widetilde d))$ is connected.
\item\label{c:cover-localization-bound} We have
\[
 E_k(a)\ge\#\mathcal S_{a,k}^\rho.
\]
\end{enumerate}
\end{lemma}
\begin{proof}
\prf{Claim~\ref{c:cover-stability-weight}} Let us check the dimension-vector and stability conditions for the support.
Definition~\ref{def:cover-dimension-classes}(1)--(2) gives
\[
 \theta^{\varpi}_{a,a+k}(\dim R)=a(a+k)-(a+k)a=0,
 \qquad\theta^{\varpi}_{a,a+k}(\dim W)=av-(a+k)u.
\]
Substituting $b=a+k$ in~\eqref{eq:stability} proves the equivalence.

\prf{Claim~\ref{c:cover-stability-connected}} Let us prove connectedness of the support by contradiction.
Suppose that $\operatorname{VecSupp}(\widetilde d)$ is disconnected.
Then, $R=R_1\oplus R_2$ for two nonzero subrepresentations.
Claim~\ref{c:cover-stability-weight} gives
\[
 0<\theta^{\varpi}_{a,a+k}(\dim R_1)+\theta^{\varpi}_{a,a+k}(\dim R_2)
 =\theta^{\varpi}_{a,a+k}(\dim R)=0,
\]
a contradiction.

\prf{Claim~\ref{c:cover-localization-bound}} Let us obtain a sum formula for $E_k(a)$ over cover dimension vectors.

Let us transfer stability to each chosen representative.
We use the stability condition in Claim~\ref{c:cover-stability-weight}
for $\cM^s_{\widetilde d}(\widetilde K_3)$.
On $\widetilde K_3$, the equivalence classes of dimension vectors in~\cite[Corollary 3.16]{Weist} are the $\operatorname{Deck}(\varpi)$-orbits.
By~\cite[Corollary 3.16]{Weist}, we have
\begin{equation}
 E_k(a)=\chi\bigl(\cM^s_{(a,a+k)}(K(3))\bigr)
 =\sum_{[\widetilde d]_{\varpi}\in\mathcal D_{a,k}}
 \chi\bigl(\cM^s_{\widetilde d}(\widetilde K_3)\bigr).
\label{eq:cover-sum}
\end{equation}

Let us restrict the sum in~\eqref{eq:cover-sum} to finitely many nonempty moduli spaces.
Claim~\ref{c:cover-stability-connected} gives connected positive support for each nonempty summand.
Lemma~\ref{lem:cover-trees-finite}
with $b=a+k$ gives only finitely many such classes.
For each such class $[\widetilde d]_{\varpi}$, let $S_{\widetilde d}=\operatorname{VecSupp}(\widetilde d)$.

Let us restrict stable representations to $S_{\widetilde d}$.
Let $\operatorname{SubQuiv}(\widetilde K_3,S_{\widetilde d})$ be the full subquiver on $S_{\widetilde d}$.
Claims~\ref{c:restriction-zero-left}--\ref{c:restriction-zero-right} of Lemma~\ref{lem:rep-restriction-zero} give inverse correspondences
\[
 [R]\longmapsto[\operatorname{Res}_{\widetilde K_3,S_{\widetilde d}}R],\qquad
 [N]\longmapsto[\operatorname{ExtZero}_{\widetilde K_3,S_{\widetilde d}}N].
\]
Claim~\ref{c:restriction-zero-subrep} of Lemma~\ref{lem:rep-restriction-zero} preserves subrepresentations and their dimensions at every vertex of $S_{\widetilde d}$, and hence their $\theta^{\varpi}_{a,a+k}$-weights.
The restriction and zero-extension formulas in Definition~\ref{def:rep-restriction-zero} apply to families of representations and induce
\begin{equation}\label{eq:cover-finite-support-moduli}
 \cM^s_{\widetilde d}(\widetilde K_3)
 \cong\cM^s_{\widetilde d|_{S_{\widetilde d}}}(\operatorname{SubQuiv}(\widetilde K_3,S_{\widetilde d})).
\end{equation}
The quiver $\operatorname{SubQuiv}(\widetilde K_3,S_{\widetilde d})$ is finite and acyclic.

Let us identify the classes in $\mathcal S_{a,k}^\rho$ with distinct cover dimension classes $[d]_{\varpi}$.
Let $M\in\mathcal H_{a,k}$.
Claim~\ref{c:flow-choices-stable} of Lemma~\ref{lem:flow-choices} gives stability of $R_M^\rho$.
Claims~\ref{c:flow-stability-support-tree}, \ref{c:flow-stability-support-colors}, and~\ref{c:flow-stability-total-dimension} of Lemma~\ref{lem:flow-stability} give a properly colored tree support and total dimension $(a,a+k)$.
Lemma~\ref{lem:cover-tree-embedding} gives a color-preserving embedding
\[
 \iota:\widetilde Q_M^\rho\longrightarrow\widetilde K_3.
\]

Let us extend the embedded representation by zero.
Let $U_M=\iota_0(\operatorname{Vert}(\widetilde Q_M^\rho))$.
An additional arrow of $\widetilde K_3$ between vertices of $U_M$ would form a cycle with the path between its endpoints in the embedded tree.
Thus, the image is the full subquiver $\operatorname{SubQuiv}(\widetilde K_3,U_M)$.
Transport $R_M^\rho$ to a representation $N_M$ of this subquiver and let
\[
 \widehat R_M^\rho=\operatorname{ExtZero}_{\widetilde K_3,U_M}N_M,\qquad
 \widehat d_M^\rho(v)=\dim_{\C}(\widehat R_M^\rho)_v
 \quad(v\in\operatorname{Vert}(\widetilde K_3)).
\]
Claim~\ref{c:restriction-zero-subrep} of Lemma~\ref{lem:rep-restriction-zero} preserves subrepresentations and their dimension vectors.
Hence, $\widehat R_M^\rho$ is stable and $[\widehat d_M^\rho]_{\varpi}\in\mathcal D^{\mathrm{st}}_{a,k}$.

Let us pass to dimension classes.
Lemma~\ref{lem:cover-dimension-equivalence} gives a well-defined injective map
\[
 \mathcal S_{a,k}^\rho\longrightarrow\mathcal D^{\mathrm{st}}_{a,k},\qquad
 [(\widetilde{\mathbb Q}_M^\rho,d_M^\rho)]\longmapsto[\widehat d_M^\rho]_{\varpi}.
\]

Let us prove that a subdimension vector of weight zero is zero or the full dimension vector.
For $\widetilde d\in\operatorname{CovDim}_{\varpi}(a,a+k)$, let $0\le e\le\widetilde d$ be a subdimension vector.
Let
\[
 u=\sum_{i\in I}e_i,\qquad v=\sum_{j\in J}e_j,
 \qquad u\in[0,a]_{\Z},\quad v\in[0,a+k]_{\Z}.
\]
By assumption~\ref{hyp:localization-coprimality}, we have $\gcd(a,a+k)=\gcd(a,k)=1$.
Definition~\ref{def:cover-dimension-classes}(2) gives
\[
 \theta^{\varpi}_{a,a+k}(e)=0
 \ \Longrightarrow\ av=(a+k)u
 \ \Longrightarrow\ (u,v)\in\{(0,0),(a,a+k)\}.
\]
If $(u,v)=(0,0)$, nonnegativity gives $e=0$.
If $(u,v)=(a,a+k)$, then
\[
 \sum_x(\widetilde d_x-e_x)=0,
 \qquad \widetilde d_x-e_x\ge0\quad(x\in\operatorname{Vert}(\widetilde K_3)).
\]
Hence, $e=\widetilde d$.
Thus, no proper nonzero subdimension vector has weight zero.

Let us bound the Euler characteristic of each nonempty summand from below.
Let $X=\cM^s_{\widetilde d}(\widetilde K_3)$.
Suppose that $X\ne\varnothing$.
Every common divisor of the entries of $\widetilde d|_{S_{\widetilde d}}$ divides $a$ and $a+k$.
Thus, $\gcd\{\widetilde d_x:x\in S_{\widetilde d}\}=1$.
Since $\theta^{\varpi}_{a,a+k}(e)\ne0$ for every subdimension vector $e$ satisfying $0\le e\le\widetilde d$ and $e\notin\{0,\widetilde d\}$, $\theta^{\varpi}_{a,a+k}$-semistability and $\theta^{\varpi}_{a,a+k}$-stability coincide for dimension vector $\widetilde d$.
By~\eqref{eq:cover-finite-support-moduli} and~\cite[Propositions 4.3 and 5.3, Remark 5.4]{King},
$X$ is a smooth projective fine moduli space for the finite acyclic support quiver.

Let us compute $\chi(X)$.
King--Walter's theorem~\cite[Theorem 3(c)]{KingWalter} gives
\[
 H^{2j+1}(X;\Q)=0\qquad(j\in\Z_{\ge0}).
\]
Hence,
\[
 \chi(X)=\sum_{j=0}^{\dim_{\C}X}\dim_{\Q}H^{2j}(X;\Q)
 \ge\dim_{\Q}H^0(X;\Q)\ge1.
\]
This proves
\begin{equation}
 \cM^s_{\widetilde d}(\widetilde K_3)\ne\varnothing
 \quad\Longrightarrow\quad
 \chi\bigl(\cM^s_{\widetilde d}(\widetilde K_3)\bigr)\ge1.
\label{eq:positive-summands}
\end{equation}
Therefore, equations~\eqref{eq:cover-sum} and~\eqref{eq:positive-summands} give
\[
 E_k(a)\ge
 \#\mathcal D^{\mathrm{st}}_{a,k},
\]
and $\#\mathcal D^{\mathrm{st}}_{a,k}\ge\#\mathcal S_{a,k}^\rho$. Proposition~\ref{prop:cover-count} bounds $\#\mathcal S_{a,k}^\rho$ from below.
\end{proof}

\Needspace{6\baselineskip}
\begin{theorem}\label{thm:central-lower}
We have the following.
\begin{enumerate}
\item\label{c:central-lower-finite} There are constants $c_*\in(0,1]$ and $q_*\in(0,1]$ with the following property.
For each $a\in\Z$ and $k\in\Z$ satisfying
\begin{hypotheses}
\item\label{hyp:lower-dimensions} $k\ge2$ and $a\ge100000k$,
\item\label{hyp:lower-coprimality} $\gcd(a,k)=1$,
\end{hypotheses}
we have
\begin{equation}\label{finite}
 E_k(a)\ge c_*\beta^a\left(\frac{k}{a}\right)^{\frac{5}{2}}q_*^{k-1}.
\end{equation}
\item\label{c:central-lower-limit} Let $((a_n,k_n))_{n\ge1}$ be a sequence in $\Z_{>0}^2$ satisfying
\[
 \begin{gathered}
 k_n\ge2,\qquad \gcd(a_n,k_n)=1\qquad(n\ge1),\\
 \lim_{n\to\infty}a_n=\infty,\qquad
 \lim_{n\to\infty}\frac{k_n}{a_n}=0.
 \end{gathered}
\]
Then, $E_{k_n}(a_n)>0$ for all sufficiently large $n$, and
\begin{equation}\label{entropy}
 \liminf_{n\to\infty}\frac{\log E_{k_n}(a_n)}{a_n}\ge\log\beta.
\end{equation}
\end{enumerate}
\end{theorem}
\begin{proof}
\prf{Claim~\ref{c:central-lower-finite}} Let us obtain the finite lower bound with constants independent of the dimensions.
For each $a\in\Z_{>0}$ and $k\in\Z_{>0}$ satisfying $k\ge2$ and $a\ge100000k$, fix a family of representatives $\rho=(\rho_M)_{M\in\mathcal H_{a,k}}$.

Let us apply the counting and localization bounds.
Choose a $(P_j,Q_j^{=},256)$-admissible pair $(A,B)$ by Lemma~\ref{lem:slot-bounds}.
The inequality $a\ge100000k$ verifies the dimension hypothesis of Proposition~\ref{prop:cover-count}.
The additional equality $\gcd(a,k)=1$ verifies the coprimality hypothesis of
Claim~\ref{c:cover-localization-bound} of Lemma~\ref{lem:cover-localization}.
Set $N=256$, $\lambda=\frac1{100}$, $x=\frac ak$, and $\gamma=1-\frac Nx$.
These results give
\[
 E_k(a)\ge\#\mathcal S_{a,k}^{\rho}
 \ge\beta^a x^{-\frac52}\frac{A^2}{8B\lambda N^2}
 \left(\frac{6AB\lambda\gamma^2}{\beta}\right)^{k-1}.
\]
Let
\[
 \gamma_0=1-\frac{N}{100000},\qquad
 c_*=\min\left\{1,\frac{A^2}{8B\lambda N^2}\right\},\qquad
 q_*=\min\left\{1,\frac{6AB\lambda\gamma_0^2}{\beta}\right\}.
\]
These constants are positive. Assumption~\ref{hyp:lower-dimensions}
gives $x=\frac ak\ge100000$. Hence, $\gamma\ge\gamma_0$. Substituting $x=\frac{a}{k}$ and $c=k-1$ in~\eqref{count}
proves~\eqref{finite}.

\prf{Claim~\ref{c:central-lower-limit}} Let us apply the finite bound to $(a_n,k_n)$.
Since $\lim_{n\to\infty}\frac{k_n}{a_n}=0$, choose $n_0\in\Z_{>0}$ such that
\[
 a_n\ge100000k_n\qquad(n\ge n_0).
\]
Together with $k_n\ge2$ and $\gcd(a_n,k_n)=1$, this verifies the hypotheses of Claim~\ref{c:central-lower-finite}.
Taking logarithms in~\eqref{finite} and dividing by $a_n$ gives
\begin{equation}\label{eq:central-lower-log-bound}
 \frac{\log E_{k_n}(a_n)}{a_n}\ge\log\beta+
 \frac{\log c_*-\frac52\log\left(\frac{a_n}{k_n}\right)+(k_n-1)\log q_*}{a_n}
 \qquad(n\ge n_0).
\end{equation}

Let us prove that each error term tends to zero as $n\to\infty$.
The constants $c_*$ and $q_*$ do not depend on $n$.
For $n\ge n_0$, we have
\[
 \begin{aligned}
 \lim_{n\to\infty}\frac{\log c_*}{a_n}&=0,\\
 0\le\frac{\log(a_n/k_n)}{a_n}&\le\frac{\log a_n}{a_n}\longrightarrow0,\\
 \left|\frac{(k_n-1)\log q_*}{a_n}\right|
 &\le\frac{k_n}{a_n}|\log q_*|\longrightarrow0.
 \end{aligned}
\]
Taking the lower limit as $n\to\infty$ in~\eqref{eq:central-lower-log-bound} proves~\eqref{entropy}.
\end{proof}

%% file: sections/central-upper-estimate-by-deletion.tex
\section{A central upper estimate by deletion}\label{sec:central-upper-deletion}
We obtain an upper bound for $\frac1a\log E_r(a)$ whose error tends to zero when $r=o(a)$.

\begin{definition}\label{sync:central-upper-deletion:4}
Let $W$ be a run. We call $W$ \textit{binary} if each of its simple runs
has vertical length one or two.
\end{definition}

\begin{lemma}\label{lem:three-delete}
Let $U$ be a nonempty binary run with $\nu$ simple runs.
Let $D$ be the number of simple runs of vertical length two.
For $m=3$, we have
\begin{equation}\label{eq:three-delete-bound}
 C_0(U)\le \binom{\nu+D}{D}64^D C_0(R_1^{\,\nu+2+8D}).
\end{equation}
\end{lemma}
\begin{proof}
Let us delete one vertical edge from every simple run of vertical length two.
We apply Claim~\ref{c:composition-bound} of Proposition~\ref{prop:composition} with $m=3$, $M=2$,
$K=A=0$, and $P=4D$.
Choose $u=v=1$ and $W_{\mathrm L}=W_{\mathrm R}=R_1^{4D}$.
Let $r_1<\cdots<r_D$ enumerate $\{j\in[\nu]_{\Z}:|\operatorname{RunFactor}_{j,1}(U)|_N=2\}$. Let $U^{(0)}=U$.
For $d\in[D]_{\Z}$, let $A_d$ be the $r_d$th simple run of $U^{(d-1)}$.
Let $v_d=\operatorname{End}_{\mathrm{NE}}(A_d)$ and define
\[
 U^{(d)}=\operatorname{Del}(U^{(d-1)},v_d).
\]
For $d\in[D]_{\Z}$, we have $v_d\in N(U^{(d-1)})$.
All intermediate simple runs have vertical length one or two.
Every added outer simple run has vertical length one.
The run $U^{(D)}$ is a translated copy of $R_1^\nu$.
The target is $R_1^{\nu+2+8D}$.
Substituting $\operatorname{Fib}_3(\nu,D)=\binom{\nu+D}{D}64^D$ into~\eqref{eq:composition} gives~\eqref{eq:three-delete-bound}.
\end{proof}

\begin{proposition}[Near-central finite comparison]\label{prop:central-delete-upper}
Let $a$ and $r$ be positive integers satisfying the following conditions.
\begin{hypotheses}
\item\label{hyp:central-coprimality} $\gcd(a,r)=1$.
\item\label{hyp:central-range} We have
\[
 2\le r,\qquad 2r\le a.
\]
\end{hypotheses}
Let $A=a+8r-6$.
Then, we have the following.

\begin{subequations}\label{eq:central-delete-upper}
\begin{enumerate}
\item\label{c:central-delete-comparison} We have
\begin{equation}
 E_r(a)\le \binom{a+r}{r-1}64^{r-1}E_1(A).
\end{equation}
\item\label{c:central-delete-exponential} We have
\begin{equation}
 \binom{a+r}{r-1}64^{r-1}E_1(A)
 \le \frac{4096}{729}\binom{a+r}{r-1}64^{r-1}\beta^{a+8r-6}.
\end{equation}
\end{enumerate}
\end{subequations}
\end{proposition}
\begin{proof}
\prf{Claim~\ref{c:central-delete-comparison}} Let us transform $\mathcal U(a,a+r)$ into a translated copy of $R_1^{a+1}$ by vertical-edge deletions.
By assumption~\ref{hyp:central-coprimality}, the pair $(a,a+r)$
is coprime. Equation~\eqref{eq:frame-word} defines the run $\mathcal{U}(a,a+r)$.
Equation~\eqref{eq:ambient} gives $\mu_1(a,a+r)\ge a+1$ and
$\Gamma(a,a+r)=\frac{a+r}{a}-\frac{1}{a\,\mu_1(a,a+r)}$.
By assumption~\ref{hyp:central-range}, we have $1<\Gamma(a,a+r)<2$.
Let $W=\operatorname{RunFactor}_{1,a}(\mathcal U(a,a+r))$.
The run $\mathcal{U}(a,a+r)=W\ccn R_1$ has $\nu=a+1$ simple runs
and $a+r$ vertical edges. Exactly $D=r-1$ simple runs have vertical length two.

Let us obtain $C_0(R_1^{A+1})=3E_1(A)$.
We apply Lemma~\ref{lem:three-delete} and~\eqref{eq:dictionary}.
Definition~\ref{def:frame} gives $\mu_1(A,A+1)=A+1$ and $\mu_2(A,A+1)=A+2$. Each of the first $A$ simple runs of $\mathcal{U}(A,A+1)$ has vertical length one. Thus, $\mathcal{U}(A,A+1)=R_1^{A+1}$. Since $A+1=\nu+2+8D$, applying~\eqref{eq:dictionary} to $(a,a+r)$ and $(A,A+1)$ in~\eqref{eq:three-delete-bound} and dividing by three gives Claim~\ref{c:central-delete-comparison}.

\prf{Claim~\ref{c:central-delete-exponential}} Let us bound the adjacent Euler characteristic exponentially.
By~\cite[Theorem 6.6]{Weist}, we have
\begin{equation}\label{eq:adjacent-three}
 E_1(A)=\frac{3}{(A+1)(2A+3)}\binom{4A+6}{A}.
\end{equation}
The binomial theorem evaluated at $\frac13$ gives
\[
 \binom{4A+6}{A}3^{-A}
 \le\sum_{j=0}^{4A+6}\binom{4A+6}{j}3^{-j}
 =\left(\frac43\right)^{4A+6}.
\]
Since $\frac{3}{(A+1)(2A+3)}\le1$, we obtain
\[
 E_1(A)\le3^A\left(\frac43\right)^{4A+6}
 =\frac{4096}{729}\left(\frac{256}{27}\right)^A
 =\frac{4096}{729}\beta^A.
\]
Multiplication by $\binom{a+r}{r-1}64^{r-1}$ gives Claim~\ref{c:central-delete-exponential}.
\end{proof}

\begin{corollary}[Upper limit at slope one]\label{cor:central-delete-limsup}
Let $((a_n,r_n))_{n\ge1}$ be a sequence in $\Z_{>0}^2$ satisfying the following conditions.
\begin{hypotheses}
\item\label{hyp:central-sequence-coprimality} $\gcd(a_n,r_n)=1$ for every $n\ge1$.
\item\label{hyp:central-sequence-size} $\displaystyle\lim_{n\to\infty}a_n=\infty$.
\item\label{hyp:central-sequence-ratio} $r_n\ge2$ for every $n\ge1$, and
$\displaystyle\lim_{n\to\infty}\frac{r_n}{a_n}=0$.
\end{hypotheses}
Then, $E_{r_n}(a_n)>0$ for all sufficiently large $n$, and
\[
 \limsup_{n\to\infty}\frac{\log E_{r_n}(a_n)}{a_n}\le\log\beta.
\]
\end{corollary}
\begin{proof}
Let us verify the hypotheses of Proposition~\ref{prop:central-delete-upper}.
Assumption~\ref{hyp:central-sequence-ratio} gives an integer $n_0\ge1$ such that
\[
 2r_n\le a_n\qquad(n\ge n_0).
\]
Together with assumption~\ref{hyp:central-sequence-coprimality}, this proves the hypotheses of that proposition for $(a,r)=(a_n,r_n)$ and $n\ge n_0$.

Let us estimate the binomial factor.
For $n\ge1$, put
\[
 \nu_n=a_n+1,\qquad D_n=r_n-1,\qquad x_n=\frac{D_n}{\nu_n}>0.
\]
For $t>0$, define $h(t)=(1+t)\log(1+t)-t\log t$, and set $h(0)=0$.
The limits $\lim_{t\downarrow0}t\log t=0$ and
$\lim_{t\downarrow0}(1+t)\log(1+t)=0$ give
\[
 \lim_{t\downarrow0}h(t)=h(0)=0.
\]
The binomial theorem gives
\[
 \binom{\nu_n+D_n}{D_n}x_n^{D_n}
 \le(1+x_n)^{\nu_n+D_n}.
\]
Taking logarithms and using $D_n=\nu_n x_n$, we obtain
\[
 0\le\frac1{a_n}\log\binom{\nu_n+D_n}{D_n}
 \le\frac{\nu_n}{a_n}h(x_n).
\]

Let us prove that $\frac{\nu_n}{a_n}h(x_n)\to0$ as $n\to\infty$.
Assumptions~\ref{hyp:central-sequence-size} and~\ref{hyp:central-sequence-ratio} give
\[
 \lim_{n\to\infty}x_n=0,\qquad
 \lim_{n\to\infty}\frac{\nu_n}{a_n}=1.
\]
Continuity of $h$ at $0$ therefore gives
\[
 \lim_{n\to\infty}\frac{\nu_n}{a_n}h(x_n)=0.
\]

Let us take logarithms of the finite comparison.
Equation~\eqref{eq:dictionary} and Lemma~\ref{lem:word-bounds} give
$E_{r_n}(a_n)>0$ for $n\ge n_0$.
Claims~\ref{c:central-delete-comparison} and~\ref{c:central-delete-exponential} of
Proposition~\ref{prop:central-delete-upper} give
\[
 \begin{aligned}
 \frac{\log E_{r_n}(a_n)}{a_n}-\log\beta
 &\le\frac{\log\left(\frac{4096}{729}\right)}{a_n}
       +\frac{\nu_n}{a_n}h(x_n)\\
 &\quad+\frac{r_n-1}{a_n}\log64
       +\frac{8r_n-6}{a_n}\log\beta
       \qquad(n\ge n_0).
 \end{aligned}
\]

Let us compute the upper limit as $n\to\infty$.
The binomial estimate and assumptions~\ref{hyp:central-sequence-size}
and~\ref{hyp:central-sequence-ratio} give
\[
 \begin{gathered}
 \lim_{n\to\infty}\frac{\log\left(\frac{4096}{729}\right)}{a_n}=0,
 \qquad \lim_{n\to\infty}\frac{\nu_n}{a_n}h(x_n)=0,\\
 \lim_{n\to\infty}\frac{r_n-1}{a_n}\log64=0,
 \qquad \lim_{n\to\infty}\frac{8r_n-6}{a_n}\log\beta=0.
 \end{gathered}
\]
Thus, $\displaystyle\limsup_{n\to\infty}\frac{\log E_{r_n}(a_n)}{a_n}\le\log\beta$.
\end{proof}

%% file: sections/three-arrow-case.tex
\section{The three-arrow case}\label{sec:three-arrow}
In this section all colored gradings use $m=3$.

\subsection{The limit at slope one}\label{sec:central}

\begin{corollary}\label{cor:central}
Let $((a_n,b_n))_{n\ge1}$ be a sequence in $\cA_3$ satisfying
\[
 \lim_{n\to\infty}a_n=\infty,\qquad
 \lim_{n\to\infty}\frac{b_n}{a_n}=1.
\]
Then,
\[
 \lim_{n\to\infty}\frac{\log E_3(a_n,b_n)}{a_n}
 =\log\left(\frac{256}{27}\right).
\]
\end{corollary}
\begin{proof}
Let us reduce the dimension pair to $(d_n,d_n+k_n)$.
For $n\ge1$, define
\[
 d_n=\min(a_n,b_n),\qquad k_n=|b_n-a_n|.
\]
Since $(a_n,b_n)\in\cA_3$, we have $\gcd(a_n,b_n)=1$.
If $a_n=b_n$, then $a_n=b_n=1$.
The condition $a_n\to\infty$ therefore gives $n_0\in\Z_{>0}$ such that $k_n\ge1$ for $n\ge n_0$.
For $n\ge n_0$, we have
\[
 \gcd(d_n,k_n)=\gcd(\min(a_n,b_n),|b_n-a_n|)=\gcd(a_n,b_n)=1.
\]
Claim~\ref{c:moduli-symmetries-identities} of Lemma~\ref{lem:moduli-symmetries} gives
\[
 E_3(a_n,b_n)=E_{k_n}(d_n).
\]

Let us compute the limits of the new dimensions.
The hypothesis $\frac{b_n}{a_n}\to1$ gives
\[
 \frac{d_n}{a_n}=\min\left(1,\frac{b_n}{a_n}\right)\longrightarrow1,
 \qquad
 \frac{k_n}{d_n}
 =\frac{\left|\frac{b_n}{a_n}-1\right|}{\min\left(1,\frac{b_n}{a_n}\right)}
 \longrightarrow0.
\]
In particular, $d_n\to\infty$.
Define the two index sets
\[
 J_1=\{n\ge n_0:k_n=1\},\qquad
 J_2=\{n\ge n_0:k_n\ge2\}.
\]

First, let us compute the limit along $J_1$ when $J_1$ is infinite.
Equation~\eqref{eq:adjacent-three} and Stirling's formula give
\[
 \log E_1(d_n)=d_n\log\beta-\frac52\log d_n+O(1)
 \qquad(n\to\infty,\ n\in J_1).
\]
Dividing by $d_n$ gives
\[
 \lim_{\substack{n\to\infty\\n\in J_1}}
 \frac{\log E_{k_n}(d_n)}{d_n}=\log\beta.
\]

Second, let us compute the limit along $J_2$ when $J_2$ is infinite.
Enumerate $J_2=\{n_1<n_2<\cdots\}$.
The sequence $((d_{n_j},k_{n_j}))_{j\ge1}$ satisfies
\[
 \begin{gathered}
 k_{n_j}\ge2,\qquad \gcd(d_{n_j},k_{n_j})=1\qquad(j\ge1),\\
 \lim_{j\to\infty}d_{n_j}=\infty,\qquad
 \lim_{j\to\infty}\frac{k_{n_j}}{d_{n_j}}=0.
 \end{gathered}
\]
Claim~\ref{c:central-lower-limit} of Theorem~\ref{thm:central-lower} and
Corollary~\ref{cor:central-delete-limsup} therefore give
\[
 \lim_{j\to\infty}\frac{\log E_{k_{n_j}}(d_{n_j})}{d_{n_j}}=\log\beta.
\]

Let us combine the two cases.
We have $J_1\sqcup J_2=\{n\ge n_0\}$.
For each infinite $J_i$, $i\in[2]_{\Z}$, the limit of $\frac{\log E_{k_n}(d_n)}{d_n}$ along $n\in J_i$ is $\log\beta$. Each finite $J_i$ contributes only finitely many terms.
Thus,
\[
 \frac{\log E_3(a_n,b_n)}{a_n}
 =\frac{d_n}{a_n}\frac{\log E_{k_n}(d_n)}{d_n}
 \longrightarrow1\cdot\log\beta=\log\left(\frac{256}{27}\right).
\]
\end{proof}

\subsection{Adjacent exchanges}
\begin{definition}[Integer lists and balance operations]\label{def:balance-operations}
Let $\operatorname{IntList}=\bigsqcup_{n\in\Z_{\ge0}}\Z^n$.
\begin{enumerate}[font=\itshape]
\item For $L=(a_1,\ldots,a_n)\in\operatorname{IntList}$, let $|L|=n$.
For $M=(b_1,\ldots,b_m)\in\operatorname{IntList}$, define
\[
 L\ccn M=(a_1,\ldots,a_n,b_1,\ldots,b_m).
\]
For $a\in\Z$ and $s\in\Z_{\ge0}$, define
\[
 \operatorname{Repeat}(a,s)=\underbrace{(a,\ldots,a)}_{s\text{ entries}},\qquad
 \operatorname{Repeat}(a,0)=().
\]
\item We define the set of \textit{balance states} by
$\operatorname{BalState}=\Z\times\operatorname{IntList}$.
\item\label{d:horizontal-balance-step} For $s\in\Z_{\ge0}$, define
\[
 \operatorname{HStep}_s:\operatorname{BalState}\longrightarrow\operatorname{BalState},\qquad
 \operatorname{HStep}_s(b,L)=(b+1,\operatorname{Repeat}(b+1,s)\ccn L).
\]
\item\label{d:vertical-balance-step} For $h\in\Z_{\ge0}$, define
\[
 \mathcal V_h=\{(b,(a)\ccn L_0)\in\operatorname{BalState}:a\le b-h\}
\]
and
\[
 \operatorname{VStep}_h:\mathcal V_h\longrightarrow\operatorname{BalState},\qquad
 \operatorname{VStep}_h(b,(a)\ccn L_0)=(b-h,L_0).
\]
We call $\operatorname{HStep}_s$ and $\operatorname{VStep}_h$ the horizontal and vertical balance operations.
\item\label{d:balance-list-order} For $\delta\in\Z_{\ge0}$ and integer lists $L$ and $L'$, define
\begin{equation}\label{eq:exchange-list-comparison}
 L'\preceq_\delta L
 \ \Longleftrightarrow\
 |L'|=|L|\ \text{ and }\ L'_i\le L_i+\delta\quad(i\in[|L|]_{\Z}).
\end{equation}
\item\label{d:boundary-balance-list} Let $Y$ be a tNE path. Let $g\in\operatorname{CtdG}(Y)$. Let $r\in[0,|Y|]_{\Z}$.
List the horizontal edges of $Y[0,r]$ as $e_1,\ldots,e_h$ in path order.
For $i\in[h]_{\Z}$, let
\[
 b_i=\operatorname{eBal}_{Y,g}(e_i),\qquad
 d_i=\operatorname{Debt}_{Y,g}(e_i,r).
\]
Define the \textit{boundary balance list} by
\[
 \operatorname{BalList}_{Y,g}(r)
 =\operatorname{Repeat}(b_h,d_h)\ccn\cdots\ccn\operatorname{Repeat}(b_1,d_1).
\]
For $h=0$, this list is $()$.
\end{enumerate}
\end{definition}

\begin{lemma}[Comparison of balance operations]\label{lem:balance-list-comparison}
Let $\nu\in\Z_{\ge0}$ and $\delta\in\Z_{\ge0}$.
For $t\in[\nu]_{\Z}$, let $\varepsilon_t\in\{E,N\}$ and $q_t\in\Z_{\ge0}$.
Assume the following.
\begin{hypotheses}
\item\label{hyp:balance-initial} $z_0=(b_0,L_0)\in\operatorname{BalState}$ and $L'_0\preceq_\delta L_0$.
Let $z'_0=(b_0+\delta,L'_0)$.
\item\label{hyp:balance-sequence} The states $z_t=(b_t,L_t)$ for $t\in[\nu]_{\Z}$ satisfy
\[
 z_t=\begin{cases}
 \operatorname{HStep}_{q_t}(z_{t-1}),&\varepsilon_t=E,\\
 \operatorname{VStep}_{q_t}(z_{t-1}),&\varepsilon_t=N.
 \end{cases}
\]
For every $t\in[\nu]_{\Z}$ with $\varepsilon_t=N$, we have $z_{t-1}\in\mathcal V_{q_t}$.
\end{hypotheses}
Starting from $z'_0$, the states $z'_t=(b'_t,L'_t)$ are well-defined by
\[
 z'_t=\begin{cases}
 \operatorname{HStep}_{q_t}(z'_{t-1}),&\varepsilon_t=E,\\
 \operatorname{VStep}_{q_t}(z'_{t-1}),&\varepsilon_t=N
 \end{cases}
 \qquad(t\in[\nu]_{\Z}).
\]
For every $t\in[0,\nu]_{\Z}$,
\[
 b'_t=b_t+\delta,\qquad L'_t\preceq_\delta L_t.
\]
In particular, $z'_{t-1}\in\mathcal V_{q_t}$ whenever $\varepsilon_t=N$.
\end{lemma}
\begin{proof}
Let us prove the assertions by induction on $t$.
The case $t=0$ is assumption~\ref{hyp:balance-initial}.
Suppose $b'_{t-1}=b_{t-1}+\delta$ and $L'_{t-1}\preceq_\delta L_{t-1}$.
If $\varepsilon_t=E$, the horizontal operation increases each balance by one and prepends $q_t$ entries.
The new entries are $b_{t-1}+1$ and $b_{t-1}+1+\delta$, respectively.
Thus, $b'_t=b_t+\delta$ and $L'_t\preceq_\delta L_t$.
If $\varepsilon_t=N$, assumption~\ref{hyp:balance-sequence} gives
$L_{t-1}=(a)\ccn M$ and $a\le b_{t-1}-q_t$.
The equal list lengths give $L'_{t-1}=(a')\ccn M'$.
By the induction hypothesis,
\[
 a'\le a+\delta\le b_{t-1}-q_t+\delta=b'_{t-1}-q_t.
\]
Hence, $z'_{t-1}\in\mathcal V_{q_t}$.
Both vertical operations delete the first entry and subtract $q_t$ from the balance.
Thus, $b'_t=b_t+\delta$ and $L'_t\preceq_\delta L_t$.
\end{proof}

\Needspace{19\baselineskip}
\begin{definition}[Adjacent exchange in a shell]\label{def:shell-adjacent-exchange}
Define $\operatorname{AdjShellExch}$ to be the set of tuples
\[
 \xi=(P,Q,o,u,v,A,B)
\]
such that $P$ and $Q$ are NE paths, $o\in\Z^2$, $u\in[2]_{\Z}$, $v\in[2]_{\Z}$, and
$(A,B)\in\{(E\ccn N,N\ccn E),(N\ccn E,E\ccn N)\}$.
We call an element of this set an \textit{adjacent exchange in a shell}.
Let $\xi=(P,Q,o,u,v,A,B)\in\operatorname{AdjShellExch}$.
\begin{enumerate}[font=\itshape]
\item\label{d:adj-shell-exch-paths} Define
\[
\begin{aligned}
 X_\xi&=\operatorname{Trans}(E\ccn P\ccn A\ccn Q,o),& X'_\xi&=\operatorname{Trans}(E\ccn P\ccn B\ccn Q,o),\\
 Y_\xi&=\operatorname{Shell}(X_\xi;u,v,0),&Y'_\xi&=\operatorname{Shell}(X'_\xi;u,v,0).
\end{aligned}
\]
\item\label{d:adj-shell-exch-edges} Let $Y=Y_\xi$ and $Y'=Y'_\xi$. Let $(e_i)_{i\in[|Y|_E]_{\Z}}$ and $(n_j)_{j\in[|Y|_N]_{\Z}}$
be the horizontal and vertical edge lists of $Y$, respectively, in path order.
Let $(e'_i)_{i\in[|Y'|_E]_{\Z}}$ and $(n'_j)_{j\in[|Y'|_N]_{\Z}}$
be the horizontal and vertical edge lists of $Y'$, respectively, in path order.
For each $j\in[|Y'_\xi|_N]_{\Z}$, let $p_j\in[|Y'_\xi|]_{\Z}$ be the unique index satisfying $(Y'_\xi)_{p_j}=n'_j$.
\item\label{d:adj-shell-edge-map} Define the bijection
\[
 \tau_\xi:\operatorname{NE}(Y_\xi)\longrightarrow\operatorname{NE}(Y'_\xi),\qquad
 \tau_\xi(e_i)=e'_i,\quad \tau_\xi(n_j)=n'_j.
\]
\item\label{d:adj-shell-exchanged-edges} Define
\[
 r_\xi=u+2+|P|,\qquad
 i_\xi=|P|_E+3,\qquad j_\xi=u+|P|_N+1,\qquad
 x_\xi=e_{i_\xi},\quad w_\xi=n_{j_\xi}.
\]
\item\label{d:adj-shell-gradings} Define
\[
 \mathcal G_\xi=
 \{f\in\operatorname{CtdG}_3(X_\xi):\operatorname{Surp}(X_\xi,f)\ge0\}.
\]
For $f\in\mathcal G_\xi$, define
\[
 g_{\xi,f}=\operatorname{GShell}_f(X_\xi;u,v,0),\qquad
 \overline g_{\xi,f}=g_{\xi,f}\circ\tau_\xi^{-1}.
\]
\item\label{d:adj-shell-counter} Let $f\in\mathcal G_\xi$ and let $s=g_{\xi,f}(x_\xi)$.
For $s>0$, define integers $q_i$ for $i\in[r_\xi+2,|Y'_\xi|]_{\Z}$ by
\[
 q_{r_\xi+2}=s-1,\qquad
 q_i=\begin{cases}
 q_{i-1}+\overline g_{\xi,f}((Y'_\xi)_i),&(Y'_\xi)_i\in E(Y'_\xi),\\
 \max(q_{i-1}-1,0),&(Y'_\xi)_i\in N(Y'_\xi)
 \end{cases}
 \quad(i>r_\xi+2).
\]
Define
\[
 \mathcal J_{\xi,f}=\{j\in[|Y'_\xi|_N]_{\Z}:p_j>r_\xi+2,\ q_{p_j-1}=0\}.
\]
For $s=0$, let $\mathcal J_{\xi,f}=\varnothing$.
\item\label{d:adj-shell-repair-index} Define the \textit{vertical repair index} by
\[
 J_\xi(f)=\begin{cases}
 \min\mathcal J_{\xi,f},&A=E\ccn N,\ s>0,\ \mathcal J_{\xi,f}\ne\varnothing,\\
 j_\xi,&\text{otherwise}.
 \end{cases}
\]
\item\label{d:adj-shell-repair-grading} Let $J=J_\xi(f)$ and define
\[
\begin{gathered}
 h=g_{\xi,f}(w_\xi),\qquad h_J=\overline g_{\xi,f}(n'_J),\\
 \vartheta=\operatorname{bBal}_{Y_\xi,g_{\xi,f}}(r_\xi)+1,\qquad
 b_J=\operatorname{bBal}_{Y'_\xi,\overline g_{\xi,f}}(p_J-1).
\end{gathered}
\]
Define the \textit{grade reduction} by
\[
 D_\xi(f)=\begin{cases}
 h-\max(h-1,0),&A=E\ccn N,\ s=0,\\
 0,&A=N\ccn E,\ s=0,\\
 h,&A=N\ccn E,\ s>0,\\
 h_J-\min\{h_J,\max(b_J-\vartheta,0)\},&A=E\ccn N,\ s>0.
 \end{cases}
\]
Define the \textit{repaired grading} $\operatorname{RepGr}_\xi(f)$ on $Y'_\xi$ by
\[
 \operatorname{RepGr}_\xi(f)(e'_i)=g_{\xi,f}(e_i),\qquad
 \operatorname{RepGr}_\xi(f)(n'_j)=
 \begin{cases}
 g_{\xi,f}(n_j)-D_\xi(f),&j=J_\xi(f),\\
 g_{\xi,f}(n_j),&j\ne J_\xi(f).
 \end{cases}
\]
\end{enumerate}
\end{definition}

\begin{lemma}[The vertical repair index]\label{lem:exchange-repair-index}
Let $\xi=(P,Q,o,u,v,A,B)\in\operatorname{AdjShellExch}$.
\begin{enumerate}
\item\label{c:repair-index-exists} If $f\in\mathcal G_\xi$, $A=E\ccn N$, and $g_{\xi,f}(x_\xi)>0$, then $\mathcal J_{\xi,f}\ne\varnothing$.
The edge $n'_{J_\xi(f)}$ lies in the copy of $X'_\xi$ in $Y'_\xi$.
\item\label{c:repair-index-horizontal} If $f\in\mathcal G_\xi$, $\widetilde f\in\mathcal G_\xi$, and
\[
 g_{\xi,f}(e_i)=g_{\xi,\widetilde f}(e_i)
 \qquad(i\in[|Y_\xi|_E]_{\Z}),
\]
then $J_\xi(f)=J_\xi(\widetilde f)$.
\end{enumerate}
\end{lemma}
\begin{proof}
\prf{Claim~\ref{c:repair-index-exists}} Let us prove $\mathcal J_{\xi,f}\ne\varnothing$.
Let $g=g_{\xi,f}$, $s=g(x_\xi)>0$, $r=r_\xi$, and $T=u+1+|X_\xi|$.
The shell grading $g=g_{\xi,f}$ belongs to $\operatorname{CtdG}_3(Y_\xi)$ by Claim~\ref{c:shell-connected} of Lemma~\ref{lem:shell}.
Since $f$ is connected, the last vertical edge of $X_\xi$ has h-precursor $(X_\xi)_1$.
Claim~\ref{c:shell-copy-precursor} of Lemma~\ref{lem:shell} and Claim~\ref{c:assigned-prefix} of Lemma~\ref{lem:precursor-prefix-counts} give
\[
 \operatorname{Debt}_{Y_\xi,g}((Y_\xi)_{u+2},r)\ge1.
\]
Since $g$ is connected, Claim~\ref{c:assigned-prefix} of Lemma~\ref{lem:precursor-prefix-counts} gives $\operatorname{Debt}_{Y_\xi,g}((Y_\xi)_1,r)\ge1$.
Let
\[
 \operatorname{BalList}_{Y_\xi,g}(r)=(a)\ccn L_0.
\]
The entry contributed by $(Y_\xi)_1$ follows an entry contributed by $(Y_\xi)_{u+2}$.
Hence, $L_0$ contains the entry contributed by $(Y_\xi)_1$, and $|L_0|\ge1$.
The edge $w_\xi$ immediately follows $x_\xi$, and $g(x_\xi)=s>0$.
Definition~\ref{def:stack}(1) gives $P_{Y_\xi,g}(w_\xi)=x_\xi$.
The node balance condition~\eqref{eq:connected-edge-balance} gives
\[
 \operatorname{eBal}_{Y_\xi,g}(x_\xi)
 \le\operatorname{eBal}_{Y_\xi,g}(w_\xi)
 =\operatorname{eBal}_{Y_\xi,g}(x_\xi)-g(w_\xi).
\]
Hence, $g(w_\xi)=0$.
Definition~\ref{def:shell-adjacent-exchange}(\ref{d:adj-shell-counter}) gives $q_t\ge0$ for $t\in[r+2,|Y'_\xi|]_{\Z}$.
By Definition~\ref{def:balance-operations}(6),
\[
 |\operatorname{BalList}_{Y_\xi,g}(r+2)|
 =s+|L_0|=q_{r+2}+1+|L_0|.
\]

Let us prove that $q_{p_j-1}=0$ for some $p_j\in[r+3,T]_{\Z}$.
Suppose, for a contradiction, that $q_{p_j-1}>0$ for every vertical position $p_j\in[r+3,T]_{\Z}$.
For $t\in[r+3,T]_{\Z}$, we have $\tau_\xi((Y_\xi)_t)=(Y'_\xi)_t$.
At a horizontal edge, both $q_t$ and $|\operatorname{BalList}_{Y_\xi,g}(t)|$ increase by $g((Y_\xi)_t)$.
At a vertical edge at position $t=p_j$, the inequality $q_{t-1}>0$ gives $q_t=q_{t-1}-1$, and Definition~\ref{def:stack}(1) decreases $|\operatorname{BalList}_{Y_\xi,g}(t)|$ by $1$.
Induction on $t$ gives
\[
 |\operatorname{BalList}_{Y_\xi,g}(t)|=q_t+1+|L_0|
 \qquad(t\in[r+2,T]_{\Z}).
\]

Let us apply the shell boundary identity at $T$.
Claim~\ref{c:shell-boundary-debts} of Lemma~\ref{lem:shell}, the equality $\operatorname{eBal}_{Y_\xi,g}((Y_\xi)_1)=1$, and Definition~\ref{def:balance-operations}(6) give $\operatorname{BalList}_{Y_\xi,g}(T)=(1)$.
Therefore,
\[
 1=|\operatorname{BalList}_{Y_\xi,g}(T)|
 =q_T+1+|L_0|\ge2,
\]
a contradiction.
Thus, some $j\in\mathcal J_{\xi,f}$ satisfies $p_j\le T$.
Definition~\ref{def:shell-adjacent-exchange}(\ref{d:adj-shell-repair-index}) gives $r+2<p_{J_\xi(f)}\le T$.
Hence, $n'_{J_\xi(f)}$ lies in the copy of $X'_\xi$ in $Y'_\xi$.

\prf{Claim~\ref{c:repair-index-horizontal}} Let us check which grades determine the repair index.
The common horizontal grades determine $s$ and every horizontal increment in the recursion for $q_i$.
The vertical recursion $q_i=\max(q_{i-1}-1,0)$ does not use a vertical grade.
Hence, $\mathcal J_{\xi,f}=\mathcal J_{\xi,\widetilde f}$ and $J_\xi(f)=J_\xi(\widetilde f)$.
\end{proof}

\begin{lemma}[The repaired grading]\label{lem:uncolored-adjacent-exchange}
Let $\xi\in\operatorname{AdjShellExch}$ and $f\in\mathcal G_\xi$.
Let
\[
 g=g_{\xi,f},\qquad g'=\operatorname{RepGr}_\xi(f),\qquad
 j=J_\xi(f),\qquad\delta=D_\xi(f).
\]
Then,
\begin{enumerate}
\item\label{c:uncolored-exchange-grades} $g'\in\operatorname{CtdG}_3(Y'_\xi)$ and $0\le\delta\le g(n_j)\le3$.
\item\label{c:uncolored-exchange-surplus} We have
\[
 \operatorname{Surp}(Y'_\xi,g')=\operatorname{Surp}(X_\xi,f)+\delta.
\]
\end{enumerate}
\end{lemma}
\begin{proof}
\prf{Claim~\ref{c:uncolored-exchange-grades}} Let us compare the h-precursors before and after exchanging the two edges.
Let $\xi=(P,Q,o,u,v,A,B)$ and $(X,X',Y,Y')=(X_\xi,X'_\xi,Y_\xi,Y'_\xi)$.
Put $c=\operatorname{Surp}(X,f)$.
Claim~\ref{c:shell-connected} of Lemma~\ref{lem:shell} gives $g\in\operatorname{CtdG}_3(Y)$.
The formulas for $D_\xi(f)$ give $0\le\delta\le g(n_j)\le3$.
Thus, $g'\in G_3(Y')$.
For $r\in[0,|Y|]_{\Z}$, let
\[
 b=\operatorname{bBal}_{Y,g}(r),\qquad L=\operatorname{BalList}_{Y,g}(r).
\]
For $r<|Y|$, Definition~\ref{def:stack}(1) and the node balance condition in
Definition~\ref{def:connected-grading}(1) give the operations $\operatorname{HStep}_{g(Y_{r+1})}$ when $Y_{r+1}\in E(Y)$ and $\operatorname{VStep}_{g(Y_{r+1})}$ when $Y_{r+1}\in N(Y)$.
Claim~\ref{c:node-boundary-balances} of Lemma~\ref{lem:node-structure}
implies that every entry of $L$ is at most $b$.
Claims~\ref{c:shell-boundary-debts} and~\ref{c:shell-boundary-balances} of Lemma~\ref{lem:shell} and Definition~\ref{def:balance-operations}(6) give
\[
\begin{aligned}
 (\operatorname{bBal}_{Y,g}(u+1),\operatorname{BalList}_{Y,g}(u+1))&=(1,(1)),\\
 (\operatorname{bBal}_{Y,g}(u+1+|X|),\operatorname{BalList}_{Y,g}(u+1+|X|))&=(c+2,(1)).
\end{aligned}
\]
The entry $1$ in $L=(1)$ is contributed by $Y_1$.
The vertical edges of $X$ remove every entry that $X$ contributes.

Let us construct states $z'_t=(b'_t,L'_t)$ for $t\in[0,|Y'|]_{\Z}$.
Let $z'_0=(0,())$ and use the recursion
\[
 z'_t=\begin{cases}
 \operatorname{HStep}_{g'(Y'_t)}(z'_{t-1}),&Y'_t\in E(Y'),\\
 \operatorname{VStep}_{g'(Y'_t)}(z'_{t-1}),&Y'_t\in N(Y').
 \end{cases}
\]

Let us verify $z'_{t-1}\in\mathcal V_{g'(Y'_t)}$ before the exchanged pair.
For $t\in[0,r_\xi]_{\Z}$, we have
\[
 z'_t=(\operatorname{bBal}_{Y,g}(t),\operatorname{BalList}_{Y,g}(t)).
\]

Let us compare the balance states at boundaries $r_\xi$ and $r_\xi+2$.
Let $\tau=\tau_\xi$, $x=x_\xi$, and $w=w_\xi$.
Let $s=g(x)$, $h=g(w)$, and $w'=\tau(w)$.
At boundary $r_\xi$ of $Y$, let $(b,(a)\ccn L_0)$ be the balance state.
The list $L_0$ is nonempty: the first edge of $Y$ contributes its entry $1$. The first edge of the copy of $X$ has positive debt before the exchanged pair.
By the definition of $g'=\operatorname{RepGr}_\xi(f)$,
\[
 g'(\tau(e_i))=g(e_i)\quad(i\in[|Y|_E]_{\Z}),\qquad
 g'(\tau(n_i))=g(n_i)\quad(i\in[|Y|_N]_{\Z}\setminus\{J_\xi(f)\}).
\]
At boundary $r_\xi+2$, let
\[
\begin{aligned}
 (b_{\mathrm{old}},L_{\mathrm{old}})&=
 (\operatorname{bBal}_{Y,g}(r_\xi+2),\operatorname{BalList}_{Y,g}(r_\xi+2)),\\
 (b_{\mathrm{new}},L_{\mathrm{new}})&=z'_{r_\xi+2}.
\end{aligned}
\]

First, let us verify the balance condition when $A=E\ccn N$ and $s=0$.
We have $J_\xi(f)=j_\xi$, $g'(w')=\max(h-1,0)$, and $\delta=h-\max(h-1,0)$.
Validity of $\operatorname{VStep}_h$ at $w$ gives $a\le b+1-h$.
We also have $a\le b$.
Hence,
\[
 a\le b-\max(h-1,0).
\]
Thus,
\[
 L_{\mathrm{new}}=L_{\mathrm{old}},\qquad
 b_{\mathrm{new}}=b_{\mathrm{old}}+\delta,\qquad \delta\in[0,1]_{\Z}.
\]

Second, let us verify the balance condition when $A=N\ccn E$ and $s=0$.
We have $\delta=0$ and $g'(w')=h$.
The inequality $a\le b-h$ implies $a\le b+1-h$.
Thus, $(b_{\mathrm{new}},L_{\mathrm{new}})=(b_{\mathrm{old}},L_{\mathrm{old}})$.

Third, let us verify the balance condition when $A=N\ccn E$ and $s>0$.
We have $J_\xi(f)=j_\xi$, $g'(w')=0$, and $\delta=h$.
The states at boundary $r_\xi+2$ satisfy
\[
 b_{\mathrm{old}}=b-h+1,\qquad
 L_{\mathrm{old}}=\operatorname{Repeat}(b-h+1,s)\ccn L_0,
\]
\[
 b_{\mathrm{new}}=b+1,\qquad
 L_{\mathrm{new}}=\operatorname{Repeat}(b+1,s-1)\ccn(a)\ccn L_0.
\]
Since $a\le b-h$, we have $L_{\mathrm{new}}\preceq_h L_{\mathrm{old}}$ by~\eqref{eq:exchange-list-comparison}.

Fourth, let us verify the balance condition when $A=E\ccn N$ and $s>0$.
Validity of $\operatorname{VStep}_h$ at $w$ immediately after $x$ gives $b+1\le b+1-h$. Thus, $h=0$.
Let $\vartheta=b+1$.
The states at boundary $r_\xi+2$ satisfy
\[
 L_{\mathrm{old}}=\operatorname{Repeat}(\vartheta,s-1)\ccn(a)\ccn L_0,
\]
\[
 L_{\mathrm{new}}=\operatorname{Repeat}(\vartheta,s)\ccn L_0,\qquad
 b_{\mathrm{new}}=b_{\mathrm{old}}=\vartheta.
\]

Let us locate the entry whose grade requires a change in the fourth case.
Mark entry $s$ of $L_{\mathrm{new}}$.
Claim~\ref{c:repair-index-exists} of Lemma~\ref{lem:exchange-repair-index} gives
\[
 r_\xi+2<p_j\le u+1+|X|.
\]
The edge $T=\tau^{-1}(n'_j)$ removes the marked entry.
At boundary $t\in[r_\xi+2,p_j-1]_{\Z}$, the mark has position $q_t+1$, where $q_t$ is the counter in Definition~\ref{def:shell-adjacent-exchange}(\ref{d:adj-shell-counter}).

Let us prove that the operations before $T$ are valid.
For $t\in[r_\xi+2,p_j-1]_{\Z}$, we prove by induction that there is a list $H_t$ satisfying
\[
\begin{aligned}
 \operatorname{BalList}_{Y,g}(t)&=H_t\ccn(a)\ccn L_0,\\
 L'_t&=H_t\ccn(\vartheta)\ccn L_0,\\
 |H_t|&=q_t,\qquad (H_t)_i\ge\vartheta\quad(i\in[q_t]_{\Z}),\\
 b'_t&=\operatorname{bBal}_{Y,g}(t)\ge\vartheta.
\end{aligned}
\]
At $t=r_\xi+2$, take $H_t=\operatorname{Repeat}(\vartheta,s-1)$.

First, let us prove the induction step when $Y'_{t+1}\in E(Y')$ and $t+1<p_j$.
Both balances increase by one.
Every entry inserted by this horizontal edge equals $b'_t+1\ge\vartheta$.
Set $H_{t+1}=\operatorname{Repeat}(b'_t+1,g'(Y'_{t+1}))\ccn H_t$.
The identity $q_{t+1}=q_t+g'(Y'_{t+1})$ verifies $|H_{t+1}|=q_{t+1}$.
Both fixed tails, $(a)\ccn L_0$ and $(\vartheta)\ccn L_0$, are unchanged.

Second, let us prove the induction step when $Y'_{t+1}\in N(Y')$ and $t+1<p_j$.
The definition of $j=J_\xi(f)$ gives $q_t>0$.
Let $a_t$ be the common first list entry.
The operation $\operatorname{VStep}_{g(Y_{t+1})}$ on $Y$ at position $t+1$ is valid, so
\[
 b'_t-g'(Y'_{t+1})
 =\operatorname{bBal}_{Y,g}(t)-g(Y_{t+1})
 \ge a_t\ge\vartheta.
\]
The inequality $a_t\le b'_t-g'(Y'_{t+1})$ verifies the condition in Definition~\ref{def:balance-operations}(\ref{d:vertical-balance-step}) for the operation at $Y'_{t+1}$.
Since $q_t>0$, deleting the first entry of $H_t$ gives $H_{t+1}$ of length $q_{t+1}=q_t-1$.
Both fixed tails are unchanged, and $b'_{t+1}=b'_t-g'(Y'_{t+1})\ge\vartheta$.

Let us verify the vertical operation at $\tau(T)$.
Let $b_T$ be the boundary balance immediately before $T$ and let $q=g(T)$.
We have $b_T\ge\vartheta$.
At $T'=\tau(T)=n'_{J_\xi(f)}$, the definition of $D_\xi(f)$ gives
\[
 g'(T')=k:=\min\{q,b_T-\vartheta\},\qquad \delta=q-k.
\]
Then, $0\le k\le q$ and $b_T-k\ge\vartheta$.
Since $q_{p_j-1}=0$, the lists immediately before $T$ and $T'$ are $(a)\ccn L_0$ and $(\vartheta)\ccn L_0$, respectively.
Deleting their first entries gives
\[
 L'_{p_j}=\operatorname{BalList}_{Y,g}(p_j)=L_0,\qquad
 b'_{p_j}=\operatorname{bBal}_{Y,g}(p_j)+\delta.
\]

Let us verify connectedness.
Lemma~\ref{lem:balance-list-comparison} gives boundary balance $c+2+\delta$ and a list of length one at $u+1+|X|$.
The entry contributed by $Y'_1$ equals $1$ and stays last in every nonempty list.
The original lists are nonempty throughout the copy of $X$.
Equality of list lengths therefore prevents this entry from being removed.
Hence,
\[
 z'_{u+1+|X|}=(c+2+\delta,(1)).
\]
The grades on the final $R_v$ are $v-1$ on its horizontal edge, zero on its first
$v-1$ vertical edges, and two on its last edge.
After its horizontal edge the state is
\[
 (c+3+\delta,\operatorname{Repeat}(c+3+\delta,v-1)\ccn(1)).
\]
The zero-grade vertical operations are valid.
The last vertical operation is valid because $1\le c+1+\delta$.
Therefore,
\[
 z'_{|Y'|}=(c+1+\delta,()).
\]

Let us define the h-precursor from the lists $L'_t$.
Assign each list entry to the horizontal edge that inserted it.
For $t\in[|Y'|]_{\Z}$ with $Y'_t\in N(Y')$, let $Y'_{r_t}$ be the horizontal edge assigned to the first entry of $L'_{t-1}$.
The list is nonempty because $z'_{t-1}\in\mathcal V_{g'(Y'_t)}$.
The order of the list gives
\[
 r_t=\max\left\{i\in[t-1]_{\Z}:
 \begin{array}{l}
 Y'_i\in E(Y'),\\
 g'(Y'_i)>\#\{j\in[t-1]_{\Z}:Y'_j\in N(Y'),\ r_j=i\}
 \end{array}\right\}.
\]
Thus, $P_{Y',g'}(Y'_t)=Y'_{r_t}$ satisfies Definition~\ref{def:stack}(1).

Let us verify the node balance and the node at $Y'_1$.
The first entry of $L'_{t-1}$ is $\operatorname{eBal}_{Y',g'}(Y'_{r_t})$.
Membership in $\mathcal V_{g'(Y'_t)}$ gives
\[
 \operatorname{eBal}_{Y',g'}(Y'_{r_t})
 \le b'_{t-1}-g'(Y'_t)
 =\operatorname{eBal}_{Y',g'}(Y'_t)
 \qquad(Y'_t\in N(Y')).
\]
Every list entry is removed because $L'_{|Y'|}=()$.
The entry from $Y'_1$ is removed at $Y'_{|Y'|}$.
Thus,
\[
 \operatorname{Node}_{Y',g'}(Y'_1)=Y',\qquad g'\in\operatorname{CtdG}_3(Y').
\]

\prf{Claim~\ref{c:uncolored-exchange-surplus}} Let us compute the surplus change.
Only the vertical grade at $n'_j$ changes by $-\delta$. Hence,
\[
 \operatorname{Surp}(Y',g')=c+\delta\in[c,c+3]_{\Z}.
\]
\end{proof}

\begin{corollary}[The adjacent-exchange repair map]\label{cor:adjacent-repair-map}
Let $\xi\in\operatorname{AdjShellExch}$.
Then, the formula
\[
 \Theta_\xi(f)=(\operatorname{RepGr}_\xi(f),J_\xi(f),D_\xi(f))
\]
defines a map
\[
 \Theta_\xi:\mathcal G_\xi\longrightarrow
 \operatorname{CtdG}_3(Y'_\xi)\times[|Y_\xi|_N]_{\Z}\times[0,3]_{\Z}.
\]
\end{corollary}
\begin{proof}
Let us verify the target of $\Theta_\xi$.
Fix $f\in\mathcal G_\xi$.
Claim~\ref{c:uncolored-exchange-grades} of Lemma~\ref{lem:uncolored-adjacent-exchange} gives
$\operatorname{RepGr}_\xi(f)\in\operatorname{CtdG}_3(Y'_\xi)$ and $D_\xi(f)\in[0,3]_{\Z}$.
The definition of $J_\xi(f)$ gives $J_\xi(f)\in[|Y_\xi|_N]_{\Z}$.
\end{proof}

\Needspace{25\baselineskip}
\begin{lemma}[An adjacent edge exchange inside a shell]\label{lem:exchange-precursors}
Let $\xi=(P,Q,o,u,v,A,B)\in\operatorname{AdjShellExch}$.
Let
\[
 (X,X',Y,Y')=(X_\xi,X'_\xi,Y_\xi,Y'_\xi).
\]
Let $e_1,\ldots,e_{|Y|_E}$ and $n_1,\ldots,n_{|Y|_N}$ be the horizontal and vertical edges of $Y$ in path order.
Let $e'_1,\ldots,e'_{|Y'|_E}$ and $n'_1,\ldots,n'_{|Y'|_N}$ be the horizontal and vertical edges of $Y'$ in path order.
There are maps
\[
\begin{aligned}
 \Psi&:\bigcup_{c\ge0}\mathcal C_c(X)\longrightarrow
             \bigcup_{c\ge0}\mathcal C_c(Y'),\\
 \kappa&:\bigcup_{c\ge0}\mathcal C_c(X)\longrightarrow
 \{(h,d):h\in[0,3]_{\Z},\ d\in[\tbinom3h]_{\Z}\}
\end{aligned}
\]
with the following properties.
For each $(f,\gamma)\in\bigcup_{c\ge0}\mathcal C_c(X)$, let $g=\operatorname{GShell}_{f}(X;u,v,0)$.
Let $\eta$ be the shell coloring in~\eqref{eq:shell-recolor} and let
$(g',\gamma')=\Psi(f,\gamma)$.
Then, we have the following.

\begin{enumerate}
\item\label{c:adjacent-exchange-surplus} We have
\[
 \operatorname{Surp}(X,f)\le\operatorname{Surp}(Y',g')
 \le\operatorname{Surp}(X,f)+3.
\]
\item\label{c:adjacent-exchange-horizontal} For $i\in[|Y|_E]_{\Z}$, we have
\[
 (g'(e'_i),\gamma'(e'_i))=(g(e_i),\eta(e_i)).
\]
\item\label{c:adjacent-exchange-vertical} We have
\[
 \#\{j\in[|Y|_N]_{\Z}:(g'(n'_j),\gamma'(n'_j))\ne(g(n_j),\eta(n_j))\}\le1.
\]
\item\label{c:adjacent-exchange-recovery} The map
\[
 (f,\gamma)\longmapsto\bigl(\Psi(f,\gamma),\kappa(f,\gamma)\bigr)
\]
is injective.
\end{enumerate}
\end{lemma}
\begin{proof}
\prf{Claim~\ref{c:adjacent-exchange-surplus}} Let us define $\Psi$ and $\kappa$.
Fix $(f,\gamma)\in\mathcal C_c(X)$.
For every $e\in\operatorname{NE}(X)$, the condition $\gamma(e)\in[\binom3{f(e)}]_{\Z}$ gives $f(e)\le3$. Thus, $f\in\mathcal G_\xi$.
Use the map $\Theta_\xi$ from Corollary~\ref{cor:adjacent-repair-map} and let $(g'_f,j_f,\delta_f)=\Theta_\xi(f)$.
Let $\tau=\tau_\xi$. Define
\[
 \gamma'=\operatorname{Recol}(Y,Y',\tau,g'_f;g,\eta),\qquad
 \Psi(f,\gamma)=(g'_f,\gamma'),\qquad
 \kappa(f,\gamma)=(g(n_{j_f}),\eta(n_{j_f})).
\]
Claim~\ref{c:uncolored-exchange-surplus} of Lemma~\ref{lem:uncolored-adjacent-exchange} gives the surplus bounds.

\prf{Claim~\ref{c:adjacent-exchange-horizontal}} Let us check the horizontal grade and color pairs.
Definition~\ref{def:shell-adjacent-exchange}(\ref{d:adj-shell-repair-grading}) keeps every horizontal grade.
Equation~\eqref{eq:color-update} therefore keeps every horizontal color.

\prf{Claim~\ref{c:adjacent-exchange-vertical}} Let us check how many vertical grade and color pairs change.
The repair in Definition~\ref{def:shell-adjacent-exchange}(\ref{d:adj-shell-repair-grading}) changes at most the vertical grade at $n'_{j_f}$.
Equation~\eqref{eq:color-update} keeps the colors at all other vertical edges.

\prf{Claim~\ref{c:adjacent-exchange-recovery}} Let us recover the input from the output and the one changed pair.
The pair $(g',\gamma')$ determines $g(e_i)=g'(e'_i)$ for every $i\in[|Y|_E]_{\Z}$.
Claim~\ref{c:repair-index-horizontal} of Lemma~\ref{lem:exchange-repair-index}
therefore determines $j_f$ from those grades.
For $j\in[|Y|_N]_{\Z}\setminus\{j_f\}$, Claim~\ref{c:adjacent-exchange-vertical} and the definition of $\gamma'$ give $(g'(n'_j),\gamma'(n'_j))=(g(n_j),\eta(n_j))$.
Recover $(g(n_{j_f}),\eta(n_{j_f}))=\kappa(f,\gamma)$.
Reverse the exchange and remove the two outer simple runs.
The result is $(f,\gamma)$. Thus, $(\Psi(f,\gamma),\kappa(f,\gamma))$ determines $(f,\gamma)$.
\end{proof}

\begin{definition}[Exchange decompositions]\label{def:exchange-decomposition}
Let $X$ be a tNE path. If $X=\emptyset$, let $\operatorname{ExcDecomp}(X)=\emptyset$.
Suppose $X\ne\emptyset$. An \textit{exchange decomposition} of $X$ is a tuple
$(B^{\mathrm L},S_1,S_2,B^{\mathrm R})$ of NE paths such that $S_1$ and $S_2$ are simple runs and
\[
 X=\operatorname{Trans}(B^{\mathrm L}\ccn S_1\ccn S_2\ccn B^{\mathrm R},\operatorname{Start}_{\mathrm{pt}}(X)),\qquad
 |S_1|_N\ne|S_2|_N.
\]
We denote the set of these tuples by $\operatorname{ExcDecomp}(X)$.
\end{definition}

\begin{lemma}\label{lem:exchange}
Let $X$ be a binary run. Let $u_0\in[2]_{\Z}$ and $v_0\in[2]_{\Z}$.
Let $(B^{\mathrm L},S_1,S_2,B^{\mathrm R})\in\operatorname{ExcDecomp}(X)$.
Let $X'=\operatorname{Exc}(X;|B^{\mathrm L}|,S_1,S_2)$.
Let $Y=R_{u_0}\ccn X\ccn R_{v_0}$ and $Y'=R_{u_0}\ccn X'\ccn R_{v_0}$.
There is a map
\[
 \Phi:\bigcup_{c\ge0}\mathcal C_c(X)\longrightarrow
       \bigcup_{c\ge0}\mathcal C_c(Y')
\]
with the following properties.
For each $(f,\gamma)\in\bigcup_{c\ge0}\mathcal C_c(X)$, let $g=\operatorname{GShell}_{f}(X;u_0,v_0,0)$
and define $\eta$ by~\eqref{eq:shell-recolor}.
Let $(g',\gamma')=\Phi(f,\gamma)$.
Label the horizontal and vertical edges of $Y$ by $e_i$ and $n_j$, respectively.
Label the horizontal and vertical edges of $Y'$ by $e'_i$ and $n'_j$, respectively, in path order.
Then, we have the following.

\begin{enumerate}
\item\label{c:exchange-connected} The grading $g'$ satisfies the following.
\begin{enumerate}[label=(\alph*),ref=\theenumi(\alph*)]
\item\label{c:exchange-connectedness} $g'\in\operatorname{CtdG}(Y')$.
\item\label{c:exchange-residual} We have
\[
 \operatorname{Surp}(X,f)\le\operatorname{Surp}(Y',g')
 \le\operatorname{Surp}(X,f)+3.
\]
\end{enumerate}
\item\label{c:exchange-horizontal} We have
\[
 (g'(e'_i),\gamma'(e'_i))=(g(e_i),\eta(e_i))
 \quad(i\in[|Y|_E]_{\Z}).
\]
\item\label{c:exchange-vertical} We have
\[
 \#\{j\in[|Y|_N]_{\Z}:(g'(n'_j),\gamma'(n'_j))\ne(g(n_j),\eta(n_j))\}\le1.
\]
\item\label{c:exchange-fibers} For every $(h,\delta)\in\bigcup_{c\ge0}\mathcal C_c(Y')$, we have
\[
 \#\Phi^{-1}(\{(h,\delta)\})\le8.
\]
\end{enumerate}
\end{lemma}
\begin{proof}
\prf{Claim~\ref{c:exchange-connectedness}} Let us specify the tuple in Definition~\ref{def:shell-adjacent-exchange}.
Define
\[
\begin{gathered}
 o=\operatorname{Start}_{\mathrm{pt}}(X),\qquad H=B^{\mathrm L}\ccn E\ccn N,\\
 P=\operatorname{Trans}(H[1,|H|],(-1,0)),\qquad Q=N\ccn B^{\mathrm R},\\
 (A,B)=
 \begin{cases}
 (E\ccn N,N\ccn E),&(S_1,S_2)=(R_1,R_2),\\
 (N\ccn E,E\ccn N),&(S_1,S_2)=(R_2,R_1).
 \end{cases}
\end{gathered}
\]
Since $X$ is binary, $H$ begins with $E$.
Thus, $P$ and $Q$ are NE paths and $E\ccn P=H$.
It follows that
\[
 \xi=(P,Q,o,u_0,v_0,A,B)\in\operatorname{AdjShellExch},\qquad
 (X_\xi,X'_\xi,Y_\xi,Y'_\xi)=(X,X',Y,Y').
\]
Apply Lemma~\ref{lem:exchange-precursors} to $\xi$ and let $\Phi=\Psi$.
The inclusion $\Psi(f,\gamma)\in\bigcup_{c\ge0}\mathcal C_c(Y')$ gives $g'\in\operatorname{CtdG}(Y')$.

\prf{Claim~\ref{c:exchange-residual}}
Claim~\ref{c:adjacent-exchange-surplus} of Lemma~\ref{lem:exchange-precursors} gives the surplus bounds.

\prf{Claim~\ref{c:exchange-horizontal}}
Claim~\ref{c:adjacent-exchange-horizontal} of Lemma~\ref{lem:exchange-precursors} preserves every horizontal grade--color pair.

\prf{Claim~\ref{c:exchange-vertical}}
Claim~\ref{c:adjacent-exchange-vertical} of Lemma~\ref{lem:exchange-precursors} bounds the number of changed vertical grade--color pairs by one.

\prf{Claim~\ref{c:exchange-fibers}} Let us bound each fiber by the number of recovery choices.
By Claim~\ref{c:adjacent-exchange-recovery} of Lemma~\ref{lem:exchange-precursors},
the map $\kappa$ is injective on each fiber of $\Phi$.
The set of possible values of $\kappa$ has cardinality
\[
 \sum_{h=0}^{3}\binom3h=8.
\]
Hence, every fiber of $\Phi$ has cardinality at most eight.
\end{proof}
\begin{definition}[Wrapped exchange]\label{def:wrapped-exchange}
Let $X$ be a binary run. Let $u_0\in[2]_{\Z}$ and $v_0\in[2]_{\Z}$.
Let $(B^{\mathrm L},S_1,S_2,B^{\mathrm R})\in\operatorname{ExcDecomp}(X)$.
Let $X'=\operatorname{Exc}(X;|B^{\mathrm L}|,S_1,S_2)$.
Let $\operatorname{WrapExch}(X,X';u_0,v_0)$ be the set of \textit{wrapped exchanges}

\[
 \Phi\colon\bigcup_{c\in\Z_{\ge0}}\mathcal C_c(X)
 \longrightarrow\bigcup_{c\in\Z_{\ge0}}\mathcal C_c(R_{u_0}\ccn X'\ccn R_{v_0})
\]
satisfying Claims~\ref{c:exchange-connected}--\ref{c:exchange-fibers}
of Lemma~\ref{lem:exchange}.
\end{definition}

\Needspace{10\baselineskip}
\subsection{Mechanical comparison by exchanges}
\begin{definition}\label{def:binary-run-count}
Let $W$ be a binary run with $n$ simple runs. Let
\[
 W=\operatorname{Trans}(R_{w_1}\ccn \cdots\ccn R_{w_n},\operatorname{Start}_{\mathrm{pt}}(W)).
\]
For $r\in[2]_{\Z}$, define
\[
 \operatorname{RunCount}_r(W)=\#\{i\in[n]_{\Z}:w_i=r\}.
\]
\end{definition}

\begin{definition}[Binary prefix discrepancies]\label{def:prefix-discrepancy}
Let $Z$ and $Y$ be binary runs with $n$ and $m$ simple runs, respectively. Let
\[
 Z=\operatorname{Trans}(R_{z_1}\ccn \cdots\ccn R_{z_n},\operatorname{Start}_{\mathrm{pt}}(Z)),\qquad
 Y=\operatorname{Trans}(R_{y_1}\ccn \cdots\ccn R_{y_m},\operatorname{Start}_{\mathrm{pt}}(Y)).
\]
For $j\in[0,\min(n,m)]_{\Z}$, define
\[
 \operatorname{PrefDiff}_j(Z,Y)=
 \#\{i\in[j]_{\Z}:z_i=2\}-\#\{i\in[j]_{\Z}:y_i=2\}.
\]
For $h\in[0,\min(n,m)]_{\Z}$, define
$\operatorname{PrefSum}_h(Z,Y)=\sum_{j=1}^h|\operatorname{PrefDiff}_j(Z,Y)|$.
\end{definition}

\begin{lemma}[Exchanges between binary runs]\label{lem:binary-prefix-exchanges}
Let $n\in\Z_{>0}$.
Let $\mathbf t\in\Z^2$ and let
\[
 Z=\operatorname{Trans}(R_{z_1}\ccn \cdots\ccn R_{z_n},\mathbf t),\qquad
 Y=\operatorname{Trans}(R_{y_1}\ccn \cdots\ccn R_{y_n},\mathbf t)
\]
be binary runs.
Assume $\operatorname{RunCount}_2(Z)=\operatorname{RunCount}_2(Y)$.
For $j\in[n]_{\Z}$, let
\[
 D_j=\operatorname{PrefDiff}_j(Z,Y).
\]
Then, there are $\ell\in\Z_{\ge0}$ and a run-edit list $\Omega=(\omega_s)_{s\in[\ell]_{\Z}}$ such that
\[
 \operatorname{RunEResult}(\Omega,Z)=Y.
\]
Every $\omega_s$ exchanges two adjacent simple runs of distinct vertical lengths, and
\begin{equation}\label{eq:binary-exchange-distance}
 \ell=\operatorname{PrefSum}_n(Z,Y).
\end{equation}
\end{lemma}
\begin{proof}
Let us decrease $\operatorname{PrefSum}_n(Z,Y)$ by adjacent exchanges to prove~\eqref{eq:binary-exchange-distance}.
If $Z=Y$, the empty run-edit list suffices. Suppose $Z\ne Y$ and let
\[
 a=\min\{i\in[n]_{\Z}:z_i\ne y_i\}.
\]
The entries before $a$ agree. Since the total counts of each vertical length agree, there is $i>a$ such that $z_i=y_a$.
Let
\[
 b=\min\{i\in[a+1,n]_{\Z}:z_i=y_a\}.
\]
For $j\in[a,b-1]_{\Z}$, minimality of $b$ gives $z_j=3-y_a$. Therefore,
\[
 D_j=\begin{cases}
 -\#\{i\in[a,j]_{\Z}:y_i=2\}\le-1,&y_a=2,\\
 \#\{i\in[a,j]_{\Z}:y_i=1\}\ge1,&y_a=1.
 \end{cases}
\]

Let us move $R_{z_b}=R_{y_a}$ to position $a$.
Set $Z^{(b)}=Z$. For $k=b-1,b-2,\ldots,a$, define
\[
 Z^{(k)}=\operatorname{Exc}\left(Z^{(k+1)};
 \sum_{i=1}^{k-1}(1+z_i),R_{3-y_a},R_{y_a}\right).
\]
For $k\in[a,b-1]_{\Z}$, the simple runs of $Z^{(k+1)}$ in positions $k$ and $k+1$ are $R_{3-y_a}$ and $R_{y_a}$, respectively, up to translation.
Their vertical lengths are distinct. Let
\[
 D_j^{(k)}=\operatorname{PrefDiff}_j(Z^{(k)},Y)
 \qquad(j\in[n]_{\Z},\ k\in[a,b]_{\Z}).
\]
Only the prefix ending at position $k$ changes in the exchange. Thus,
\[
 D_j^{(k)}=\begin{cases}
 D_k^{(k+1)}+1,&j=k,\ y_a=2,\\
 D_k^{(k+1)}-1,&j=k,\ y_a=1,\\
 D_j^{(k+1)},&j\ne k.
 \end{cases}
\]
The exchanges at positions greater than $k$ do not change this prefix. Thus, $D_k^{(k+1)}=D_k$.
The bounds $D_k\le-1$ when $y_a=2$ and $D_k\ge1$ when $y_a=1$ give
\[
 |D_k^{(k)}|=|D_k^{(k+1)}|-1,\qquad
 \operatorname{PrefSum}_n(Z^{(k)},Y)
 =\operatorname{PrefSum}_n(Z^{(k+1)},Y)-1.
\]

Let us prove that the sequence of exchanges terminates.
After these $b-a$ exchanges, $Z^{(a)}$ agrees with $Y$ through position $a$.
If $Z^{(a)}\ne Y$, use $Z^{(a)}$ as the new input $Z$.
Recompute its vertical lengths $z_1,\ldots,z_n$, choose $a$ and $b$ by the defining minima, and perform the exchanges for $k=b-1,\ldots,a$.
Each application increases the number of initial simple runs that agree with $Y$. Thus, the process ends after finitely many applications.
Consequently, the resulting run-edit list ends at $Y$ and has length
\[
 \ell=\operatorname{PrefSum}_n(Z,Y)-\operatorname{PrefSum}_n(Y,Y)
 =\operatorname{PrefSum}_n(Z,Y).
\]
\end{proof}

\begin{lemma}[Completion of a binary run]\label{lem:binary-completion}
Let $L\in\Z_{>0}$ and $C\in\Z_{\ge0}$.
Let $\mathbf t\in\Z^2$ and let
\[
 X=\operatorname{Trans}(R_{x_1}\ccn \cdots\ccn R_{x_L},\mathbf t),\qquad
 Y=\operatorname{Trans}(R_{y_1}\ccn \cdots\ccn R_{y_{L+C}},\mathbf t)
\]
be binary runs.
For $r\in[2]_{\Z}$, assume
\[
 b_r=\operatorname{RunCount}_r(Y)-\operatorname{RunCount}_r(X)\ge0.
\]
Let $Z=X\ccn R_1^{b_1}\ccn R_2^{b_2}$.
Then, $b_1+b_2=C$, and there are $\ell\in\Z_{\ge0}$ and a run-edit list $\Omega=(\omega_s)_{s\in[\ell]_{\Z}}$ such that
\[
 \operatorname{RunEResult}(\Omega,Z)=Y.
\]
Every $\omega_s$ exchanges two adjacent simple runs of distinct vertical lengths, and
\begin{equation}\label{eq:binary-completion-bound}
 \ell\le\operatorname{PrefSum}_L(X,Y)+\frac{C(C-1)}2.
\end{equation}
\end{lemma}
\begin{proof}
Let us check the counts after appending the missing simple runs.
The two run counts sum to the number of simple runs. Thus, $b_1+b_2=C$.
Since $Z=X\ccn R_1^{b_1}\ccn R_2^{b_2}$, we have
$\operatorname{RunCount}_r(Z)=\operatorname{RunCount}_r(X)+b_r=\operatorname{RunCount}_r(Y)$ for $r\in[2]_{\Z}$.
Apply Lemma~\ref{lem:binary-prefix-exchanges}.
For $j\in[L+C]_{\Z}$, let
\[
 D_j=\operatorname{PrefDiff}_j(Z,Y).
\]
For $j\le L$, the first $j$ simple runs of $Z$ belong to $X$.
For $j=L+r$ with $r\in[C]_{\Z}$, equality of the total counts gives
\[
 |D_{L+r}|\le C-r,\qquad
 \sum_{r=1}^{C}|D_{L+r}|\le\sum_{r=1}^{C}(C-r)=\frac{C(C-1)}2.
\]
Adding $\sum_{j=1}^L|D_j|$ gives~\eqref{eq:binary-completion-bound}.
\end{proof}

\begin{definition}[Binary run-edit lists]\label{def:binary-run-edit-lists}
Let $U$ and $V$ be binary runs and let $K\in\Z_{\ge0}$ and $k\in\Z_{\ge0}$.
For a run-edit list $\Omega=(\omega_s)_{s\in[K+k]_{\Z}}$ from $U$ to $V$, let
\[
 U_s=\operatorname{RunEResult}((\omega_1,\ldots,\omega_s),U)
 \qquad(s\in[0,K+k]_{\Z}).
\]
We define membership in $\operatorname{BinEdit}_{K,k}(U,V)$ as follows.
We have $\Omega\in\operatorname{BinEdit}_{K,k}(U,V)$ if and only if every $U_s$ is binary for $s\in[0,K+k]_{\Z}$ and there is a partition
\[
 [K+k]_{\Z}=I_R\sqcup I_X,\qquad \#I_R=K,\qquad \#I_X=k,
\]
with the following conditions.
\begin{enumerate}[font=\itshape]
\item For $s\in I_R$, there are $j_s\in[2]_{\Z}$ and $i_s\in[0,|U_{s-1}|]_{\Z}$ such that
\[
 \omega_s(U_{s-1})=\operatorname{Ins}(U_{s-1},i_s,R_{j_s}),\qquad
 i_s=|U_{s-1}|\ \text{or}\ (U_{s-1})_{i_s+1}\in E(U_{s-1}).
\]
\item For $s\in I_X$, the run-edit $\omega_s$ exchanges adjacent simple runs $A$ and $B$ such that
\[
 (|A|_N,|B|_N)\in\{(1,2),(2,1)\}.
\]
\end{enumerate}
\end{definition}

\begin{lemma}[Comparison of mechanical factors by exchanges]\label{lem:binary-factor-comparison}
Let $r\in(1,2)$ be irrational and $0<\eta<\frac14$. Let
\[
 \tau=\frac12\min(r-1,2-r),\qquad C=\ceil{\frac4\tau}.
\]
Then, there are $L\in\Z_{>0}$ and $\delta>0$ depending only on $r$ and $\eta$ with the following properties.
\begin{enumerate}
\item\label{c:binary-factor-constants} We have
\[
 C\le\eta L,\qquad \delta<\min\left\{\tau,\frac{\eta}{4L}\right\}.
\]
\item\label{c:binary-factor-uniform} For each $\sigma_1\in(r-\delta,r+\delta)$, $\sigma_2\in(r-\delta,r+\delta)$,
$u_1\in[0,1)$, $u_2\in[0,1)$, $p\in\Z_{\ge0}$, $V=\operatorname{MechRun}(\sigma_2,u_2)$, and translated copy $X$ of
$\operatorname{RunFactor}_{1,L}(\operatorname{MechRun}(\sigma_1,u_1))$,
there are $q\in\Z_{\ge0}$, $k\in\Z_{\ge0}$, $\zeta\in\Z^2$, and a run-edit list $\Omega$ with the following properties.
\begin{enumerate}[label=(\alph*),ref=\theenumi(\alph*)]
\item\label{c:binary-factor-budget} We have $C+q\le\eta L$.
\item\label{c:binary-factor-target} We have
\[
 \Omega\in\operatorname{BinEdit}_{C+q,k}
 \bigl(X,\operatorname{Trans}(\operatorname{RunFactor}_{p+1,L+C+q}(V),\zeta)\bigr).
\]
\item\label{c:binary-factor-exchanges} We have $k\le7\eta L$.
\end{enumerate}
\end{enumerate}
\end{lemma}
\begin{proof}
\prf{Claim~\ref{c:binary-factor-constants}} Let us choose a rational approximation common to both slopes.
Choose relatively prime $p_0\in\Z_{>0}$ and $q_0\in\Z_{>0}$ satisfying~\eqref{eq:mechanical-approximation}.
Choose $L\in\Z_{>0}$ and $\delta>0$ such that
\[
 q_0+C\le\eta L,\qquad C(C+1)\le\eta L,\qquad
 \delta<\min\left\{\tau,q_0^{-2},\frac{\eta}{4L}\right\}.
\]
These constants depend only on $r$ and $\eta$. They satisfy $C\le q_0+C\le\eta L$ and $\delta<\min\{\tau,\frac\eta{4L}\}$.

\prf{Claim~\ref{c:binary-factor-budget}} Let us bound the count of each simple-run length.
For $i\in[2]_{\Z}$, we have $1+\tau<\sigma_i<2-\tau$ and $|\sigma_i-\frac{p_0}{q_0}|<2q_0^{-2}$.
Let us choose $q$ and the simple runs inserted before $X$.
Apply Claim~\ref{c:phase-selection} of Lemma~\ref{lem:phase-prefix-discrepancy}
with offset $b=p$. Let $q=\min J$ and $v_q=\{u_2+(p+q)\sigma_2\}$ as in Claim~\ref{c:phase-selection} of Lemma~\ref{lem:phase-prefix-discrepancy}.
Then, $0\le q<q_0$ and $|u_1-v_q|<\eta$.
If $q>0$, insert $\operatorname{RunFactor}_{p+1,q}(V)$ immediately before $X$.
We have $q+C\le q_0+C\le\eta L$.

\prf{Claim~\ref{c:binary-factor-target}} Let us compare the completed factor with its target.
Let $Y$ be the translated copy of $\operatorname{RunFactor}_{p+q+1,L+C}(V)$
that starts at $\operatorname{Start}_{\mathrm{pt}}(X)$.

Let us specify the simple runs appended to $X$.
A mechanical factor of $a$ simple runs and slope $\sigma$ has
$a(\sigma-1)+e$ simple runs of length two and $a(2-\sigma)-e$ of length one
for some $|e|<1$.
For $j\in[2]_{\Z}$, define
\[
 b_j=\operatorname{RunCount}_j(Y)-\operatorname{RunCount}_j(X).
\]
The inequalities $\sigma_2-1>\tau$, $2-\sigma_2>\tau$, and
$|\sigma_2-\sigma_1|<2\delta$ give
\[
 b_j>\tau C-2-2\delta L>0,\qquad b_1+b_2=C.
\]
Append $b_1$ copies of $R_1$ and $b_2$ copies of $R_2$ immediately after $X$.

Let us construct the exchanges from $X\ccn R_1^{b_1}\ccn R_2^{b_2}$ to $Y$.
Lemma~\ref{lem:binary-completion} supplies a run-edit list $\mathcal S=(S_1,\ldots,S_\nu)$ from $X\ccn R_1^{b_1}\ccn R_2^{b_2}$ to $Y$
consisting of exchanges of adjacent simple runs of distinct vertical lengths.
Let $H$ be the prefix of $q$ simple runs inserted before $X$.
When a local run-edit acts on the current factor $Z$ after $H$, replace its edge index $i$ by $|H|+i$.
Replace $\operatorname{Exc}(Z;i,R_a,R_b)$ by
\[
 \operatorname{Exc}(H\ccn Z;|H|+i,R_a,R_b).
\]

Let us identify the path obtained by the complete run-edit list.
The $q$ prefix insertions followed by the $C$ appended insertions and these exchanges give $\Omega$.
Let $k=\nu$ and
\[
 \zeta=\operatorname{Start}_{\mathrm{pt}}(\operatorname{RunEResult}(\Omega,X))
       -(p,\floor{u_2+p\sigma_2}).
\]
The inserted and exchanged simple runs have vertical lengths in $[2]_{\Z}$. We have
\[
 \operatorname{RunEResult}(\Omega,X)
 =\operatorname{Trans}(\operatorname{RunFactor}_{p+1,L+C+q}(V),\zeta).
\]
Each insertion is between simple runs or at an end.
Thus, Definition~\ref{def:binary-run-edit-lists} gives the membership assertion in Claim~\ref{c:binary-factor-target}.

\prf{Claim~\ref{c:binary-factor-exchanges}} Let us bound the number of exchanges in $\mathcal S$.
For $i\in[L]_{\Z}$, the difference of the length-two prefix counts is
\[
 D_i=\operatorname{PrefDiff}_i(X,Y)=\floor{u_1+i\sigma_1}-\floor{v_q+i\sigma_2}.
\]
Claim~\ref{c:phase-prefix-mass} of Lemma~\ref{lem:phase-prefix-discrepancy} gives
\[
 \operatorname{PrefSum}_L(X,Y)=\sum_{i=1}^L|D_i|\le\frac{11}{2}\eta L.
\]
The bound in Lemma~\ref{lem:binary-completion} gives
\[
 k=\nu\le\sum_{i=1}^L|D_i|+\frac{C(C-1)}2
 \le\frac{13}{2}\eta L\le7\eta L.
\]
\end{proof}

\begin{lemma}[Completion with one exceptional simple run]\label{lem:binary-final-completion}
Let $t\in\Z_{>0}$, $C\in\Z_{>0}$, $p\in\Z_{\ge0}$, $\sigma_1\in(1,2)$, $\sigma_2\in(1,2)$, $u_1\in[0,1)$, and $u_2\in[0,1)$.
Assume
\[
 C\min(\sigma_2-1,2-\sigma_2)>3+t|\sigma_1-\sigma_2|.
\]
Let $A_i$ be simple runs satisfying $|A_i|_N\in[2]_{\Z}$ for $i\in[t]_{\Z}$.
Let $X=A_1\ccn \cdots\ccn A_t$.
Assume
\[
 |A_i|_N=\floor{u_1+i\sigma_1}-\floor{u_1+(i-1)\sigma_1}
 \qquad(i\in[t-1]_{\Z}).
\]
Let
\[
 Y=\operatorname{Trans}\bigl(\operatorname{RunFactor}_{p+1,t+C}(\operatorname{MechRun}(\sigma_2,u_2)),
  \operatorname{Start}_{\mathrm{pt}}(X)-(p,\floor{u_2+p\sigma_2})\bigr).
\]
Then, there are $k\in\Z_{\ge0}$ and a run-edit list $\Omega$ satisfying the following.
\begin{enumerate}
\item\label{c:binary-final-target} We have $\Omega\in\operatorname{BinEdit}_{C,k}(X,Y)$.
\item\label{c:binary-final-exchanges} We have $k\le(t+C)^2$.
\end{enumerate}
\end{lemma}
\begin{proof}
\prf{Claim~\ref{c:binary-final-target}} Let us specify the simple runs appended to $X$.
For $j\in[2]_{\Z}$, put $b_j=\operatorname{RunCount}_j(Y)-\operatorname{RunCount}_j(X)$.
Replacing the final mechanical simple run by $A_t$ changes either run count by at most one.
Definition~\ref{def:mechanical-paths}(\ref{d:mechanical-sequence}) gives
\[
\begin{gathered}
 \operatorname{RunCount}_2(X)=t(\sigma_1-1)+e_1,\qquad
 \operatorname{RunCount}_2(Y)=(t+C)(\sigma_2-1)+e_2,\\
 |e_1|<2,\qquad |e_2|<1.
\end{gathered}
\]
We have
\[
 \operatorname{RunCount}_1(X)=t-\operatorname{RunCount}_2(X),\qquad
 \operatorname{RunCount}_1(Y)=t+C-\operatorname{RunCount}_2(Y).
\] Hence,
\[
 b_j>C\min(\sigma_2-1,2-\sigma_2)-3-t|\sigma_1-\sigma_2|>0,
 \qquad b_1+b_2=C.
\]
Append $R_1^{b_1}\ccn R_2^{b_2}$ to $X$.
Lemma~\ref{lem:binary-prefix-exchanges} supplies exchanges from $X\ccn R_1^{b_1}\ccn R_2^{b_2}$ to $Y$.
Let $\Omega$ consist of these $C$ insertions followed by the exchanges. Let $k$ be the number of exchanges.
Every intermediate path is binary. Hence, $\Omega\in\operatorname{BinEdit}_{C,k}(X,Y)$.

Let us bound the exchanges in the completed final factor.
\prf{Claim~\ref{c:binary-final-exchanges}} Let $Z=X\ccn R_1^{b_1}\ccn R_2^{b_2}$.
For $i\in[t+C]_{\Z}$, we have $|\operatorname{PrefDiff}_i(Z,Y)|\le t+C$.
Lemma~\ref{lem:binary-prefix-exchanges} gives
\[
 k\le\operatorname{PrefSum}_{t+C}(Z,Y)\le(t+C)^2.
\]
\end{proof}

\Needspace{23\baselineskip}
\begin{lemma}\label{lem:binary-mechanical}
Let $r\in(1,2)$ be irrational and let $0<\eta<1$.
Then, there are constants $0<\delta<\min(r-1,2-r)$ and $B>0$
depending only on $r$ and $\eta$ with the following property.
For each $(2,r,\eta,\delta,B)$-admissible tuple
$(l,T,\sigma_1,\sigma_2,u_1,u_2,(A_1,\ldots,A_l))$ and
\[
 U=A_1\ccn\cdots\ccn A_l,\qquad V=\operatorname{MechRun}(\sigma_2,u_2),
\]
there are $K\in\Z_{\ge0}$, $k\in\Z_{\ge0}$, $\zeta\in\Z^2$, and a run-edit list $\Omega$ such that we have the following.
\begin{enumerate}
\item\label{c:binary-mechanical-target} $\Omega\in\operatorname{BinEdit}_{K,k}
 \bigl(U,\operatorname{Trans}(\operatorname{RunFactor}_{1,T}(V),\zeta)\bigr)$.
\item\label{c:binary-mechanical-insertions} $K=T-l$.
\item\label{c:binary-mechanical-exchanges} $k\le7\eta l+B$.
\end{enumerate}
\end{lemma}
\begin{proof}
\prf{Claim~\ref{c:binary-mechanical-target}} Let us choose the parameters for the binary factor comparison.
Let $\eta_0=\min(\eta,\frac18)$. Then, $0<\eta_0<\frac14$ and $\eta_0\le\eta$.
Choose $\tau$, $C$, $L$, and $\delta$ from Lemma~\ref{lem:binary-factor-comparison} with parameter $\eta_0$. Let
$B=C+(L+C)^2+1$.
The inequalities $|\sigma_i-r|<\delta<\min(r-1,2-r)$ give $1<\sigma_i<2$ for $i\in[2]_{\Z}$.
Thus, every simple run of $U$ and $V$ has vertical length in $[2]_{\Z}$.
Let
\[
\begin{gathered}
 h=\floor{\frac{l-1}{L}},\qquad t=l-hL\in[L]_{\Z},\\
 X_s=A_{(s-1)L+1}\ccn \cdots\ccn A_{sL}\quad(s\in[h]_{\Z}),\qquad X_*=A_{hL+1}\ccn \cdots\ccn A_l.
\end{gathered}
\]

Let us apply Claims~\ref{c:binary-factor-budget}--\ref{c:binary-factor-exchanges} of Lemma~\ref{lem:binary-factor-comparison} to $X_1,\ldots,X_h$.
Let $p_0=0$.
For $s\in[h]_{\Z}$, use the source intercept $\{u_1+(s-1)L\sigma_1\}$ and the target index $p_{s-1}$.
Claim~\ref{c:binary-factor-target} of Lemma~\ref{lem:binary-factor-comparison} gives $q_s\in\Z_{\ge0}$ and a run-edit list
$\Omega_s=(\omega_{s,1},\ldots,\omega_{s,n_s})$ from $X_s$ to a translated copy of
$\operatorname{RunFactor}_{p_{s-1}+1,L+C+q_s}(V)$.
Define $p_s=p_{s-1}+L+C+q_s$.

Let us bound the number of run-edits in each $\Omega_s$.
For $s\in[h]_{\Z}$, let $K_s$ and $k_s$ be its numbers of simple-run insertions and adjacent exchanges, respectively.
Claims~\ref{c:binary-factor-budget} and~\ref{c:binary-factor-exchanges} of Lemma~\ref{lem:binary-factor-comparison} give
\[
 K_s\le\eta_0L,\qquad k_s\le7\eta_0L\qquad(s\in[h]_{\Z}).
\]

Let us apply Claims~\ref{c:binary-final-target} and~\ref{c:binary-final-exchanges} of Lemma~\ref{lem:binary-final-completion} to $X_*$.
Use the source intercept $\{u_1+hL\sigma_1\}$ and target index $p_h$.
Since $\tau C\ge4$, $\min(\sigma_2-1,2-\sigma_2)>\tau$, and $2\delta t<\frac{\eta_0}2<1$, we have
\[
 C\min(\sigma_2-1,2-\sigma_2)>3+t|\sigma_1-\sigma_2|.
\]
Claim~\ref{c:binary-final-target} of Lemma~\ref{lem:binary-final-completion} gives a run-edit list $\Omega_*=(\omega_{*,1},\ldots,\omega_{*,n_*})$ from $X_*$
to a translated copy of $\operatorname{RunFactor}_{p_h+1,t+C}(V)$.
Its insertion count is $C$. Claim~\ref{c:binary-final-exchanges} of Lemma~\ref{lem:binary-final-completion} bounds its exchange count by $(t+C)^2$.

Let us assemble the run-edits on the full path.
Apply $\Omega_1,\ldots,\Omega_h$ and $\Omega_*$ successively to their factors in the whole current path.
Let $B^{\mathrm L}$ and $B^{\mathrm R}$ be the prefix and suffix paths of the current factor $Z$, so the current path is $B^{\mathrm L}\ccn Z\ccn B^{\mathrm R}$.
Replace a local insertion $\operatorname{Ins}(Z,i,R_j)$ by
\[
 \operatorname{Ins}(B^{\mathrm L}\ccn Z\ccn B^{\mathrm R},|B^{\mathrm L}|+i,R_j).
\]
Replace a local exchange $\operatorname{Exc}(Z;i,R_a,R_b)$ by
\[
 \operatorname{Exc}(B^{\mathrm L}\ccn Z\ccn B^{\mathrm R};|B^{\mathrm L}|+i,R_a,R_b).
\]
Induction on the local run-edits gives translated copies of $\operatorname{RunFactor}_{p_{s-1}+1,L+C+q_s}(V)$ for $s=1,\ldots,h$, followed by $\operatorname{RunFactor}_{p_h+1,t+C}(V)$.
Every intermediate path is binary and the insertion and exchange counts are unchanged.

Let us count the simple-run insertions and exchanges.
Let $Q$ be the resulting path and let $T_0=p_h+t+C$.
Then,
\[
\begin{gathered}
 Q=\operatorname{Trans}(\operatorname{RunFactor}_{1,T_0}(V),\operatorname{Start}_{\mathrm{pt}}(Q)),\\
 T_0=l+\sum_{s=1}^h(C+q_s)+C
 \le l+h\eta_0 L+C\le(1+\eta)l+B\le T.
\end{gathered}
\]
If $T_0<T$, append $\operatorname{RunFactor}_{T_0+1,T-T_0}(V)$.
Concatenating the run-edit lists on the whole path gives $\Omega$.
Let $\zeta=\operatorname{Start}_{\mathrm{pt}}(Q)$ and write $\Omega=(\omega_1,\ldots,\omega_n)$.
Define
\[
\begin{aligned}
 I_R&=\{s\in[n]_{\Z}:\omega_s\text{ is a simple-run insertion}\},\\
 I_X&=\{s\in[n]_{\Z}:\omega_s\text{ exchanges adjacent simple runs}\},\\
 K&=\#I_R,\qquad k=\#I_X.
\end{aligned}
\]
Then, $[n]_{\Z}=I_R\sqcup I_X$.
Every intermediate path is binary. Each exchange has the unequal-length property in Definition~\ref{def:binary-run-edit-lists}.
The endpoint construction and Definition~\ref{def:binary-run-edit-lists} give
\[
 \Omega\in\operatorname{BinEdit}_{K,k}
 \bigl(U,\operatorname{Trans}(\operatorname{RunFactor}_{1,T}(V),\zeta)\bigr).
\]

\prf{Claim~\ref{c:binary-mechanical-insertions}}
The simple-run insertion count is
\[
 K=(T_0-l)+(T-T_0)=T-l.
\]

\prf{Claim~\ref{c:binary-mechanical-exchanges}}
The exchange counts of $\Omega_1,\ldots,\Omega_h$ and $\Omega_*$ give
\[
 k\le7\eta_0 hL+(t+C)^2\le7\eta l+B.
\]
\end{proof}

\begin{proposition}\label{prop:exchange-composition}
Let $m=3$. Let $K\in\Z_{\ge0}$, $k\in\Z_{\ge0}$, and $P\in\Z_{\ge0}$ satisfy $P\ge K+4k$.
Let $U$ and $V$ be binary runs and fix $\Omega\in\operatorname{BinEdit}_{K,k}(U,V)$.
Let $u\in[2]_{\Z}$ and $v\in[2]_{\Z}$. Let $W_{\mathrm L}$ and $W_{\mathrm R}$ be binary runs with $P$ simple runs each.
Let $V_*=W_{\mathrm L}\ccn R_u\ccn V\ccn R_v\ccn W_{\mathrm R}$.
Then, the following statements hold.
\begin{enumerate}
\item\label{c:exchange-composition-map} There is a map
\[
 \Phi:\mathcal C_0(U)\longrightarrow\mathcal C_0(V_*)
\]
such that
\begin{equation}\label{eq:exchange-composition-fibers}
 \#\Phi^{-1}(\{\xi\})\le8^k\qquad(\xi\in\mathcal C_0(V_*)).
\end{equation}
\item\label{c:exchange-composition-count} We have
\begin{equation}\label{eq:exchange-composition}
 C_0(V_*)\ge\frac{C_0(U)}{8^k}.
\end{equation}
\end{enumerate}
\end{proposition}
\begin{proof}
\prf{Claim~\ref{c:exchange-composition-map}} Let us compose the maps for the individual run-edits.
Let $\Omega=(\omega_s)_{s\in[K+k]_{\Z}}$.
For $s\in[0,K+k]_{\Z}$, let $U_s=\operatorname{RunEResult}((\omega_1,\ldots,\omega_s),U)$.
Choose $I_R$, $I_X$, and $(i_s,j_s)_{s\in I_R}$ satisfying the conditions for $\Omega\in\operatorname{BinEdit}_{K,k}(U,V)$ in Definition~\ref{def:binary-run-edit-lists}.
Let us apply $\omega_1,\ldots,\omega_{K+k}$ to the copy of $U$ in $\operatorname{Shell}(U;u,v,0)$ and bound the surplus.
Define $a_j\in[2]_{\Z}$ and $b_j\in[2]_{\Z}$ for $j\in[P]_{\Z}$ by
\[
\begin{aligned}
 W_{\mathrm L}&=\operatorname{Trans}(R_{a_P}\ccn \cdots\ccn R_{a_1},\operatorname{Start}_{\mathrm{pt}}(W_{\mathrm L})),\\
 W_{\mathrm R}&=\operatorname{Trans}(R_{b_1}\ccn \cdots\ccn R_{b_P},\operatorname{Start}_{\mathrm{pt}}(W_{\mathrm R})).
\end{aligned}
\]
For $(f,\gamma)\in\mathcal C_0(U)$, let
\[
 Q_0=\operatorname{Shell}(U;u,v,0),\qquad
 (h_0,\eta_0)=\Phi_{\mathrm{sh}}(f,\gamma).
\]
Here $\Phi_{\mathrm{sh}}$ is the map in Claim~\ref{c:shell-recolor} of Lemma~\ref{lem:shell} with $j=1$ and $\epsilon=0$.
Claims~\ref{c:shell-connected} and~\ref{c:shell-residual} of Lemma~\ref{lem:shell} give $(h_0,\eta_0)\in\mathcal C_0(Q_0)$.
For $s\in I_R$, let $L\ccn U_{s-1}\ccn R$ be the current path and apply Lemma~\ref{lem:simple-run-insertion} at $|L|+i_s$.
Equation~\eqref{eq:run-insertion-surplus} shows that this run-edit adds one to the surplus.
For $s\in I_X$, let $e_s=1+\#(I_X\cap[s-1]_{\Z})$.
Let $Z=L\ccn U_{s-1}\ccn R$ be the current path and let $Z'=L\ccn U_s\ccn R$.
Choose
\[
 \Phi_s\in\operatorname{WrapExch}(Z,Z';a_{e_s},b_{e_s}).
\]
Lemma~\ref{lem:exchange} gives a member of this set because $Z'$ exchanges the two simple runs specified by $\omega_s$ inside $Z$.
The choice of $\Phi_s$ is independent of the input colored grading.
Apply $\Phi_s$ to the colored grading at step $s-1$.
Claim~\ref{c:exchange-residual} of Lemma~\ref{lem:exchange} bounds the increase in surplus by three.
Thus, the colored grading $(h,\eta)$ on the resulting path $Q$ satisfies
\[
 0\le\operatorname{Surp}(Q,h)\le K+3k\le P-k.
\]

Let us reduce the surplus to zero.
Apply Lemma~\ref{lem:surplus-reduction} with $n=P-k$ and
\[
 \alpha'=(a_{k+j})_{j\in[P-k]_{\Z}},\qquad
 \mu'=(b_{k+j})_{j\in[P-k]_{\Z}}.
\]
Lemma~\ref{lem:word-reindexing} translates the final colored grading to $V_*$.
This defines a map $\Phi\colon\mathcal C_0(U)\to\mathcal C_0(V_*)$.

Let us bound the fibers of $\Phi$.
The map $\Delta_{Q,\alpha',\mu'}$ from Lemma~\ref{lem:surplus-reduction} is injective.
Restriction recovers the colored grading before each simple-run insertion
by~\eqref{eq:run-insertion-grades} and~\eqref{eq:run-insertion-recolor}.

Let us count the choices at the $k$ exchanges.
Claim~\ref{c:exchange-fibers} of Lemma~\ref{lem:exchange} gives at most eight
inputs at each exchange.
Claim~\ref{c:shell-recolor} of Lemma~\ref{lem:shell} recovers $(f,\gamma)$ by restriction from $(h_0,\eta_0)$.
Recovering the inputs in the reverse order of $\Omega$ gives
\[
 \#\Phi^{-1}(\{\xi\})\le8^k\qquad(\xi\in\mathcal C_0(V_*)).
\]
This proves~\eqref{eq:exchange-composition-fibers}.

\prf{Claim~\ref{c:exchange-composition-count}} Let us sum the fiber bounds to compare the grading counts.
Summing~\eqref{eq:exchange-composition-fibers} over $\xi\in\mathcal C_0(V_*)$ gives
\[
 C_0(U)=\sum_{\xi\in\mathcal C_0(V_*)}\#\Phi^{-1}(\{\xi\})
 \le8^k C_0(V_*).
\]
Since $8^k>0$, division by $8^k$ proves~\eqref{eq:exchange-composition}.
\end{proof}

\begin{proposition}\label{prop:three-irrational}
We have the following.
\begin{enumerate}
\item\label{c:three-irrational-limit} For every irrational $r\in(1,2)$,
there is $\Upsilon_3(r)\in(0,\infty)$ such that
\[
 \lim_{\substack{a\to\infty,\ (a,b)\in\cA_3\\\frac{b}{a}\to r}}
 \frac{\log E_3(a,b)}a=\Upsilon_3(r).
\]
\item\label{c:three-irrational-continuity} Let $\Upsilon_3(1)=\log\left(\frac{256}{27}\right)$.
Use $\Upsilon_3(r)$ from~\eqref{eq:rational-period} for rational $r\in(1,\frac32]$ and from Claim~\ref{c:three-irrational-limit} for irrational $r\in(1,\frac32)$.
Then, $\Upsilon_3|_{[1,\frac32]}$ takes values in $(0,\infty)$ and satisfies
\[
 \lim_{\substack{s\to r\\s\in[1,\frac32]}}\Upsilon_3(s)=\Upsilon_3(r)
 \qquad\left(r\in[1,\frac32]\right).
\]
\end{enumerate}
\end{proposition}
\begin{proof}
\prf{Claim~\ref{c:three-irrational-limit}} Let us choose the comparison parameters for Lemma~\ref{lem:binary-mechanical} and Claim~\ref{c:exchange-composition-count} of Proposition~\ref{prop:exchange-composition}.
Consider $\eta$ with $0<\eta<\frac{1}{1024}$. Let $d=1-144\eta$.
Choose $\delta$ and $B$ from Lemma~\ref{lem:binary-mechanical} for $r$ and $\eta$.
Choose relatively prime $p\in\Z_{>0}$ and $q\in\Z_{>0}$ such that $q\ge2$,
$\alpha=\frac{p}{q}\in(1,2)$, and $|\alpha-r|<\delta$.
Let $V_*$ have $\nu$ simple runs and let $P=\floor{64\eta\nu}$.
Let us translate $V_*$ so that
\[
 V_*=R_{a_P}\ccn \cdots\ccn R_{a_1}\ccn R_u\ccn V\ccn R_v\ccn R_{b_1}\ccn \cdots\ccn R_{b_P},
 \qquad a_j,b_j,u,v\in[2]_{\Z}.
\]
Thus, $V$ has $T=\nu-2P-2$ simple runs.
For $l=d\nu+O(1)$, we have
\[
 T-(1+\eta)l=(15\eta+144\eta^2)\nu+O(1)\longrightarrow\infty.
\]
Hence, $T\ge(1+\eta)l+B$ for all sufficiently large $\nu$.

Let us specify the two applications of Lemma~\ref{lem:binary-mechanical}.
Consider $(a,b)\in\cA_3$ as $a\to\infty$ and $\frac ba\to r$.
\begin{enumerate}[label=(\roman*)]
\item For the lower bound, let
\[
\begin{gathered}
 \nu=a+1,\quad j=\floor{\frac{d\nu}{q}},\quad l=jq,\quad
 U=W_{p,q}^j,\quad V_*=\mathcal U(a,b),\\
 (\sigma_1,\sigma_2,u_1,u_2)
 =\bigl(\alpha,\Gamma(a,b),0,\{(P+1)\Gamma(a,b)\}\bigr).
\end{gathered}
\]
\item For the upper bound, let
\[
\begin{gathered}
 l=a+1,\quad N=\ceil{\frac{l}{dq}},\quad \nu=Nq,\quad
 U=\mathcal U(a,b),\quad V_*=W_{p,q}^N,\\
 (\sigma_1,\sigma_2,u_1,u_2)
 =\bigl(\Gamma(a,b),\alpha,0,\{(P+1)\alpha\}\bigr).
\end{gathered}
\]
\end{enumerate}

Let us verify admissibility for each of the two comparisons.
In each case, let $A_1,\ldots,A_l$ be the simple runs of $U$ in path order.
Equations~\eqref{eq:ambient}, \eqref{eq:frame-word}, and~\eqref{eq:rational-period} show that
$(l,T,\sigma_1,\sigma_2,u_1,u_2,(A_1,\ldots,A_l))$ is
$(2,r,\eta,\delta,B)$-admissible for all sufficiently large $a$.
The path $V$ is a translated copy of
$\operatorname{RunFactor}_{1,T}(\operatorname{MechRun}(\sigma_2,u_2))$.
In case (i), the final $R_1$ of $\mathcal U(a,b)$ belongs to $W_{\mathrm R}$.
In case (ii), we have $A_l=R_1$ and $|A_l|_N=1\in[2]_{\Z}$, as permitted by condition~\ref{hyp:mechanical-source} of Definition~\ref{def:admissible-mechanical}.

Let us bound the insertions and exchanges in each comparison.
Lemma~\ref{lem:binary-mechanical} gives an insertion count $K=T-l$ and an exchange count $k$ satisfying
\[
 \lim_{\nu\to\infty}\frac K\nu=16\eta,\qquad k\le7\eta l+B\le8\eta\nu,
 \qquad K+4k\le P
\]
for all sufficiently large $\nu$.
Apply Claim~\ref{c:exchange-composition-count} of Proposition~\ref{prop:exchange-composition} and Lemma~\ref{lem:word-reindexing} in each case.

First, let us obtain the lower bound from case (i), keeping $\eta$, $p$, and $q$ fixed.
Equation~\eqref{eq:dictionary} gives
\[
 \frac{\log E_3(a,b)}a
 \ge\frac ja\frac{\log C_0(W_{p,q}^j)}j
       -\frac ka\log8-\frac{\log3}a,
 \qquad
 \frac ja\longrightarrow\frac dq,\quad \frac\nu a\longrightarrow1.
\]

Second, let us obtain the upper bound from case (ii).
Equation~\eqref{eq:dictionary} gives
\[
 \frac{\log E_3(a,b)}a
 \le\frac Na\frac{\log C_0(W_{p,q}^N)}N
       +\frac ka\log8-\frac{\log3}a,
 \qquad
 \frac Na\longrightarrow\frac1{dq},\quad \frac\nu a\longrightarrow\frac1d.
\]

Let us apply the periodic limit to both bounds.
We have $j\to\infty$ in case (i) and $N\to\infty$ in case (ii).
Use $k\le8\eta\nu$ and Claim~\ref{c:periodic-limit} of Proposition~\ref{prop:periodic} to obtain
\begin{equation}\label{eq:three-irrational-comparison}
 \begin{aligned}
 \liminf_{\substack{a\to\infty,\ (a,b)\in\cA_3\\\frac{b}{a}\to r}}
 \frac{\log E_3(a,b)}a
 &\ge d\Upsilon_3(\alpha)-8\eta\log8,\\[1ex]
 \limsup_{\substack{a\to\infty,\ (a,b)\in\cA_3\\\frac{b}{a}\to r}}
 \frac{\log E_3(a,b)}a
  &\le\frac{\Upsilon_3(\alpha)+8\eta\log8}{d}.
 \end{aligned}
\end{equation}
For each fixed $\eta$, we choose $\alpha=\frac{p}{q}$, $\delta$, and $B$
before letting $a\to\infty$ and $\frac ba\to r$ through $(a,b)\in\cA_3$.

Let us bound the difference between the upper and lower bounds in~\eqref{eq:three-irrational-comparison}. Lemma~\ref{lem:word-bounds} with $m=3$ and $L=2$, together with Claim~\ref{c:periodic-limit} of Proposition~\ref{prop:periodic} and~\eqref{eq:rational-period}, gives
\[
 \Upsilon_3(\alpha)\le3\log8.
\]
For each fixed $\eta$, equation~\eqref{eq:three-irrational-comparison} bounds the lower and upper limits for every sequence $((a_n,b_n))_{n\ge1}$ in $\cA_3$ satisfying
\[
 \lim_{n\to\infty}a_n=\infty,\qquad
 \lim_{n\to\infty}\frac{b_n}{a_n}=r.
\]

Let us prove that the upper bound minus the lower bound in~\eqref{eq:three-irrational-comparison} tends to zero as $\eta\downarrow0$.
The difference is at most
\[
 3\log8\left(\frac1d-d\right)
 +8\eta\log8\left(1+\frac1d\right)
 \longrightarrow0.
\]
Since $r\in(1,2)\subset I_3$, Lemma~\ref{lem:approaching-sequence} gives a sequence $((a_n,b_n))_{n\ge1}$ in $\cA_3$ satisfying
\[
 \lim_{n\to\infty}a_n=\infty,\qquad
 \lim_{n\to\infty}\frac{b_n}{a_n}=r.
\]
For any two such sequences, denote their lower and upper limits by $L_i$ and $U_i$, $i\in[2]_{\Z}$.
They are finite by Lemma~\ref{lem:word-bounds} and~\eqref{eq:dictionary}.
The common interval in~\eqref{eq:three-irrational-comparison} gives
\[
 U_i-L_j\le3\log8\left(\frac1d-d\right)
             +8\eta\log8\left(1+\frac1d\right)
 \qquad(i,j\in[2]_{\Z}).
\]
Letting $\eta\downarrow0$ gives $U_i\le L_j$ for all $i,j$.
Since $L_i\le U_i$, both sequences converge and their limits are equal.
Define $\Upsilon_3(r)$ to be this common limit.

Let us check positivity of the limit.
The run $\mathcal U(a,b)$ has $a+1$ simple runs of vertical length at most two for all sufficiently large $a$.
Lemma~\ref{lem:word-bounds} and~\eqref{eq:dictionary} give
\[
 \frac{\log E_3(a,b)}a
 \ge\frac{\lceil\frac{a+1}{2}\rceil-1}{a}\log3
 \longrightarrow\frac12\log3.
\]
Hence, $\Upsilon_3(r)\ge\frac12\log3>0$.

\prf{Claim~\ref{c:three-irrational-continuity}} Let us combine the three kinds of slopes before applying the continuity criterion.
Claim~\ref{c:three-irrational-limit}, Theorem~\ref{thm:rational}, and Corollary~\ref{cor:central} give, respectively, the irrational, noninteger rational $r>1$, and $r=1$ cases of
\[
 \lim_{\substack{a\to\infty,\ (a,b)\in\cA_3\\\frac ba\to r}}
 \frac{\log E_3(a,b)}a=\Upsilon_3(r)
 \qquad\left(r\in[1,\frac32]\right).
\]
For rational $r\in(1,\frac32]$, Claim~\ref{c:periodic-limit} of Proposition~\ref{prop:periodic} and~\eqref{eq:rational-period} give $\Upsilon_3(r)>0$.
At $r=1$, we have $\Upsilon_3(1)=\log\frac{256}{27}>0$.
Lemma~\ref{lem:approaching-sequence} gives hypothesis~\ref{hyp:approaching-sequences} of Lemma~\ref{lem:continuity} for $J=[1,\frac32]$.
Claim~\ref{c:continuity-function} of Lemma~\ref{lem:continuity} therefore gives continuity on $[1,\frac32]$.
Claim~\ref{c:continuity-endpoint-function} of Lemma~\ref{lem:continuity} specifies the endpoint limits:
\[
 \lim_{s\downarrow1}\Upsilon_3(s)=\Upsilon_3(1),\qquad
 \lim_{s\uparrow\frac32}\Upsilon_3(s)=\Upsilon_3\left(\frac32\right)
 \qquad\left(s\in\left[1,\frac32\right]\right).
\]
\end{proof}

%% file: sections/reflection-endpoints-and-the-central-value.tex
\section{Reflection, endpoints, and the central value}\label{sec:endpoints}
\subsection{Limits on the full interior interval}\label{sec:interior-transport}
\Needspace{18\baselineskip}
\begin{lemma}[A specified reduction word and its scale factor]\label{lem:specified-reduction-word}
Let $m\in\Z_{\ge3}$. Let
\[
 \lambda_\pm=\frac{m\pm\sqrt{m^2-4}}2,\qquad
 I=(\lambda_-,\lambda_+),\qquad
 T(t)=\frac1{m-t},\qquad
 z(t)=\frac{t-\lambda_-}{\lambda_+-t}\quad(t\in I).
\]
Define the iterates by $T^0(t)=t$ and $T^{j+1}(t)=T(T^j(t))$ for $j\in\Z_{\ge0}$. Then, the following statements hold.

\begin{enumerate}
\item\label{c:reduction-cross-ratio} The maps $t\mapsto\frac1t$, $t\mapsto m-t$, and $T$ preserve $I$. The function $z$ is positive and strictly increasing on $I$. We have
\[
z(1)=\lambda_-,\qquad z(m-1)=\lambda_+,\qquad
z(T^j(t))=\lambda_+^{-2j}z(t)
\quad(t\in I,\ j\in\Z_{\ge0}).
\]
In particular,
\begin{equation}\label{eq:reduction-iterate}
T^j(t)=
\frac{\lambda_-+\lambda_+^{1-2j}z(t)}{1+\lambda_+^{-2j}z(t)}.
\end{equation}

\item\label{c:reduction-first-index} For $r\in I$, define
\[
\nu=\begin{cases}1,&r<1,\\0,&r\ge1,\end{cases}
\qquad
s=\begin{cases}\frac1r,&\nu=1,\\r,&\nu=0,\end{cases}
\]
\[
A=\frac{\log(\frac{z(s)}{\lambda_+})}{2\log\lambda_+},\qquad
N=\max\{0,\lceil A\rceil\},\qquad x=T^N(s),
\]
\[
\varepsilon=\begin{cases}1,&x>\frac m2,\\0,&x\le \frac m2,\end{cases}
\qquad
r_*=\begin{cases}m-x,&\varepsilon=1,\\x,&\varepsilon=0.\end{cases}
\]
Then, $N\in\Z_{\ge0}$, $1\le x\le m-1$, and $1\le r_*\le \frac m2$. Moreover,
\[
N=\min\{j\in\Z_{\ge0}:T^j(s)\le m-1\}.
\]

\item\label{c:reduction-word-scale} Use $\nu$, $N$, $\varepsilon$, and $r_*$ from Claim~\ref{c:reduction-first-index}. Define the integer matrices
\[
S=\begin{pmatrix}0&1\\1&0\end{pmatrix},\qquad
R=\begin{pmatrix}1&0\\m&-1\end{pmatrix},\qquad
G=R^\varepsilon(SR)^N S^\nu=(G_{ij})_{i,j\in[2]_{\Z}}.
\]
For $t\in I$, define
\begin{equation}\label{eq:reduction-scale-map}
c_G(t)=G_{11}+G_{12}t,\qquad
\rho_G(t)=\frac{G_{21}+G_{22}t}{G_{11}+G_{12}t}.
\end{equation}
We have $G_{11}+G_{12}t>0$, $\rho_G(t)\in I$, and $\rho_G(r)=r_*$. Both functions are continuous on $I$. If $(a,b)\in\cA_m$ and $t=\frac ba\in I$, define
\[
\binom{\widetilde a}{\widetilde b}=G\binom ab.
\]
Then,
\[
(\widetilde a,\widetilde b)\in\cA_m,\qquad
E_m(\widetilde a,\widetilde b)=E_m(a,b),\qquad
\frac{\widetilde a}{a}=c_G(t),\qquad
\frac{\widetilde b}{\widetilde a}=\rho_G(t).
\]
In particular, $c_m(r):=c_G(r)>0$ depends only on $m$ and $r$. If $r\in[1,\frac m2]$, then $G$ is the identity matrix, $c_m(r)=1$, and $r_*=r$.
\end{enumerate}
\end{lemma}
\begin{proof}
\prf{Claim~\ref{c:reduction-cross-ratio}} Let us use $\lambda_\pm=\frac{m\pm\sqrt{m^2-4}}2$. We have
\[
 \lambda_-+\lambda_+=m,\qquad \lambda_-\lambda_+=1,
 \qquad 0<\lambda_-<1<m-1<\lambda_+<m.
\]
For the strict inequality $m-1<\lambda_+$, evaluate the polynomial $t^2-mt+1$ at $m-1$: its value is $2-m<0$. Thus, $m-1$ lies between the two roots. The other inequalities follow from their positivity, product one, and sum $m>2$. Consequently, if $t\in I$, then
\[
 \lambda_-<\frac1t<\lambda_+,\qquad
 \lambda_-<m-t<\lambda_+.
\]
Their composition $T(t)=\frac1{m-t}$ also belongs to $I$. We have $m-t>0$.

Let us check that the change of variable is strictly increasing.
For $\lambda_-<u<v<\lambda_+$,
\[
 z(v)-z(u)=
 \frac{(\lambda_+-\lambda_-)(v-u)}{(\lambda_+-v)(\lambda_+-u)}>0.
\]
For $t\in I$, we have $t-\lambda_->0$ and $\lambda_+-t>0$. The identities $\lambda_-\lambda_+=1$ and $\lambda_-+\lambda_+=m$ give
\[
 z(1)=\frac{1-\lambda_-}{\lambda_+-1}=\lambda_-,\qquad
 z(m-1)=\frac{\lambda_+-1}{1-\lambda_-}=\lambda_+.
\]
They also give
\[
 \begin{aligned}
 z(T(t))
 &=\frac{1-(m-t)\lambda_-}{(m-t)\lambda_+-1}\\
 &=\frac{\lambda_-(t-\lambda_-)}{\lambda_+(\lambda_+-t)}
 =\lambda_+^{-2}z(t).
 \end{aligned}
\]
Induction on $j$ gives $z(T^j(t))=\lambda_+^{-2j}z(t)$. Solving $z(u)=\frac{u-\lambda_-}{\lambda_+-u}$ gives
\[
 u=\frac{\lambda_-+\lambda_+z(u)}{1+z(u)}.
\]
Substituting $u=T^j(t)$ and $z(T^j(t))=\lambda_+^{-2j}z(t)$ proves~\eqref{eq:reduction-iterate}.

\prf{Claim~\ref{c:reduction-first-index}} Let us verify the specified stopping index. Claim~\ref{c:reduction-cross-ratio} gives $s\in I$ and $s\ge1$. Since $z(s)>0$ and $\lambda_+>1$, the real number $A$ is defined. For every integer $j\ge0$, Claim~\ref{c:reduction-cross-ratio} and strict monotonicity of $z$ give
\[
 \begin{aligned}
 T^j(s)\le m-1
 &\Longleftrightarrow\lambda_+^{-2j}z(s)\le\lambda_+\\
 &\Longleftrightarrow j\ge A.
 \end{aligned}
\]
Therefore, $N=\max\{0,\lceil A\rceil\}$ is exactly the least nonnegative integer with $T^N(s)\le m-1$. If $N=0$, then $x=s\ge1$. If $N\ge1$, minimality gives $y:=T^{N-1}(s)>m-1$, while $y\in I$ gives $y<m$. Hence,
\[
 0<m-y<1,\qquad x=\frac1{m-y}>1.
\]
This proves $x\in[1,m-1]$. If $x\le \frac m2$, then $r_*=x\in[1,\frac m2]$. If $x>\frac m2$, then
\[
 1\le m-x<\frac m2,
\]
so $r_*=m-x\in[1,\frac m2]$.

\prf{Claim~\ref{c:reduction-word-scale}} Let us factor $G=R^\varepsilon(SR)^N S^\nu$. Put $d=\nu+2N+\varepsilon$. Define matrices $L_i$, $i\in[d]_{\Z}$, in application order by
\begin{equation}\label{eq:reduction-factor-indices}
 \begin{gathered}
 L_1=S\quad\text{if }\nu=1,\\
 L_{\nu+2j-1}=R,\quad L_{\nu+2j}=S\quad(j\in[N]_{\Z}),\\
 L_d=R\quad\text{if }\varepsilon=1.
 \end{gathered}
\end{equation}
The index ranges in~\eqref{eq:reduction-factor-indices} are disjoint and cover $[d]_{\Z}$. We have $G=L_d\cdots L_1$. For $d=0$, this product is the identity. For $t\in I$, put $t_0=t$ and, for $i\in[d]_{\Z}$, define
\[
 t_i=\begin{cases}\frac1{t_{i-1}},&L_i=S,\\m-t_{i-1},&L_i=R.\end{cases}
\]

Let us compute the first coordinate after applying $G$.
Claim~\ref{c:reduction-cross-ratio} gives $t_i\in I$ for $i\in[0,d]_{\Z}$.
For $i\in[d]_{\Z}$ and $a>0$, we have
\[
 S\binom{a}{at_{i-1}}=at_{i-1}\binom{1}{1/t_{i-1}},\qquad
 R\binom{a}{at_{i-1}}=a\binom{1}{m-t_{i-1}}.
\]
Both first coordinates are positive since $t_{i-1}\in I\subset(0,m)$.
Starting with $(1,t)^{\mathsf T}$ and applying $L_1,\ldots,L_d$ gives
\[
 c_G(t)=\prod_{\substack{i\in[d]_{\Z}\\L_i=S}}t_{i-1}>0,
 \qquad \rho_G(t)=t_d\in I.
\]
The empty product equals one.
The construction in Claim~\ref{c:reduction-first-index} gives $t_d=r_*$ when $t=r$.

Let us prove continuity of $c_G$ and $\rho_G$ at an arbitrary $t\in I$.
For any sequence $(t_n)_{n\ge1}$ in $I$ satisfying $\lim_{n\to\infty}t_n=t$, we have
\[
 \lim_{n\to\infty}(G_{11}+G_{12}t_n)=G_{11}+G_{12}t=c_G(t)>0.
\]
Equation~\eqref{eq:reduction-scale-map} and the quotient limit give
\[
 \lim_{n\to\infty}c_G(t_n)=c_G(t),\qquad
 \lim_{n\to\infty}\rho_G(t_n)
 =\frac{G_{21}+G_{22}t}{G_{11}+G_{12}t}=\rho_G(t).
\]

Let us identify the corresponding transformations of dimension pairs.
On integer dimension pairs the two factors act by
\[
 S\binom ab=\binom ba,\qquad
 R\binom ab=\binom a{ma-b}.
\]
Each intermediate pair is positive because its slope belongs to $I$. Claims~\ref{c:moduli-symmetries-identities}, \ref{c:moduli-coprimality}, and~\ref{c:moduli-quadratic} of Lemma~\ref{lem:moduli-symmetries} and~\eqref{eq:nonempty-domain} preserve membership in $\cA_m$ and the Euler characteristic at each step. More explicitly,
\[
 \gcd(b,a)=\gcd(a,b),\qquad \gcd(a,ma-b)=\gcd(a,b),
\]
\[
 m ba-b^2-a^2=mab-a^2-b^2,
\]
\[
 ma(ma-b)-a^2-(ma-b)^2=mab-a^2-b^2.
\]
The Euler identities are those of Claim~\ref{c:moduli-symmetries-identities} of Lemma~\ref{lem:moduli-symmetries}. Multiplication by $G$ gives
\[
 \frac{\widetilde a}{a}=G_{11}+G_{12}t=c_G(t),\qquad
 \frac{\widetilde b}{\widetilde a}
 =\frac{G_{21}+G_{22}t}{G_{11}+G_{12}t}=\rho_G(t).
\] If $1\le r\le \frac m2$, then $\nu=0$, $r\le m-1$, $N=0$, and $\varepsilon=0$. Thus, $G=\left(\begin{smallmatrix}1&0\\0&1\end{smallmatrix}\right)$, $c_m(r)=1$, and $r_*=r$.
\end{proof}

\begin{lemma}[Endpoint folding]\label{lem:reduction-endpoint-fold}
Let $m\in\Z_{\ge3}$, put $I=(\lambda_-,\lambda_+)$, and let $\rho\in[1,\frac m2]$. Define
\[
 \delta_\rho=
 \begin{cases}
 \frac12\min\{\rho-1,\frac m2-\rho\},&1<\rho<\frac m2,\\
 \frac{m-2}{2m},&\rho=1,\\
 \frac{m-2}{4},&\rho=\frac m2.
 \end{cases}
\]
The cases are disjoint because $m>2$. We have $\delta_\rho>0$. For $(a,b)\in\cA_m$, $t=\frac ba\in I$, and $|t-\rho|<\delta_\rho$, define
\[
 \operatorname{Fold}_\rho(a,b)=(a',b')=
 \begin{cases}
 (a,b),&1<\rho<\frac m2,\\
 (\min(a,b),\max(a,b)),&\rho=1,\\
 (a,\min(b,ma-b)),&\rho=\frac m2.
 \end{cases}
\]
Then,
\[
 (a',b')\in\cA_m,\qquad E_m(a',b')=E_m(a,b),\qquad
 1\le \frac{b'}{a'}\le \frac m2.
\]
Moreover, the ratios are exactly
\[
 \left(\frac{a'}a,\frac{b'}{a'}\right)=
 \begin{cases}
 (1,t),&1<\rho<\frac m2,\\
 (\min(1,t),\max(t,\frac1t)),&\rho=1,\\
 (1,\min(t,m-t)),&\rho=\frac m2.
 \end{cases}
\]
Consequently, let $((a_n,b_n))_{n\ge1}$ be a sequence in $\cA_m$ satisfying
\[
 \lim_{n\to\infty}\frac{b_n}{a_n}=\rho.
\]
Then, there is $N\in\Z_{>0}$ such that $\operatorname{Fold}_\rho(a_n,b_n)=(a'_n,b'_n)$ is defined for every $n\ge N$ and satisfies
\[
 (a'_n,b'_n)\in\cA_m,\qquad E_m(a'_n,b'_n)=E_m(a_n,b_n),\qquad 1\le\frac{b'_n}{a'_n}\le\frac m2.
\]
We have
\[
 \lim_{n\to\infty}\frac{a'_n}{a_n}=1,\qquad
 \lim_{n\to\infty}\frac{b'_n}{a'_n}=\rho.
\]
\end{lemma}
\begin{proof}
First, suppose $1<\rho<\frac m2$. Then,
\[
 1<\rho-\delta_\rho<t<\rho+\delta_\rho<\frac m2.
\]
The pair is unchanged. Thus, all conclusions follow.

Second, suppose $\rho=1$. The chosen radius satisfies
\[
 1-\delta_1=\frac{m+2}{2m}>\frac2m,
 \qquad
 1+\delta_1=\frac{3m-2}{2m}<\frac m2,
\]
where $\frac{3m-2}{2m}<\frac m2$ follows from $(m-1)(m-2)>0$. Therefore,
\[
 \frac2m<t<\frac m2,
 \qquad 1\le\max(t,\frac1t)<\frac m2.
\]
If $t\ge1$, then $\operatorname{Fold}_1(a,b)=(a,b)$. If $t<1$, then $\operatorname{Fold}_1(a,b)=(b,a)$. These two cases give
\[
 \frac{a'}a=\min(1,t),\qquad \frac{b'}{a'}=\max(t,\frac1t).
\]
Vertex exchange preserves membership in $\cA_m$ and the Euler characteristic by Claims~\ref{c:moduli-symmetries-identities}, \ref{c:moduli-coprimality}, and~\ref{c:moduli-quadratic} of Lemma~\ref{lem:moduli-symmetries} and~\eqref{eq:nonempty-domain}.

Third, suppose $\rho=\frac m2$. In this case
\[
 1<\frac{m+2}{4}=\rho-\delta_\rho<t
 <\rho+\delta_\rho=\frac{3m-2}{4}<m-1.
\]
Therefore,
\[
 1<\min(t,m-t)\le \frac m2.
\]
For $t\le \frac m2$, the pair $\operatorname{Fold}_{\frac m2}(a,b)$ equals $(a,b)$.
For $t>\frac m2$, it equals $(a,ma-b)$.
Thus, $\frac{a'}a=1$ and $\frac{b'}{a'}=\min(t,m-t)$.
Claims~\ref{c:moduli-symmetries-identities}, \ref{c:moduli-coprimality}, and~\ref{c:moduli-quadratic} of Lemma~\ref{lem:moduli-symmetries} and~\eqref{eq:nonempty-domain} give $(a',b')\in\cA_m$ and $E_m(a',b')=E_m(a,b)$.

Let us pass to the limit in each case. We have $[1,\frac m2]\subset I$. Convergence to $\rho$ puts $t_n=\frac{b_n}{a_n}$ in $I\cap(\rho-\delta_\rho,\rho+\delta_\rho)$ for all sufficiently large $n$. The limits $\frac{a'_n}{a_n}\to1$ and $\frac{b'_n}{a'_n}\to\rho$ follow from continuity of $\min$, $\max$, and inversion at a positive number. At $\rho=1$, these limits are $\min(1,1)=1$ and $\max(1,1)=1$. At $\rho=\frac m2$, the slope limit is $\min(\frac m2,m-\frac m2)=\frac m2$.
\end{proof}

\begin{lemma}[Reduction of sequences of dimension pairs]\label{lem:approaching-dimension-reduction}
Let $m\in\Z_{\ge3}$ and $r\in(\lambda_-,\lambda_+)$.
Then, there are numbers $r_*\in[1,\frac m2]$ and $c_m(r)>0$ depending only on $m$ and $r$ with the following property.
For each sequence $((a_n,b_n))_{n\ge1}$ in $\cA_m$ such that
\begin{equation}\label{eq:approaching-reduction-input}
 \lim_{n\to\infty}a_n=\infty,\qquad
 \lim_{n\to\infty}\frac{b_n}{a_n}=r,
\end{equation}
there are an integer $N\ge1$ and pairs $(a'_n,b'_n)\in\cA_m$ for $n\ge N$ such that
\begin{equation}\label{eq:approaching-reduction-pairs}
 1\le\frac{b'_n}{a'_n}\le\frac m2,\qquad
 E_m(a'_n,b'_n)=E_m(a_n,b_n)\qquad(n\ge N),\qquad\text{and}
\end{equation}
\begin{equation}\label{eq:approaching-reduction-limits}
 \frac{a'_n}{a_n}\longrightarrow c_m(r),\qquad
 \frac{b'_n}{a'_n}\longrightarrow r_*\qquad(n\longrightarrow\infty).
\end{equation}
\end{lemma}
\begin{proof}
Let us apply one fixed reduction word to the sequence.
Put $I=(\lambda_-,\lambda_+)$. Let $G$, $r_*$, and $c_m(r)=G_{11}+G_{12}r>0$ be the objects specified in Lemma~\ref{lem:specified-reduction-word}. By~\eqref{eq:approaching-reduction-input}, the slopes $\frac{b_n}{a_n}$ belong to $I$ for all sufficiently large $n$. For these indices, define
\[
 \binom{\widetilde a_n}{\widetilde b_n}=G\binom{a_n}{b_n}.
\]
Claim~\ref{c:reduction-word-scale} of Lemma~\ref{lem:specified-reduction-word} gives
\[
 \begin{gathered}
 (\widetilde a_n,\widetilde b_n)\in\cA_m,
 \qquad E_m(\widetilde a_n,\widetilde b_n)=E_m(a_n,b_n),\\
 \frac{\widetilde a_n}{a_n}
 =G_{11}+G_{12}\frac{b_n}{a_n}\longrightarrow c_m(r),\\
 \frac{\widetilde b_n}{\widetilde a_n}
 =\frac{G_{21}+G_{22}(\frac{b_n}{a_n})}{G_{11}+G_{12}(\frac{b_n}{a_n})}
 \longrightarrow r_*.
 \end{gathered}
\]
Apply Lemma~\ref{lem:reduction-endpoint-fold} with $\rho=r_*$. Define
\[
 (a'_n,b'_n)=\operatorname{Fold}_{r_*}(\widetilde a_n,\widetilde b_n)
\]
for $n\ge N$, where $N\ge1$ is large enough for Lemma~\ref{lem:reduction-endpoint-fold} to apply. That lemma gives $(a'_n,b'_n)\in\cA_m$, $E_m(a'_n,b'_n)=E_m(\widetilde a_n,\widetilde b_n)$, and $1\le\frac{b'_n}{a'_n}\le\frac m2$. Thus~\eqref{eq:approaching-reduction-pairs} holds. Lemma~\ref{lem:reduction-endpoint-fold} gives
\[
 \frac{a'_n}{a_n}
 =\frac{a'_n}{\widetilde a_n}\frac{\widetilde a_n}{a_n}
 \longrightarrow1\cdot c_m(r)=c_m(r),
 \qquad \frac{b'_n}{a'_n}\longrightarrow r_*.
\]
This proves~\eqref{eq:approaching-reduction-limits}. The definitions of $G$, $r_*$, and $c_m(r)$ use only $m$ and $r$.
\end{proof}

\begin{proposition}\label{prop:interior-limit}
For every integer $m\ge3$, there is a function
$F_m:(\lambda_-,\lambda_+)\to(0,\infty)$ such that the following hold.
\begin{enumerate}
\item\label{c:interior-limit} We have
\[
 \lim_{\substack{a\to\infty,\ (a,b)\in\cA_m\\\frac{b}{a}\to r}}
 \frac{\log E_m(a,b)}a=F_m(r)
 \qquad\bigl(r\in(\lambda_-,\lambda_+)\bigr).
\]
\item\label{c:interior-continuity} The function $F_m$ is continuous.
\end{enumerate}
\end{proposition}
\begin{proof}
\prf{Claim~\ref{c:interior-limit}} Let us first obtain the limit for $r\in[1,\frac m2]$.
Fix such an $r$ and consider $(a,b)\in\cA_m$ such that
\[
 a\longrightarrow\infty,\qquad \frac ba\longrightarrow r,\qquad
 1\le\frac ba\le\frac m2.
\]

First, let us obtain the limit when $m\ge4$.
Define $\Upsilon_m(1)=\ell(R_1)$ and $\Upsilon_m(m-2)=\ell(R_{m-2})$.
For all sufficiently large $a$, we have $b>a\ge2$ and
$1<\Gamma(a,b)<\frac ba\le \frac{m}{2}\le m-2$ by~\eqref{eq:ambient}.
The floor differences in~\eqref{eq:frame-word} and the appended simple run have vertical lengths in $[m-2]_{\Z}$.
Claims~\ref{c:integer-limit} and~\ref{c:integer-endpoints} of Proposition~\ref{prop:integer},
Theorem~\ref{thm:rational}, and Theorem~\ref{thm:irrational-large-m} give the limits on $[1,\frac{m}{2}]$.
Lemma~\ref{lem:word-bounds} and~\eqref{eq:dictionary} give positivity.

Second, let us obtain the limit when $m=3$.
Theorem~\ref{thm:rational}, Claim~\ref{c:three-irrational-limit} of
Proposition~\ref{prop:three-irrational}, and Corollary~\ref{cor:central} give the limits on $[1,\frac32]$.
Claim~\ref{c:three-irrational-continuity} of Proposition~\ref{prop:three-irrational} gives positivity.
Thus, the positive limits on $[1,\frac m2]$ equal $\Upsilon_m(r)$.

Let us obtain the limit at each $r\in(\lambda_-,\lambda_+)$.
Let $((a_n,b_n))_{n\ge1}$ be any sequence in $\cA_m$ satisfying
\[
 \lim_{n\to\infty}a_n=\infty,\qquad
 \lim_{n\to\infty}\frac{b_n}{a_n}=r.
\]
Apply Lemma~\ref{lem:approaching-dimension-reduction} to this sequence.
Since $\frac{a'_n}{a_n}\to c_m(r)>0$, we have $a'_n\to\infty$.
Hence,
\[
 \frac{\log E_m(a_n,b_n)}{a_n}
 =\frac{a'_n}{a_n}\frac{\log E_m(a'_n,b'_n)}{a'_n}
 \longrightarrow c_m(r)\Upsilon_m(r_*)>0.
\]

Let us define $F_m$ on $(\lambda_-,\lambda_+)$.
The numbers $r_*$ and $c_m(r)$ depend only on $m$ and $r$.
Lemma~\ref{lem:approaching-sequence} supplies at least one sequence approaching $r$.
Define $F_m(r)$ to be the common limit of $\frac1{a_n}\log E_m(a_n,b_n)$ along all such sequences.
Then,
\[
 F_m(r)=c_m(r)\Upsilon_m(r_*)
 \qquad\bigl(r\in(\lambda_-,\lambda_+)\bigr).
\]
For $r\in[1,\frac{m}{2}]$, Claim~\ref{c:reduction-word-scale} of Lemma~\ref{lem:specified-reduction-word}
gives $r_*=r$ and $c_m(r)=1$.
Thus, Claim~\ref{c:interior-limit} holds at every interior slope.

Let us identify $F_m(r)$ with $\Upsilon_m(r)$ on the entire interval $(1,M_{\mathrm{run}}(m))$.
Fix $r$ in this interval and choose a sequence $((a_n,b_n))_{n\ge1}$ as in~\eqref{eq:approaching-sequence}.

If $m\ge4$ and $r$ is an integer, equations~\eqref{eq:ambient} and~\eqref{eq:frame-word} give $(a_n,b_n)\in\mathcal S_m$ for all sufficiently large $n$.
Claim~\ref{c:integer-limit} of Proposition~\ref{prop:integer} and~\eqref{eq:rational-period} then give the limit $\Upsilon_m(r)=\ell(R_r)$.

For noninteger rational $r$, Theorem~\ref{thm:rational} and~\eqref{eq:rational-period} give the limit $\Upsilon_m(r)$.
For irrational $r$, Theorem~\ref{thm:irrational-large-m} applies when $m\ge4$, and Claim~\ref{c:three-irrational-limit} of Proposition~\ref{prop:three-irrational} applies when $m=3$.

Comparing each limit with Claim~\ref{c:interior-limit} gives
\[
 F_m(r)=\Upsilon_m(r)\qquad\bigl(1<r<M_{\mathrm{run}}(m)\bigr).
\]

Let $r\in(\lambda_-,\lambda_+)$.
First, let us prove $F_m(r)=rF_m(1/r)$.
By Lemma~\ref{lem:approaching-sequence}, choose a sequence $((a_n,b_n))_{n\ge1}$ in $\cA_m$ satisfying
\[
 \lim_{n\to\infty}a_n=\infty,\qquad
 \lim_{n\to\infty}\frac{b_n}{a_n}=r.
\]
Since $r>0$, we have $b_n\to\infty$ and $\frac{a_n}{b_n}\to\frac1r$.
Claim~\ref{c:moduli-symmetries-identities} of Lemma~\ref{lem:moduli-symmetries} gives
\[
 \frac{\log E_m(a_n,b_n)}{a_n}
 =\frac{b_n}{a_n}\frac{\log E_m(b_n,a_n)}{b_n}
 \longrightarrow rF_m\left(\frac1r\right).
\]

Second, let us prove $F_m(r)=F_m(m-r)$.
Since $m-r>0$, choose $n_0\in\Z_{>0}$ such that
\[
 m-\frac{b_n}{a_n}>\frac{m-r}{2}>0\qquad(n\ge n_0).
\]
Thus, $ma_n-b_n>0$ for $n\ge n_0$, and
\[
 \frac{ma_n-b_n}{a_n}\longrightarrow m-r,\qquad
 \frac{\log E_m(a_n,b_n)}{a_n}
 =\frac{\log E_m(a_n,ma_n-b_n)}{a_n}
 \longrightarrow F_m(m-r).
\]
Both transformed sequences belong to $\cA_m$ by Lemma~\ref{lem:moduli-symmetries}
and~\eqref{eq:nonempty-domain}. Comparing these limits with Claim~\ref{c:interior-limit} gives
\begin{equation}\label{eq:entropy-symmetries}
 F_m(r)=rF_m\left(\frac1r\right),\qquad F_m(r)=F_m(m-r).
\end{equation}

\prf{Claim~\ref{c:interior-continuity}} Let us prove continuity at each interior slope.
Let $r\in(\lambda_-,\lambda_+)$ and choose
$\lambda_-<h_1<r<h_2<\lambda_+$.
Claim~\ref{c:interior-limit} and Lemma~\ref{lem:approaching-sequence} verify the hypotheses of
Lemma~\ref{lem:continuity} for $F_m|_{[h_1,h_2]}$.
Claim~\ref{c:continuity-function} of Lemma~\ref{lem:continuity} gives continuity at $r$.
\end{proof}

\subsection{\texorpdfstring{Proof of Theorem~\ref{thm:main} and upper bound}{Proof of Theorem 1.1 and upper bound}}
\begin{proof}[Proof of Theorem~\ref{thm:main}]
Let us prove the limit of $\frac1a\log E_m(a,b)$ at both endpoints.
Claim~\ref{c:interior-limit} of Proposition~\ref{prop:interior-limit} gives positivity and~\eqref{eq:main-limit}
on $(\lambda_-,\lambda_+)$. Let $F_m(\lambda_\pm)=0$.
Fix $r\in\{\lambda_-,\lambda_+\}$ and any sequence $((a_n,b_n))_{n\ge1}$ in $\cA_m$ such that
$a_n\to\infty$ and $r_n=\frac{b_n}{a_n}\to r$.
Since $mr-1-r^2=0$, we have
\[
 \max(mr_n-1-r_n^2,0)\longrightarrow0,
 \qquad
 \frac{\log(\frac{2^m}{m})}{a_n}\longrightarrow0
 \qquad(n\to\infty).
\]
The bound~\eqref{eq:finite-endpoint} gives
\[
 0\le\frac{\log E_m(a_n,b_n)}{a_n}
 \le m\log2\sqrt{\frac{m+2}{m-2}}\sqrt{\max(mr_n-1-r_n^2,0)}
   +\frac{\log(\frac{2^m}{m})}{a_n}\longrightarrow0.
\]
Thus, we obtain
\eqref{eq:main-limit} for $r=\lambda_-$ and $r=\lambda_+$.

Let us prove continuity on the closed interval $I_m$.
For every $r\in I_m$, Lemma~\ref{lem:approaching-sequence} supplies a sequence satisfying~\eqref{eq:approaching-sequence}.
Apply Claim~\ref{c:continuity-function} of Lemma~\ref{lem:continuity} with $J=I_m=[\lambda_-,\lambda_+]$.
\end{proof}

\begin{corollary}[Quantitative upper bound]\label{cor:upper-bound}
For each integer $m\ge3$, the function $F_m$ in
Theorem~\ref{thm:main} satisfies
\begin{equation}\label{eq:main-bound}
 F_m(r)\le m\log2\sqrt{\frac{m+2}{m-2}}
       \sqrt{(\lambda_+-r)(r-\lambda_-)}
       \qquad(r\in I_m).
\end{equation}
\end{corollary}
\begin{proof}
Let us obtain~\eqref{eq:main-bound} from~\eqref{eq:finite-endpoint}.
For $r\in(\lambda_-,\lambda_+)$, choose a sequence $((a_\nu,b_\nu))_{\nu\ge1}$ in $\cA_m$ satisfying~\eqref{eq:approaching-sequence}.
Put $r_\nu=\frac{b_\nu}{a_\nu}$. Equations~\eqref{eq:main-limit} and~\eqref{eq:finite-endpoint} give
\[
 \begin{aligned}
 F_m(r)
 &=\lim_{\nu\to\infty}\frac{\log E_m(a_\nu,b_\nu)}{a_\nu}\\
 &\le\lim_{\nu\to\infty}\left(
 m\log2\sqrt{\frac{m+2}{m-2}}\sqrt{\max(mr_\nu-1-r_\nu^2,0)}
 +\frac{\log(\frac{2^m}{m})}{a_\nu}\right)\\
 &=m\log2\sqrt{\frac{m+2}{m-2}}\sqrt{mr-1-r^2}.
 \end{aligned}
\]
The identity $mr-1-r^2=(\lambda_+-r)(r-\lambda_-)$ proves~\eqref{eq:main-bound}.
For $r=\lambda_\pm$, both sides are zero.
\end{proof}

\subsection{The value at slope one and some figures for \texorpdfstring{$F_m$}{Fm}}
\begin{corollary}\label{cor:central-value}
For every $m\ge3$,
\begin{equation}\label{eq:central-value}
 \begin{aligned}
 \lim_{\substack{a\to\infty,\ (a,b)\in\cA_m\\\frac{b}{a}\to1}}
 \frac{\log E_m(a,b)}a
 &=F_m(1)\\
 &=(m-1)^2\log((m-1)^2)-m(m-2)\log(m(m-2)).
 \end{aligned}
\end{equation}
\end{corollary}
\begin{proof}
Let us obtain the limit along dimensions $(a,a+1)$.
The adjacent-dimension formula of~\cite[Theorem 6.6]{Weist} is
\[
 E_m(a,a+1)=\frac{m}{(a+1)((m-1)a+m)}
          \binom{(m-1)^2a+(m-1)m}{a}.
\]
Fix $m$ and let $c=(m-1)^2$ and $d=m(m-1)$.
Stirling's formula $\log(n!)=n\log n-n+O(\log n)$ gives, as $a\to\infty$ through positive integers,
\[
 \begin{aligned}
 \log\binom{ca+d}{a}
 &=(ca+d)\log(ca+d)-a\log a\\
 &\quad-((c-1)a+d)\log((c-1)a+d)+O(\log a)\\
 &=a\bigl(c\log c-(c-1)\log(c-1)\bigr)+O(\log a).
 \end{aligned}
\]
The constants in the error terms depend only on $m$.
For the prefactor, we have
\[
 \log\frac{m}{(a+1)((m-1)a+m)}=-2\log a+O(1).
\]
Dividing by $a$ gives
\[
 \begin{aligned}
 \frac1a\log\binom{ca+(m-1)m}{a}
 &\longrightarrow c\log c-(c-1)\log(c-1),\\
 \frac1a\log\frac{m}{(a+1)((m-1)a+m)}&\longrightarrow0.
 \end{aligned}
\]
Since $c-1=m(m-2)$, equation~\eqref{eq:main-limit} gives~\eqref{eq:central-value}.
\end{proof}

\begin{figure}[H]
\centering
\includegraphics[width=\textwidth]{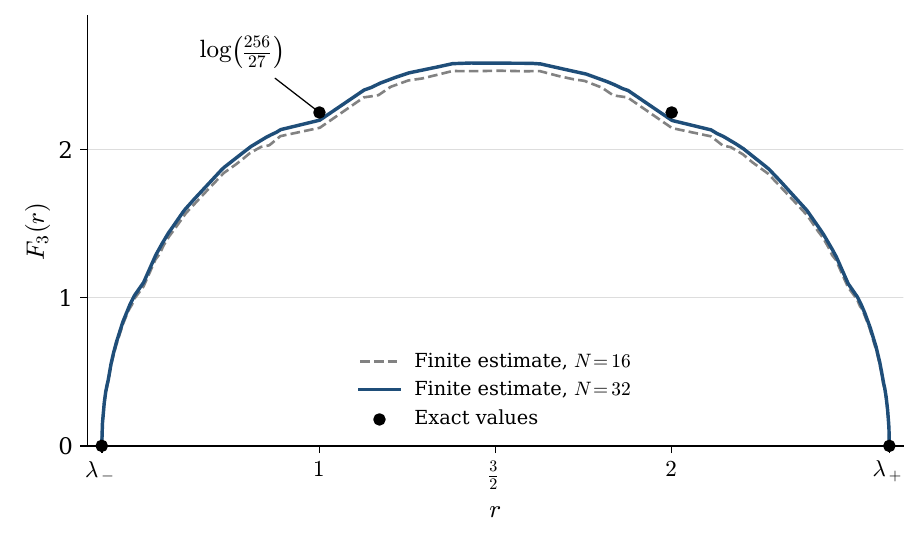}
\caption[Numerical illustration of $F_3$]{Numerical illustration of $F_3$ on $I_3$.
Here $\lambda_\pm=\frac{3\pm\sqrt5}{2}$.
For each reduced $\alpha=\frac pq\in[1,\frac32]$ with $q\le8$,
let $N\in\{16,32\}$ and $j=\floor{\frac Nq}$. We use
$\widehat F_N(\alpha)=
\frac{\log C_0(W_{p,q}^{2j})-\log C_0(W_{p,q}^j)}{jq}$.
We interpolate these values linearly and extend them by
\eqref{eq:entropy-symmetries}. The dots mark the exact values
$F_3(\lambda_\pm)=0$ and
$F_3(1)=F_3(2)=\log\left(\frac{256}{27}\right)$.
The curves are finite approximations; we do not assert an error bound.}
\label{fig:f3}
\end{figure}

\begin{figure}[H]
\centering
\includegraphics[width=\textwidth]{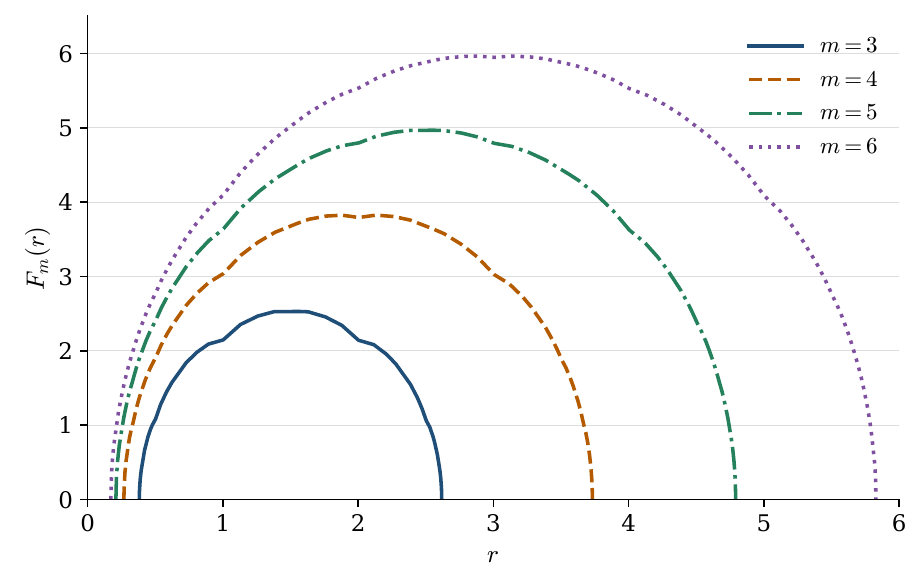}
\caption[Finite approximations for $m=3$, $4$, $5$, and $6$]{Finite approximations to
$F_m(r)$ for $m=3$, $4$, $5$, and $6$, each on its interval $I_m$.
For $\alpha=1+\frac{k}{8}$ with $k\in[0,4(m-2)]_{\Z}$, let
$\alpha=\frac pq$ in lowest terms and let $j=\frac{16}{q}$.
We evaluate $C_0$ for the indicated $m$ and use
$\widehat F_m(\alpha)=
\frac{\log C_0(W_{p,q}^{2j})-\log C_0(W_{p,q}^j)}{jq}$.
We interpolate these values linearly on $[1,\frac m2]$ and extend them by
\eqref{eq:entropy-symmetries}. We use the exact value zero at both endpoints.
The curves are finite approximations; we do not assert an error bound.}
\label{fig:fm-comparison}
\end{figure}

%% file: sections/conjectural-log-concavity.tex
\Needspace{0.72\textheight}
\section{Conjectural log-concavity}\label{sec:log-concavity-conjecture}

We state the following conjecture on Euler characteristics near slope one.
See \cite{CDW, Oko} for related studies and
\cite[Section 1.12]{Oka2026} for a study of log-concavity
in mathematics and statistical mechanics.

\begin{conjecture}[Euler-characteristic ratios near slope one]\label{conj:central-ratio-log-concavity}
Let $a\in\Z_{\ge5}$ satisfy $\gcd(a,6)=1$. Then,
\[
 \frac{E_3(a,a+4)}{E_3(a,a+3)}
 \le
 \frac{E_3(a,a+3)}{E_3(a,a+2)}.
\]
\end{conjecture}

We provide the following figure on Conjecture~\ref{conj:central-ratio-log-concavity}.

\begin{figure}[H]
\centering
\includegraphics[width=0.9\textwidth]{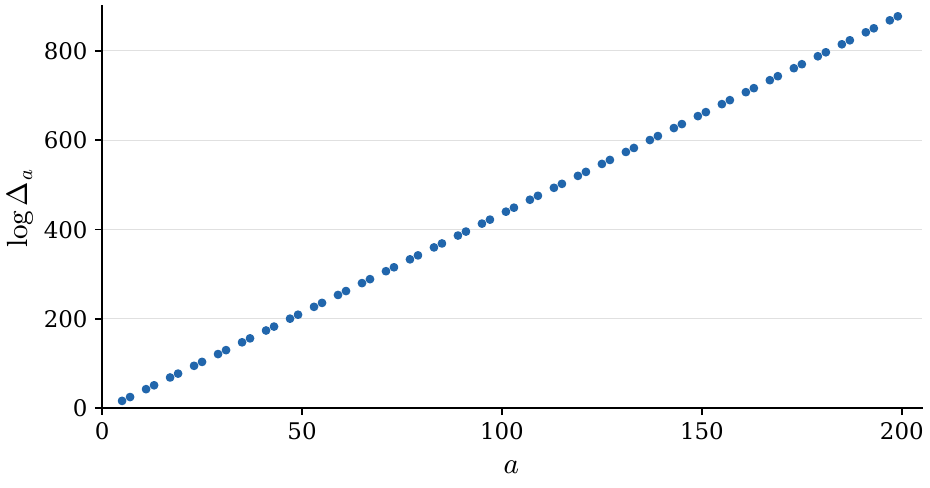}
\caption[Logarithms of the log-concavity differences]{Let
$\Delta_a=E_3(a,a+3)^2-E_3(a,a+2)E_3(a,a+4)$.
The points plot $\log\Delta_a$ by explicit computation for
$a\in[5,200]_{\Z}$ satisfying $\gcd(a,6)=1$.
All 66 differences are positive.}
\label{fig:conjecture-delta-log}
\end{figure}

%% file: sections/acknowledgements.tex
\section*{Acknowledgements}
The author formulated the initial questions, conducted the mathematical
discussions, prepared the manuscript, and organized the Lean
formalizations, using ChatGPT~6 Astra and Claude Fable~5.1 with MathBox.
The author is responsible for the mathematical content and its
presentation.